\documentclass[a4,12pt,oneside]{book}

\usepackage[utf8x]{inputenc}
\usepackage[english]{babel}
\usepackage{url}
\usepackage{amssymb,amsmath,amsthm,amsfonts}
\usepackage{mathtools}

\usepackage{fancyhdr}
\usepackage{blindtext} 

\fancypagestyle{plain}{%
\fancyhf{} 
\fancyfoot[C]{\thepage} 
}

\title{NOTES ON RIDGE FUNCTIONS AND NEURAL NETWORKS}
\author{Vugar E. Ismailov}

\begin{document}
\maketitle

\newpage \vspace*{3cm}
\thispagestyle{empty}
\rightline{\large \emph{To the Memory of My Parents}}

\chapter*{Preface}

These notes are about \textit{ridge functions}. Recent years have
witnessed a flurry of interest in these functions. Ridge functions
appear in various fields and under various guises. They appear in fields as diverse as
partial differential equations (where they are called \textit{plane waves}),
computerized tomography and statistics. These functions are also
the underpinnings of many central models in neural networks.

We are interested in ridge functions from the point of view of approximation theory.
The basic goal in approximation theory is to approximate complicated objects by simpler objects.
Among many classes of multivariate functions, linear combinations of ridge functions are a class of simpler functions.
These notes study some problems of approximation of multivariate functions by linear combinations of ridge functions.
We present here various properties of these functions.
The questions we ask are as follows. When can a multivariate function be expressed as a linear combination of ridge functions from a certain class? When do such linear combinations represent each multivariate function? If a precise representation is not possible, can one approximate arbitrarily well? If well approximation fails, how can one compute/estimate the error of approximation, know that a best approximation exists? How can one characterize and construct best approximations? If a smooth function is a sum of arbitrarily behaved ridge functions, is it true that it can be expressed as a sum of smooth ridge functions?
We also study properties of generalized ridge functions, which are very much related to linear superpositions and Kolmogorov's famous superposition theorem. These notes end with a few applications of ridge functions to the problem of approximation by single and two hidden layer neural networks with a restricted set of weights.

We hope that these notes will be useful and interesting to both researchers and graduate students.

\newpage

\tableofcontents

\newpage

\chapter*{Introduction}

\addcontentsline{toc}{chapter}{Introduction}

Recent years have seen a growing interest in the study of special
multivariate functions called ridge functions. A \textit{ridge function}, in
its simplest format, is a multivariate function of the form $g\left( \mathbf{%
a}\cdot \mathbf{x}\right) $, where $g:\mathbb{R}\rightarrow \mathbb{R}$, $%
\mathbf{a}=\left( a_{1},...,a_{d}\right) $ is a fixed vector (direction) in $%
\mathbb{R}^{d}\backslash \left\{ \mathbf{0}\right\} $, $\mathbf{x}=\left(
x_{1},...,x_{d}\right) $ is the variable and $\mathbf{a}\cdot \mathbf{x}$ is
the standard inner product. In other words, a ridge function is a
multivariate function constant on the parallel hyperplanes $\mathbf{a}\cdot
\mathbf{x}=c$, $c\in \mathbb{R}$. These functions arise naturally in various
fields. They arise in computerized tomography (see, e.g., \cite%
{72,73,74,97,106,111}), statistics (see, e.g., \cite{13,14,27,33,42}) and
neural networks (see, e.g., \cite{22,Is,58,94,100,119,123}). These
functions are also used in modern approximation theory as an effective and
convenient tool for approximating complicated multivariate functions (see,
e.g., \cite{38,66,57,Kr,101,114,118,137}).

It should be remarked that long before the appearance of the name
\textquotedblleft ridge", these functions were used in PDE theory under the name of \textit{plane waves}. 
For example, see the book by F. John \cite{69}. In general, sums of ridge functions with
fixed directions occur in the study of hyperbolic constant coefficient partial
differential equations. As an example, assume that $(\alpha _{i},\beta
_{i}),~i=1,...,r,$ are pairwise linearly independent vectors in $\mathbb{R}%
^{2}$. Then the general solution to the homogeneous partial differential
equation
\begin{equation*}
\prod\limits_{i=1}^{r}\left( \alpha {_{i}{\frac{\partial }{\partial {x}}}%
+\beta _{i}{\frac{\partial }{\partial {y}}}}\right) {u}\left( {x,y}\right) =0
\end{equation*}%
are all functions of the form
\begin{equation*}
u(x,y)=\sum\limits_{i=1}^{r}g_{i}\left( \beta {_{i}x-\alpha _{i}y}\right)
\end{equation*}%
for arbitrary continuous univariate functions $g_{i}$, $i=1,...,r$. Here the
derivatives are understood in the sense of distributions.

The term \textquotedblleft ridge function" was coined by Logan and Shepp in
their seminal paper \cite{97} devoted to the basic mathematical problem of
computerized tomography. This problem consists of reconstructing a given
multivariate function from values of its integrals along certain straight
lines in the plane. The integrals along parallel lines can be considered as
a ridge function. Thus, the problem is to reconstruct $f$ from some set of
ridge functions generated by the function $f$ itself. In practice, one can
consider only a finite number of directions along which the above integrals
are taken. Obviously, reconstruction from such data needs some additional
conditions to be unique, since there are many functions $g$ having the same
integrals. For uniqueness, Logan and Shepp \cite{97} used the criterion of
minimizing the $L_{2}$ norm of $g$. That is, they found a function $g(x,y)$
with the minimum $L_{2}$ norm among all functions, which has the same
integrals as $f$. More precisely, let $D$ be the unit disk in the plane and
an unknown function $f(x,y)$ be square integrable and supported on $D.$ We
are given projections $P_{f}(t,\theta )$ (integrals of $f$ along the lines $%
x\cos \theta +y\sin \theta =t$) and looking for a function $g=g(x,y)$ of
minimum $L_{2}$ norm, which has the same projections as $f:$ $P_{g}(t,\theta
_{j})=P_{f}(t,\theta _{j}),$ $j=0,1,...,n-1$, where the angles $\theta _{j}$
generate equally spaced directions, i.e. $\theta _{j}=\frac{j\pi }{n},$ $%
j=0,1,...,n-1.$ The authors of \cite{97} showed that this problem of
tomography is equivalent to the problem of $L_{2}$-approximation of the
function $f$ by sums of ridge functions with the equally spaced directions $%
(\cos \theta _{j},\sin \theta _{j})$, $j=0,1,...,n-1.$ They gave a
closed-form expression for the unique function $g(x,y)$ and showed that the
unique polynomial $P(x,y)$ of degree $n-1$ which best approximates $f$ in $%
L_{2}(D)$ is determined from the above $n$ projections of $f$ and can be
represented as a sum of $n$ ridge functions.

Kazantsev \cite{72} solved the above problem of tomography without requiring
that the considered directions are equally spaced. Marr \cite{106}
considered the problem of finding a polynomial of degree $n-2$, whose
projections along lines joining each pair of $n$ equally spaced points on
the circumference of $D$ best matches the given projections of $f$ in the
sense of minimizing the sum of squares of the differences. Thus we see that
the problems of tomography give rise to an independent study of
approximation theoretic properties of the following set of linear
combinations of ridge functions:
\begin{equation*}
\mathcal{R}\left( \mathbf{a}^{1},...,\mathbf{a}^{r}\right) =\left\{
\sum\limits_{i=1}^{r}g_{i}\left( \mathbf{a}^{i}\cdot \mathbf{x}\right)
:g_{i}:\mathbb{R}\rightarrow \mathbb{R},i=1,...,r\right\} ,
\end{equation*}%
where directions $\mathbf{a}^{1},...,\mathbf{a}^{r}$ are fixed and belong to
the $d$-dimensional Euclidean space. Note that the set $\mathcal{R}\left(
\mathbf{a}^{1},...,\mathbf{a}^{r}\right) $ is a linear space.

Ridge function approximation also appears in statistics in \textit{%
Projection Pursuit}. This term was introduced by Friedman and Tukey \cite{32}
to name a technique for the explanatory analysis of large and multivariate
data sets. This technique seeks out ``interesting" linear projections of the
multivariate data onto a line or a plane. Projection Pursuit algorithms
approximate a multivariate function $f$ by sums of ridge functions with
variable directions, that is, by functions from the set
\begin{equation*}
\mathcal{R}_{r}=\left\{ \sum\limits_{i=1}^{r}g_{i}\left( \mathbf{a}^{i}\cdot
\mathbf{x}\right) :\mathbf{a}^{i}\in \mathbb{R}^{d}\setminus \{\mathbf{0}%
\},\ g_{i}:\mathbb{R}\rightarrow \mathbb{R},i=1,...,r\right\} .
\end{equation*}%
Here $r$ is the only fixed parameter, directions $\mathbf{a}^{1},...,\mathbf{%
a}^{r}$ and functions $g_{1},...,g_{r}$ are free to choose. The first method
of such approximation was developed by Friedman and Stuetzle \cite{33}.
Their approximation process called \textit{Projection Pursuit Regression}
(PPR) operates in a stepwise and greedy fashion. The process does not find a
best approximation from $\mathcal{R}_{r}$, it algorithmically constructs
functions $g_{r}\in \mathcal{R}_{r},$ such that $\left\Vert
g_{r}-f\right\Vert _{L_{2}}\rightarrow 0,$ as $r\rightarrow \infty $. At
stage $m$, PPR looks for a univariate function $g_{m}$ and direction $%
\mathbf{a}^{m}$ such that the ridge function $g_{m}\left( \mathbf{a}%
^{m}\cdot \mathbf{x}\right) $ best approximates the residual $%
f(x)-\sum\limits_{j=1}^{m-1}g_{j}\left( \mathbf{a}^{j}\cdot \mathbf{x}%
\right) $. Projection pursuit regression has been proposed as an approach to
bypass the curse of dimensionality and now is applied to prediction in
applied sciences. In \cite{13,14}, Candes developed a new approach based not
on stepwise construction of approximation but on a new transform called the
\textit{ridgelet transform}. The ridgelet transform represents general
functions as integrals of \textit{ridgelets} -- specifically chosen ridge
functions.

The significance of approximation by ridge functions is well understood from
its role in the theory of \textit{neural networks}. Ridge functions appear
in the definitions of many central neural network models. It is a broad
knowledge that neural networks are being successfully applied across an
extraordinary range of problem domains, in fields as diverse as finance,
medicine, engineering, geology and physics. Generally speaking, neural
networks are being introduced anywhere that there are problems of
prediction, classification or control. Thus not surprisingly, there is a
great interest to this powerful and very popular area of research (see,
e.g., \cite{119} and a great deal of references therein). An artificial
neural network is a way to perform computations using networks of
interconnected computational units vaguely analogous to neurons simulating
how our brain solves them. An artificial neuron, which forms the basis for
designing neural networks, is a device with $d$ real inputs and an output.
This output is generally a ridge function of the given inputs. In
mathematical terms, a neuron may be described as
\begin{equation*}
y=\sigma (\mathbf{w\cdot x}-\theta ),
\end{equation*}%
where $\mathbf{x=(x}_{1},...,x_{d})\in \mathbb{R}^{d}$ are the input
signals, $w=(w_{1},...,w_{d})\in \mathbb{R}^{d}$ are the synaptic weights, $%
\theta \in \mathbb{R}$ is the bias, $\sigma $ is the activation function and
$y$ is the output signal of the neuron. In a layered neural network the
neurons are organized in the form of layers. We have at least two layers: an
input and an output layer. The layers between the input and the output
layers (if any) are called hidden layers, whose computation nodes are
correspondingly called hidden neurons or hidden units. The output signals of
the first layer are used as inputs to the second layer, the output signals
of the second layer are used as inputs to the third layer, and so on for the
rest of the network. Neural networks with this kind of architecture is
called a \textit{Multilayer Feedforward Perceptron} (MLP). This is the most
popular model among other neural network models. In this model, a neural
network with a single hidden layer and one output represents a function of
the form
\begin{equation*}
\sum_{i=1}^{r}c_{i}\sigma (\mathbf{w}^{i}\mathbf{\cdot x}-\theta _{i}).
\end{equation*}
Here the weights $\mathbf{w}^{i}$ are vectors in $\mathbb{R}^{d}$, the
thresholds $\theta _{i}$ and the coefficients $c_{i}$ are real numbers and \
the activation function $\sigma $ is a univariate function. We fix only $%
\sigma $ and $r$. Note that the functions $\sigma (\mathbf{w}^{i}\mathbf{%
\cdot x}-\theta _{i})$ are ridge functions. Thus it is not surprising that
some approximation theoretic problems related to neural networks have strong
association with the corresponding problems of approximation by ridge
functions.

It is clear that in the special case, linear combinations of ridge functions
turn into sums of univariate functions. This is also the simplest case. The
simplicity of the approximation apparatus itself guarantees its utility in
applications where multivariate functions are constant obstacles. In
mathematics, this type of approximation has arisen, for example, in
connection with the classical functional equations \cite{11}, the numerical
solution of certain PDE boundary value problems \cite{9}, dimension theory
\cite{132,133}, etc. In computer science, it arises in connection with the
efficient storage of data in computer databases (see, e.g., \cite{140}).
There is an interesting interconnection between the theory of approximation
by univariate functions and problems of equilibrium construction in
economics (see \cite{136}).

Linear combinations of ridge functions with fixed directions allow a natural
generalization to functions of the form $g(\alpha _{1}(x_{1})+\cdot \cdot
\cdot +\alpha _{d}(x_{d}))$, where $\alpha _{i}(x_{i})$, $i=\overline{1,d},$
are real univariate functions. Such a generalization has a strong
association with linear superpositions. A linear superposition is a function
expressed as the sum%
\begin{equation*}
\sum\limits_{i=1}^{r}g_{i}(h_{i}(x)), \; x \in X,
\end{equation*}%
where $X$ is any set (in particular, a subset of $\mathbb{R}^{d}$), $%
h_{i}:X\rightarrow {{\mathbb{R}}},~i=1,...,r,$ are arbitrarily fixed
functions, and $g_{i}:\mathbb{R}\rightarrow \mathbb{R},~i=1,...,r.$ Note
that here we deal with more complicated composition than the composition of
a univariate function with the inner product. A starting point in the study
of linear superpositions was the well known superposition theorem of
Kolmogorov \cite{83} (see also the paper on Kolmogorov's works by Tikhomirov
\cite{139}). This theorem states that for the unit cube $\mathbb{I}^{d},~%
\mathbb{I}=[0,1],~d\geq 2,$ there exist $2d+1$ functions $%
\{s_{q}\}_{q=1}^{2d+1}\subset C(\mathbb{I}^{d})$ of the form
\begin{equation*}
s_{q}(x_{1},...,x_{d})=\sum_{p=1}^{d}\varphi _{pq}(x_{p}),~\varphi _{pq}\in
C(\mathbb{I}),~p=1,...,d,~q=1,...,2d+1
\end{equation*}%
such that each function $f\in C(\mathbb{I}^{d})$ admits the representation
\begin{equation*}
f(x)=\sum_{q=1}^{2d+1}g_{q}(s_{q}(x)),~x=(x_{1},...,x_{d})\in \mathbb{I}%
^{d},~g_{q}\in C({{\mathbb{R)}}}.
\end{equation*}%
Thus, any continuous function on the unit cube can be represented as a
linear superposition with the fixed inner functions $s_{1},...,s_{2d+1}$. In
literature, these functions are called universal functions or the Kolmogorov
functions. Note that all the functions $g_{q}(s_{q}(x))$ in the Kolmogorov
superposition formula are generalized ridge functions, since each $s_{q}$ is a
sum of univariate functions.

In these notes, we consider some problems of approximation and/or representation of multivariate
functions by linear combinations of ridge functions,
generalized ridge functions and feedforward neural networks. The notes consist
of five chapters.

Chapter 1 is devoted to the approximation from some sets of ridge functions
with arbitrarily fixed directions in $C$ and $L_{2}$ metrics. First, we
study problems of representation of multivariate functions by linear
combinations of ridge functions. Then, in case of two fixed directions and
under suitable conditions, we give complete solutions to three basic
problems of uniform approximation, namely, problems on existence,
characterization, and construction of a best approximation. We also study
problems of well approximation (approximation with arbitrary accuracy)
and representation of continuous multivariate functions by sums
of two continuous ridge functions. The reader will see the main difficulties
and remained open problems in the uniform approximation by sums of more than
two ridge functions. For $L_{2}$ approximation, a number of summands does
not play such an essential role as it plays in the uniform approximation. In
this case, it is known that a best approximation always exists and unique.
For some special domains in $\mathbb{R}^{d}$, we characterize and then
construct the best approximation. We also give an explicit formula for the
approximation error.

Chapter 2 explores the following open problem raised in Buhmann and
Pinkus \cite{12}, and Pinkus \cite[p. 14]{117}. Assume we are given a
function $f(\mathbf{x})=f(x_{1},...,x_{n})$ of the form
\begin{equation*}
f(\mathbf{x})=\sum_{i=1}^{k}f_{i}(\mathbf{a}^{i}\cdot \mathbf{x}),
\end{equation*}%
where the $\mathbf{a}^{i},$ $i=1,...,k,$ are pairwise linearly independent
vectors (directions) in $\mathbb{R}^{d}$, $f_{i}$ are arbitrarily behaved
univariate functions and $\mathbf{a}^{i}\cdot \mathbf{x}$ are standard inner
products. Assume, in addition, that $f$ is of a certain smoothness class,
that is, $f\in C^{s}(\mathbb{R}^{d})$, where $s\geq 0$ (with the convention
that $C^{0}(\mathbb{R}^{d})=C(\mathbb{R}^{d})$). Is it true that there will
always exist $g_{i}\in C^{s}(\mathbb{R})$ such that
\begin{equation*}
f(\mathbf{x})=\sum_{i=1}^{k}g_{i}(\mathbf{a}^{i}\cdot \mathbf{x})\text{ ?}%
\end{equation*}
In this chapter, we solve this problem up to some multivariate polynomial.
We find various conditions on the directions $%
\mathbf{a}^{i}$ allowing to express this polynomial as a sum of smooth ridge functions with these directions. We also
consider the question of constructing $g_{i}$ using the information about the known functions $f_{i}$.

Chapter 3 is devoted to the simplest type of ridge functions -- univariate
functions. Note that a ridge function depends only on one variable if its
direction coincides with the coordinate direction. Thus, in case of
coincidence of all given directions with the coordinate directions, the
problem of ridge function approximation turns into the problem of
approximation of multivariate functions by sums of univariate functions. In
this chapter, we first consider the approximation of a bivariate function $%
f(x,y)$ by sums $\varphi (x)+\psi (y)$ on a rectangular domain $R$. We
construct special classes of continuous functions depending on a numerical
parameter and characterize each class in terms of the approximation error
calculation formulas. This parameter will show which points of $R$ the
calculation formula involves. We will also construct a best approximating
sum $\varphi _{0}(x)+\psi _{0}(y)$ to a function from constructed classes.
Then we develop a method for obtaining explicit formulas for the error of
approximation of bivariate functions, defined on a union of rectangles, by
sums of univariate functions. It should be remarked that formulas of such
type were known only for functions defined on a rectangle with sides
parallel to the coordinate axes. Our method, based on a maximization process
over certain objects, called \textquotedblleft closed bolts", allows the
consideration of functions defined on hexagons, octagons and stairlike
polygons with sides parallel to the coordinate axes. At the end of this
chapter we discuss one important result from Golomb's paper \cite{37}. This
paper, published in 1959, made a start of a systematic study of
approximation of multivariate functions by various compositions, including
sums of univariate functions. In \cite{37}, along with many other results,
Golomb obtained a duality formula for the error of approximation to a
multivariate function from the set of sums of univariate functions.
Unfortunately, his proof had a gap, which was 24 years later pointed out by
Marshall and O'Farrell \cite{107}. But the question if Golomb's formula was
correct, remained unsolved. In Chapter 3, we show that Golomb's formula is
correct, and moreover it holds in a stronger form.

Chapter 4 tells us about some problems concerning generalized ridge functions $%
g(\alpha _{1}(x_{1})+\cdot \cdot \cdot +\alpha _{d}(x_{d}))$ and linear
superpositions. We consider the problem of representation of general
functions by linear superpositions. We show that if some representation by
linear superpositions, in particular by linear combinations of generalized ridge
functions, holds for continuous functions, then it holds for all functions.
This leads us to extensions of many superpositions theorems (such as the
well-known Kolmogorov superposition theorem, Ostrand's superposition
theorem, etc.) from continuous to arbitrarily behaved multivariate
functions. Concerning generalized ridge functions, we see that every
multivariate function can be written as a generalized ridge function or as a sum
of finitely many such functions. We also study the uniqueness of representation of
functions by linear superpositions.

Chapter 5 is about neural network approximation. The analysis in this chapter is
based on properties of ordinary and generalized ridge functions.
We consider a single and two hidden layer feedforward neural network models with a restricted set of
weights. Such network models are important from the point of view of
practical applications. We study approximation properties of single hidden
layer neural networks with weights varying on a finite set of directions and
straight lines. We give several necessary and sufficient conditions for well
approximation by such networks. For a set of weights consisting of two
directions (and two straight lines), we show that there is a geometrically
explicit solution to the problem. Regarding
two hidden layer feedforward neural networks, we prove that two hidden layer
neural networks with $d$ inputs, $d$ neurons in the first hidden layer, $2d+2
$ neurons in the second hidden layer and with a specifically constructed
sigmoidal, infinitely differentiable and almost monotone activation function can approximate
any continuous multivariate function with arbitrary precision. We show that
for this approximation only a finite number of fixed weights (precisely, $d$
fixed weights) suffice.

There are topics related to ridge functions that are not presented here. The
glaring omission is that of interpolation at points and on straight lines by
ridge functions. We also do not address, for example, questions of linear
independence and spanning by linear combinations of ridge monomials in the
spaces of homogeneous and algebraic polynomials of a fixed degree, integral
representations of functions where the kernel is a ridge function,
approximation algorithms for finding best approximations from spaces of
linear combinations of ridge functions. These and similar topics may be
found in the monograph by Pinkus \cite{117}. The reader may also consult
the survey articles \cite{57,87,118}.

\newpage

\chapter{Properties of linear combinations of ridge functions}

In this chapter, we consider approximation-theoretic problems arising in
ridge function approximation. First we briefly review some results on approximation
by sums of ridge functions with both fixed and variable directions. Then we
analyze the problem of representability of an arbitrary multivariate
function by linear combinations of ridge functions with fixed directions. In
the special case of two fixed directions, we characterize a best uniform
approximation from the set of sums of ridge functions with these directions.
For a class of bivariate functions we use this result to construct
explicitly a best approximation. Questions on existence of a best
approximation are also studied. We also study problems of well approximation
(approximation with arbitrary accuracy) and representation of continuous
multivariate functions by sums of two continuous ridge functions. The reader
will see the main difficulties and remained open problems in the uniform
approximation by sums of more than two ridge functions. For $L_{2}$
approximation, a number of summands does not play such an essential role as
it plays in the uniform approximation. In this case, it is known that a best
approximation always exists and unique. For some special domains in $%
\mathbb{R}^{d}$, we characterize and then construct the best approximation.
We also give an explicit formula for the approximation error.

\bigskip

\section{A brief excursion into the approximation theory of ridge functions}

In this section we briefly review some results on approximation properties of the sets $%
\mathcal{R}\left( \mathbf{a}^{1},...,\mathbf{a}^{r}\right) $ and $\mathcal{R}%
_{r}$. These results are presented without proofs but with discussions and complete references.
We hope this section will whet the reader's appetite for the rest of these notes,
where a more comprehensive study of concrete mathematical problems is provided.

\subsection{$\mathcal{R}\left( \mathbf{a}^{1},...,\mathbf{a}^{r}\right) $ --
ridge functions with fixed directions}

It is clear that well approximation of a multivariate function $f:X\rightarrow \mathbb{R}$ from
some normed space by using elements of the set $\mathcal{R}\left( \mathbf{a}%
^{1},...,\mathbf{a}^{r}\right)$ is not always possible. The value of the
approximation error depends not only on the approximated function $f$ but
also on geometrical structure of the given set $X$. This poses challenging research problems on
computing the error of approximation and constructing best approximations from
$\mathcal{R}\left( \mathbf{a}^{1},...,\mathbf{a}^{r}\right)$. Serious difficulties arise
when one attempts to solve these problems in continuous function spaces endowed with
the uniform norm. For example, let us consider the algorithm for finding
best approximations, called the \textit{Diliberto-Straus algorithm} (see \cite{18}). The essence of
this algorithm is as follows. Let $X$ be a compact subset of $\mathbb{R}^{d}$
and $A_{i}$ be a best approximation operator from the space of continuous
functions $C(X)$ to the subspace of ridge functions $G_{i}=\{g_{i}\left(
\mathbf{a}^{i}\cdot \mathbf{x}\right) :~g_{i}\in C(\mathbb{R)},~\mathbf{x}%
\in X\}$, $i=1,...,r.$ That is, for each function $f$ $\in C(X)$, the
function $A_{i}f$ is a best approximation to $f$ from $G_{i}.$ Set
\begin{equation*}
Tf=(I-A_{r})(I-A_{r-1})\cdot \cdot \cdot (I-A_{1})f,
\end{equation*}%
where $I$ is the identity operator. It is clear that
\begin{equation*}
Tf=f-g_{1}-g_{2}-\cdot \cdot \cdot -g_{r},
\end{equation*}%
where $g_{k}$ is a best approximation from $G_{k}$ to the function $%
f-g_{1}-g_{2}-\cdot \cdot \cdot -g_{k-1}$, $k=1,...,r.$ Consider powers of
the operator $T$: $T^{2},T^{3}$ and so on. Is the sequence $%
\{T^{n}f\}_{n=1}^{\infty }$ convergent? In case of an affirmative answer,
which function is the limit of $T^{n}f,$ as $n\rightarrow \infty$? One may
expect that the sequence $\{T^{n}f\}_{n=1}^{\infty }$ converges to $%
f-g^{\ast },$ where $g^{\ast }$ is a best approximation from $\mathcal{R}%
\left( \mathbf{a}^{1},...,\mathbf{a}^{r}\right) $ to $f$. This conjecture
was first stated by Diliberto and Straus \cite{26} in 1951 for the
uniform approximation of a multivariate function, defined on the unit cube,
by sums of univariate functions (that is, sums of ridge functions with the
coordinate directions). But later it was shown by Aumann \cite{4} that the sequence generated by this
algorithm may not converge if $r>2$. For $r=2$ and certain convex
compact sets $X$, the sequence $\{\|T^{n}f\|\}_{n=1}^{\infty }$
converges to the approximation error $\|f-g_{0}\|$,
where $g_{0}$ is a best approximation from $\mathcal{R}%
\left( \mathbf{a}^{1},\mathbf{a}^{2}\right)$ (see \cite{61,117}).
However, it is not yet clear whether $\|T^{n}f-(f-g_{0})\|$ converges to zero
as $n\rightarrow \infty$. In the case $r>2$ no efficient algorithm is known for finding a best uniform
approximation from $\mathcal{R}\left( \mathbf{a}^{1},...,\mathbf{a}%
^{r}\right)$. Note that in the $L_{2}$ metric, the Diliberto-Straus algorithm
converges as desired for an arbitrary number of distinct
directions. This also holds in the $L_{p}$ space setting, provided that $p>1$ and
$\mathcal{R}\left( \mathbf{a}^{1},...,\mathbf{a}^{r}\right)$ is closed (see \cite{118}). But in the $L_{1}$ space setting,
the alternating algorithm does not work even in the case of two directions (see \cite{117}).

One of the basic problems concerning the approximation by sums of ridge
functions with fixed directions is the problem of verifying if a given
function $f$ belongs to the space $\mathcal{R}\left( \mathbf{a}^{1},...,%
\mathbf{a}^{r}\right) $. This problem has a simple solution if the space
dimension $d=2$ and a given function $f(x,y)$ has partial derivatives up to $%
r$-th order. For the representation of $f(x,y)$ in the form%
\begin{equation*}
f(x,y)=\sum_{i=1}^{r}g_{i}(a_{i}x+b_{i}y),
\end{equation*}%
it is necessary and sufficient that

\begin{equation*}
\prod\limits_{i=1}^{r}\left( b_{i}\frac{\partial }{\partial x}-a_{i}\frac{%
\partial }{\partial y}\right) f=0.\eqno(1.1)
\end{equation*}%
This recipe is also valid for continuous bivariate functions provided that
the derivatives are understood in the sense of distributions.

Unfortunately such a simple characterization does not carry over to the
case of more than two variables. Below we provide two results
concerning the general case of arbitrarily many variables.

\bigskip

\textbf{Proposition 1.1 }(Diaconis, Shahshahani \cite{25}). \textit{Let $%
\mathbf{a}^{1},...,\mathbf{a}^{r}$ be pairwise linearly independent vectors
in $\mathbb{R}^{d}. $ Let for $i=1,2,...,r$, $H^{i}$ denote the hyperplane $%
\{\mathbf{c}\in \mathbb{R}^{d}$: $\mathbf{c\cdot a}^{i}=0\}.$ Then a
function $f\in C^{r}(\mathbb{R}^{d})$ can be represented in the form}
\begin{equation*}
f(\mathbf{x})=\sum\limits_{i=1}^{r}g_{i}\left( \mathbf{a}^{i}\cdot \mathbf{x}%
\right) +P(\mathbf{x}),
\end{equation*}%
\textit{where $P(\mathbf{x})$ is a polynomial of degree not more than $r$,
if and only if}
\begin{equation*}
\prod\limits_{i=1}^{r}\sum_{s=1}^{d}c_{s}^{i}\frac{\partial f}{\partial x_{s}%
}=0,
\end{equation*}%
\textit{for all vectors $\mathbf{c}^{i}=(c_{1}^{i},c_{2}^{i},...,c_{d}^{i})%
\in H^{i}, $ $i=1,2,...,r.$}

\bigskip

There are examples showing that one cannot simply dispense with the
polynomial $P(\mathbf{x})$ in the above proposition (see \cite{25}). In
fact, a polynomial term appears in the sufficiency part of the proof of this
proposition.

Lin and Pinkus \cite{95} obtained more general result on the representation
by sums of ridge functions with fixed directions.
We need some notation to present their result. Each
polynomial $p(x_{1},...,x_{d})$ generates the differential operator $p(\frac{%
\partial }{\partial x_{1}},...,\frac{\partial }{\partial x_{d}}).$ Let $P(%
\mathbf{a}^{1},...,\mathbf{a}^{r})$ denote the set of polynomials which
vanish on all the lines $\{\lambda \mathbf{a}^{i},\lambda \in \mathbb{R}\},$
$i=1,...,r.$ Obviously, this is an ideal in the ring of all polynomials. Let
$Q$ be the set of polynomials $q=q(x_{1},...,x_{d})$ such that $p(\frac{%
\partial }{\partial x_{1}},...,\frac{\partial }{\partial x_{d}})q=0$, for
all $p(x_{1},...,x_{d})\in P(\mathbf{a}^{1},...,\mathbf{a}^{r}).$

\bigskip

\textbf{Proposition 1.2 }(Lin, Pinkus \cite{95}). \textit{Let $\mathbf{a}%
^{1},...,\mathbf{a}^{r}$ be pairwise linearly independent vectors in $%
\mathbb{R}^{d}.$ A function $f\in C(\mathbb{R}^{d})$ can be expressed in the
form}
\begin{equation*}
f(\mathbf{x})=\sum\limits_{i=1}^{r}g_{i}(\mathbf{a}^{i}\cdot \mathbf{x)},
\end{equation*}%
\textit{if and only if $f$ belongs to the closure of the linear span of $Q.$}

\bigskip

In \cite{120}, A.Pinkus considered the problems of smoothness and uniqueness
in ridge function representation. For a given function $f$ $\in $ $\mathcal{R%
}\left( \mathbf{a}^{1},...,\mathbf{a}^{r}\right) $, he posed and answered the
following questions. If $f$ belongs to some smoothness class, what can we
say about the smoothness of the functions $g_{i}$? How many different ways
can we write $f$ as a linear combination of ridge functions? These and similar problems
will be extensively discussed in Chapter 2.

The above problem of representation of fixed functions by sums of
ridge functions gives rise to the problem of representation of some classes
of functions by such sums. For example, one may consider the following
problem. Let $X$ be a subset of the $d$-dimensional Euclidean space. Let $%
C(X),$ $B(X),$ $T(X)$ denote the set of continuous, bounded and all real
functions defined on $X$, respectively. In the first case, we additionally
suppose that $X$ is a compact set. Let $\mathcal{R}_{c}\left( \mathbf{a}%
^{1},...,\mathbf{a}^{r}\right) $ and $\mathcal{R}_{b}\left( \mathbf{a}%
^{1},...,\mathbf{a}^{r}\right) $ denote the subspaces of $\mathcal{R}\left(
\mathbf{a}^{1},...,\mathbf{a}^{r}\right) $ comprising only sums with
continuous and bounded terms $g_{i}\left( \mathbf{a}^{i}\cdot \mathbf{x}%
\right) $, $i=1,...,r$, respectively. The following questions naturally
arise: For which sets $X$,

$(1)$ $C(X)=\mathcal{R}_{c}\left(
\mathbf{a}^{1},...,\mathbf{a}^{r}\right)$?

$(2)$ $B(X)=\mathcal{R}_{b}\left(
\mathbf{a}^{1},...,\mathbf{a}^{r}\right)$?

$(3)$ $T(X)=\mathcal{R}\left(
\mathbf{a}^{1},...,\mathbf{a}^{r}\right)$?

The first two questions in a more
general setting were answered in Sternfeld \cite{132,131}.
The third question will be answered in the next section. Let us briefly
discuss some results of Sternfeld concerning ridge function representation. These results
have been mostly overlooked in the corresponding ridge function literature,
as they have to do with more general superpositions of functions and do not directly mention ridge functions.
Assume we are given directions $%
\mathbf{a}^{1},...,\mathbf{a}^{r}\in \mathbb{R}^{d}\backslash \{\mathbf{0}\}$
and a set $X\subseteq \mathbb{R}^{d}.$ Following Sternfeld, we say that a family $F=\{\mathbf{a}^{1},...,%
\mathbf{a}^{r}\}$ \textit{uniformly separates points} of $X$ if there exists
a number $0<\lambda \leq 1$ such that for each pair $\{\mathbf{x}%
_{j}\}_{j=1}^{m}$, $\{\mathbf{z}_{j}\}_{j=1}^{m}$ of disjoint finite
sequences in $X$, there exists some direction $\mathbf{a}^{k}\in F$ so that
if from the two sequences $\{\mathbf{a}^{k}\cdot \mathbf{x}_{j}\}_{j=1}^{m}$
and $\{\mathbf{a}^{k}\cdot \mathbf{z}_{j}\}_{j=1}^{m}$ we remove a maximal
number of pairs of points $\mathbf{a}^{k}\cdot \mathbf{x}_{j_{1}}$ and $%
\mathbf{a}^{k}\cdot \mathbf{z}_{j_{2}}$ with $\mathbf{a}^{k}\cdot \mathbf{x}%
_{j_{1}}=\mathbf{a}^{k}\cdot \mathbf{z}_{j_{2}},$ then there remains at
least $\lambda m$ points in each sequence (or, equivalently, at most $%
(1-\lambda )m$ pairs can be removed). Sternfeld \cite{132}, in particular,
proved that a family of directions $F=\{\mathbf{a}^{1},...,\mathbf{a}^{r}\}$
uniformly separates points of $X$ if and only if $\mathcal{R}_{b}\left(
\mathbf{a}^{1},...,\mathbf{a}^{r}\right) =B(X)$. In \cite{132}, he also
obtained a practically convenient sufficient condition for the equality $%
\mathcal{R}_{b}\left( \mathbf{a}^{1},...,\mathbf{a}^{r}\right) =B(X).$ To
describe his condition, define the set functions
\begin{equation*}
\tau _{i}(Z)=\{\mathbf{x}\in Z:~|p_{i}^{-1}(p_{i}(\mathbf{x}))\bigcap Z|\geq
2\},
\end{equation*}%
where $Z\subset X,~p_{i}(\mathbf{x})=\mathbf{a}^{i}\cdot \mathbf{x}$, \ $%
i=1,\ldots ,r,$ and $|Y|$ denotes the cardinality of a set $Y$.
Define $\tau (Z)$ to be $\bigcap_{i=1}^{k}\tau _{i}(Z)$ and define $\tau
^{2}(Z)=\tau (\tau (Z))$, $\tau ^{3}(Z)=\tau (\tau ^{2}(Z))$ and so on
inductively.

\bigskip

\textbf{Proposition 1.3 }(Sternfeld \cite{132}). \textit{If $\tau
^{n}(X)=\emptyset $ for some $n$, then $\mathcal{R}_{b}\left( \mathbf{a}%
^{1},...,\mathbf{a}^{r}\right) =B(X)$. If $X$ is a compact subset of $%
\mathbb{R}^{d}$, and $\tau ^{n}(X)=\emptyset $ for some $n$, then $\mathcal{R%
}_{c}\left( \mathbf{a}^{1},...,\mathbf{a}^{r}\right) =C(X)$.}

\bigskip

If $r=2$, the sufficient condition \textquotedblleft $\tau ^{n}(X)=\emptyset $ for
some $n$" turns out to be also necessary. In this case,
the equality $\mathcal{R}_{b}\left( \mathbf{a}^{1},\mathbf{a}^{2}\right)
=B(X)$ is equivalent to the equality $\mathcal{R}_{c}\left( \mathbf{a}^{1},%
\mathbf{a}^{2}\right) =C(X)$. In another work \cite{131}, \ Sternfeld
obtained a measure-theoretic necessary and sufficient condition for the
equality $\mathcal{R}_{c}\left( \mathbf{a}^{1},...,\mathbf{a}^{r}\right)
=C(X)$. Let $p_{i}(\mathbf{x})=\mathbf{a}^{i}\cdot \mathbf{x}$, \ $%
i=1,\ldots ,r$, $X$ be a compact set in $\mathbb{R}^{d}$ and $M(X)$ be a
class of measures defined on some field of subsets of $X$.
Following Sternfeld, we say that a family $F=\{%
\mathbf{a}^{1},...,\mathbf{a}^{r}\}$ \textit{uniformly separates measures}
of the class $M(X)$ if there exists a number $0<\lambda \leq 1$ such that
for each measure $\mu $ in $M(X)$ the equality $\left\Vert \mu \circ
p_{k}^{-1}\right\Vert \geq \lambda \left\Vert \mu \right\Vert $ holds for
some direction $\mathbf{a}^{k}\in F$. Sternfeld \cite{134,131}, in
particular, proved that the equality $\mathcal{R}_{c}\left( \mathbf{a}%
^{1},...,\mathbf{a}^{r}\right) =C(X)$ holds if and only if the family of
directions $\{\mathbf{a}^{1},...,\mathbf{a}^{r}\}$ uniformly separates
measures of the class $C(X)^{\ast }$ (that is, the class of regular Borel
measures). In addition, he proved that $\mathcal{R}_{b}\left( \mathbf{a}^{1},...,%
\mathbf{a}^{r}\right) =B(X)$ if and only if the family of directions $\{%
\mathbf{a}^{1},...,\mathbf{a}^{r}\}$ uniformly separates measures of the
class $l_{1}(X)$ (that is, the class of finite measures defined on countable
subsets of $X$). Since $l_{1}(X)\subset C(X)^{\ast },$ the first equality $%
\mathcal{R}_{c}\left( \mathbf{a}^{1},...,\mathbf{a}^{r}\right) =C(X)$
implies the second equality $\mathcal{R}_{b}\left( \mathbf{a}^{1},...,%
\mathbf{a}^{r}\right) =B(X).$ The inverse is not true (see \cite{131}).
We emphasize again that the above results of Sternfeld were obtained for
more general functions, than linear combinations of ridge functions, namely
for functions of the form $\sum_{i=1}^{r}g_{i}(h_{i}(x))$, where $h_{i}$
arbitrarily fixed functions (bounded or continuous) defined on $X.$ Such functions will
be discussed in Chapter 4.

\bigskip

\subsection{$\mathcal{R}_{r}$ -- ridge functions with variable directions}

Obviously, the\ set $\mathcal{R}_{c}\left( \mathbf{a}^{1},...,\mathbf{a}%
^{r}\right) $ is not dense in $C(\mathbb{R}^{d})$ in the topology of uniform
convergence on compact subsets of $\mathbb{R}^{d}.$ Density here does not
hold because the number of considered directions is finite. If consider all
the possible directions, then the set $\mathcal{R}=span\{g(\mathbf{a}\cdot
\mathbf{x)}:~g\in C(\mathbb{R)},~\mathbf{a}\in \mathbb{R}^{d}\backslash \{%
\mathbf{0}\}\}$ will certainly be dense in the space $C(\mathbb{R}^{d})$ in
the above mentioned topology. In order to be sure, it is enough to consider
only the functions $e^{\mathbf{a}\cdot \mathbf{x}}\in \mathcal{R}$, the
linear span of which is dense in $C(\mathbb{R}^{d})$ by the
Stone-Weierstrass theorem. In fact, for density it is not necessary to
comprise all directions. The following theorem shows how many directions in
totality satisfy the density requirements.

\bigskip

\textbf{Proposition 1.4 }(Vostrecov and Kreines \cite{142}, Lin and Pinkus
\cite{95}). \textit{For density of the set}
\begin{equation*}
\mathcal{R(A)}=span\{g(\mathbf{a}\cdot \mathbf{x)}:~g\in C(\mathbb{R)},~%
\mathbf{a}\in \mathcal{A}\subset \mathbb{R}^{d}\}
\end{equation*}%
\textit{in $C(\mathbb{R}^{d})$ (in the topology of uniform convergence on
compact sets) it is necessary and sufficient that the only homogeneous
polynomial which vanishes identically on $\mathcal{A}$ is the zero
polynomial.}

\bigskip

Since in the definition of $\mathcal{R(A)}$ we vary over all univariate
functions$~g,$ allowing one direction $\mathbf{a}$ is equivalent to allowing
all directions $k\mathbf{a}$ for every real $k$. Thus it is sufficient to
consider only the set $\mathcal{A}$ of directions normalized to the unit
sphere $S^{n-1}.$ For example, if $\mathcal{A}$ is a subset of the sphere $%
S^{n-1},$ which contains an interior point (interior point with respect to
the induced topology on $S^{n-1}$), then $\mathcal{R(A)}$ is dense in the
space $C(\mathbb{R}^{d}).$ The proof of Proposition 1.4 highlights an
important fact that the set $\mathcal{R(A)}$ is dense in $C(\mathbb{R}^{d})$
in the topology of uniform convergence on compact subsets if and only if $%
\mathcal{R(A)}$ contains all the polynomials (see \cite{95}).

Representability of polynomials by sums of ridge functions is a building
block for many results. In many works (see, e.g., \cite{119}), the following
fact is fundamental: Every multivariate polynomial $h(\mathbf{x}%
)=h(x_{1},...,x_{d})$ of degree $k $ can be represented in the form
\begin{equation*}
h(\mathbf{x})=\sum\limits_{i=1}^{l}p_{i}(\mathbf{a}^{i}\cdot \mathbf{x),}
\end{equation*}%
where $p_{i}$ is a univariate polynomial, $\mathbf{a}^{i}\in \mathbb{R}^{d}$%
, and $l=$ $\binom{d-1+k}{k}$.

For example, for the representation of a bivariate polynomial of degree $k$,
it is needed $k+1$ univariate polynomials and $k+1$ directions (see \cite{97}%
). The proof of this fact is organized so that the directions $\mathbf{a}^{i}
$, $i=1,...,k+1$, are chosen once for all multivariate polynomials of $k$-th
degree. At one of the seminars in the Technion -- Israel Institute of
Technology in 2007, A. Pinkus posed two problems:

1) Can every multivariate polynomial of degree $k$ be represented by less
than $l$ ridge functions?

2) How large is the set of polynomials represented by $l-1,$ $l-2,...$ ridge
functions?

Note that for bivariate polynomials the 1-st problem is solved positively,
that is, the number $l=k+1$ can be reduced. Indeed, for a bivariate
polynomial $P(x,y)$ of $k$-th degree, there exist many combinations of real
numbers $c_{0},...,c_{k}$ such that
\begin{equation*}
\sum_{i=0}^{k}c_{i}\frac{\partial ^{k}}{\partial x^{i}\partial y^{k-i}}%
P(x,y)=0.
\end{equation*}%
Further the numbers $c_{i}$, $i=0,...,k$, can be selected to enjoy the
property that the polynomial $\sum_{i=0}^{k}c_{i}t^{i}$ has distinct real
zeros. Then it is not difficult to verify that the differential operator $%
\sum_{i=0}^{k}c_{i}\frac{\partial ^{k}}{\partial x^{i}\partial y^{k-i}}$ can
be written in the form
\begin{equation*}
\prod\limits_{i=1}^{k}\left( b_{i}\frac{\partial }{\partial x}-a_{i}\frac{%
\partial }{\partial y}\right) ,
\end{equation*}%
for some pairwise linearly independent vectors $(a_{i},b_{i})$, $i=1,...,k$.
Now from the above criterion (1.1) we obtain that the polynomial $P(x,y)$
can be represented as a sum of $k$ ridge functions. Note that the problem of
representation of a multivariate algebraic polynomial $P(\mathbf{x})$ in the
form $\sum_{i=0}^{r}g_{i}(\mathbf{a}^{i}\cdot \mathbf{x})$ with minimal $r$
was extensively studied in the monograph by Pinkus \cite{117}.

In connection with the 2-nd problem of Pinkus, V. Maiorov \cite{103} studied
certain geometrical properties of the manifold $\mathcal{R}_{r}$. Namely, he
estimated the $\varepsilon $-entropy numbers in terms of smaller $%
\varepsilon $-covering numbers of the compact class formed by the
intersection of the class $\mathcal{R}_{r}$ with the unit ball in the space
of polynomials of degree at most $s$ on $\mathbb{R}^{d}$. Let $E$ be a
Banach space and let for $x\in E$ and $\delta >0,$ $S(x,\delta )$ denote the
ball of radius $\delta $ centered at the point $x$. For any positive number $%
\varepsilon $, the $\varepsilon $-covering number of a set $F$ in the space $%
E$ represents the quantity
\begin{equation*}
L_{\varepsilon }(F,E)=\min \left\{ N:~\exists x_{1},...,x_{N}\in F\text{
such that }F\subset \bigcup_{i=1}^{N}S(x_{i},\varepsilon )\text{ }\right\}.
\end{equation*}

The $\varepsilon $-entropy of $F$ is defined as the number $H_{\varepsilon
}(F,E)\overset{def}{=}$ $\log _{2}L_{\varepsilon }(F,E)$. The notion of $%
\varepsilon $-entropy has been devised by A.N.Kolmogorov (see \cite{84,85})
to classify compact metric sets according to their massivity.

In order to formulate Maiorov's result, let $\mathcal{P}_{s}^{d}$ be the
space of all polynomials of degree at most $s$ on $\mathbb{R}^{d}$, $%
L_{q}=L_{q}(I)$, $1\leq q\leq \infty $, be the space of $q$-integrable
functions on the unit cube $I=[0,1]^{d}$ with the norm $\left\Vert
f\right\Vert _{q}=\left( \int_{I}\left\vert f(x)\right\vert ^{q}dx\right)
^{1/q}$, $BL_{q}$ be the unit ball in the space $L_{q},$ and $B_{q}\mathcal{P%
}_{s}^{d}=$ $BL_{q}\cap $ $\mathcal{P}_{s}^{d}$ be the unit ball in the
space $\mathcal{P}_{s}^{d}$ equipped with the $L_{q}$ metric.

\bigskip

\textbf{Proposition 1.5 }(Maiorov \cite{103}). \textit{Let $r,s\in \mathbb{N}
$, $1\leq q\leq \infty $, $0 < \varepsilon < 1$. The $\varepsilon $-entropy
of the class $B_{q}\mathcal{P}_{s}^{d}\cap \mathcal{R}_{r}$ in the space $%
L_{q}$ satisfies the inequalities}

1)
\begin{equation*}
c_{1}rs\leq \frac{H_{\varepsilon }(B_{q}\mathcal{P}_{s}^{d}\cap \mathcal{R}%
_{r},L_{q})}{\log _{2}\frac{1}{\varepsilon }}\leq c_{2}rs\log _{2}\frac{%
2es^{d-1}}{r},
\end{equation*}
\textit{for $r\leq s^{d-1}.$}

2)
\begin{equation*}
c_{1}^{^{\prime }}s^{d}\leq \frac{H_{\varepsilon }(B_{q}\mathcal{P}%
_{s}^{d}\cap \mathcal{R}_{r},L_{q})}{\log _{2}\frac{1}{\varepsilon }}\leq
c_{2}^{^{\prime }}s^{d},
\end{equation*}
\textit{for $r>s^{d-1}.$ In these inequalities $c_{1},c_{2},c_{1}^{^{\prime
}},c_{2}^{^{\prime }}$ are constants depending only on $d$.}

\bigskip

Let us consider $\mathcal{R}_{r}$ as a subspace of some normed linear space $%
X$ endowed with the norm $\left\Vert \cdot \right\Vert _{X}.$ The error of
approximation of a given function $f\in X$ by functions $g\in \mathcal{R}%
_{r} $ is defined as follows
\begin{equation*}
E(f,\mathcal{R}_{r},X)\overset{def}{=}\underset{g\in \mathcal{R}_{r}}{\inf }%
\left\Vert f-g\right\Vert _{X}.
\end{equation*}

Let $B^{d}$ denote the unit ball in the space $\mathbb{R}^{d}.$ Besides, let
$\mathbb{Z}_{+}^{d}$ denote the lattice of nonnegative multi-integers in $%
\mathbb{R}^{d}.$ For $k=(k_{1},...,k_{d})\in \mathbb{Z}_{+}^{d},$ set $%
\left\vert k\right\vert =k_{1}+\cdot \cdot \cdot +k_{d}$, $\mathbf{x}^{%
\mathbf{k}}=x_{1}^{k_{1}}\cdot \cdot \cdot x_{d}^{k_{d}}$ and

\begin{equation*}
D^{\mathbf{k}}=\frac{\partial ^{\left\vert k\right\vert }}{\partial
^{k_{1}}x_{1}\cdot \cdot \cdot \partial ^{k_{d}}x_{d}}
\end{equation*}

The Sobolev space $W_{p}^{m}(B^{d})$ is the space of functions defined on $%
B^{d}$ with the norm
\begin{equation*}
\left\Vert f\right\Vert _{m,p}=\left\{
\begin{array}{c}
\left( \sum_{0\leq \left\vert \mathbf{k}\right\vert \leq m}\left\Vert D^{%
\mathbf{k}}f\right\Vert _{p}^{p}\right) ^{1/p},\text{ if }1\leq p<\infty \\
\max_{0\leq \left\vert \mathbf{k}\right\vert \leq m}\left\Vert D^{\mathbf{k}%
}f\right\Vert _{\infty },\text{ if }p=\infty .%
\end{array}%
\right.
\end{equation*}

Here
\begin{equation*}
\left\Vert h(\mathbf{x})\right\Vert _{p}=\left\{
\begin{array}{c}
\left( \int_{B^{n}}\left\vert h(\mathbf{x})\right\vert ^{p}d\mathbf{x}%
\right) ^{1/p},\text{ if }1\leq p<\infty \\
ess\sup_{\mathbf{x}\in B^{d}}\left\vert h(\mathbf{x})\right\vert ,\text{ if }%
p=\infty .%
\end{array}%
\right.
\end{equation*}

Let $S_{p}^{m}(B^{d})$ be the unit ball in $W_{p}^{m}(B^{d})$:
\begin{equation*}
S_{p}^{m}(B^{d})=\{f\in W_{p}^{m}(B^{d}):\left\Vert f\right\Vert _{m,p}\leq
1~\}.
\end{equation*}

In 1999, Maiorov \cite{102} proved the following result

\bigskip

\textbf{Proposition 1.6 }(Maiorov \cite{102}). \textit{Assume $m\geq 1$ and $%
d\geq 2$. Then for each $r\in \mathbb{N}$ there exists a function $f\in
S_{2}^{m}(B^{d})$ such that}

\begin{equation*}
E(f,\mathcal{R}_{r},L_{2})\geq Cr^{-m/(d-1)},\eqno(1.2)
\end{equation*}
\textit{where $C$ is a constant independent of $f$ and $r.$}

\bigskip

For $d=2,$ this inequality was proved by Oskolkov \cite{114}. In \cite{102},
Maiorov also proved that for each function $f\in S_{2}^{m}(B^{d})$
\begin{equation*}
E(f,\mathcal{R}_{r},L_{2})\leq Cr^{-m/(d-1)}.\eqno(1.3)
\end{equation*}
Thus he established the following order for the error of approximation to
functions in $S_{2}^{m}(B^{d})$ from the class $\mathcal{R}_{r}$:
\begin{equation*}
E(S_{2}^{m}(B^{d}),\mathcal{R}_{r},L_{2})\overset{def}{=}\sup_{f\in
S_{2}^{m}(B^{d})}E(f,\mathcal{R}_{r},L_{2})\asymp r^{-m/(d-1)}.
\end{equation*}

Pinkus \cite{119} revealed that the upper bound (1.3) is also valid in the $L_{p}
$ metric ($1\leq p\leq \infty $). In other words, for every function $f\in
S_{p}^{m}(B^{d})$
\begin{equation*}
E(f,\mathcal{R}_{r},L_{p})\leq Cr^{-m/(d-1)}.
\end{equation*}

These inequalities were successfully applied to some problems of
approximation of multivariate functions by neural networks with a single
hidden layer. Recall that such networks are given by the formula $%
\sum_{i=1}^{r}c_{i}\sigma (\mathbf{w}^{i}\mathbf{\cdot x}-\theta _{i}).$ By $%
\mathcal{M}_{r}(\sigma )$ let us denote the set of all single hidden layer
networks with the activation function $\sigma $. That is,
\begin{equation*}
\mathcal{M}_{r}(\sigma )=\left\{ \sum_{i=1}^{r}c_{i}\sigma (\mathbf{w}^{i}%
\mathbf{\cdot x}-\theta _{i}):~c_{i},\theta _{i}\in \mathbb{R},~\mathbf{w}%
^{i}\in \mathbb{R}^{d}\right\} .
\end{equation*}

The above results on ridge approximation from $\mathcal{R}_{r}$ enable us to
estimate the rate with which the approximation error $E(f,\mathcal{M}%
_{r}(\sigma ),L_{2})$ tends to zero. First note that $\mathcal{M}_{r}(\sigma
)\subset \mathcal{R}_{r},$ since each function of the form $\sigma (\mathbf{%
w\cdot x}-\theta )$ is a ridge function with the direction $\mathbf{w}$.
Thus the lower bound (1.2) holds also for the set $\mathcal{M}_{r}(\sigma )$%
: there exists a function $f\in S_{2}^{m}(B^{d})$ for which
\begin{equation*}
E(f,\mathcal{M}_{r}(\sigma ),L_{2})\geq Cr^{-m/(d-1)}.
\end{equation*}

It remains to see whether the upper bound (1.3) is valid for $\mathcal{M}%
_{r}(\sigma )$. Clearly, it cannot be valid if $\sigma $ is an arbitrary
continuous function. Here we are dealing with the question if there exists a
function $\sigma ^{\ast }\in C(\mathbb{R})$, for which
\begin{equation*}
E(f,\mathcal{M}_{r}(\sigma ^{\ast }),L_{2})\leq Cr^{-m/(d-1)}.
\end{equation*}

This question is answered affirmatively by the following result.

\bigskip

\textbf{Proposition 1.7 }(Maiorov, Pinkus \cite{99}). \textit{There exists a
function $\sigma ^{\ast }\in C(\mathbb{R})$ with the following properties}

1) \textit{$\sigma ^{\ast }$ is infinitely differentiable and strictly
increasing;}

2) \textit{$\lim_{t\rightarrow \infty }\sigma^{\ast } (t)=1$ and $\lim_{t\rightarrow
-\infty }\sigma^{\ast } (t)=0;$}

3) \textit{for every $g\in \mathcal{R}_{r}$ and $\varepsilon >0$ there exist
$c _{i},\theta _{i}\in \mathbb{R}$ and $\mathbf{w}^{i}\in \mathbb{R}^{d}$
satisfying}
\begin{equation*}
\sup_{\mathbf{x}\in B^{d}}\left\vert g(\mathbf{x})-\sum_{i=1}^{r+d+1}c_{i}%
\sigma^{\ast } (\mathbf{w}^{i}\mathbf{\cdot x}-\theta _{i})\right\vert <\varepsilon .
\end{equation*}

\bigskip

Temlyakov \cite{138} considered the approximation from some certain subclass
of $\mathcal{R}_{r}$ in $L_{2}$ metric. More precisely, he considered the
approximation of a function $f\in L_{2}(D),$ where $D$ is the unit disk in $%
\mathbb{R}^{2}$, by functions $\sum\limits_{i=1}^{r}g_{i}(\mathbf{a}%
^{i}\cdot \mathbf{x)}\in \mathcal{R}_{r}\cap L_{2}(D)$, which satisfy the
additional condition $\left\Vert g_{i}(\mathbf{a}^{i}\cdot \mathbf{x)}%
\right\Vert _{2}\leq B\left\Vert f\right\Vert _{2},$ $i=1,...,r$ ($B$ is a
given positive number). Let $\sigma _{r}^{B}(f)$ be the error of this
approximation. For this approximation error, the author of \cite{138}
obtained upper and lower bounds. Let, for $\alpha >0,$ $H^{\alpha }(D)$
denote the set of all functions $f\in L_{2}(D)$, which can be represented in
the form
\begin{equation*}
f=\sum_{n=1}^{\infty }P_{n},
\end{equation*}%
where $P_{n}$ are bivariate algebraic polynomials of total degree $2^{n}-1$
satisfying the inequalities
\begin{equation*}
\left\Vert P_{n}\right\Vert _{2}\leq 2^{-\alpha n},\text{ }n=1,2,...
\end{equation*}

\bigskip

\textbf{Proposition 1.8 }(Temlyakov \cite{138}). \textit{1) For every $f\in
H^{\alpha }(D) $, we have}
\begin{equation*}
\sigma _{r}^{1}(f)\leq C(\alpha )r^{-\alpha }.
\end{equation*}

\textit{2) For any given $\alpha >0$, $B>0$, $r>1$, there exists a function $%
f\in H^{\alpha }(D)$ such that}
\begin{equation*}
\sigma _{r}^{B}(f)\geq C(\alpha ,B)(r\ln r)^{-\alpha }.
\end{equation*}

\bigskip

Petrushev \cite{116} proved the following interesting result: Let $X_{k}$ be
the $k$ dimensional linear space of univariate functions in $L_{2}[-1,1],$ $%
k=1,2,...$. Besides, let $B^{d}$ and $S^{d-1}$ denote correspondingly the
unit ball and unit sphere in the space $\mathbb{R}^{d}$. If $X_{k}$ provides
order of approximation $O(k^{-m})$ for univariate functions with $m$
derivatives in$\ L_{2}[-1,1]$ and $\Omega _{k}$ are appropriately chosen
finite sets of directions distributed on $S^{d-1}$, then the space $%
Y_{k}=span\{p_{k}(\mathbf{a}\cdot \mathbf{x}):~p_{k}\in X_{k},~\mathbf{a}\in
\Omega _{k}\}$ will provide approximation of order $O(k^{-m-d/2+1/2})$ for
every function $f\in L_{2}(B^{d})$ with smoothness of order $m+d/2-1/2$.
Thus, Petrushev showed that the above form of ridge approximation has the
same efficiency of approximation as the traditional multivariate polynomial
approximation.

Many other results concerning the approximation of multivariate functions by
functions from the set $\mathcal{R}_{r}$ and their applications in neural
network theory may be found in \cite{43,94,99,119,123}.

\bigskip

\section{Representation of multivariate functions by linear combinations of
ridge functions}

In this section we develop a technique for verifying if a multivariate function can be expressed as a sum of ridge functions with given directions. We also obtain a necessary and sufficient condition for the representation of all multivariate functions on a subset $X$ of $\mathbb{R}^{d}$ by sums of ridge functions with fixed directions.

\subsection{Two representation problems}

Let $X$ be a subset of ${{\mathbb{R}}}^{d}$ and $\{\mathbf{a}%
^{i}\}_{i=1}^{r} $ be arbitrarily fixed nonzero directions (vectors) in ${{%
\mathbb{R}}}^{d}$. Consider the following set of linear combinations of
ridge functions.
\begin{equation*}
\mathcal{R}(\mathbf{a}^{1},...,\mathbf{a}^{r};X)=\left\{
\sum\limits_{i=1}^{r}g_{i}(\mathbf{a}^{i}\cdot \mathbf{x}),~\mathbf{x}\in
X,~g_{i}:\mathbb{R}\rightarrow \mathbb{R},~i=1,...,r\right\}
\end{equation*}%
In this section, we are going to deal with the following two problems:

\bigskip

\textbf{Problem 1.} \textit{What conditions imposed on $f:X\rightarrow
\mathbb{R}$ are necessary and sufficient for the inclusion $f\in \mathcal{R}(%
\mathbf{a}^{1},...,\mathbf{a}^{r};X)$?}

\bigskip

\textbf{Problem 2.} \textit{What conditions imposed on $X$ are necessary and
sufficient that every function defined on $X$ belongs to the space $\mathcal{%
R}(\mathbf{a}^{1},...,\mathbf{a}^{r};X)$?}

\bigskip

As noticed in Section 1.1, Problem 1 was considered for continuous functions
in \cite{95} and a theoretical result was obtained. It was also noticed
there that the similar problem of representation of $f$ in the form $%
\sum_{i=1}^{r}g_{i}(\mathbf{a}^{i}\cdot \mathbf{x})+P(\mathbf{x})$ with
polynomial $P(\mathbf{x})$ was solved for continuously differentiable
functions in \cite{25}. Problem 2 was solved in \cite{10} for finite subsets
$X$ of ${{\mathbb{R}}}^{d}$ and in \cite{81} for the case when $r=d$ and $%
\mathbf{a}^{i}$ are the coordinate directions.

Here we consider both Problem 1 and Problem 2 without imposing on $X$, $f$
and $r$ any conditions. In fact, we solve these problems for more general,
than $\mathcal{R}(\mathbf{a}^{1},...,\mathbf{a}^{r};X)$, set of functions.
Namely, we solve them for the set
\begin{equation*}
\mathcal{B}(X)=\mathcal{B}(h_{1},...,h_{r};X)=\left\{
\sum\limits_{i=1}^{r}g_{i}(h_{i}(x)),~x\in X,~g_{i}:\mathbb{R}\rightarrow
\mathbb{R},~i=1,...,r\right\} ,
\end{equation*}%
where $h_{i}:X\rightarrow {{\mathbb{R}}},~i=1,...,r,$ are arbitrarily fixed
functions. In particular, the functions $h_{i},~i=1,...,r$, may be equal to
scalar products of the variable $\mathbf{x}$ with some vectors $\mathbf{a}%
^{i}$, $i=1,...,r$. Only in this special case, we have $\mathcal{B}%
(h_{1},...,h_{r};X)=\mathcal{R}(\mathbf{a}^{1},...,\mathbf{a}^{r};X).$

\bigskip

\subsection{Cycles}

The main idea leading to solutions of the above problems is in using new
objects called \textit{cycles} with respect to $r$ functions $%
h_{i}:X\rightarrow \mathbb{R},~i=1,...,r$ (and in particular, with respect
to $r$ directions $\mathbf{a}^{1},...,\mathbf{a}^{r}$). In the sequel, by $%
\delta _{A}$ we will denote the characteristic function of a set $\ A\subset
\mathbb{R}.$ That is,

\begin{equation*}
\delta _{A}(y)=\left\{
\begin{array}{c}
1,~if~y\in A \\
0,~if~y\notin A.%
\end{array}%
\right.
\end{equation*}

\bigskip

\textbf{Definition 1.1.} \textit{Given a subset $X\subset \mathbb{R}^{d}$\
and functions $h_{i}:X\rightarrow \mathbb{R},~i=1,...,r$. A set of points $%
\{x_{1},...,x_{n}\}\subset X$ is called a cycle with respect to the
functions $h_{1},...,h_{r}$ (or, concisely, a cycle if there is no
confusion), if there exists a vector $\lambda =(\lambda _{1},...,\lambda
_{n})$ with the nonzero real coordinates $\lambda _{i},~i=1,...,n,$ such that%
}
\begin{equation*}
\sum_{j=1}^{n}\lambda _{j}\delta _{h_{i}(x_{j})}=0,~i=1,...,r.\eqno(1.4)
\end{equation*}

\bigskip

If $h_{i}=\mathbf{a}^{i}\cdot \mathbf{x}$, $i=1,...,r$, where $\mathbf{a}%
^{1},...,\mathbf{a}^{r}$ are some directions in $\mathbb{R}^{d}$, a cycle,
with respect to the functions $h_{1},...,h_{r}$, is called a cycle with
respect to the directions\textit{\ }$\mathbf{a}^{1},...,\mathbf{a}^{r}.$

Let for $i=1,...,r,$ the set $\{h_{i}(x_{j}),~j=1,...,n\}$ have $k_{i}$
different values. Then it is not difficult to see that Eq. (1.4) stands for
a system of $\sum_{i=1}^{r}k_{i}$ homogeneous linear equations in unknowns $%
\lambda _{1},...,\lambda _{n}.$ If this system has any solution with the
nonzero components, then the given set $\{x_{1},...,x_{n}\}$ is a cycle. In
the last case, the system has also a solution $m=(m_{1},...,m_{n})$ with the
nonzero integer components $m_{i},~i=1,...,n.$ Thus, in Definition 1.1, the
vector $\lambda =(\lambda _{1},...,\lambda _{n})$ can be replaced with a
vector $m=(m_{1},...,m_{n})$ with $m_{i}\in \mathbb{Z}\backslash \{0\}.$

For example, the set $l=\{(0,0,0),~(0,0,1),~(0,1,0),~(1,0,0),~(1,1,1)\}$ is
a cycle in $\mathbb{R}^{3}$ with respect to the functions $%
h_{i}(z_{1},z_{2},z_{3})=z_{i},~i=1,2,3.$ The vector $\lambda $ in
Definition 1.1 can be taken as $(2,1,1,1,-1).$

In case $r=2,$ the picture of cycles becomes more clear. Let, for example, $%
h_{1}$ and $h_{2}$ be the coordinate functions on $\mathbb{R}^{2}.$ In this
case, a cycle is the union of some sets $A_{k}$ with the property: each $%
A_{k}$ consists of vertices of a closed broken line with the sides parallel
to the coordinate axis. These objects (sets $A_{k}$) have been exploited in
practically all works devoted to the approximation of bivariate functions by
univariate functions, although under various different names (see ``bolt of lightning"
in Section 1.3). If the functions $h_{1}$ and $h_{2}$ are
arbitrary, the sets $A_{k}$ can be described as a trace of some point
traveling alternatively in the level sets of $h_{1}$ and $h_{2},$ and then
returning to its primary position. It should be remarked that in the case $%
r>2,$ cycles do not admit such a simple geometric description. We refer the
reader to Braess and Pinkus \cite{10} for the description of cycles when $r=3
$ and $h_{i}(\mathbf{x})=\mathbf{a}^{i}\cdot \mathbf{x},$ $\mathbf{x}\in
\mathbb{R}^{2},~\mathbf{a}^{i}\in \mathbb{R}^{2}\backslash \{\mathbf{0}%
\},~i=1,2,3.$

Let $T(X)$ denote the set of all functions on $X.$ With each pair $%
\left\langle p,\lambda \right\rangle ,$ where $p=\{x_{1},...,x_{n}\}$ is a
cycle in $X$ and $\lambda =(\lambda _{1},...,\lambda _{n})$ is a vector
known from Definition 1.1, we associate the functional
\begin{equation*}
G_{p,\lambda }:T(X)\rightarrow \mathbb{R},~~G_{p,\lambda
}(f)=\sum_{j=1}^{n}\lambda _{j}f(x_{j}).
\end{equation*}%
In the following, such pairs $\left\langle p,\lambda \right\rangle $ will be
called \textit{cycle-vector pairs} of $X.$ It is clear that the functional $%
G_{p,\lambda }$ is linear and $G_{p,\lambda }(g)=0$ for all functions $g\in
\mathcal{B}(h_{1},...,h_{r};X).$

\bigskip

\textbf{Lemma 1.1.} \textit{Let $X$ have cycles and $h_{i}(X)\cap
h_{j}(X)=\varnothing ,$ for all $i,j\in \{1,...,r\},~i\neq j.$ Then a
function $f:X\rightarrow \mathbb{R}$ belongs to the set $\mathcal{B}%
(h_{1},...,h_{r};X)$ if and only if $G_{p,\lambda }(f)=0$ for any
cycle-vector pair $\left\langle p,\lambda \right\rangle $ of $X.$}

\bigskip

\begin{proof} The necessity is obvious, since the functional $G_{p,\lambda }$
annihilates all members of
$\mathcal{B}(h_{1},...,h_{r};X)$. Let us prove the sufficiency. Introduce
the notation
\begin{eqnarray*}
Y_{i} &=&h_{i}(X),~i=1,...,r; \\
\Omega &=&Y_{1}\cup ...\cup Y_{r}.
\end{eqnarray*}

Consider the following set.
\begin{equation*}
\mathcal{L}=\{Y=\{y_{1},...,y_{r}\}:\text{if there exists }x\in X\text{ such
that }h_{i}(x)=y_{i},~i=1,...,r\}\eqno(1.5)
\end{equation*}

Note that $\mathcal{L}$ is not a subset of $\Omega $. It is a set
of some certain subsets of $\Omega .$ Each element of $\mathcal{L}$ is a set
$Y=\{y_{1},...,y_{r}\}\subset \Omega $ with the property that there exists $%
x\in X$ such that $h_{i}(x)=y_{i},~i=1,...,r.$

In what follows, all the points $x$ associated with $Y$ by (1.5) will be
called $(\ast )$-points of $Y.$ It is clear that the number of such points
depends on $Y$ as well as on the functions $h_{1},...,h_{r}$, and may be
greater than 1. But note that if any two points $x_{1}$ and $x_{2}$ are $%
(\ast )$-points of $Y$, then the set $\{x_{1}$, $x_{2}\}$ necessarily forms
a cycle with the associated vector $\lambda _{0}=(1;-1).$ Indeed, if $x_{1}$
and $x_{2}$ are $(\ast )$-points of $Y$, then $h_{i}(x_{1})=h_{i}(x_{2})$, $%
i=1,...,r,$ whence
\begin{equation*}
1\cdot \delta _{h_{i}(x_{1})}+(-1)\cdot \delta _{h_{i}(x_{2})}\equiv
0,~i=1,...,r.
\end{equation*}

The last identity means that the set $p_{0}=\{x_{1},$ $x_{2}\}$ forms a
cycle and $\lambda _{0}=(1;-1)$ is an associated vector. Then by the
the sufficiency condition, $G_{p_{0},\lambda _{0}}(f)=0$, whcih yields
that $f(x_{1})=f(x_{2})$.

Let now $Y^{\ast }$ be the set of all $(\ast )$-points of $Y.$ Since we have
already known that $f(Y^{\ast })$ is a single number, we can define the
function
\begin{equation*}
t:\mathcal{L}\rightarrow \mathbb{R},~t(Y)=f(Y^{\ast }).
\end{equation*}%
Or, equivalently, $t(Y)=f(x),$ where $x$ is an arbitrary $(\ast )$-point of $%
Y$.

Consider now a class $\mathcal{S}$ of functions of the form $%
\sum_{j=1}^{k}r_{j}\delta _{D_{j}},$ where $k$ is a positive integer, $r_{j}$
are real numbers and $D_{j}$ are elements of $\mathcal{L},~j=1,...,k.$ We
fix neither the numbers $\ k,~r_{j},$ nor the sets $D_{j}.$ Clearly, $%
\mathcal{S\ }$is a linear space. Over $\mathcal{S}$, we define the functional

\begin{equation*}
F:\mathcal{S}\rightarrow \mathbb{R},~F\left( \sum_{j=1}^{k}r_{j}\delta
_{D_{j}}\right) =\sum_{j=1}^{k}r_{j}t(D_{j}).
\end{equation*}

First of all, we must show that this functional is well defined. That is,
the equality
\begin{equation*}
\sum_{j=1}^{k_{1}}r_{j}^{\prime }\delta _{D_{j}^{\prime
}}=\sum_{j=1}^{k_{2}}r_{j}^{\prime \prime }\delta _{D_{j}^{\prime \prime }}
\end{equation*}%
always implies the equality

\begin{equation*}
\sum_{j=1}^{k_{1}}r_{j}^{\prime }t(D_{j}^{\prime
})=\sum_{j=1}^{k_{2}}r_{j}^{\prime \prime }t(D_{j}^{\prime \prime }).
\end{equation*}%
In fact, this is equivalent to the implication
\begin{equation*}
\sum_{j=1}^{k}r_{j}\delta _{D_{j}}=0\Longrightarrow
\sum_{j=1}^{k}r_{j}t(D_{j})=0,~\text{for all }k\in \mathbb{N}\text{, }%
r_{j}\in \mathbb{R}\text{, }D_{j}\subset \mathcal{L}\text{.}\eqno(1.6)
\end{equation*}

Suppose that the left-hand side of the implication (1.6) be satisfied. Each
set $D_{j}$ consists of $r$ real numbers $y_{1}^{j},...,y_{r}^{j}$, $%
j=1,...,k.$ By the hypothesis of the lemma, all these numbers are different.
Therefore,

\begin{equation*}
\delta _{D_{j}}=\sum_{i=1}^{r}\delta _{y_{i}^{j}},~j=1,...,k.\eqno(1.7)
\end{equation*}%
Eq. (1.7) together with the left-hand side of (1.6) gives

\begin{equation*}
\sum_{i=1}^{r}\sum_{j=1}^{k}r_{j}\delta _{y_{i}^{j}}=0.\eqno(1.8)
\end{equation*}%
Since the sets $\{y_{i}^{1},y_{i}^{2},...,y_{i}^{k}\}$, $i=1,...,r,$ are
pairwise disjoint, we obtain from (1.8) that
\begin{equation*}
\sum_{j=1}^{k}r_{j}\delta _{y_{i}^{j}}=0,\text{ }i=1,...,r.\eqno(1.9)
\end{equation*}

Let now $x_{1},...,x_{k}$ be some $(\ast )$-points of the sets $%
D_{1},...,D_{k}$ respectively. Since by (1.5), $y_{i}^{j}=h_{i}(x_{j})$, for
$i=1,...,r$ and $j=1,...,k,$ it follows from (1.9) that the set $%
\{x_{1},...,x_{k}\}$ is a cycle. Then by the condition of the sufficiency, $%
\sum_{j=1}^{k}r_{j}f(x_{j})=0.$ Hence $\sum_{j=1}^{k}r_{j}t(D_{j})=0.$ We
have proved the implication (1.6) and hence the functional $F$ is well
defined. Note that the functional $F$ is linear (this can be easily seen
from its definition).

Consider now the following space:

\begin{equation*}
\mathcal{S}^{\prime }=\left\{ \sum_{j=1}^{k}r_{j}\delta _{\omega
_{j}}\right\} ,
\end{equation*}%
where $k\in \mathbb{N}$, $r_{j}\in \mathbb{R}$, $\omega _{j}\subset \Omega .$
As above, we do not fix the parameters $k$, $r_{j}$ and $\omega _{j}.$
Clearly, the space $\mathcal{S}^{\prime }$ is larger than $\mathcal{S}$. Let
us prove that the functional $F$ can be linearly extended to the space $%
\mathcal{S}^{\prime }$. So, we must prove that there exists a linear
functional $F^{\prime }:\mathcal{S}^{\prime }\rightarrow \mathbb{R}$ such
that $F^{\prime }(x)=F(x)$, for all $x\in \mathcal{S}$. Let $H$ denote the
set of all linear extensions of $F$ to subspaces of $\mathcal{S}^{\prime }$
containing $\mathcal{S}$. The set $H$ is not empty, since it contains a
functional $F.$ For each functional $v\in H$, let $dom(v)$ denote the domain
of $v$. Consider the following partial order in $H$: $v_{1}\leq v_{2}$, if $%
v_{2}$ is a linear extension of $v_{1}$ from the space $dom(v_{1})$ to the
space $dom(v_{2}).$ Let now $P$ be any chain (linearly ordered subset) in $H$%
. Consider the following functional $u$ defined on the union of domains of
all functionals $p\in P$:
\begin{equation*}
u:\bigcup\limits_{p\in P}dom(p)\rightarrow \mathbb{R},~u(x)=p(x),\text{ if }%
x\in dom(p)
\end{equation*}

Obviously, this functional is well defined and linear. Besides, the
functional $u$ provides an upper bound for $P.$ We see that the arbitrarily
chosen chain $P$ has an upper bound. Then by Zorn's lemma, there is a
maximal element $F^{\prime }\in H$. We claim that the functional $F^{\prime
} $ must be defined on the whole space $\mathcal{S}^{\prime }$. Indeed, if $%
F^{\prime }$ is defined on a proper subspace $\mathcal{D\subset }$ $\mathcal{%
S}^{\prime }$, then it can be linearly extended to a space larger than $%
\mathcal{D}$ by the following way: take any point $x\in \mathcal{S}^{\prime
}\backslash \mathcal{D}$ and consider the linear space $\mathcal{D}^{\prime
}=\{\mathcal{D}+\alpha x\}$, where $\alpha $ runs through all real numbers.
For an arbitrary point $y+\alpha x\in \mathcal{D}^{\prime }$, set $%
F^{^{\prime \prime }}(y+\alpha x)=F^{\prime }(y)+\alpha b$, where $b$ is any
real number considered as the value of $F^{^{\prime \prime }}$ at $x$. Thus,
we constructed a linear functional $F^{^{\prime \prime }}\in H$ satisfying $%
F^{\prime }\leq F^{^{\prime \prime }}.$ The last contradicts the maximality
of $F^{\prime }.$ This means that the functional $F^{\prime }$ is defined on
the whole $\mathcal{S}^{\prime }$ and $F\leq F^{\prime }$ ($F^{\prime }$ is
a linear extension of $F$).

Define the following functions by means of the functional $F^{\prime }$:
\begin{equation*}
g_{i}:Y_{i}\rightarrow \mathbb{R},\text{ }g_{i}(y_{i})\overset{def}{=}%
F^{\prime }(\delta _{y_{i}}),\text{ }i=1,...,r.
\end{equation*}%
Let $x$ be an arbitrary point in $X.$ Obviously, $x$ is a $(\ast )$-point of
some set $Y=\{y_{1},...,y_{r}\}\subset \mathcal{L}.$ Thus,
\begin{eqnarray*}
f(x) &=&t(Y)=F(\delta _{Y})=F\left( \sum_{i=1}^{r}\delta _{y_{i}}\right)
=F^{\prime }\left( \sum_{i=1}^{r}\delta _{y_{i}}\right) = \\
\sum_{i=1}^{r}F^{\prime }(\delta _{y_{i}})
&=&\sum_{i=1}^{r}g_{i}(y_{i})=\sum_{i=1}^{r}g_{i}(h_{i}(x)).
\end{eqnarray*}%
\end{proof}

\bigskip

\subsection{Minimal cycles and the main results}

\textbf{Definition 1.2.} \textit{A cycle $p=\{x_{1},...,x_{n}\}$ is said to
be minimal if $p$ does not contain any cycle as its proper subset.}

\bigskip

For example, the set $l=\{(0,0,0),~(0,0,1),~(0,1,0),~(1,0,0),~(1,1,1)\}$
considered above is a minimal cycle with respect to the functions $%
h_{i}(z_{1},z_{2},z_{3})=z_{i},~i=1,2,3.$ Adding the point $(0,1,1)$ to $l$,
we will have a cycle, but not minimal. The vector $\lambda $ associated with
$l\cup \{(0,1,1)\}$ can be taken as $(3,-1,-1,-2,2,-1).$

A minimal cycle $p=\{x_{1},...,x_{n}\}$ has the following obvious properties:

\begin{description}
\item[(a)] \textit{The vector $\lambda $ associated with $p$ through Eq.
(1.4) is unique up to multiplication by a constant;}

\item[(b)] \textit{If in (1.4), $\sum_{j=1}^{n}\left\vert \lambda
_{j}\right\vert =1,$ then all the numbers $\lambda _{j},~j=1,...,n,$ are
rational.}
\end{description}

Thus, a minimal cycle $p$ uniquely (up to a sign) defines the functional

\begin{equation*}
~G_{p}(f)=\sum_{j=1}^{n}\lambda _{j}f(x_{j}),\text{ \ }\sum_{j=1}^{n}\left%
\vert \lambda _{j}\right\vert =1.
\end{equation*}

\bigskip

\textbf{Lemma 1.2.} \textit{The functional $G_{p,\lambda }$ is a linear
combination of functionals $G_{p_{1}},...,G_{p_{k}},$ where $p_{1},...,p_{k}$
are minimal cycles in $p.$}

\bigskip

\begin{proof} Let $\left\langle p,\lambda \right\rangle $ be a cycle-vector
pair of $X$, where $p=\{x_{1},...,x_{n}\}$ and $\lambda =(\lambda
_{1},...,\lambda _{n})$. Let $p_{1}=$ $\{y_{1}^{1},...,y_{s_{1}}^{1}\},$ $%
s_{1}<n$, be a minimal cycle in $p$ and

\begin{equation*}
G_{p_{1}}(f)=\sum_{j=1}^{s_{1}}\nu _{j}^{1}f(y_{j}^{1}),\text{ }%
\sum_{j=1}^{s_{1}}\left\vert \nu _{j}^{1}\right\vert =1.
\end{equation*}

Without loss of generality, we may assume that $y_{1}^{1}=x_{1}.$ Put

\begin{equation*}
t_{1}=\frac{\lambda _{1}}{\nu _{1}^{1}}.
\end{equation*}%
Then the functional $G_{p,\lambda }-t_{1}G_{p_{1}}$ has the form

\begin{equation*}
G_{p,\lambda }-t_{1}G_{p_{1}}=\sum_{j=1}^{n_{1}}\lambda _{j}^{1}f(x_{j}^{1}),
\end{equation*}%
where $x_{j}^{1}\in p$, $\lambda _{j}^{1}\neq 0$, $j=1,...,n_{1}$. Clearly,
the set $l_{1}=\{x_{1}^{1},...,x_{n_{1}}^{1}\}$ is a cycle in $p$ with the
associated vector $\lambda ^{1}=(\lambda _{1}^{1},...,\lambda _{n_{1}}^{1})$%
. Besides, $x_{1}\notin l_{1}$. Thus, $n_{1}<n$ and $G_{l_{1},\lambda ^{1}}=$
$G_{p,\lambda }-t_{1}G_{p_{1}}$. If $l_{1}$ is minimal, then the proof is
completed. Assume $l_{1}$ is not minimal. Let $p_{1}=$ $%
\{y_{1}^{2},...,y_{s_{2}}^{2}\},$ $s_{2}<n_{1},$ be a minimal cycle in $%
l_{1} $ and

\begin{equation*}
G_{p_{2}}(f)=\sum_{j=1}^{s_{2}}\nu _{j}^{2}f(y_{j}^{2}),\text{ }%
\sum_{j=1}^{s_{2}}\left\vert \nu _{j}^{2}\right\vert =1.
\end{equation*}

Without loss of generality, we may assume that $y_{1}^{2}=x_{1}^{1}.$ Put

\begin{equation*}
t_{2}=\frac{\lambda _{1}^{1}}{\nu _{1}^{2}}.
\end{equation*}%
Then the functional $G_{l_{1},\lambda ^{1}}-t_{2}G_{p_{2}}$ has the form

\begin{equation*}
G_{l_{1},\lambda ^{1}}-t_{2}G_{p_{2}}=\sum_{j=1}^{n_{2}}\lambda
_{j}^{2}f(x_{j}^{2}),
\end{equation*}%
where $x_{j}^{2}\in l_{1}$, $\lambda _{j}^{2}\neq 0$, $j=1,...,n_{2}$.
Clearly, the set $l_{2}=\{x_{1}^{2},...,x_{n_{2}}^{2}\}$ is a cycle in $%
l_{1} $ with the associated vector $\lambda ^{2}=(\lambda
_{1}^{2},...,\lambda _{n_{2}}^{2})$. Besides, $x_{1}^{1}\notin l_{2}$. Thus,
$n_{2}<n_{1}$ and $G_{l_{2},\lambda ^{2}}=$ $G_{l_{1},\lambda
^{1}}-t_{2}G_{p_{2}}.$ If $l_{2}$ is minimal, then the proof is completed.
Let $l_{2}$ be not minimal. Repeating the above process for $l_{2}$, then
for $l_{3}$, etc., after some $k-1$ steps we will come to a minimal cycle $%
l_{k-1}$ and the functional

\begin{equation*}
G_{l_{k-1},\lambda ^{k-1}}=G_{l_{k-2},\lambda
^{k-2}}-t_{k-1}G_{p_{k-1}}=\sum_{j=1}^{n_{k-1}}\lambda
_{j}^{k-1}f(x_{j}^{k-1}).
\end{equation*}%
Since the cycle $l_{k-1}$ is minimal,

\begin{equation*}
G_{l_{k-1},\lambda ^{k-1}}=t_{k}G_{l_{k-1}},\text{ \ where }%
t_{k}=\sum_{j=1}^{n_{k-1}}\left\vert \lambda _{j}^{k-1}\right\vert .
\end{equation*}%
Now putting $p_{k}=l_{k-1}$ and considering the above chain relations
between the functionals $G_{l_{i},\lambda ^{i}}$, $i=1,...,k-1,$ we obtain
that
\begin{equation*}
G_{p,\lambda }=\sum_{i=1}^{k}t_{i}G_{p_{i}}.
\end{equation*}
\end{proof}

\textbf{Theorem 1.1.} \textit{Assume $X\subset \mathbb{R}^{d}$ and $%
h_{1},...,h_{r}$ are arbitrarily fixed real functions on $X.$ The following
assertions are valid.}

\textit{1) Let $X$ have cycles with respect to the functions $%
h_{1},...,h_{r} $. A function $f:X\rightarrow \mathbb{R}$ belongs to the
space $\mathcal{B}(h_{1},...,h_{r};X)$ if and only if $G_{p}(f)=0$ for any
minimal cycle $p\subset X$.}

\textit{2) Let $X$ have no cycles. Then $\mathcal{B}%
(h_{1},...,h_{r};X)=T(X). $}

\bigskip

\begin{proof} 1) The necessity is clear. Let us prove the sufficiency. On the
strength of Lemma 1.2, it is enough to prove that if $G_{p,\lambda }(f)=0$
for any cycle-vector pair $\left\langle p,\lambda \right\rangle $ of $X$,
then $f\in \mathcal{B}(X).$

Consider a system of intervals $\{(a_{i},b_{i})\subset \mathbb{R}%
\}_{i=1}^{r} $ such that $(a_{i},b_{i})\cap (a_{j},b_{j})=\varnothing $ for
all the indices $i,j\in \{1,...,r\}$, $~i\neq j.$ For $i=1,...,r$, let $\tau
_{i}$ be one-to-one mappings of $\mathbb{R}$\ onto $(a_{i},b_{i}).$
Introduce the following functions on $X$:
\begin{equation*}
h_{i}^{^{\prime }}(x)=\tau _{i}(h_{i}(x)),\text{ }i=1,...,r.
\end{equation*}

It is clear that any cycle with respect to the functions $h_{1},...,h_{r}$
is also a cycle with respect to the functions $h_{1}^{^{\prime
}},...,h_{r}^{^{\prime }}$, and vice versa. Besides, $h_{i}^{\prime }(X)\cap
h_{j}^{\prime }(X)=\varnothing ,$ for all $i,j\in \{1,...,r\},~i\neq j.$
Then by Lemma 1.1,
\begin{equation*}
f(x)=g_{1}^{\prime }(h_{1}^{\prime }(x))+\cdots +g_{r}^{\prime
}(h_{r}^{\prime }(x)),
\end{equation*}%
where $g_{1}^{\prime },...,g_{r}^{\prime }$ are univariate functions
depending on $f$. From the last equality we obtain that
\begin{equation*}
f(x)=g_{1}^{\prime }(\tau _{1}(h_{1}(x)))+\cdots +g_{r}^{\prime }(\tau
_{r}(h_{r}(x)))=g_{1}(h_{1}(x))+\cdots +g_{r}(h_{r}(x)).
\end{equation*}%
That is, $f\in \mathcal{B}(X)$.

2) Let $f:X\rightarrow \mathbb{R}$ be an arbitrary function. First suppose
that $h_{i}(X)\cap h_{j}(X)=\varnothing ,$ for all $i,j\in \{1,...,r\}$,$%
~i\neq j.$ In this case, the proof is similar to and even simpler than that
of Lemma 1.1. Indeed, the set of all $(\ast )$-points of $Y$ consists of a
single point, since otherwise we would have a cycle with two points, which
contradicts the hypothesis of the 2-nd part of the theorem. Further,
well definition of the functional $F$ becomes obvious, since the left-hand
side of (1.6) also contradicts the nonexistence of cycles. Thus, as in the
proof of Lemma 1.1, we can extend $F$ to the space $\mathcal{S}^{\prime }$
and then obtain the desired representation for the function $f$. Since $f$
is arbitrary, $T(X)=\mathcal{B}(X).$

Using the techniques from the proof of the 1-st part of the theorem, one can
easily generalize the above argument to the case when the functions $%
h_{1},...,h_{r}$ have arbitrary ranges.
\end{proof}

\bigskip

\textbf{Theorem 1.2.} \textit{$\mathcal{B}(h_{1},...,h_{r};X)=T(X)$ if and
only if $X$ has no cycles with respect to the functions }$h_{1},...,h_{r}$%
\textit{.}

\bigskip

\begin{proof} The sufficiency immediately follows from Theorem 1.1. To prove
the necessity, assume that $X$ has a cycle $p=\{x_{1},...,x_{n}\}$. Let $%
\lambda =(\lambda _{1},...,\lambda _{n})$ be a vector associated with $p$ by
Eq. (1.4). Consider a function $f_{0}$ on $X$ with the property: $%
f_{0}(x_{i})=1,$ for indices $i$ such that $\lambda _{i}\,>0$ and $%
f_{0}(x_{i})=-1,$ for indices $i$ such that $\lambda _{i}\,<0$. For this
function, $G_{p,\lambda }(f_{0})\neq 0$. Then by Theorem 1.1, $f_{0}\notin
\mathcal{B}(X)$. Hence $\mathcal{B}(X)\neq T(X)$. The contradiction shows
that $X$ does not admit cycles.
\end{proof}

\subsection{Corollaries}

From Theorems 1.1 and 1.2 we obtain the following corollaries for the ridge
function representation.

\bigskip

\textbf{Corollary 1.1.} \textit{Assume $X\subset \mathbb{R}^{d}$ and
$\mathbf{a}^{1},...,\mathbf{a}^{r}\in \mathbb{R}^{d}\backslash \{\mathbf{0}%
\}$. The following assertions are valid.}

\textit{1) Let $X$ have cycles with respect to the directions $%
\mathbf{a}^{1},...,\mathbf{a}^{r}$. A function $f:X\rightarrow
\mathbb{R}$ belongs to the space $\mathcal{R}(\mathbf{a}^{1},...,%
\mathbf{a}^{r};X)$ if and only if $G_{p}(f)=0$ for any
minimal cycle $p\subset X$.}

\textit{2) Let $X$ have no cycles. Then every function $%
f:X\rightarrow \mathbb{R}$ belongs to the space $\mathcal{R}(%
\mathbf{a}^{1},...,\mathbf{a}^{r};X)$.}

\bigskip

\textbf{Corollary 1.2.} \textit{$\mathcal{R}(\mathbf{a}^{1},...,\mathbf{a}^{r};X)=T(X)$
if and only if $X$ has no cycles with respect to the
directions $\mathbf{a}^{1},...,\mathbf{a}^{r}$.}

\bigskip

Note that solutions to Problems 1 and 2 are given by Corollaries 1.1 and
1.2, correspondingly. Although it is not always easy to find all cycles of a
given set $X$ and even to know if $X$ possesses a single cycle, Corollaries
1.1 and 1.2 are of more practical than theoretical character. Particular
cases of Problems 1 and 2 evidence in favor of our opinion. For example, for
the problem of representation by sums of two ridge functions, the picture of
cycles is completely describable (see the beginning of this section). The
interpretation of cycles with respect to three directions in the plane can
be found in Braess and Pinkus \cite{10}. A geometric description of cycles
with respect to 4 and more directions is quite complicated and requires deep
techniques from geometry and graph theory. This is not within the aim of our
study.

From the last corollary, it follows that if representation by sums of ridge
functions with fixed directions $\mathbf{a}^{1},...,\mathbf{a}^{r}$ is valid
in the class of continuous functions (or in the class of bounded functions),
then such representation is valid in the class of all functions. For a rigid
mathematical formulation of this result, let us introduce the notation:
\begin{equation*}
\mathcal{R}_{c}(\mathbf{a}^{1},...,\mathbf{a}^{r};X)=\left\{
\sum\limits_{i=1}^{r}g_{i}(\mathbf{a}^{i}\cdot \mathbf{x}),~\mathbf{x}\in
X,~g_{i}(\mathbf{a}^{i}\cdot \mathbf{x})\in C(X\mathbb{)},~i=1,...,r\right\}
\end{equation*}
and
\begin{equation*}
\mathcal{R}_{b}(\mathbf{a}^{1},...,\mathbf{a}^{r};X)=\left\{
\sum\limits_{i=1}^{r}g_{i}(\mathbf{a}^{i}\cdot \mathbf{x}),~\mathbf{x}\in
X,~g_{i}(\mathbf{a}^{i}\cdot \mathbf{x})\in B(X\mathbb{)},~i=1,...,r\right\}
\end{equation*}

Here $C(X)$ and $B(X)$ denote the spaces of continuous and bounded functions
defined on $X\subset \mathbb{R}^{d}$ correspondingly (for the first space,
the set $X$ is supposed to be compact). As we know (see Section 1.1) from
the results of Sternfeld it follows that the equality $\mathcal{R}_{c}(%
\mathbf{a}^{1},...,\mathbf{a}^{r};X)=C(X)$ implies the equality $\mathcal{R}%
_{b}(\mathbf{a}^{1},...,\mathbf{a}^{r};X)=B(X).$ In other words, if every
continuous function is represented by sums of ridge functions (with fixed
directions!), then every bounded function also obeys such representation
(with bounded summands). Corollaries 1.1 and 1.2 allow us to obtain the
following result.

\bigskip

\textbf{Corollary 1.3.} \textit{Let $X$ be a compact subset of $\mathbb{R}%
^{d}$ and $\mathbf{a}^{1},...,\mathbf{a}^{r}$ be given directions in $%
\mathbb{R}^{d}\backslash \{\mathbf{0}\}$. If $\mathcal{R}_{c}(\mathbf{a}%
^{1},...,\mathbf{a}^{r};X)=C(X),$ then $\mathcal{R}(\mathbf{a}^{1},...,%
\mathbf{a}^{r};X)=T(X).$}

\bigskip

\begin{proof} If every continuous function defined on $X\subset \mathbb{R}^{d}$ is
represented by sums of ridge functions with the directions $\mathbf{a}%
^{1},...,\mathbf{a}^{r}$, then it can be shown by applying the same idea (as
in the proof of Theorem 1.2) that the set $X$ has no cycles with respect to
the given directions. Only, because of continuity, Urysohn's great lemma
should be taken into account. That is, it should be taken into account that,
by assuming the existence of a cycle $p_{0}=\{x_{1},...,x_{n}\}$ with an
associated vector $\lambda _{0}=(\lambda _{1},...,\lambda _{n})$, we can
deduce from Urysohn's great lemma the existence of a continuous
function $u:X\rightarrow \mathbb{R}$ satisfying

1) $u(x_{i})=1,$ for indices $i$ such that $\lambda _{i}\,>0$

2) $u(x_{j})=-1,$ for indices $j$ such that $\lambda _{j}\,<0$,

3) $-1<u(x)<1,$ for all $x\in X\backslash p_{0}.$

These properties mean that $G_{p_{0},\lambda _{0}}(u)\neq
0\Longrightarrow u\notin \mathcal{R}_{c}(\mathbf{a}^{1},...,\mathbf{a}%
^{r};X)\Longrightarrow \mathcal{R}_{c}(\mathbf{a}^{1},...,\mathbf{a}%
^{r};X)\neq C(X).$

But if $X$ has no cycles with respect to the directions $\mathbf{a}^{1},...,%
\mathbf{a}^{r}$, then by Corollary 1.2, $\mathcal{R}(\mathbf{a}^{1},...,%
\mathbf{a}^{r};X)=T(X).$
\end{proof}

Let us now give some examples of sets over which the representation by
linear combinations of ridge functions is possible.

\begin{description}
\item[(1)] Let $r=2$ and $X$ be the union of two parallel lines not
perpendicular to the directions $\mathbf{a}^{1}$ and $\mathbf{a}^{2}$.
Then $X$ has no cycles with respect to $\{\mathbf{a}^{1},\mathbf{a}^{2}\}$.
Therefore, by Corollary 1.2, $\mathcal{R}\left( \mathbf{a}^{1},\mathbf{a}%
^{2};X\right) =T(X).$

\item[(2)] Let $r=2,$ $\mathbf{a}^{1}=(1,1)$, $\mathbf{a}^{2}=(1,-1)$ and $X$
be the graph of the function $y=\arcsin (\sin x)$. Then $X$ has no cycles
and hence $\mathcal{R}\left( \mathbf{a}^{1},\mathbf{a}^{2};X\right) =T(X).$

\item[(3)] Assume now we are given $r$ directions $\{\mathbf{a}^{j}\}_{j=1}^{r}$ and $%
r+1$ points $\{\mathbf{x}^{i}\}_{i=1}^{r+1}\subset \mathbb{R}^{d}$ such that%
\begin{eqnarray*}
\mathbf{a}^{1}\cdot \mathbf{x}^{i} &=&\mathbf{a}^{1}\cdot \mathbf{x}^{j}\neq
\mathbf{a}^{1}\cdot \mathbf{x}^{2}\text{, \ for }1\leq i,j\leq r+1\text{, }%
i,j\neq 2 \\
\mathbf{a}^{2}\cdot \mathbf{x}^{i} &=&\mathbf{a}^{2}\cdot \mathbf{x}^{j}\neq
\mathbf{a}^{2}\cdot \mathbf{x}^{3}\text{, \ for }1\leq i,j\leq r+1\text{, }%
i,j\neq 3 \\
&&\mathbf{......................................} \\
\mathbf{a}^{r}\cdot \mathbf{x}^{i} &=&\mathbf{a}^{r}\cdot \mathbf{x}^{j}\neq
\mathbf{a}^{r}\cdot \mathbf{x}^{r+1}\text{, \ for }1\leq i,j\leq r.
\end{eqnarray*}%
The simplest data realizing these equations are the basis directions in $%
\mathbb{R}^{d}$ and the points $(0,0,...,0)$, $(1,0,...,0)$, $(0,1,...,0)$%
,..., $(0,0,...,1)$. From the first equation we obtain that $\mathbf{x}^{2}$
cannot be a point of any cycle in $X=\{\mathbf{x}^{1},...,\mathbf{x}^{r+1}\}$%
. Sequentially, from the second, third, ..., $r$-th equations it follows
that the points $\mathbf{x}^{3},\mathbf{x}^{4},...,\mathbf{x}^{r+1}$ also
cannot be points of cycles in $X$, respectively. Thus the set $X$ does not
contain cycles at all. By Corollary 1.2, $\mathcal{R}\left( \mathbf{a}%
^{1},...,\mathbf{a}^{r};X\right) =T(X).$

\item[(4)] Assume we are given directions $\{\mathbf{a}^{j}\}_{j=1}^{r}$ and a curve $%
\gamma $ in $\mathbb{R}^{d}$ such that for any $c\in \mathbb{R}$, $\gamma $
has at most one common point with at least one of the hyperplanes $\mathbf{a}%
^{j}\cdot \mathbf{x}=c$, $j=1,...,r.$ Clearly, the curve $\gamma $ has no
cycles and hence $\mathcal{R}\left( \mathbf{a}^{1},...,\mathbf{a}^{r};\gamma
\right) =T(\gamma ).$
\end{description}

Braess and Pinkus \cite{10} considered the partial case of Problem 2:
characterize a set of points $\left( \mathbf{x}^{1},...,\mathbf{x}%
^{k}\right) \subset \mathbb{R}^{d}$ such that for any data $\{\alpha
_{1},...,\alpha _{k}\}\subset \mathbb{R}$ there exists a function $g\in
\mathcal{R}\left( \mathbf{a}^{1},...,\mathbf{a}^{r};\mathbb{R}^{d}\right) $
satisfying $g(\mathbf{x}^{i})=\alpha _{i},$ $i=1,...,k$. In connection with
this problem, they introduced the notion of the \textit{NI}-property (non
interpolation property) and \textit{MNI}-property (minimal non interpolation
property) of a finite set of points as follows:

Given directions $\{\mathbf{a}^{j}\}_{j=1}^{r}\subset \mathbb{R}%
^{d}\backslash \{\mathbf{0}\}$, we say that a set of points $\{\mathbf{x}%
^{i}\}_{i=1}^{k}\subset \mathbb{R}^{d}$ has the \textit{NI}-property with
respect to $\{\mathbf{a}^{j}\}_{j=1}^{r}$, if there exists $\{\alpha
_{i}\}_{i=1}^{k}\subset \mathbb{R}$ such that we cannot find a function $%
g\in \mathcal{R}\left( \mathbf{a}^{1},...,\mathbf{a}^{r};\mathbb{R}%
^{d}\right) $ satisfying $g(\mathbf{x}^{i})=\alpha _{i},$ $i=1,...,k$. We
say that \ the set $\{\mathbf{x}^{i}\}_{i=1}^{k}\subset \mathbb{R}^{d}$ has
the \textit{MNI}-property with respect to $\{\mathbf{a}^{j}\}_{j=1}^{r}$, if
$\{\mathbf{x}^{i}\}_{i=1}^{k}$ but no proper subset thereof has the \textit{%
NI}-property.

It follows from Corollary 1.2 that a set $\{\mathbf{x}^{i}\}_{i=1}^{k}$ has
the \textit{NI}-property if and only if $\{\mathbf{x}^{i}\}_{i=1}^{k}$
contains a cycle with respect to the functions $h_{i}=\mathbf{a}^{i}\cdot
\mathbf{x},$ $i=1,...,r$ (or, simply, to the directions $\mathbf{a}^{i},$ $%
i=1,...,r$) and the \textit{MNI}-property if and only if the set $\{\mathbf{x%
}^{i}\}_{i=1}^{k}$ itself is a minimal cycle with respect to the given
directions. Taking into account this argument and Definitions 1.1 and 1.2,
we obtain that the set $\{\mathbf{x}^{i}\}_{i=1}^{k}$ has the \textit{NI}%
-property if and only if there is a vector $\mathbf{m}=(m_{1},...,m_{k})\in
\mathbb{Z}^{k}\backslash \{\mathbf{0}\}$ such that
\begin{equation*}
\sum_{j=1}^{k}m_{j}g(\mathbf{a}^{i}\cdot \mathbf{x}^{j})=0,
\end{equation*}%
for $i=1,...,r$ and all functions $g:\mathbb{R\rightarrow R}$. This set has
the \textit{MNI}-property if and only if the vector $\mathbf{m}$ has the
additional properties: it is unique up to multiplication by a constant and
all its components are different from zero. This special consequence of
Corollary 1.2 was proved in \cite{10}.

\bigskip

\section{Characterization of an extremal sum of ridge functions}

The approximation problem considered in this section is to approximate a
continuous multivariate function $f\left( \mathbf{x}\right) =f\left( {%
x_{1},...,x_{d}}\right) $ by sums of two ridge functions in the uniform
norm. We give a necessary and sufficient condition for a sum of two ridge
functions to be a best approximation to $f\left( \mathbf{x}\right) .$ This
main result is next used in a special case to obtain an explicit formula for
the approximation error and to construct a best approximation. The problem
of well approximation by such sums is also considered.

\subsection{Exposition of the problem}

Consider the following set of sums of ridge functions
\begin{equation*}
\mathcal{R}=\mathcal{R}\left( \mathbf{a}^{1},\mathbf{a}^{2}\right) ={\left\{ {%
g_{1}\left( \mathbf{a}^{1}{\cdot }\mathbf{x}\right) +g_{2}\left( \mathbf{a}^{2}{%
\cdot }\mathbf{x}\right) :g}_{i}{\in C\left( {\mathbb{R}}\right) ,i=1,2}%
\right\} }.
\end{equation*}
That is, we fix directions $\mathbf{a}^{1}$ and $\mathbf{a}^{2}$ and consider linear
combinations of ridge functions with these directions.

Assume $f\left( \mathbf{x}\right) $ is a continuous function on a
compact subset $Q$ of $\mathbb{R}^{d}$. We want to find conditions that are
necessary and sufficient for a function $g_{_{0}}\in \mathcal{R}\left(
\mathbf{a}^{1},\mathbf{a}^{2}\right) $ to be an extremal element (or a best
approximation) to $f$. In other words, we want to characterize such sums $%
g_{0}\left( \mathbf{x}\right) =g_{1}\left( \mathbf{a}^{1}{\cdot }\mathbf{x}%
\right) +g_{2}\left( \mathbf{a}^{2}{\cdot }\mathbf{x}\right) $ of ridge
functions that
\begin{equation*}
{\left\Vert {f-g_{0}}\right\Vert }={\max\limits_{{\mathbf{x}\in Q}}}{%
\left\vert {f\left( \mathbf{x}\right) -g}_{{0}}{\left( \mathbf{x}\right) }%
\right\vert }=E\left( {f}\right) ,
\end{equation*}%
where
\begin{equation*}
E\left( {f}\right) =E(f,\mathcal{R})\overset{def}{=}{\inf_{g \in \mathcal{R}\left(
\mathbf{a}^{1},\mathbf{a}^{2}\right)}}{\left\Vert {f-g}%
\right\Vert }
\end{equation*}%
is the error in approximating from $\mathcal{R}\left( \mathbf{a}^{1},\mathbf{a}^{2}%
\right) .$ The other related problem is how to construct these sums of ridge
functions. We also want to know if we can approximate well, i.e. for which
compact sets $Q,$ $\mathcal{R}\left( \mathbf{a}^{1},\mathbf{a}^{2}\right) $ is
dense in $C\left( {Q}\right) $ in the topology of uniform convergence. It
should be remarked that solutions to these problems may be useful in
connection with the study of partial differential equations. For example,
assume that $\left( {a_{1},b_{1}}\right) $ and $\left( {a_{2},b_{2}}\right) $
are linearly independent vectors in $\mathbb{R}^{2}.$ Then the general
solution to the homogeneous partial differential equation
\begin{equation*}
\left( {a_{1}{\frac{\partial }{\partial {x}}}+b_{1}{\frac{\partial }{%
\partial {y}}}}\right) \left( {a_{2}{\frac{\partial }{\partial {x}}}+b_{2}{%
\frac{\partial }{\partial {y}}}}\right) {u}\left( {x,y}\right) =0\eqno(1.10)
\end{equation*}%
are all functions of the form
\begin{equation*}
u\left( {x,y}\right) =g_{1}\left( {b_{1}x-a_{1}y}\right) +g_{_{2}}\left( {%
b_{2}x-a_{2}y}\right) \eqno(1.11)
\end{equation*}%
for arbitrary $g_{1}$ and $g_{2}.$ In \cite{36}, Golitschek and Light
described an algorithm that computes the error of approximation of a
continuous function $f\left( {x,y}\right) $ by solutions of
equation (1.10), provided that $a_{1}=b_{2}=1$, $a_{2}=b_{1}=0.$ Using our
result (see Theorem 1.3), one can characterize those solutions (1.11) that are
extremal to a given function $f(x,y)$. For a certain class of functions $f(x,y)
$, one can also easily calculate the approximation error and construct an
extremal solution (see Theorems 1.5 and 1.6 below).

The problem of approximating by functions from the set $\mathcal{R}\left(
\mathbf{a}^{1},\mathbf{a}^{2}\right) $ arises in other contexts too. Buck \cite{11}
studied the classical functional equation: given $\beta (t)\in C[0,1]$, $%
0\leq \beta (t)\leq 1$, for which $u\in C[0,1]$ does there exist $\varphi
\in C[0,1]$ such that
\begin{equation*}
\varphi (t)=\varphi \left( \beta (t)\right) +u(t)?
\end{equation*}%
He proved that the set of all $u$ satisfying this condition is dense in the
set
\begin{equation*}
\{v\in C[0,1]:\ v(t)=0\ \mbox{whenever}\ \beta (t)=t\}
\end{equation*}%
if and only if $\mathcal{R}\left( \mathbf{a}^{1},\mathbf{a}^{2}\right) $ with the
unit directions $\mathbf{a}^{1}=(1;0)$ and $\mathbf{a}^{2}=(0,1)$ is dense in $C(K)$%
, where $K=\{(x,y):y=x\ \mbox{or}\ y=\beta (x),\ 0\leq x\leq 1\}$.

Although there are enough reasons to consider approximation problems
associated with the set $\mathcal{R}\left( \mathbf{a}^{1},\mathbf{a}^{2}\right) $
in an independent way, one may ask why sums of only two ridge functions are
considered instead of sums with an arbitrary number of terms. We will try to
answer this fair question in Section 1.3.4.

\

\subsection{The characterization theorem}

Let $Q$ be a compact subset of $\mathbb{R}^{d}$ and $\mathbf{a}^{1},\mathbf{a}^{2}%
\in \mathbb{R}^{d}\backslash {\left\{ \mathbf{0}\right\} }.$

\bigskip

\textbf{Definition 1.3.} \textit{A finite or infinite ordered set $p=\left(
\mathbf{p}{_{1},\mathbf{p}_{2},...}\right) \subset Q$ with $\mathbf{p}%
_{i}\neq \mathbf{p}_{i+1},$ and either $\mathbf{a}^{1}\cdot \mathbf{p}_{1}=%
\mathbf{a}^{1}\cdot \mathbf{p}_{2},\mathbf{a}^{2}\cdot \mathbf{p}_{2}=\mathbf{a}^{2}%
\cdot \mathbf{p}_{3},\mathbf{a}^{1}\cdot \mathbf{p}_{3}=\mathbf{a}^{1}\cdot \mathbf{p%
}_{4},...$ or $\mathbf{a}^{2}\cdot \mathbf{p}_{1}=\mathbf{a}^{2}\cdot \mathbf{p}%
_{2},~\mathbf{a}^{1}\cdot \mathbf{p}_{2}=\mathbf{a}^{1}\cdot \mathbf{p}_{3},
\mathbf{a}^{2}\cdot \mathbf{p}_{3}=\mathbf{a}^{2}\cdot \mathbf{p}_{4},...$is called a path
with respect to the directions $\mathbf{a}^{1}$ and $\mathbf{a}^{2}$.}

\bigskip

This notion (in the two-dimensional case) was introduced by Braess and Pinkus
\cite{10}. They showed that paths give geometric means of deciding if a set
of points ${\left\{ {\mathbf{x}}^{i}\right\} }_{i=1}^{m}\subset \mathbb{R}%
^{2}$ has the \textit{NI} property (see Section 1.2.4). Ismailov and Pinkus
\cite{63} used these objects to study the problem of interpolation on
straight lines by linear combinations of a finite number of ridge functions
with fixed directions. In \cite{51,53,44} paths were generalized to those
with respect to two functions. The last objects turned out to be useful in
problems of approximation and representation by sums of compositions of
fixed multivariate functions with univariate functions.

If $\mathbf{a}^{1}$ and $\mathbf{a}^{2}$ are the coordinate vectors in $\mathbb{R}%
^{2}$, then Definition 1.3 defines a \textit{bolt of lightning}. The idea of
bolts was first introduced in Diliberto and Straus \cite{26}, where these
objects are called \textit{permissible lines}. They appeared further in a number
of papers, although under several different names (see, e.g., \cite%
{29,34,36,55,56,79,78,76,82,93,108,107,113}). Note that the term
\textquotedblleft bolt of lightning" is due to Arnold \cite{3}.

For the sake of brevity, we use the term ``path" instead of the long expression
``path with respect to the directions $\mathbf{a}^{1}$ and $\mathbf{a}^{2}$".

The length of a path is the number of its points. A single point is a path
of the unit length. A finite path $\left( \mathbf{p}_{1} ,\mathbf{p}_{2}
,...,\mathbf{p}_{2n} \right)$ is said to be closed if $\left(\mathbf{p}_{1} ,%
\mathbf{p}_{2} ,...,\mathbf{p}_{2n}, \mathbf{p}_{1}\right)$ is a path.

We associate each closed path $p=\left(\mathbf{p}_{1}, \mathbf{p}_{2} ,...,%
\mathbf{p}_{2n} \right) $ with the functional
\begin{equation*}
G_{p} (f)=\frac{1}{2n} \sum\limits_{k=1}^{2n}(-1)^{k+1} f(\mathbf{p}_{k}).
\end{equation*}

This functional has the following obvious properties:

(a) If $g\in \mathcal{R}\left( \mathbf{a}^{1},\mathbf{a}^{2}\right)$, then $G_{p}
(g)=0$.

(b) $\left\| G_{p} \right\| \leq 1$ and if $\mathbf{p}_{i} \neq \mathbf{p}_{j}$ for all $%
i\neq j,$ $1\leq i,j\leq 2n$ , then $\left\| G_{p} \right\| =1$.

\bigskip

\textbf{Lemma 1.3.} \textit{Let a compact set $Q$ have closed paths. Then
\begin{equation*}
\sup\limits_{p\subset Q}\left\vert G_{p}(f)\right\vert \leq E\left( f\right)
,\eqno(1.12)
\end{equation*}%
where the sup is taken over all closed paths. Moreover, inequality (1.12) is
sharp, i.e. there exist functions for which (1.12) turns into equality.}

\bigskip

\begin{proof} Let $p$ be a closed path in $Q$ and $g$ be any function from $%
\mathcal{R}\left( \mathbf{a}^{1},\mathbf{a}^{2}\right) $. By the linearity of
$G_{p}$ and properties (a) and (b),
\begin{equation*}
\left\vert G_{p}(f)\right\vert =\left\vert G_{p}(f-g)\right\vert \leq
\left\Vert f-g\right\Vert .\eqno(1.13)
\end{equation*}%
Since the left-hand and the right-hand sides of (1.13) do not depend on $g$
and $p$ respectively, it follows from (1.13) that
\begin{equation*}
\sup_{p\subset Q}\left\vert G_{p}(f)\right\vert \leq \inf_{g \in \mathcal{R}\left(
\mathbf{a}^{1},\mathbf{a}^{2}\right) }\left\Vert f-g\right\Vert .\eqno(1.14)
\end{equation*}

Now we prove the sharpness of (1.12). By assumption $Q$ has closed paths.
Then $Q$ has a closed path $p^{\prime }=\left( \mathbf{p}_{1}^{\prime },...,%
\mathbf{p}_{2m}^{\prime }\right)$ with distinct points $\mathbf{p}_{1}^{\prime },...,%
\mathbf{p}_{2m}^{\prime}$. In fact, such a special path can be obtained
from any closed path $p=\left( \mathbf{p}_{1},...,\mathbf{p}_{2n}\right) $
by the following simple algorithm: if the points of the path $p$ are not all
distinct, let $i$ and $k>0$ be the minimal indices such that $\mathbf{p}_{i}=%
\mathbf{p}_{i+2k}$; delete from $p$ the subsequence $\mathbf{p}_{i+1},...,%
\mathbf{p}_{i+2k}$ and call $p$ the obtained path; repeat the above step
until all points of $p$ are all distinct; set $p^{\prime }:=p$. On the other
hand there exist continuous functions $h=h(\mathbf{x})$ on $Q$ such that $h(%
\mathbf{p}_{i}^{\prime })=1$, $i=1,3,...,2m-1$, $h(\mathbf{p}_{i}^{\prime
})=-1$, $i=2,4,...,2m$ and $-1<h(\mathbf{x})<1$ elsewhere. For such
functions we have
\begin{equation*}
G_{p^{\prime }}(h)=\Vert h\Vert =1\eqno(1.15)
\end{equation*}%
and
\begin{equation*}
E(h)\leq \Vert h\Vert ,\eqno(1.16)
\end{equation*}%
where the last inequality follows from the fact that $0\in \mathcal{R}\left(
\mathbf{a}^{1},\mathbf{a}^{2}\right) .$ From (1.14)-(1.16) it follows that
\begin{equation*}
\sup_{p\subset Q}\left\vert G_{p}(h)\right\vert =E\left( h\right) .
\end{equation*}%
\end{proof}

\textbf{Lemma 1.4.} \textit{Let $Q$ be a convex compact subset of $%
\mathbb{R}^{d}$ and $f \in C(Q)$. For a vector $\mathbf{e}\in
\mathbb{R}^{d}\backslash \mathbf{\{0\}}$ and a real number $t$ set
\begin{equation*}
Q_{t}=\left\{ {\mathbf{x}}\in Q:\mathbf{e}\cdot \mathbf{x}=t\right\} ,\ \ \
T_{h}=\left\{ t\in \mathbb{R}:Q_{t}\neq \emptyset \right\} .
\end{equation*}%
The functions
\begin{equation*}
g_{1} (t)=\max_{\mathbf{x} \in Q_t} f(\mathbf{x} ),\ \ t\in T_{h}\ \ %
\mbox{and}\ \ g_{2} (t)=\min\limits_{\mathbf{x} \in Q_t} f(\mathbf{x} ),\ \
\ t\in T_{h}
\end{equation*}
are defined and continuous on $T_{h} $.}

\bigskip

The proof of this lemma is not difficult and can be obtained by the
well-known elementary methods of mathematical analysis.

\bigskip

\textbf{Definition 1.4.} \textit{A finite or infinite path $(\mathbf{p}_{1},%
\mathbf{p}_{2},...)$ is said to be extremal for a function $u \in
C(Q)$ if $u(\mathbf{p}_{i})=(-1)^{i}\left\Vert u\right\Vert ,i=1,2,...$ or $%
u(\mathbf{p}_{i})=(-1)^{i+1}\left\Vert u\right\Vert ,$ $i=1,2,...$}

\bigskip

\textbf{Theorem 1.3.} \textit{Let $Q\subset \mathbb{R}^{d}$ be a convex
compact set satisfying the following condition}

\textit{Condition (A): For any path $q=(\mathbf{q}_{1},\mathbf{q}_{2},...,%
\mathbf{q}_{n})\subset Q$ there exist points $\mathbf{q}_{n+1},\mathbf{q}%
_{n+2},...,\mathbf{q}_{n+s}\in Q$ such that $(\mathbf{q}_{1},\mathbf{q}%
_{2},...,\mathbf{q}_{n+s})$ is a closed path and $s$ is not more than some
positive integer $N_{0}$ independent of $q$.}

\textit{Then a necessary and sufficient condition for a function $g_{0}\in
\mathcal{R}\left( \mathbf{a}^{1},\mathbf{a}^{2}\right) $ to be an extremal element
to the given function $f \in C(Q)$ is the existence of a closed or infinite
path $l=(\mathbf{p}_{1},\mathbf{p}_{2},...)$ extremal for the function $%
f_{1}=f-g_{0}$.}

\bigskip

It should be remarked that the above condition (A) strongly
depends on the fixed directions $\mathbf{a}^{1}$ and $\mathbf{a}^{2}$. For example,
in the familiar case of a square $S\subset \mathbb{R}^2$ there are many
directions which are not allowed. If it is possible to reach a corner of $S$
with not more than one of the two directions orthogonal to $\mathbf{a}^{1}
$ and $\mathbf{a}^{2}$, respectively (we don't differentiate between directions $%
\mathbf{c}$ and $-\mathbf{c}$), the triple $(S,\mathbf{a}^{1}, \mathbf{a}^{2})$ does
not satisfy condition (A) of the theorem. Here are simple examples: Let $%
S=[0;1]^2$, $\mathbf{a}^{1}=(1;0)$, $\mathbf{a}^{2}=(1;1)$. Then the ordered set $%
\{(0;1), (1;0), (1;1)\}$ is a path in $S$ which can not be made closed. In
this case, $(1;1)$ is not reached with the direction orthogonal to $\mathbf{b%
}$. Let now $\mathbf{a}^{1}=\left(1;\frac{1}{2}\right)$, $\mathbf{a}^{2}=(1;1)$.
Then the corner $(1;1)$ is reached with none of the directions orthogonal to
$\mathbf{a}^{1}$ and $\mathbf{a}^{2}$ respectively. In this case, for any positive
integer $N_0$ and any point $\mathbf{q}_0$ in $S$ one can chose a point $%
\mathbf{q}_1\in S$ from a sufficiently small neighborhood of the corner $%
(1;1)$ so that any path containing $\mathbf{q}_0$ and $\mathbf{q}_1$ has the
length more than $N_0$. These examples and a little geometry show that if a convex
compact set $Q\subset\mathbb{R}^2$ satisfies condition (A) of Theorem
1.3, then any point in the boundary of $Q$ must be reached with each of the
two directions orthogonal to $\mathbf{a}^{1}$ and $\mathbf{a}^{2}$ respectively. If $%
Q\subset \mathbb{R}^d, \mathbf{a}^{1}, \mathbf{a}^{2}\in \mathbb{R}^d\backslash\{%
\mathbf{0}\}$, $d>2$, there are many directions orthogonal to $\mathbf{a}^{1}$
and $\mathbf{a}^{2} $. In this case, condition (A) requires that any point in
the boundary of $Q$ should be reached with at least two directions
orthogonal to $\mathbf{a}^{1}$ and $\mathbf{a}^{2}$, respectively.

\begin{proof} \textit{Necessity.} Let $g_{0}(\mathbf{x}) =g_{1,0} \left( \mathbf{a}^{1}{\cdot
}\mathbf{x}\right) +g_{2,0} \left(\mathbf{a}^{2}{\cdot }\mathbf{x}\right)$ be an
extremal element from $\mathcal{R}\left( \mathbf{a}^{1},\mathbf{a}^{2} \right)$ to
$f$. We must show that if there is not a closed path extremal for $f_{1} $,
then there exists a path extremal for $f_{1} $ with the infinite length
(number of points). Suppose the contrary. Suppose that there exists a
positive integer $N$ such that the length of each path extremal for $f_{1} $
is not more than $N$. Set the following functions:
\begin{equation*}
f_{n} =f_{n-1} -g_{1,n-1} -g_{2,n-1} ,\ \ n=2, 3, ...,
\end{equation*}
where
\begin{equation*}
g_{1,n-1} =g_{1,n-1} \left(\mathbf{a}^{1}{\cdot }\mathbf{x}\right)=\frac{1}{2}
\left(\max\limits_{\substack{ \mathbf{y}\in Q  \\ \mathbf{a}^{1}{\cdot }\mathbf{y%
}=\mathbf{a}^{1}{\cdot }\mathbf{x}}} f_{n-1}(\mathbf{y}) +\min\limits_{\substack{
\mathbf{y}\in Q  \\ \mathbf{a}^{1}{\cdot }\mathbf{y}=\mathbf{a}^{1}{\cdot }\mathbf{x}
}} f_{n-1}(\mathbf{y})\right)
\end{equation*}
\begin{equation*}
g_{2,n-1} =g_{2,n-1} (\mathbf{a}^{2}{\cdot }\mathbf{x})=\frac{1}{2} \left( \max
_{\substack{ \mathbf{y}\in Q  \\ \mathbf{a}^{2}{\cdot }\mathbf{y}=\mathbf{a}^{2}{%
\cdot }\mathbf{x}}} \left( f_{n-1} (\mathbf{y})-g_{1,n-1}(\mathbf{a}^{1}{\cdot }%
\mathbf{y})\right)\right.
\end{equation*}
\begin{equation*}
\left.+\min\limits_{\substack{ \mathbf{y}\in Q  \\ \mathbf{a}^{2}{\cdot }\mathbf{%
y}=\mathbf{a}^{2}{\cdot }\mathbf{x}}} \left( f_{n-1}(\mathbf{y}) -g_{1,n-1}(%
\mathbf{a}^{1}{\cdot }\mathbf{y}) \right) \right).
\end{equation*}

By Lemma 1.4, all the functions $f_{n}(\mathbf{x}),$ $n=2,3,...,$ are
continuous on $Q$. By assumption $g_{0}$ is a best approximation to $f$.
Hence $\left\Vert f_{1}\right\Vert =E\left( f\right) $. Now let us show that $%
\left\Vert f_{2}\right\Vert =E\left( f\right) $. Indeed, for any $\mathbf{x}%
\in Q$
\begin{equation*}
f_{1}(\mathbf{x})-g_{1,1}(\mathbf{a}^{1}{\cdot }\mathbf{x})\leq \frac{1}{2}%
\left( \max\limits_{\substack{ \mathbf{y}\in Q  \\ \mathbf{a}^{1}{\cdot }\mathbf{%
y}=\mathbf{a}^{1}{\cdot }\mathbf{x}}}f_{1}(\mathbf{y})-\min\limits_{\substack{
\mathbf{y}\in Q  \\ \mathbf{a}^{1}{\cdot }\mathbf{y}=\mathbf{a}^{1}{\cdot }\mathbf{x}
}}f_{1}(\mathbf{y})\right) \leq E(f)\eqno(1.17)
\end{equation*}%
and
\begin{equation*}
f_{1}(\mathbf{x})-g_{1,1}(\mathbf{a}^{1}{\cdot }\mathbf{x})\geq \frac{1}{2}%
\left( \min\limits_{\substack{ \mathbf{y}\in Q  \\ \mathbf{a}^{1}{\cdot }\mathbf{%
y}=\mathbf{a}^{1}{\cdot }\mathbf{x}}}f_{1}(\mathbf{y})-\max\limits_{\substack{
\mathbf{y}\in Q  \\ \mathbf{a}^{1}{\cdot }\mathbf{y}=\mathbf{a}^{1}{\cdot }\mathbf{x}
}}f_{1}(\mathbf{y})\right) \geq -E(f).\eqno(1.18)
\end{equation*}

Using the definition of $g_{2,1}(\mathbf{a}^{2}\cdot \mathbf{x})$, for any $%
\mathbf{x}\in Q$ we have
\begin{equation*}
f_{1}(\mathbf{x})-g_{1,1}(\mathbf{a}^{1}\cdot \mathbf{x})-g_{2,1}(\mathbf{a}^{2}%
\cdot \mathbf{x})
\end{equation*}
\begin{equation*}
\leq \frac{1}{2}\left( \max\limits_{\substack{ \mathbf{y}\in Q  \\ \mathbf{a}^{2}%
\cdot \mathbf{y}=\mathbf{a}^{2}\cdot \mathbf{x}}}\left( f_{1}(\mathbf{y}%
)-g_{1,1}(\mathbf{a}^{1}\cdot \mathbf{y})\right) -\min\limits_{ _{\substack{
\mathbf{y}\in Q  \\ \mathbf{a}^{2}\cdot \mathbf{y}=\mathbf{a}^{2}\cdot \mathbf{x}}}%
}\left( f_{1}(\mathbf{y})-g_{1,1}(\mathbf{a}^{1}\cdot \mathbf{y})\right) \right)
\end{equation*}
and
\begin{equation*}
f_{1}(\mathbf{x})-g_{1,1}(\mathbf{a}^{1}\cdot \mathbf{x})-g_{2,1}(\mathbf{a}^{2}%
\cdot \mathbf{x})
\end{equation*}
\begin{equation*}
\leq \frac{1}{2}\left( \min\limits_{_{\substack{ \mathbf{y}\in Q  \\ \mathbf{%
a}^{2}\cdot \mathbf{y}=\mathbf{a}^{2}\cdot \mathbf{x}}}}\left( f_{1}(\mathbf{y}%
)-g_{1,1}(\mathbf{a}^{1}\cdot \mathbf{y})\right) -\max\limits_{\substack{
\mathbf{y}\in Q  \\ \mathbf{a}^{2}\cdot \mathbf{y}=\mathbf{a}^{2}\cdot \mathbf{x}}}%
\left( f_{1}(\mathbf{y})-g_{1,1}(\mathbf{a}^{1}\cdot \mathbf{y})\right) \right).
\end{equation*}

Using (1.17) and (1.18) in the last two inequalities, we obtain that for any
$\mathbf{x}\in Q$
\begin{equation*}
-E(f)\leq f_{2}(\mathbf{x})=f_{1}(\mathbf{x})-g_{1,1}(\mathbf{a}^{1}{\cdot }%
\mathbf{x})-g_{2,1}(\mathbf{a}^{2}{\cdot }\mathbf{x})\leq E(f).
\end{equation*}

Therefore,
\begin{equation*}
\left\Vert f_{2}\right\Vert \leq E(f).\eqno(1.19)
\end{equation*}

Since $f_{2}(\mathbf{x})-f(\mathbf{x})$ belongs to $\mathcal{R}\left(
\mathbf{a}^{1},\mathbf{a}^{2}\right) $, we deduce from (1.19) that
\begin{equation*}
\left\Vert f_{2}\right\Vert =E(f).
\end{equation*}

By the same way, one can show that $\|f_3\|=E(f)$, $\|f_4\|=E(f)$, and so
on. Thus we can write
\begin{equation*}
\left\| f_{n} \right\| =E(f), \ \mbox{for any}\ n .
\end{equation*}

Let us now prove the implications
\begin{equation*}
f_{1}(\mathbf{p}_{0})<E(f)\Rightarrow f_{2}(\mathbf{p}_{0})<E(f)\eqno(1.20)
\end{equation*}%
and
\begin{equation*}
f_{1}(\mathbf{p}_{0})>-E(f)\Rightarrow f_{2}(\mathbf{p}_{0})>-E(f),\eqno%
(1.21)
\end{equation*}%
where $\mathbf{p}_{0}\in Q$. First, we are going to prove the implication
\begin{equation*}
f_{1}(\mathbf{p}_{0})<E(f)\Rightarrow f_{1}(\mathbf{p}_{0})-g_{1,1}(\mathbf{a%
}\cdot \mathbf{p}_{0})<E(f).\eqno(1.22)
\end{equation*}

There are two possible cases.

1) $\max\limits_{\substack{ \mathbf{y}\in Q  \\ \mathbf{a}^{1}{\cdot }\mathbf{y}=%
\mathbf{a}^{1}{\cdot }\mathbf{p}_0}} f_{1} (\mathbf{y} )=E(f)$ and $\min\limits
_{\substack{ \mathbf{y}\in Q  \\ \mathbf{a}^{1}{\cdot }\mathbf{y}=\mathbf{a}^{1}{%
\cdot }\mathbf{p}_0}} f_{1} (\mathbf{y}) = -E(f). $ In this case, $g_{1,1}(%
\mathbf{a}^{1}\cdot \mathbf{p}_0)=0$. Hence
\begin{equation*}
f_1(\mathbf{p}_0)-g_{1,1}(\mathbf{a}^{1}\cdot \mathbf{p}_0)< E(f).
\end{equation*}

2) $\max\limits_{\substack{ \mathbf{y}\in Q  \\ \mathbf{a}^{1}{\cdot }\mathbf{y}=%
\mathbf{a}^{1}{\cdot }\mathbf{p}_{0}}}f_{1}(\mathbf{y})=E(f)-\varepsilon _{1}$
and $\min\limits_{\substack{ \mathbf{y}\in Q  \\ \mathbf{a}^{1}{\cdot }\mathbf{y}%
=\mathbf{a}^{1}{\cdot }\mathbf{p}_{0}}}f_{1}(\mathbf{y})=-E(f)+\varepsilon _{2}$,%
\newline
where $\varepsilon _{1}$, $\varepsilon _{2}$ are nonnegative real numbers
with the sum $\varepsilon _{1}+\varepsilon _{2}\not=0$. In this case,
\begin{eqnarray*}
f_{1}(\mathbf{p}_{0})-g_{1,1}(\mathbf{a}^{1}{\cdot }\mathbf{p}_{0}) &\leq
&\max\limits_{\substack{ \mathbf{y}\in Q  \\ \mathbf{a}^{1}{\cdot }\mathbf{y}=%
\mathbf{a}^{1}{\cdot }\mathbf{p}_{0}}}f_{1}(\mathbf{y})-g_{1,1}(\mathbf{a}^{1}{\cdot
}\mathbf{p}_{0})= \\
&=&\frac{1}{2}\left( \max\limits_{\substack{ \mathbf{y}\in Q  \\ \mathbf{a}^{1}{%
\cdot }\mathbf{y}=\mathbf{a}^{1}{\cdot }\mathbf{p}_{0}}}f_{1}(\mathbf{y}%
)-\min\limits_{\substack{ \mathbf{y}\in Q  \\ \mathbf{a}^{1}{\cdot }\mathbf{y}=%
\mathbf{a}^{1}{\cdot }\mathbf{p}_{0}}}f_{1}(\mathbf{y})\right) =
\end{eqnarray*}%
\begin{equation*}
=E(f)-\frac{\varepsilon _{1}+\varepsilon _{2}}{2}<E(f).
\end{equation*}

Thus we have proved (1.22). Using this method, we can also prove that
\begin{equation*}
f_{1}(\mathbf{p}_{0})-g_{1,1}(\mathbf{a}^{1}{\cdot }\mathbf{p}%
_{0})<E(f)\Rightarrow f_{1}(\mathbf{p}_{0})-g_{1,1}(\mathbf{a}^{1}{\cdot }%
\mathbf{p}_{0})-g_{2,1}(\mathbf{a}^{2}{\cdot }\mathbf{p}_{0})<E(f).\eqno(1.23)
\end{equation*}%
Now (1.20) follows from (1.22) and (1.23). By the same way we can prove
(1.21). It follows from implications (1.20) and (1.21) that if $f_{2}(%
\mathbf{p}_{0})=E(f)$, then $f_{1}(\mathbf{p}_{0})=E(f)$ and if $f_{2}(%
\mathbf{p}_{0})=-E(f)$, then $f_{1}(\mathbf{p}_{0})=-E(f)$. This simply
means that each path extremal for $f_{2}$ will be extremal for $f_{1}$.

Now we show that if any path extremal for $f_{1} $ has the length not more
than $N$, then any path extremal for $f_{2} $ has the length not more than $%
N-1$. Suppose the contrary. Suppose that there is a path extremal for $f_{2}
$ with the length equal to $N$. Denote it by $q=(\mathbf{q}_{1} ,\mathbf{q}%
_{2} ,...,\mathbf{q}_{N} )$. Without loss of generality we may assume that $%
\mathbf{a}^{2}\cdot \mathbf{q}_{N-1} =\mathbf{a}^{2}\cdot \mathbf{q}_{N}$. As it has
been shown above, the path $q$ is also extremal for $f_{1}$. Assume that $%
f_{1} (\mathbf{q}_{N} )=E(f)$. Then there is not any $\mathbf{q}_{0} \in Q$
such that $\mathbf{q}_{0} \neq \mathbf{q}_{N} $, $\mathbf{a}^{1}\cdot \mathbf{q}%
_{0} =\mathbf{a}^{1}\cdot \mathbf{q}_{N}$ and $f_{1} (\mathbf{q}_{0} )=-E(f)$.
Indeed, if there was such $\mathbf{q}_{0}$ and $\mathbf{q}_{0} \not\in q$,
the path $(\mathbf{q}_{1} ,\mathbf{q}_{2} ,...,\mathbf{q}_{N} ,\mathbf{q}%
_{0} )$ would be extremal for $f_{1} $. But this would contradict our
assumption that any path extremal for $f_{1} $ has the length not more than $%
N$. Besides, if there was such $\mathbf{q}_{0} $ and $\mathbf{q}_{0} \in q$,
we could form some closed path extremal for $f_{1} $. This also would
contradict our assumption that there does not exist a closed path extremal
for $f_{1} $.

Hence
\begin{equation*}
\max\limits_{\substack{ \mathbf{y}\in Q  \\ \mathbf{a}^{1}{\cdot }\mathbf{y}=%
\mathbf{a}^{1}{\cdot }\mathbf{q}_N}} f_{1} (\mathbf{y} )=E(f),\ \ \min\limits
_{\substack{ \mathbf{y}\in Q  \\ \mathbf{a}^{1}{\cdot }\mathbf{y}=\mathbf{a}^{1}{%
\cdot }\mathbf{q}_N}} f_{1} (\mathbf{y})>-E(f).
\end{equation*}
Therefore,
\begin{equation*}
\left| f_{1} (\mathbf{q}_{N} )-g_{1,1} (\mathbf{a}^{1}{\cdot }\mathbf{q}%
_N)\right| <E(f).
\end{equation*}

From the last inequality it is easy to obtain that (see the proof of
implications (1.20) and (1.21))
\begin{equation*}
\left\vert f_{2}(\mathbf{q}_{N})\right\vert <E(f).
\end{equation*}%
This means, on the contrary to our assumption, that the path $(\mathbf{q}%
_{1},\mathbf{q}_{2},...,\mathbf{q}_{N})$ can not be extremal for $f_{2}$.
Hence any path extremal for $f_{2}$ has the length not more than $N-1$.

By the same way, it can be shown that any path extremal for $f_{3}$ has the
length not more than $N-2$, any path extremal for $f_{4}$ has the length not
more than $N-3$ and so on. Finally, we will obtain that there is not a path
extremal for $f_{N+1}$. Hence there is not a point $\mathbf{p}_{0}\in Q$
such that $|f_{N+1}(\mathbf{p}_{0})|=\Vert f_{N+1}\Vert $. But by Lemma 1.4,
all the functions $f_{2}$, $f_{3},...,f_{N+1}$ are continuous on the compact
set $Q$; hence the norm $\Vert f_{N+1}\Vert $ must be attained. This
contradiction means that there exists a path extremal for $f_{1}$ with the
infinite length.

\bigskip

\textit{Sufficiency.} Let a path $p=(\mathbf{p}_{1},\mathbf{p}_{2},...,%
\mathbf{p}_{2n})$ be closed and extremal for $f_{1}$. Then
\begin{equation*}
\left\vert G_{p}(f)\right\vert =\left\Vert f-g_{0}\right\Vert .\eqno(1.24)
\end{equation*}

By Lemma 1.3,
\begin{equation*}
\left\vert G_{p}(f)\right\vert \leq E(f).\eqno(1.25)
\end{equation*}

It follows from (1.24) and (1.25) that $g_{0}$ is a best approximation.

Let now a path $p=(\mathbf{p}_{1},\mathbf{p}_{2},...,\mathbf{p}_{n},...)$ be
infinite and extremal for $f_{1}$. Consider the sequence $p_{n}=(\mathbf{p}%
_{1},\mathbf{p}_{2},...,\mathbf{p}_{n})$, $n=1,2,...,$ of finite paths. By
condition (A) of the theorem, for each $p_{n}$
there exists a closed path $p_{n}^{m_{n}}=(\mathbf{p}_{1},\mathbf{p}_{2},...,%
\mathbf{p}_{n},\mathbf{q}_{n+1},...,\mathbf{q}_{n+m_{n}})$, where $m_{n}\leq
N_{0}$. Then for any positive integer $n$,
\begin{equation*}
\left\vert G_{p_{n}^{m_{n}}}(f)\right\vert =\left\vert
G_{p_{n}^{m_{n}}}(f-g_{0})\right\vert \leq \frac{n\left\Vert
f-g_{0}\right\Vert +m_{n}\left\Vert f-g_{0}\right\Vert }{n+m_{n}}=\left\Vert
f-g_{0}\right\Vert
\end{equation*}%
and
\begin{equation*}
\left\vert G_{p_{n}^{m_{n}}}(f)\right\vert \geq \frac{n\left\Vert
f-g_{0}\right\Vert -m_{n}\left\Vert f-g_{0}\right\Vert }{n+m_{n}}=\frac{%
n-m_{n}}{n+m_{n}}\left\Vert f-g_{0}\right\Vert .
\end{equation*}

It follows from the above two inequalities for $\left\vert G_{p_{n}^{m_{n}}}(f)\right\vert$ that
\begin{equation*}
\sup_{p_{n}^{m_{n}}}\left\vert G_{p_{n}^{m_{n}}}(f)\right\vert =\left\Vert
f-g_{0}\right\Vert.
\end{equation*}
This together with Lemma 1.3 give that
\begin{equation*}
\Vert f-g_{0}\Vert \leq E(f).
\end{equation*}%
Hence $g_{0}$ is a best approximation.
\end{proof}

\bigskip

Theorem 1.3 has been proved by using only methods of classical analysis. By
implementing more deep techniques from functional analysis we will see below
that condition (A) and the convexity assumption on a compact set $Q$ can be
dropped.

\bigskip

\textbf{Theorem 1.4.}\ \textit{\ Assume $Q$ is a compact subset of $\mathbb{R%
}^{d}$. A function $g_{0}\in \mathcal{R}$ is a best approximation to a
function $f\in C(Q)$ if and only if there exists a closed or infinite path $%
p=(\mathbf{p}_{1},\mathbf{p}_{2},...)$ extremal for the function $f-g_{0}$.}

\bigskip

\begin{proof} \textit{Sufficiency.} There are two possible cases. The first case
happens when there exists a closed path $(\mathbf{p}_{1},...,\mathbf{p}%
_{2n}) $ extremal for the function $f-g_{0}.$ Let us check that in this
case, $f-g_{0}$ is a best approximation. Indeed, on the one hand, the
following equalities are valid
\begin{equation*}
\left\vert \sum_{i=1}^{2n}(-1)^{i}f(\mathbf{p}_{i})\right\vert =\left\vert
\sum_{i=1}^{2n}(-1)^{i}\left[ f-g_{0}\right] (\mathbf{p}_{i})\right\vert
=2n\left\Vert f-g_{0}\right\Vert .
\end{equation*}%
On the other hand, for any function $g\in \mathcal{R}$, we have
\begin{equation*}
\left\vert \sum_{i=1}^{2n}(-1)^{i}f(\mathbf{p}_{i})\right\vert =\left\vert
\sum_{i=1}^{2n}(-1)^{i}\left[ f-g\right] (\mathbf{p}_{i})\right\vert \leq
2n\left\Vert f-g\right\Vert .
\end{equation*}%
Therefore, $\left\Vert f-g_{0}\right\Vert \leq \left\Vert f-g\right\Vert $
for any $g\in \mathcal{R}$. That is, $g_{0}$ is a best approximation.

The second case happens when we do not have closed paths extremal for $%
f-g_{0}$, but there exists an infinite path $(\mathbf{p}_{1},\mathbf{p}%
_{2},...)$ extremal for $f-g_{0}$. To analyze this case, consider the
following linear functional
\begin{equation*}
L_{q}:C(Q)\rightarrow \mathbb{R}\text{, \ }L_{q}(F)=\frac{1}{n}%
\sum_{i=1}^{n}(-1)^{i}F(\mathbf{q}_{i}),
\end{equation*}%
where $q=\{\mathbf{q}_{1},...,\mathbf{q}_{n}\}$ is a finite path in $Q$. It
is easy to see that the norm $\left\Vert L_{q}\right\Vert \leq 1$ and $%
\left\Vert L_{q}\right\Vert =1$ if and only if the set of points of $q$ with
odd indices $O=\{\mathbf{q}_{i}\in q:$ $i$ \textit{is an odd number}$\}$ do
not intersect with the set of points of $q$ with even indices $E=\{\mathbf{q}%
_{i}\in q:$ $i$ \textit{is an even number}$\}$. Indeed, from the definition
of $L_{q}$ it follows that $\left\vert L_{q}(F)\right\vert \leq \left\Vert
F\right\Vert $ for all functions $F\in C(Q)$, whence $\left\Vert
L_{q}\right\Vert \leq 1.$ If $O\cap E=\varnothing $, then for a function $%
F_{0}$ with the property $F_{0}(\mathbf{q}_{i})=-1$ if $i$ is odd, $F_{0}(%
\mathbf{q}_{i})=1$ if $i$ is even and $-1<F_{0}(x)<1$ elsewhere on $Q,$ we
have $\left\vert L_{q}(F_{0})\right\vert =\left\Vert F_{0}\right\Vert .$
Hence, $\left\Vert L_{q}\right\Vert =1$. Recall that such a function $F_{0}$
exists on the basis of Urysohn's great lemma.

Note that if $q$ is a closed path, then $L_{q}$ annihilates all members of
the class $\mathcal{R}$. But in general, when $q$ is not closed, we do not
have the equality $L_{q}(g)=0,$ for all members $g\in \mathcal{R}$.
Nonetheless, this functional has the important property that
\begin{equation*}
\left\vert L_{q}(g_{1}+g_{2})\right\vert \leq \frac{2}{n}(\left\Vert
g_{1}\right\Vert +\left\Vert g_{2}\right\Vert ),\eqno(1.26)
\end{equation*}%
where $g_{1}$ and $g_{2}$ are ridge functions with the directions $\mathbf{a}%
_{1}$ and $\mathbf{a}_{2}$, respectively, that is, $g_{1}=g_{1}(\mathbf{a}%
_{1}\cdot \mathbf{x})$ and $g_{2}=g_{2}(\mathbf{a}_{2}\cdot \mathbf{x}).$
This property is important in the sense that if $n$ is sufficiently large,
then the functional $L_{q}$ is close to an annihilating functional. To prove
(1.26), note that $\left\vert L_{q}(g_{1})\right\vert \leq \frac{2}{n}%
\left\Vert g_{1}\right\Vert $ and $\left\vert L_{q}(g_{2})\right\vert \leq
\frac{2}{n}\left\Vert g_{2}\right\Vert $. These estimates become obvious if
consider the chain of equalities $g_{1}(\mathbf{a}_{1}\cdot \mathbf{q}%
_{1})=g_{1}(\mathbf{a}_{1}\cdot \mathbf{q}_{2}),$ $g_{1}(\mathbf{a}_{1}\cdot
\mathbf{q}_{3})=g_{1}(\mathbf{a}_{1}\cdot \mathbf{q}_{4}),...$(or $g_{1}(%
\mathbf{a}_{1}\cdot \mathbf{q}_{2})=g_{1}(\mathbf{a}_{1}\cdot \mathbf{q}%
_{3}),$ $g_{1}(\mathbf{a}_{1}\cdot \mathbf{q}_{4})=g_{1}(\mathbf{a}_{1}\cdot
\mathbf{q}_{5}),...$) for $g_{1}(\mathbf{a}_{1}\cdot \mathbf{x})$ and the
corresponding chain of equalities for $g_{2}(\mathbf{a}_{2}\cdot \mathbf{x})$%
.

Now consider the infinite path $p=(\mathbf{p}_{1},\mathbf{p}_{2},...)$ and
form the finite paths $p_{k}=(\mathbf{p}_{1},...,\mathbf{p}_{k}),$ $%
k=1,2,... $. For ease of notation, let us set $L_{k}=L_{p_{k}}.$ The
sequence $\{L_{_{k}}\}_{k=1}^{\infty }$ is a subset of the unit ball of the
conjugate space $C^{\ast }(Q).$ By the Banach-Alaoglu theorem, the unit ball
is weak$^{\text{*}}$ compact in the weak$^{\text{*}}$ topology of $C^{\ast
}(Q)$ (see \cite[p.68]{122}). It follows from this theorem that
the sequence $\{L_{_{k}}\}_{k=1}^{\infty }$ must have weak$^{\text{*}}$
cluster points. Suppose $L^{\ast }$ denotes one of them. Without loss of
generality we may assume that $L_{k}\overset{weak^{\ast }}{\longrightarrow }%
L^{\ast },$ as $k\rightarrow \infty .$ From (1.26) it follows that $L^{\ast
}(g_{1}+g_{2})=0.$ That is, $L^{\ast }\in \mathcal{R}^{\bot },$ where the
symbol $\mathcal{R}^{\bot }$ stands for the annihilator of $\mathcal{R}$.
Since in addition $\left\Vert L^{\ast }\right\Vert \leq 1,$ we can write that
\begin{equation*}
\left\vert L^{\ast }(f)\right\vert =\left\vert L^{\ast }(f-g)\right\vert
\leq \left\Vert f-g\right\Vert ,\eqno(1.27)
\end{equation*}%
for all functions $g\in \mathcal{R}.$ On the other hand, since the infinite
bolt $p$ is extremal for $f-g_{0}$
\begin{equation*}
\left\vert L_{k}(f-g_{0})\right\vert =\left\Vert f-g_{0}\right\Vert ,\text{ }%
k=1,2,...
\end{equation*}%
Therefore,
\begin{equation*}
\left\vert L^{\ast }(f)\right\vert =\left\vert L^{\ast }(f-g_{0})\right\vert
=\left\Vert f-g_{0}\right\Vert .\eqno(1.28)
\end{equation*}%
From (1.27) and (1.28) we conclude that
\begin{equation*}
\left\Vert f-g_{0}\right\Vert \leq \left\Vert f-g\right\Vert ,
\end{equation*}%
for all $g\in \mathcal{R}.$ In other words, $g_{0}$ is a best approximation
to $f$. We proved the sufficiency of the theorem.

\textit{Necessity.} The proof of this part is mainly based on the following
result of Singer \cite{S}: Let $X$ be a compact space, $U$
be a linear subspace of $C(X)$, $f\in C(X)\backslash U$ and $u_{0}\in U.$
Then $u_{0}$ is a best approximation to $f$ if and only if there exists a
regular Borel measure $\mu $ on $X$ such that

(1) The total variation $\left\Vert \mu \right\Vert =1$;

(2) $\mu $ is orthogonal to the subspace $U$, that is, $\int_{X}ud\mu =0$
for all $u\in U$;

(3) For the Jordan decomposition $\mu =\mu ^{+}-\mu ^{-}$,
\begin{equation*}
f(x)-u_{0}(x)=\left\{
\begin{array}{c}
\left\Vert f-u_{0}\right\Vert \text{ for }x\in S^{+}, \\
-\left\Vert f-u_{0}\right\Vert \text{ for }x\in S^{-},%
\end{array}%
\right.
\end{equation*}%
where $S^{+}$ and $S^{-}$ are closed supports of the positive measures $\mu
^{+}$ and $\mu ^{-}$, respectively.

Let us show how we use this theorem in the proof of necessity part of our
theorem. Assume $g_{0}\in \mathcal{R}$ is a best approximation. For the
subspace $\mathcal{R},$ the existence of a measure $\mu $ satisfying the
conditions (1)-(3) is a direct consequence of Singer's result. Let $\mathbf{x}%
_{0}$ be any point in $S^{+}.$\ Consider the point $y_{0}=\mathbf{a}%
_{1}\cdot \mathbf{x}_{0}$ and a $\delta $-neighborhood of $y_{0}$. That is,
choose an arbitrary $\delta >0$ and consider the set $I_{\delta
}=(y_{0}-\delta ,y_{0}+\delta )\cap \mathbf{a}_{1}\cdot Q.$ Here, $\mathbf{a}%
_{1}\cdot Q=\{\mathbf{a}_{1}\cdot \mathbf{x}:$ $\mathbf{x}\in Q\}.$ For any
subset $E\subset \mathbb{R}$, put
\begin{equation*}
E^{i}=\{\mathbf{x}\in Q:\mathbf{a}_{i}\cdot \mathbf{x}\in E\},\text{ }i=1,2.%
\text{ }
\end{equation*}

Clearly, for some sets $E,$ one or both the sets $E^{i}$ may be empty. Since
$I_{\delta }^{1}\cap S^{+}$ is not empty (note that $\mathbf{x}_{0}\in
I_{\delta }^{1}$), it follows that $\mu ^{+}(I_{\delta }^{1})>0.$ At the
same time $\mu (I_{\delta }^{1})=0,$ since $\mu $ is orthogonal to all
functions $g_{1}(\mathbf{a}_{1}\cdot \mathbf{x}).$ Therefore, $\mu
^{-}(I_{\delta }^{1})>0.$ We conclude that $I_{\delta }^{1}\cap S^{-}$ is
not empty. Denote this intersection by $A_{\delta }.$ Tending $\delta $ to $%
0,$ we obtain a set $A$ which is a subset of $S^{-}$ and has the property
that for each $\mathbf{x}\in A,$ we have $\mathbf{a}_{1}\cdot \mathbf{x}=%
\mathbf{a}_{1}\cdot \mathbf{x}_{0}.$ Fix any point $\mathbf{x}_{1}\in A$.
Changing $\mathbf{a}_{1}$, $\mu ^{+}$, $S^{+}$ to $\mathbf{a}_{2}$, $\mu
^{-} $ and $S^{-}$ correspondingly, repeat the above process with the point $%
y_{1}=\mathbf{a}_{2}\cdot \mathbf{x}_{1}$ and a $\delta $-neighborhood of $%
y_{1}$. Then we obtain a point $\mathbf{x}_{2}\in S^{+}$ such that $\mathbf{a%
}_{2}\cdot \mathbf{x}_{2}=\mathbf{a}_{2}\cdot \mathbf{x}_{1}.$ Continuing
this process, one can construct points $\mathbf{x}_{3}$, $\mathbf{x}_{4}$,
and so on. Note that \ the set of all constructed points $\mathbf{x}_{i}$, $%
i=0,1,...,$ forms a path. By Singer's above result, this path is extremal for the
function $f-g_{0}$. We have proved the necessity and hence Theorem 1.4.
\end{proof}

Theorem 1.4, in a more general setting, was proven in Pinkus \cite[p.99]{117}
under additional assumption that $Q$ is convex. Convexity assumption was
made to guarantee continuity of the following functions
\begin{equation*}
g_{1,i}(t)=\max_{\substack{ \mathbf{x}\in Q  \\ \mathbf{a}_{i}\cdot \mathbf{x%
}=t}}F(\mathbf{x})\ \ \text{and }\ g_{2,i}(t)=\min\limits_{\substack{
\mathbf{x}\in Q  \\ \mathbf{a}_{i}\cdot \mathbf{x}=t}}F(\mathbf{x}),\text{ }%
i=1,2,
\end{equation*}%
where $F$ is an arbitrary continuous function on $Q$. Note that in the proof
of Theorem 1.4 we did not need continuity of these functions.

\bigskip

It is well known that characterization theorems of this type are very
essential in approximation theory. Chebyshev was the first to prove a
similar result for polynomial approximation. Khavinson \cite{79}
characterized extremal elements in the special case of the problem
considered here. His case allows the approximation of a continuous bivariate
function $f\left( {x,y}\right) $ by functions of the form $\varphi \left( {x}%
\right) +\psi \left( {y}\right)$.

\bigskip

\subsection{Construction of an extremal element}

In 1951, Diliberto and Straus \cite{26} established a formula for the error
of approximation of a bivariate function by sums of univariate functions. Their
formula contains the supremum over all closed bolts (see Section 3.3.1).
Although the mentioned formula is
valid for all continuous functions, it is not easily calculable. Therefore,
it cannot give a desired effect if one is interested in the precise
value of the approximation error. After this general result some authors
started to seek easily calculable formulas for the approximation error by
considering not the whole space of continuous functions, but some subsets thereof
(see, for example, \cite{4,7,55,56,79,121}). These subsets were chosen so
that they could provide precise and easy computation of the approximation
error. Since the set of ridge functions contains univariate functions,
one may ask for explicit formulas for the error of
approximation of a multivariate function by sums of ridge functions.

In this section, we see how with the use of Theorem 1.3 (or 1.4) it is
possible to find the approximation error and construct an extremal element in the problem of approximation
by sums of ridge functions. We restrict ourselves to $%
\mathbb{R}^{2}.$ To make the problem more precise, let $\Omega $ be a
compact set in $\mathbb{R}^{2},$ $f \in C\left( {%
\Omega }\right)$, $\mathbf{a}=\left( {a_{1},a_{2}}\right)$ and $\mathbf{b}%
=\left( {b_{1},b_{2}}\right)$ be linearly independent vectors.
Consider the approximation of $f$ by functions from $\mathcal{R}=\mathcal{R}\left(
\mathbf{a},\mathbf{b}\right)$. We want, under
some suitable conditions on $f\;$and $\Omega $, to establish a formula for an easy
and direct computation of the approximation error $E\left(f,\mathcal{R}\right)$.

\bigskip

\textbf{Theorem 1.5.} \textit{Let
\begin{equation*}
\Omega =\left\{ \mathbf{x}\in \mathbb{R}^{2}:c_{1}\leq \mathbf{a}\cdot
\mathbf{x}\leq d_{1},\ \ c_{2}\leq \mathbf{b}\cdot \mathbf{x}\leq
d_{2}\right\} ,
\end{equation*}%
where $c_{1}<d_{1}$ and $c_{2}<d_{2}$. Let a function $f(\mathbf{x})\in
C(\Omega )$ have the continuous partial derivatives $\frac{\partial ^{2}f}{%
\partial x_{1}^{2}},\frac{\partial ^{2}f}{\partial x_{1}\partial x_{2}},%
\frac{\partial ^{2}f}{\partial x_{2}^{2}}$ and for any $\mathbf{x}\in \Omega
$
\begin{equation*}
\frac{\partial ^{2}f}{\partial x_{1}\partial x_{2}}\left(
a_{1}b_{2}+a_{2}b_{1}\right) -\frac{\partial ^{2}f}{\partial x_{1}^{2}}%
a_{2}b_{2}-\frac{\partial ^{2}f}{\partial x_{2}^{2}}a_{1}b_{1}\geq 0.
\end{equation*}%
Then
\begin{equation*}
E\left(f,\mathcal{R}\right)=\frac{1}{4}\left(
f_{1}(c_{1},c_{2})+f_{1}(d_{1},d_{2})-f_{1}(c_{1},d_{2})-f_{1}(d_{1},c_{2})%
\right) ,
\end{equation*}%
where
\begin{equation*}
f_{1}(y_{1},y_{2})=f\left( \frac{y_{1}b_{2}-y_{2}a_{2}}{a_{1}b_{2}-a_{2}b_{1}%
},\frac{y_{2}a_{1}-y_{1}b_{1}}{a_{1}b_{2}-a_{2}b_{1}}\right) .\eqno(1.29)
\end{equation*}%
}

\bigskip

\begin{proof} Introduce the new variables
\begin{equation*}
y_{1}=a_{1}x_{1}+a_{2}x_{2},\ \ y_{2}=b_{1}x_{1}+b_{2}x_{2}.\eqno(1.30)
\end{equation*}

Since the vectors $(a_{1},a_{2})$ and $(b_{1},b_{2})$ are linearly
independent, for any $(y_{1},y_{2})\in Y$, where $Y=[c_{1},d_{1}]\times
\lbrack c_{2},d_{2}]$, there exists only one solution $(x_{1},x_{2})\in
\Omega $ of the system (1.30). The coordinates of this solution are
\begin{equation*}
x_{1}=\frac{y_{1}b_{2}-y_{2}a_{2}}{a_{1}b_{2}-a_{2}b_{1}},\qquad \ x_{2}=%
\frac{y_{2}a_{1}-y_{1}b_{1}}{a_{1}b_{2}-a_{2}b_{1}}.\eqno(1.31)
\end{equation*}

The linear transformation (1.31) transforms the function $f(x_{1},x_{2})$ to
the function $f_{1}(y_{1},y_{2})$. Consider the approximation of $%
f_{1}(y_{1},y_{2})$ from the set
\begin{equation*}
\mathcal{Z}=\left\{ z_{1}(y_{1})+z_{2}(y_{2}):z_{i}\in C(\mathbb{R}),\
i=1,2\right\} .
\end{equation*}

It is easy to see that
\begin{equation*}
E\left( f,\mathcal{R}\right) =E\left( f_{1},\mathcal{Z}\right) .\eqno(1.32)
\end{equation*}

With each rectangle $S=\lbrack u_{1} ,v_{1} ]\times \lbrack u_{2} ,v_{2}
]\subset Y$ we associate the functional
\begin{equation*}
L\left(h, S\right) =\frac{1}{4} \left(h(u_{1} ,u_{2} )+h(v_{1} ,v_{2}
)-h(u_{1} ,v_{2} )-h (v_{1} ,u_{2})\right),\ \ h\in C(Y).
\end{equation*}

This functional has the following obvious properties:

(i) $L(z,S)=0$ for any $z\in \mathcal{Z}$ and $S\subset Y$.

(ii) For any point $(y_{1} ,y_{2} )\in Y$, $L(f_1,
Y)=\sum\limits_{i=1}^{4}L(f_1, S_{i} ) $, where $S_{1} =[c_{1} ,y_{1}]\times
[c_{2} ,y_{2}],$ $S_{2} =[y_{1} ,d_{1}]\times [y_{2} ,d_{2}],$ $S_{3}
=[c_{1} ,y_{1}]\times [y_{2} ,d_{2}],$ $S_{4} =[y_{1} ,d_{1}]\times [c_{2}
,y_{2}]$.

By the conditions of the theorem, it is not difficult to verify that
\begin{equation*}
\frac{\partial ^{2} f_1}{\partial y_{1} \partial y_{2} } \geq 0\ \
\mbox{for
any}\ \ (y_{1} ,y_{2} )\in Y.
\end{equation*}

Integrating both sides of the last inequality over arbitrary rectangle $%
S=[u_{1},v_{1}]\times \lbrack u_{2},v_{2}]\subset Y$, we obtain that
\begin{equation*}
L\left( f_{1},S\right) \geq 0.\eqno(1.33)
\end{equation*}%
Set the function
\begin{equation*}
f_{2}(y_{1},y_{2})=L\left( f_{1},S_{1}\right) +L\left( f_{1},S_{2}\right)
-L\left( f_{1},S_{3}\right) -L\left( f_{1},S_{4}\right) .\eqno(1.34)
\end{equation*}

It is not difficult to verify that the function $f_{1}-f_{2}$ belongs to $%
\mathcal{Z}$. Hence
\begin{equation*}
E\left( f_{1},\mathcal{Z}\right) =E\left( f_{2},\mathcal{Z}\right) .\eqno%
(1.35)
\end{equation*}

Calculate the norm $\left\Vert f_{2}\right\Vert $. From the property (ii),
it follows that
\begin{equation*}
f_{2}(y_{1},y_{2})=L(f_{1},Y)-2(L(f_{1},S_{3})+L(f_{1},S_{4}))
\end{equation*}%
and
\begin{equation*}
f_{2}(y_{1},y_{2})=2\left( L\left( f_{1},S_{1}\right) +L\left(
f_{1},S_{2}\right) \right) -L\left( f_{1},Y\right) .
\end{equation*}%
From the last equalities and (1.33), we obtain that
\begin{equation*}
\left\vert f_{2}(y_{1},y_{2})\right\vert \leq L\left( f_{1},Y\right) ,\ %
\mbox{for any}\ (y_{1},y_{2})\in Y.
\end{equation*}%
On the other hand, one can check that
\begin{equation*}
f_{2}(c_{1},c_{2})=f_{2}(d_{1},d_{2})=L\left( f_{1},Y\right) \eqno(1.36)
\end{equation*}%
and
\begin{equation*}
f_{2}(c_{1},d_{2})=f_{2}(d_{1},c_{2})=-L\left( f_{1},Y\right) .\eqno(1.37)
\end{equation*}%
Therefore,
\begin{equation*}
\left\Vert f_{2}\right\Vert =L\left( f_{1},Y\right) .\eqno(1.38)
\end{equation*}%
Note that the points $%
(c_{1},c_{2}),(c_{1},d_{2}),(d_{1},d_{2}),(d_{1},c_{2}) $ in the given order
form a closed path with respect to the directions $(0;1) $ and $(1;0)$. We
conclude from (1.36)-(1.38) that this path is extremal for $f_{2}$. By
Theorem 1.3, $z_{0}=0$ is a best approximation to $f_{2}$. Hence
\begin{equation*}
E\left( f_{2},\mathcal{Z}\right) =L\left( f_{1},Y\right) .\eqno(1.39)
\end{equation*}

Now from (1.32),(1.35) and (1.39) we finally conclude that
\begin{equation*}
E\left( f,\mathcal{R}\right) =L\left( f_{1},Y\right) =\frac{1}{4}\left(
f_{1}(c_{1},c_{2})+f_{1}(d_{1},d_{2})-f_{1}(c_{1},d_{2})-f_{1}(d_{1},c_{2})%
\right) ,
\end{equation*}%
which is the desired result. \end{proof}

\textbf{Corollary 1.4.}\ \textit{Let all the conditions of Theorem 1.5 hold
and $f_{1}(y_{1},y_{2})$ is the function defined in (1.29). Then the
function $g_{0}(y_{1},y_{2})=g_{1,0}(y_{1})+g_{2,0}(y_{2})$, where
\begin{equation*}
g_{1,0}(y_{1})=\frac{1}{2}f_{1}(y_{1},c_{2})+\frac{1}{2}f_{1}(y_{1},d_{2})-%
\frac{1}{4}f_{1}(c_{1},c_{2})-\frac{1}{4}f_{1}(d_{1},d_{2}),
\end{equation*}%
\begin{equation*}
g_{2,0}(y_{2})=\frac{1}{2}f_{1}(c_{1},y_{2})+\frac{1}{2}f_{1}(d_{1},y_{2})-%
\frac{1}{4}f_{1}(c_{1},d_{2})-\frac{1}{4}f_{1}(d_{1},c_{2})
\end{equation*}%
and $y_{1}=a_{1}x_{1}+a_{2}x_{2}$, $y_{2}=b_{1}x_{1}+b_{2}x_{2}$, is a best
approximation from the set $\mathcal{R}(a,b)$ to the function $f$.}

\bigskip

\begin{proof} It is not difficult to verify that the function $%
f_{2}(y_{1},y_{2})$ defined in (1.34) has the form
\begin{equation*}
f_{2}(y_{1},y_{2})=f_{1}(y_{1},y_{2})-g_{1,0}(y_{1})-g_{2,0}(y_{2}).
\end{equation*}

On the other hand, we know from the proof of Theorem 1.5 that
\begin{equation*}
E(f_{1},\mathcal{Z})=\left\Vert f_{2}\right\Vert .
\end{equation*}%
Therefore, the function $g_{1,0}(y_{1})+g_{2,0}(y_{2})$ is a best
approximation to $f_{1}$. Then the function $g_{1,0}(\mathbf{a}\cdot \mathbf{%
x})+g_{2,0}(\mathbf{b}\cdot \mathbf{x})$ is an extremal element from $%
\mathcal{R}(\mathbf{a},\mathbf{b})$ to $f(\mathbf{x})$. \end{proof}

\textbf{Remark 1.1.} Rivlin and Sibner \cite{121}, and Babaev \cite{7}
proved Theorem 1.5 for the case in which $\mathbf{a}${\ and\ }$\mathbf{b}$
are the coordinate vectors. Our proof of Theorem 1.5 is different, short and
elementary. Moreover, it has turned out to be useful in constructing an
extremal element (see the proof of Corollary 1.4).

\bigskip

\subsection{Density of ridge functions and some problems}

Let $\mathbf{a}^{1}${\ and }$\mathbf{a}^{2}$ be nonzero directions in $%
\mathbb{R}^{d}$. One may ask the following question: are there cases in
which the set $\mathcal{R}\left( \mathbf{a}^{1},\mathbf{a}^{2}\right) $ is
dense in the space of all continuous functions? Undoubtedly, a positive
answer depends on the geometrical structure of compact sets over which all
the considered functions are defined. This problem may be interesting in the
theory of partial differential equations. Take, for example, equation
(1.10). A positive answer to the problem means that for any continuous
function $f$ there exist solutions of the given equation uniformly
converging to $f$.

It should be remarked that our problem is a special case of the problem
considered by Marshall and O'Farrell. In \cite{107}, they obtained a
necessary and sufficient condition for a sum $A_{1}+A_{2}$ of two
subalgebras to be dense in $C(U)$, where $C(U)$ denotes the space of
real-valued continuous functions on a compact Hausdorff space $U$. Below we
describe Marshall and O' Farrell's result for sums of ridge functions.

Let $X$ be a compact subset of $\mathbb{R}^{d}.$ The relation on $X$,
defined by setting $\mathbf{x}\approx \mathbf{y}$ if $\mathbf{x}$ and $%
\mathbf{y}$ belong to some path in $X$, is an equivalence relation. The
equivalence classes we call orbits. \bigskip

\textbf{Theorem 1.6.} \textit{Let $X$ be a compact subset of $\mathbb{R}^{d}$
with all its orbits closed. The set $\mathcal{\ R}\left( \mathbf{a}^{1},%
\mathbf{a}^{2}\right) $ is dense in $C(X)$ if and only if $X$ contains no
closed path with respect to the directions $\mathbf{a}^{1}${\ and }$\mathbf{a%
}^{2}$.}

\bigskip

The proof immediately follows from proposition 2 in \cite{108} established
for the sum of two algebras. Since that proposition was given without proof,
for completeness of the exposition we give the proof of Theorem 1.6.

\begin{proof}
\textit{Necessity}. If $X$ has closed paths, then $X$ has a closed path $%
p^{\prime }=\left( \mathbf{p}_{1}^{\prime },...,\mathbf{p}_{2m}^{\prime
}\right) $ such that all points $\mathbf{p}_{1}^{\prime },...,\mathbf{p}%
_{2m}^{\prime }$ are distinct. In fact, such a special path can be obtained
from any closed path $p=\left( \mathbf{p}_{1},...,\mathbf{p}_{2n}\right) $
by the following simple algorithm: if the points of the path $p$ are not all
distinct, let $i$ and $k>0$ be the minimal indices such that $\mathbf{p}_{i}=%
\mathbf{p}_{i+2k}$; delete from $p$ the subsequence $\mathbf{p}_{i+1},...,%
\mathbf{p}_{i+2k}$ and call $p$ the obtained path; repeat the above step
until all points of $p$ are all distinct; set $p^{\prime }:=p$. By Urysohn's
great lemma, there exist continuous functions $h=h(\mathbf{x})$ on $X$ such
that $h(\mathbf{p}_{i}^{\prime })=1$, $i=1,3,...,2m-1$, $h(\mathbf{p}%
_{i}^{\prime })=-1$, $i=2,4,...,2m$ and $-1<h(\mathbf{x})<1$ elsewhere.
Consider the measure

\begin{equation*}
\mu _{p^{\prime }}=\frac{1}{2m}\sum_{i=1}^{2m}(-1)^{i-1}\delta _{\mathbf{p}%
_{i}^{\prime }}\text{ ,}
\end{equation*}%
where $\delta _{\mathbf{p}_{i}^{\prime }}$ is a point mass at $\mathbf{p}%
_{i}^{\prime }$. For this measure, $\int\limits_{X}hd\mu _{p^{\prime }}=1$
and $\int\limits_{X}gd\mu _{p^{\prime }}=0$ for all functions $g\in \mathcal{%
R}\left( \mathbf{a}^{1},\mathbf{a}^{2}\right) $. Thus the set $\mathcal{R}%
\left( \mathbf{a}^{1},\mathbf{a}^{2}\right) $ cannot be dense in\textit{\ }%
$C(X)$.

\textit{Sufficiency}. We are going to prove that the only annihilating
regular Borel measure for $\mathcal{R}\left( \mathbf{a}^{1},\mathbf{a}%
^{2}\right) $ is the zero measure. Suppose, contrary to this assumption,
there exists a nonzero annihilating measure on $X$ for $\mathcal{R}\left(
\mathbf{a}^{1},\mathbf{a}^{2}\right) $. The class of such measures with
total variation not more than $1$ we denote by $S.$ Clearly, $S$ is weak-*
compact and convex. By the Krein-Milman theorem, there exists an extreme
measure $\mu $ in $S.$ Since the orbits are closed, $\mu $ must be supported
on a single orbit. Denote this orbit by $T.$

For $i=1,2,$ let $X_{i}$ be the quotient space of $X$ obtained by
identifying the points $\mathbf{y}$ and $\mathbf{z}$ whenever $\mathbf{a}^{i}%
{\cdot }\mathbf{y}=\mathbf{a}^{i}{\cdot }\mathbf{z}$. Let $\pi _{i}$ be the
natural projection of $X$ onto $X_{i}$. For a fixed point $t\in X$ set $%
T_{1}=\{t\}$, $T_{2}=\pi _{1}^{-1}\left( \pi _{1}T_{1}\right) $, $T_{3}=\pi
_{2}^{-1}\left( \pi _{2}T_{2}\right) $, $T_{4}=\pi _{1}^{-1}\left( \pi
_{1}T_{3}\right) $, $...$ Obviously, $T_{1}\subset T_{2}\subset T_{3}\subset
\cdot \cdot \cdot $ . Therefore, for some $k\in \mathbb{N}$, $\left\vert \mu
\right\vert (T_{2k})>0$, where $\left\vert \mu \right\vert $ is a total
variation measure of $\mu$. Since $\mu $ is orthogonal to every continuous
function of the form ${g\left( \mathbf{a}^{1}{\cdot }\mathbf{x}\right) }$, $%
\mu (T_{2k})=0$. From the Haar decomposition $\mu (T_{2k})=\mu
^{+}(T_{2k})-\mu ^{-}(T_{2k})$ it follows that $\mu ^{+}(T_{2k})=\mu
^{-}(T_{2k})>0$. Fix a Borel subset $S_{0}\subset T_{2k}$ such that $\mu
^{+}(S_{0})>0$ and $\mu ^{-}(S_{0})=0$. Since $\mu $ is orthogonal to every
continuous function of the form ${g\left( \mathbf{a}^{2}{\cdot }\mathbf{x}%
\right) }$, $\mu (\pi _{2}^{-1}\left( \pi _{2}S_{0}\right) )=0.$ Therefore,
one can chose a Borel set $S_{1}$ such that $S_{1}\subset \pi
_{2}^{-1}\left( \pi _{2}S_{0}\right) \subset T_{2k+1}$, $S_{1}\cap
S_{0}=\varnothing $, $\mu ^{+}(S_{1})=0$, $\mu ^{-}(S_{1})\geqslant \mu
^{+}(S_{0})$. By the same way one can chose a Borel set $S_{2}$ such that $%
S_{2}\subset \pi _{1}^{-1}\left( \pi _{1}S_{1}\right) \subset T_{2k+2}$, $%
S_{2}\cap S_{1}=\varnothing $, $\mu ^{-}(S_{2})=0$, $\mu
^{+}(S_{2})\geqslant \mu ^{-}(S_{1})$, and so on.

The sets $S_{0},S_{1},S_{2},...$are pairwise disjoint. For otherwise, there
would exist positive integers $n$ and $m,$ with $n<m$ and a path $%
(y_{n},y_{n+1},...,y_{m})$ such that $y_{i}\in S_{i}$ for $i=n,...,m$ and $%
y_{m}\in S_{m}\cap S_{n}$. But then there would exist paths $%
(z_{1},z_{2},...,z_{n-1},y_{n})$ and $(z_{1},z_{2}^{^{\prime
}},...,z_{n-1}^{^{\prime }},y_{m})$ with $z_{i}$ and $z_{i}^{^{\prime }}$ in
$T_{i}$ for $i=2,...,n-1.$ Hence, the set
\begin{equation*}
\{z_{1},z_{2},...,z_{n-1},y_{n},y_{n+1},...,y_{m},z_{n-1}^{^{\prime
}},...,z_{2}^{^{\prime }},z_{1}\}
\end{equation*}%
would contain a closed path. This would contradict our assumption on $X.$

Now, since the sets $S_{0},S_{1},S_{2},...,$ are pairwise disjoint and $%
\left\vert \mu \right\vert (S_{i})\geqslant \mu ^{+}(S_{0})>0$ for each $%
i=1,2,...,$ it follows that the total variation of $\mu $ is infinite. This
contradiction completes the proof.
\end{proof}

The following corollary concerns the problem considered by Colitschek and
Light \cite{36}.

\bigskip \textbf{Corollary 1.5.} \textit{Let $D$ be a compact subset of \ $%
\mathbb{R}^{2}$ with all its orbits closed. Let $W$ denote the set of all
solutions of the wave equation
\begin{equation*}
\frac{\partial ^{2}w}{\partial s\partial t}(s,t)=0,\;\ \ \ \ (s,t)\in D.
\end{equation*}%
Then
\begin{equation*}
\inf\limits_{w\in W}\left\Vert f-w\right\Vert =0
\end{equation*}%
for any continuous function $f(s,t)$ on $D$ if and only if $D$ contains no
closed bolt of lightning.}

\bigskip

\begin{proof} Let $\pi _{1}$ and $\pi _{2}$ denote the usual
coordinate projections, viz: $\pi _{1}(s,t)=s$ and $\pi _{2}(s,t)=t$, $%
(s,t)\in \mathbb{R}^{2}$. Set $S=\pi _{1}(D)$ and $T=\pi _{2}(D)$. It is
easy to see that
\begin{equation*}
W=\left\{ w\in C(D):w(s,t)=x(s)+y(t),\;\ \ x\in C^{2}(S),\;\ y\in
C^{2}(T)\right\} .
\end{equation*}
Set
\begin{equation*}
\widetilde{W}=\left\{ w\in C(D):w(s,t)=x(s)+y(t),\;\ \ x\in C(S),\;\ y\in
C(T)\right\} .
\end{equation*}
Since the set $W$ is dense in $\widetilde{W},$
\begin{equation*}
\inf\limits_{w\in W}\left\Vert f-w\right\Vert =\inf\limits_{w\in \widetilde{W%
}}\left\Vert f-w\right\Vert.
\end{equation*}
But by Theorem 1.6, the equality
\begin{equation*}
\inf\limits_{w\in \widetilde{W}}\left\Vert f-w\right\Vert =0
\end{equation*}%
holds for any $f\in C(D)$ if and only if $D$ contains no closed bolt of lightning.
\end{proof}

Let us discuss some difficulties that arise when studying sums of more than two ridge functions.
Consider the set
\begin{equation*}
\mathcal{R}\left( \mathbf{a}{^{1},...,}\mathbf{a}{^{r}}\right) ={\left\{
\sum\limits_{i=1}^{r}{g}_{i}{{{\left( \mathbf{a}{^{i}\cdot }\mathbf{x}%
\right) ,g}}}_{{i}}{{{\ \in C\left( \mathbb{R}\right) ,i=1,...,r}}}\right\} }%
,
\end{equation*}%
where $\mathbf{a}{^{1},...,}\mathbf{a}{^{r}}$ are pairwise linearly
independent vectors in $\mathbb{R}^{d}\backslash \{\mathbf{0}\}$. Let $r\geq
3$. How can one define a path in this general case? Recall that in the case when $r=2$, a path is
an ordered set of points $\left( \mathbf{p}_{1},\mathbf{p}_{2},...,\mathbf{p}%
_{n}\right) $ in $\mathbb{R}^{d}$ with edges $\mathbf{p}_{i}\mathbf{p}_{i+1}$
in alternating hyperplanes. The first, the third, the fifth,... hyperplanes
(also the second, the fourth, the sixth,... hyperplanes) are parallel. If
not differentiate between parallel hyperplanes, the path $\left( \mathbf{p}%
_{1},\mathbf{p}_{2},...,\mathbf{p}_{n}\right) $ can be considered as a trace
of some point traveling in two alternating hyperplanes. In this case, if the
point starts and stops at the same location (i.e., if $\mathbf{p}_{n}=%
\mathbf{p}_{1})$ and $n$ is an odd number, then the path functional
\begin{equation*}
G(f)=\frac{1}{n-1}\sum\limits_{i=1}^{n-1}(-1)^{i+1}f(\mathbf{p}_{i}),
\end{equation*}%
annihilates sums of ridge functions with the corresponding two fixed directions. The
picture becomes quite different and more complicated when the number of directions more than
two. The simple generalization of the above-mentioned arguments demands a
point traveling in three or more alternating hyperplanes. But in this case
the appropriate generalization of the functional $G$ does not annihilate
functions from $\mathcal{R}\left( \mathbf{a}{^{1},...,}\mathbf{a}{^{r}}%
\right) $.

There were several attempts to fill this gap in the special case when $r=d$
and $\mathbf{a}{^{1},...,}\mathbf{a}{^{r}}$ are the coordinate vectors.
Unfortunately, all these attempts failed (see, for example, the attempts in
\cite{26,37} and the refutations in \cite{4,109}).

At the end of this subsection we want to draw the readers attention to the following problems.
All these problems are open and cannot be solved by the methods presented here.

Let $Q$ be a compact subset of $\mathbb{R}^{d}$. Consider the approximation
of a continuous function defined on $Q$ by functions from $\mathcal{R}\left(
\mathbf{a}{^{1},...,}\mathbf{a}{^{r}}\right)$. Let $r\ge 3$.

\textbf{Problem 3.} \textit{Characterize those functions from $\mathcal{R}%
\left( \mathbf{a}{^{1},...,}\mathbf{a}{^{r}}\right) $ that are extremal to a
given continuous function.}

\textbf{Problem 4.} \textit{Establish explicit formulas for the error in
approximating from $\mathcal{R}\left( \mathbf{a}{^{1},...,}\mathbf{a}{^{r}}%
\right) $ \ and construct a best approximation.}

\textbf{Problem 5.} \textit{Find necessary and sufficient geometrical
conditions for the set $\mathcal{R}\left( \mathbf{a}{^{1},...,}\mathbf{a}{%
^{r}}\right) $ to be dense in $C(Q)$.}

It should be remarked that in \cite{108}, Problem 5 was set up for the sum
of $r $ subalgebras of $C(Q)$. Lin and Pinkus \cite{95} proved that the set $%
\mathcal{R}\left( \mathbf{a}{^{1},...,}\mathbf{a}{^{r}}\right) $ ($r$ may be
very large) is not dense in $C(\mathbb{R}^{d})$ in the topology of uniform
convergence on compact subsets of $\mathbb{R}^{d}$. That is, there are
compact sets $Q\subset \mathbb{R}^{d}$ such that $\mathcal{R}\left( \mathbf{a%
}{^{1},...,}\mathbf{a}{^{r}}\right) $ is not dense in $C(Q)$. In the case $%
r=2$, Theorem 1.6 complements this result, by describing compact sets $%
Q\subset \mathbb{R}^{2}$, for which $\mathcal{R}\left( \mathbf{a}{^{1},}%
\mathbf{a}{^{2}}\right) $ is dense in $C(Q)$.

\bigskip

\section{Sums of continuous ridge functions}

In this section, we find geometric means of deciding if any continuous
multivariate function can be represented by a sum of two continuous ridge
functions.

\subsection{Exposition of the problem}

In this section, we will consider the following representation problem
associated with the set $\mathcal{R}\left( \mathbf{a}^{1},...,\mathbf{a}%
^{r}\right) .$

\bigskip

\textbf{Problem 6.}\ \textit{Let $X$ be a compact subset of \ $\mathbb{R}%
^{d}.$ Give geometrical conditions that are necessary and sufficient for
\begin{equation*}
\mathcal{R}\left( \mathbf{a}^{1},...,\mathbf{a}^{r}\right)=C\left( X\right)
,
\end{equation*}%
where $C\left( X\right) $ is the space of continuous functions on $X$
furnished with the uniform norm.}

\bigskip

We solve this problem for $r=2$. Problem 6, like Problems 3--5 from the previous section, is open in
the case $r\geq 3$. Geometrical characterization of compact sets $X \subset \mathbb{R}^{d}$ with the property
$\mathcal{R}\left(\mathbf{a}^{1},...,\mathbf{a}^{r}\right)=C\left(X\right)$, $r\geq 3$,
seems to be beyond the scope of the methods discussed herein. Nevertheless, recall that this
problem in a quite abstract form, which involves regular Borel measures on $X$,
was solved by Sternfeld (see Section 1.1.1).

In the sequel, we will use the notation
\begin{equation*}
H_{1}=H_{1}\left( X\right) =\left\{ g_{1}\left( \mathbf{a}^{1}\cdot \mathbf{x%
}\right) :g_{1}\in C\left( \mathbb{R}\right) \right\} ,
\end{equation*}
\begin{equation*}
H_{2}=H_{2}\left( X\right) =\left\{ g_{2}\left( \mathbf{a}^{2}\cdot \mathbf{x%
}\right) :g_{2}\in C\left( \mathbb{R}\right) \right\} .
\end{equation*}
Note that by this notation, $\mathcal{R}\left( \mathbf{a}^{1},\mathbf{a}%
^{2}\right) =H_{1}+H_{2}.$

At the end of this section, we generalize the obtained result from $%
H_{1}+H_{2}$ to the set of sums $g_{1}\left( h_{1}\left( \mathbf{x}\right)
\right) +g_{2}\left( h_{2}\left( \mathbf{x}\right) \right)$, where $h_{1},
h_{2}$ are fixed continuous functions on $X$.

\bigskip

\subsection{The representation theorem}

\textbf{Theorem 1.7.} \textit{Let $X$ be a compact subset of $\mathbb{R}%
^{d} $. The equality
\begin{equation*}
H_{1}\left( X\right) +H_{2}\left( X\right) =C\left( X\right)
\end{equation*}%
holds if and only if $X$ contains no closed path and there exists a positive
integer $n_{0}$ such that the lengths of paths in $X$ are bounded by $n_{0}$.%
}

\bigskip

\begin{proof}
\textit{Necessity.} Let $H_{1}+H_{2}=C\left( X\right) $. Consider the linear
operator
\begin{equation*}
A:H_{1}\times H_{2}\rightarrow C\left( X\right), ~~~ A\left[ \left(
g_{1},g_{2}\right) \right] =g_{1}+g_{2},
\end{equation*}%
where $g_{1}\in H_{1},g_{2}\in H_{2}.$ The norm on $H_{1}\times H_{2}$ we
define as
\begin{equation*}
\left\Vert \left( g_{1},g_{2}\right) \right\Vert =\left\Vert
g_{1}\right\Vert +\left\Vert g_{2}\right\Vert .
\end{equation*}%
It is obvious that the operator $A$ is continuous with respect to this norm.
Besides, since $C\left( X\right) =H_{1}+H_{2},$ $A$ is a surjection.
Consider the conjugate operator
\begin{equation*}
A^*:C\left( X\right) ^{\ast }\rightarrow \left[ H_{1}\times H_{2}\right]
^{\ast }, ~~~ A^{\ast }\left[ G\right] =\left( G_{1},G_{2}\right) ,
\end{equation*}%
where the functionals $G_{1}$ and $G_{2}$ are defined as follows
\begin{equation*}
G_{1}\left( g_{1}\right) =G\left( g_{1}\right) ,g_{1}\in H_{1}; ~~~
G_{2}\left( g_{2}\right) =G\left( g_{2}\right) ,g_{2}\in H_{2}.
\end{equation*}

An element $\left( G_{1},G_{2}\right) $ from $\left[ H_{1}\times H_{2}\right]
^{\ast }$ has the norm
\begin{equation*}
\left\Vert \left( G_{1},G_{2}\right) \right\Vert =\max \left\{ \left\Vert
G_{1}\right\Vert ,\left\Vert G_{2}\right\Vert \right\} .\eqno(1.40)
\end{equation*}

Let now $p=\left( p_{1},...,p_{m}\right) $ be any path with different
points: $p_{i}\neq p_{j}$ for any $i\neq j$, $1\leq i,~j\leq m$. We
associate with $p$ the following functional over $C\left( X\right) $
\begin{equation*}
L\left[ f\right] =\frac{1}{m}\sum\limits_{i=1}^{m}\left( -1\right)
^{i-1}f\left( p_{i}\right) .
\end{equation*}%
Since $\left\vert L(f)\right\vert \leq \left\Vert f\right\Vert $ and $%
\left\vert L(g)\right\vert =\left\Vert g\right\Vert $ for a continuous
function $g(\mathbf{x})$ such that $g(p_{i})=1,\ $for odd indices $i,\
g(p_{j})=-1,$ for even\ indices$\ j\ $and $-1<g(\mathbf{x})<1$ elsewhere, we
obtain that $\left\Vert L\right\Vert =1$. Let $A^{\ast }\left[ L\right]
=\left( L_{1},L_{2}\right) $. One can easily verify that
\begin{equation*}
\left\Vert L_{i}\right\Vert \leq \frac{2}{m},i=1,2.
\end{equation*}%
Therefore, from (1.40) we obtain that
\begin{equation*}
\left\Vert A^{\ast }\left[ L\right] \right\Vert \leq \frac{2}{m}.\eqno(1.41)
\end{equation*}%
Since $A$ is a surjection, there exists $\delta >0$ such that
\begin{equation*}
\left\Vert A^{\ast }\left[ G\right] \right\Vert \geq \delta \left\Vert
G\right\Vert ~~~~\;\mbox{for any functional}\;\ G\in C\left( X\right) ^{\ast
}
\end{equation*}%
Hence
\begin{equation*}
\left\Vert A^{\ast }\left[ L\right] \right\Vert \geq \delta .\eqno(1.42)
\end{equation*}%
Now from (1.41) and (1.42) we conclude that
\begin{equation*}
m\leq \frac{2}{\delta }.
\end{equation*}

This means that for a path with different points, $n_{0}$ can be chosen as $%
\left[ \frac{2}{\delta }\right] +1$.

Let now $p=\left( p_{1},...,p_{m}\right) $ be a path with at least two
coinciding points. Then we can form a closed path with different points.
This may be done by the following way: let $i\ $and $j\ $be indices such
that $p_{i}=\ p_{j}\ $and $j-i\ $takes its minimal value. Note that in this
case all the points $p_{i},p_{i+1},...,p_{j-1}\ $are distinct. Now if $j-i\ $%
is an even number, then the path $(p_{i},p_{i+1},...,p_{j-1})\ $, and if $\
j-i\ $is an odd number, then the path $(p_{i+1},...,p_{j-1})$ is a closed
path with different points. It remains to show that $X$ can not possess
closed paths with different points. Indeed, if $q=\left(
q_{1},...,q_{2k}\right) $ is a path of this type, then the functional $L,$
associated with $q,$ annihilates all functions from $H_{1}+H_{2}$. On the
other hand, $L\left[ f\right] =1$ for a continuous function $f$ on $X$
satisfying the conditions $f\left( t\right) =1$ if $t\in \left\{
q_{1},q_{3},...,q_{2k-1}\right\} ;$ $f\left( t\right) =-1$ if $t\in \left\{
q_{2},q_{4},...,q_{2k}\right\} ;$ $f\left( t\right) \in \left( -1;1\right) $
if $t\in X\backslash q$ . This implies on the contrary to our assumption
that $H_{1}+H_{2}\neq C\left( X\right) $. The necessity has been proved.

\textit{Sufficiency.} Let $X$ contains no closed path and the lengths of all
paths are bounded by some positive integer $n_{0}$. We may suppose that any
path has different points. Indeed, in other case we can form a closed path,
which contradicts our assumption.

For $i=1,2,$ let $X_{i}$ be the quotient space of $X$ obtained by
identifying the points $a$ and $b$ whenever $g\left( a\right) =g\left(
b\right) $ for each $g$ in $H_{i}$. Let $\pi _{i}$ be the natural projection
of $X$ onto $X_{i}$. For a point $t\in X$ set $T_{1}=\pi _{1}^{-1}\left( \pi
_{1}t\right) ,T_{2}=\pi _{2}^{-1}\left( \pi _{2}T_{1}\right) ,\ldots .$ By $%
O\left( t\right) $ denote the orbit of $X$ containing $t.$ Since the length
of any path in $X$ is not more than $n_{0}$, we conclude that $O\left(
t\right) =T_{n_{0}}$. Since $X\ $is compact, the sets $%
T_{1},T_{2},...,T_{n_{0}},\ $hence $O(t),$ are compact. By Theorem 1.6,
$\overline{H_{1}+H_{2}}=C\left( X\right)$.

Now let us show that $H_{1}+H_{2}$ is closed in $C\left( X\right) $. Set
\begin{equation*}
H_{3}=H_{1}\cap H_{2}.
\end{equation*}

Let $X_{3}$ and $\pi _{3}$ be the associated quotient space and projection.
Fix some $a\in X_{3}$. Show, within conditions of our theorem, that if $t\in
\pi _{3}^{-1}\left( a\right) ,$ then $O\left( t\right) =\pi _{3}^{-1}\left(
a\right) $. The inclusion $O\left( t\right) \subset \pi _{3}^{-1}\left(
a\right) $ is obvious. Suppose that there exists a point $t_{1}\in \pi
_{3}^{-1}\left( a\right) $ such that $t_{1}\notin O\left( t\right)
$. Then $O\left( t\right) \cap O\left( t_{1}\right) =\emptyset $. By $X|O$
denote the factor space generated by orbits of $X$. $X|O$ is a normal
topological space with its natural factor topology. Hence we can construct a
continuous function $u\in C\left( X|O\right) $ such that $u\left( O\left(
t\right) \right) =0,$ $u\left( O\left( t_{1}\right) \right) =1$. The
function $\upsilon \left( x\right) =u\left( O\left( x\right) \right) ,\;\
x\in X,$ is continuous on $X$ and belongs to $H_{3}$ as a function being
constant on each orbit. But, since $O\left( t\right) \subset \pi
_{3}^{-1}\left( a\right) $ and $O\left( t_{1}\right) \subset \pi
_{3}^{-1}\left( a\right) $, the function $\upsilon \left( x\right) $ can not
take different values on $O\left( t\right) $ and $O\left( t_{1}\right) $.
This contradiction means that there is not a point $t_{1}\in \pi
_{3}^{-1}\left( a\right) $ such that $t_{1}\notin O\left( t\right) $. Thus,
\begin{equation*}
O\left( t\right) =\pi _{3}^{-1}\left( a\right) \eqno(1.43)
\end{equation*}%
for any $a\in X_{3}$ and $t\in \pi _{3}^{-1}\left( a\right) $.

Now prove that there exists a positive real number $c$ such that
\begin{equation*}
\sup\limits_{z\in X_{3}}\underset{\pi _{3}^{-1}\left( z\right) }{var}f\leq
c\sup\limits_{y\in X_{2}}\underset{\pi _{2}^{-1}\left( y\right) }{var}f\eqno%
(1.44)
\end{equation*}%
for all $f$ in $H_{1}$. Note that for $Y\subset X,\ \;\underset{Y}{var}f$ is
the variation of $f$ on the set $Y.$ That is, $\;$%
\begin{equation*}
\underset{Y}{var}f=\sup\limits_{x,y\in Y}\left\vert f\left( x\right)
-f\left( y\right) \right\vert .
\end{equation*}

Due to (1.43), inequality (1.44) can be written in the following form
\begin{equation*}
\sup_{t\in X}\underset{O\left( t\right) }{var}f\leq c\sup_{t\in X}\underset{%
\pi _{2}^{-1}\left( \pi _{2}\left( t\right) \right) }{var}f\eqno(1.45)
\end{equation*}%
for all $f\in H_{1}$.

Let $t\in X$ and $t_{1},t_{2}$ be arbitrary points of $O\left( t\right) $.
Then there is a path $\left( b_{1},b_{2},...,b_{m}\right) $ with $%
b_{1}=t_{1} $ and $b_{m}=t_{2}$. Besides, by the condition, $m\leq n_{0}$ .
Let first $\mathbf{a}^{2}\cdot b_{1}=\mathbf{a}^{2}\cdot b_{2},$ $\mathbf{a}%
^{1}\cdot b_{2}=\mathbf{a}^{1}\cdot b_{3},...,\mathbf{a}^{2}\cdot b_{m-1}=%
\mathbf{a}^{2}\cdot b_{m}$. Then for any function $f\in H_{1}$
\begin{equation*}
\left\vert f\left( t_{1}\right) -f\left( t_{2}\right) \right\vert
=\left\vert f\left( b_{1}\right) -f\left( b_{2}\right) +...-f\left(
b_{m}\right) \right\vert \leq
\end{equation*}%
\begin{equation*}
\leq \left\vert f\left( b_{1}\right) -f\left( b_{2}\right) \right\vert
+...+\left\vert f\left( b_{m-1}\right) -f\left( b_{m}\right) \right\vert
\leq \frac{n_{o}}{2}\sup_{t\in X}\underset{\pi _{2}^{-1}\left( \pi
_{2}\left( t\right) \right) }{var}f.\eqno(1.46)
\end{equation*}

It is not difficult to verify that inequality (1.46) holds in all other
possible cases of the path $\left( b_{1},...,b_{m}\right) $. Now from (1.46)
we obtain (1.45), hence (1.44), where $c=\frac{n_{0}}{2}$. In \cite{108},
Marshall and O'Farrell proved the following result (see
\cite[Proposition 4]{108}): Let $A_{1}\ $and $A_{2}\ $be closed subalgebras of $C(X)\ $that
contain the constants. Let $(X_{1},\pi _{1}),\ (X_{2},\pi _{2})\ $and $%
(X_{3},\pi _{3})\ $be the quotient spaces and projections associated with
the algebras $A_{1},$ $A_{2}\ $and $A_{3}=A_{1}\cap A_{2}\ $respectively.
Then $A_{1}+A_{2}\ $is closed in $C(X)\ $if and only if there exists a
positive real number $c$ such that
\begin{equation*}
\sup\limits_{z\in X_{3}}\underset{\pi _{3}^{-1}\left( z\right) }{var}f\leq
c\sup\limits_{y\in X_{2}}\underset{\pi _{2}^{-1}\left( y\right) }{var}f
\end{equation*}%
for all $f\ $in $A_{1}.$

By this proposition, (1.44) implies that $H_{1}+H_{2}$ is closed in $C\left(
X\right) $. Thus we finally obtain that $H_{1}+H_{2}=C\left( X\right) $.
\end{proof}

Paths with respect to two directions are explicit objects and give geometric
means of deciding if $H_{1}+H_{2}=C\left( X\right) $. Let us show this in
the example of the bivariate ridge functions $g_{1}=x_{1}+x_{2}\ $and $%
g_{2}=x_{1}-x_{2}.$ If $X$ is the union of two parallel line segments in $%
\mathbb{R}^{2},$ not parallel to any of the lines $x_{1}+x_{2}=0$ and $%
x_{1}-x_{2}=0,\ $then Theorem 1.7 holds. If $X$ is any bounded part of the
graph of the function $x_{2}=\arcsin (\sin x_{1}),$ then Theorem 1.7 also
holds. Let now $X\ $be the set

\begin{equation*}
\begin{array}{c}
\{(0,0),(1,-1),(0,-2),(-1\frac{1}{2},-\frac{1}{2}),(0,1),(\frac{3}{4},\frac{1%
}{4}),(0,-\frac{1}{2}), \\
(-\frac{3}{8},-\frac{1}{8}),(0,\frac{1}{4}),(\frac{3}{16},\frac{1}{16}%
),...\}.%
\end{array}%
\end{equation*}

In this case, there is no positive integer bounding lengths of all paths.
Thus Theorem 1.7 fails. Note that since orbits of all paths are closed,
Theorem 1.6 from the previous section shows $H_{1}+H_{2}$ is dense in $%
C\left( X\right) .$

If $X$ is any set with interior points, then both Theorem 1.6 and Theorem
1.7 fail, since any such set contains the vertices of some parallelogram
with sides parallel to the directions $\mathbf{a}^{1}$ and $\mathbf{a}^{2}$,
that is a closed path.

Theorem 1.7 admits a direct generalization to the representation by sums $%
g_{1}\left( h_{1}\left( \mathbf{x}\right) \right) +g_{2}\left( h_{2}\left(
\mathbf{x}\right) \right) $, where $h_{1}\left( \mathbf{x}\right) $ and $%
h_{2}\left( \mathbf{x}\right) $ are fixed continuous functions on $X$. This
generalization needs consideration of new objects -- paths with respect to
two continuous functions.

\bigskip

\textbf{Definition 1.5.} \textit{Let $X$ be a compact set in $\mathbb{R}%
^{d}$ and $h_{1},h_{2}\in C\left( X\right)$. A finite ordered subset $%
\left( p_{1},p_{2},...,p_{m}\right)$ of $X$ with $p_{i}\neq p_{i+1}\left(
i=1,...,m-1\right)$, and either $h_{1}\left( p_{1}\right) =h_{1}\left(
p_{2}\right)$, $h_{2}\left( p_{2}\right) =h_{2}\left( p_{3}\right)$, $%
h_{1}\left( p_{3}\right) =h_{1}\left( p_{4}\right),...,$ or $h_{2}\left(
p_{1}\right) =h_{2}\left( p_{2}\right)$, $h_{1}\left( p_{2}\right)
=h_{1}\left( p_{3}\right)$, $h_{2}\left( p_{3}\right) =h_{2}\left(
p_{4}\right),...,$ is called a path with respect to the functions $h_{1}$
and $h_{2}$ or, shortly, an $h_{1}$-$h_{2}$ path.}

\bigskip

\textbf{Theorem 1.8.} \textit{Let $X$ be a compact subset of $\mathbb{R}^{d}$%
. All functions $f \in C(X)$ admit a representation
\begin{equation*}
f(\mathbf{x})=g_{1}\left( h_{1}\left( \mathbf{x}\right) \right) +g_{2}\left(
h_{2}\left( \mathbf{x}\right) \right) ,~g_{1},g_{2}\in C(\mathit{\mathbb{R}})
\end{equation*}%
if and only if the set $X$ contains no closed $h_{1}$-$h_{2}$ path and there
exists a positive integer $n_{0}$ such that the lengths of $h_{1}$-$h_{2}$
paths in $X$ are bounded by $n_{0}$.}

\bigskip

The proof can be carried out by the same arguments as above.

It should be noted that Theorem 1.8 was first proved by Khavinson in his
monograph \cite{76}. Khavinson's proof (see \cite[p.87]{76}) used theorems
of Sternfeld \cite{132} and Medvedev \cite[Theorem 2.2]{76}, whereas our
proof, which generalizes the ideas of Khavinson, was based on the above
proposition of Marshall and O'Farrell.

\bigskip

\section{On the proximinality of ridge functions}

In this section, using two results of Garkavi, Medvedev and Khavinson \cite%
{35}, we give sufficient conditions for proximinality of sums of two ridge
functions with bounded and continuous summands in the spaces of bounded and
continuous multivariate functions, respectively. In the first case, we give
an example which shows that the corresponding sufficient condition cannot be
made weaker for certain subsets of $\mathbb{R}^{n}$. In the second case, we
obtain also a necessary condition for proximinality. All the results are
furnished with plenty of examples. The results, examples and following
discussions naturally lead us to a conjecture on the proximinality of the
considered class of ridge functions.

\subsection{Problem statement}

Let $E$ be a normed linear space and $F$ be its subspace. We say that $F$ is
proximinal in $E$ if for any element $e\in E$ there exists at least one
element $f_{0}\in F$ such that
\begin{equation*}
\left\Vert e-f_{0}\right\Vert =\inf_{f\in F}\left\Vert e-f\right\Vert .
\end{equation*}

In this case, the element $f_{0}$ is said to be extremal to $e$.

We are interested in the problem of proximinality of the set of linear
combinations of ridge functions in the spaces of bounded and continuous
functions respectively. This problem will be considered in the simplest case
when the class of approximating functions is the set
\begin{equation*}
\mathcal{R}=\mathcal{R}\left( \mathbf{a}^{1},\mathbf{a}^{2}\right) ={\
\left\{ {g_{1}\left( \mathbf{a}^{1}{\cdot }\mathbf{x}\right) +g_{2}\left(
\mathbf{a}^{2}{\cdot }\mathbf{x}\right) :g}_{i}:{{\mathbb{R\rightarrow R}}
,i=1,2}\right\} }.
\end{equation*}
Here $\mathbf{a}^{1}${and }$\mathbf{a}^{2}$ are fixed directions and we vary
over ${g}_{i}$. It is clear that this is a linear space. Consider the
following three subspaces of $\mathcal{R}$. The first is obtained by taking
only bounded sums ${g_{1}\left( \mathbf{a}^{1}{\cdot }\mathbf{x}\right)
+g_{2}\left( \mathbf{a}^{2}{\cdot }\mathbf{x}\right) }$ over some set $X$ in
$\mathbb{R}^{n}.$ We denote this subspace by $\mathcal{R}_{a}(X)$. The
second and the third are subspaces of $\mathcal{R}$ with bounded and
continuous summands $g_{i}\left( \mathbf{a}^{i}\cdot \mathbf{x}\right)
,~i=1,2,$ on $X$ respectively. These subspaces will be denoted by $\mathcal{%
R }_{b}(X)$ and $\mathcal{R}_{c}(X).$ In the case of $\mathcal{R}_{c}(X),$
the set $X$ is considered to be compact.

Let $B(X)$ and $C(X)$ be the spaces of bounded and continuous multivariate
functions over $X$ respectively. What conditions must one impose on $X$ in
order that the sets $\mathcal{R}_{a}(X)$ and $\mathcal{R}_{b}(X)$ be
proximinal in $B(X)$ and the set $\mathcal{R}_{c}(X)$ be proximinal in $C(X)$%
? We are also interested in necessary conditions for proximinality. It
follows from one result of Garkavi, Medvedev and Khavinson (see
\cite[Theorem 1]{35}) that $\mathcal{R}_{a}(X)$ is proximinal in $B(X)$ for all subsets
$X$ of $\mathbb{R}^{n}$. There is also an answer (see \cite[Theorem 2]{35})
for proximinality of $\mathcal{R}_{b}(X)$ in $B(X)$. This will be discussed
in Section 1.5.2. Is the set $\mathcal{R}_{b}(X)$ always proximinal in $B(X)$%
? There is an an example of a set $X\subset \mathbb{R}^{n}$ and a bounded
function $f$ on $X$ for which there does not exist an extremal element in $%
\mathcal{R}_{b}(X)$.

In Section 1.5.3, we will obtain sufficient conditions for the existence of
extremal elements from $\mathcal{R}_{c}(X)$ to an arbitrary function $f$ $%
\in $ $C(X)$. Based on one result of Marshall and O'Farrell \cite{108}, we
will also give a necessary condition for proximinality of $\mathcal{R}_{c}(X)
$ in $C(X)$. All the theorems, following discussions and examples of the
paper will lead us naturally to a conjecture on the proximinality of the
subspaces $\mathcal{R}_{b}(X)$ and $\mathcal{R}_{c}(X)$ in the spaces $B(X)$
and $C(X)$ respectively.

The reader may also be interested in the more general case with the set $%
\mathcal{R}=\mathcal{R}\left( \mathbf{a}^{1},...,\mathbf{a}^{r}\right)$. In
this case, the corresponding sets $\mathcal{R}_{a}(X)$, $\mathcal{R}_{b}(X)$
and $\mathcal{R}_{c}(X)$ are defined similarly. Using the results of \cite%
{35}, one can obtain sufficient (but not necessary) conditions for
proximinality of these sets. This needs, besides paths, the consideration of
some additional and more complicated relations between points of $X$. Here
we will not consider the case $r\geq 3$, since our main purpose is to draw
the reader's attention to the arisen problems of proximinality in the
simplest case of approximation. For the existing open problems connected
with the set $\mathcal{R}\left( \mathbf{a}^{1},...,\mathbf{a}^{r}\right) $,
where $r\geq 3$, see \cite{53} and \cite{118}.

\bigskip

\subsection{Proximinality of $\mathcal{R}_{b}(X)$ in $B(X)$}

Let $\mathbf{a}^{1}$ and $\mathbf{a}^{2}$ be two different directions in $%
\mathbb{R}^{n}$. In the sequel, we will use paths with respect to the
directions $\mathbf{a}^{1}$ and $\mathbf{a}^{2}$. Recall that a length of a
path is the number of its points and can be equal to $\infty $ if the path
is infinite. A singleton is a path of the unit length. We say that a path $%
\left( \mathbf{x}^{1},...,\mathbf{x}^{m}\right) $ belonging to some subset $X
$ of $\mathbb{R}^{n}$ is irreducible if there is not another path $\left(
\mathbf{y}^{1},...,\mathbf{y}^{l}\right) \subset X$ with $\mathbf{y}^{1}=%
\mathbf{x}^{1},~\mathbf{y}^{l}=\mathbf{x}^{m}$ and $l<m$.

The following theorem follows from \cite[Theorem 2]{35}.

\bigskip

\textbf{Theorem 1.9.} \textit{Let $X\subset $ $\mathbb{R}^{n}$ and the
lengths of all irreducible paths in $X$ be uniformly bounded by some
positive integer. Then each function in $B(X)$ has an extremal element in $%
\mathcal{R}_{b}(X)$.}

\bigskip

There are a large number of sets in $\mathbb{R}^{n}$ satisfying the
hypothesis of this theorem. For example, if a set $X$ has a cross section
according to one of the directions $\mathbf{a}^{1}$ or $\mathbf{a}^{2}$,
then the set $X$ satisfies the hypothesis of Theorem 1.9. By a cross section
according to the direction $\mathbf{a}^{1}$ we mean any set $X_{\mathbf{a}%
^{1}}=\{x\in X:\ \mathbf{a}^{1}\cdot \mathbf{x}=c\}$, $c\in \mathbb{R}$, with
the property: for any $\mathbf{y}\in X$ there exists a point $\mathbf{y}%
^{1}\in X_{\mathbf{a}^{1}}$ such that $\mathbf{a}^{2}\cdot \mathbf{y}=%
\mathbf{a}^{2}\cdot \mathbf{y}^{1}$. By the similar way, one can define a
cross section according to the direction $\mathbf{a}^{2}$. For more on cross
sections in problems of proximinality of sums of univariate functions see
\cite{34,77}. Regarding Theorem 1.9 one may ask if the condition of the
theorem is necessary for proximinality of $\mathcal{R}_{b}(X)$ in $B(X)$.
While we do not know a complete answer to this question, we are going to
give an example of a set $X $ for which Theorem 1.9 fails. Let $\mathbf{a}%
^{1}=(1;-1),\ \mathbf{a}^{2}=(1;1).$ Consider the set
\begin{eqnarray*}
X &=&\{(2;\frac{2}{3}),(\frac{2}{3};-\frac{2}{3}),(0;0),(1;1),(1+\frac{1}{2}%
;1-\frac{1}{2}),(1+\frac{1}{2}+\frac{1}{4};1-\frac{1}{2}+\frac{1}{4}), \\
&&(1+\frac{1}{2}+\frac{1}{4}+\frac{1}{8};1-\frac{1}{2}+\frac{1}{4}-\frac{1}{8%
}),...\}.
\end{eqnarray*}%
In what follows, the elements of $X$ in the given order will be denoted by $%
\mathbf{x}^{0},\mathbf{x}^{1},\mathbf{x}^{2},...$ . It is clear that $X$ is
a path of the infinite length and $\mathbf{x}^{n}\rightarrow \mathbf{x}^{0}$
as $n\rightarrow \infty $. Let $\sum_{n=1}^{\infty }c_{n}$ be any
divergent series with the terms $c_{n}>0$ and $c_{n}\rightarrow 0$ as $%
n\rightarrow \infty $. Besides let $f_{0}$ be a function vanishing at the
points $\mathbf{x}^{0},\mathbf{x}^{2},\mathbf{x}^{4},...,$ and taking values
$c_{1},c_{2},c_{3},...$ at the points $\mathbf{x}^{1},\mathbf{x}^{3},\mathbf{%
\ x}^{5},...$, respectively. It is obvious that $f_{0}$ is continuous on $X$.
The set $X$ is compact and satisfies all the conditions of Theorem 1.6. By
that theorem, $\overline{\mathcal{R}_{c}(X)}=C(X).$ Therefore, for any
continuous function on $X$, thus for $f_{0}$,
\begin{equation*}
\inf_{g\in \mathcal{R}_{c}(X)}\left\Vert f_{0}-g\right\Vert _{C(X)}=0.\eqno%
(1.47)
\end{equation*}
Since $\mathcal{R}_{c}(X)\subset \mathcal{R}_{b}(X),$ we obtain from (1.47)
that
\begin{equation*}
\inf_{g\in \mathcal{R}_{b}(X)}\left\Vert f_{0}-g\right\Vert _{B(X)}=0.\eqno%
(1.48)
\end{equation*}

Suppose that $f_{0}$ has an extremal element ${g_{1}^{0}\left( \mathbf{a}^{1}%
{\cdot }\mathbf{x}\right) +g_{2}^{0}\left( \mathbf{a}^{2}{\ \cdot }\mathbf{x}%
\right) }$ in $\mathcal{R}_{b}(X).$ By the definition of $\mathcal{R}_{b}(X)$%
, the ridge functions ${g_{i}^{0},i=1,2}$, are bounded on $X.$ From (1.48)\
it follows that $f_{0}={g_{1}^{0}\left( \mathbf{a}^{1}{\ \cdot }\mathbf{x}%
\right) +g_{2}^{0}\left( \mathbf{a}^{2}{\cdot }\mathbf{x}\right) .}$ Since $%
\mathbf{a}^{1}\cdot \mathbf{x}^{2n}=\mathbf{a}^{1}\cdot \mathbf{x}^{2n+1}$
and $\mathbf{a}^{2}\cdot \mathbf{x}^{2n+1}=\mathbf{a}^{2}\cdot \mathbf{x}%
^{2n+2},$ for $n=0,1,...,$ we can write that
\begin{equation*}
\sum_{n=0}^{k}c_{n+1}=\sum_{n=0}^{k}\left[ f(\mathbf{x}^{2n+1})-f(\mathbf{x}%
^{2n})\right]
\end{equation*}
\begin{equation*}
=\sum_{n=0}^{k}\left[ {g_{2}^{0}}(\mathbf{x}^{2n+1})-{g_{2}^{0}}(\mathbf{x}%
^{2n})\right] ={g_{2}^{0}(}\mathbf{a}^{2}\cdot \mathbf{x}^{2k+1})-{g_{2}^{0}(%
}\mathbf{a}^{2}\cdot \mathbf{x}^{0}).\eqno(1.49)
\end{equation*}
Since $\sum_{n=1}^{\infty }c_{n}=\infty ,$ we deduce from (1.49)\ that the
function ${g_{2}^{0}\left( \mathbf{a}^{2}{\cdot }\mathbf{x}\right) }$ is not
bounded on $X.$ This contradiction means that the function $f_{0}$ does not
have an extremal element in $\mathcal{R}_{b}(X).$ Therefore, the space $%
\mathcal{R}_{b}(X)$ is not proximinal in $B(X).$

\bigskip

\subsection{Proximinality of $\mathcal{R}_{c}(X)$ in $C(X)$}

In this section, we give a sufficient condition and also a necessary condition
for proximinality of $\mathcal{R}_{c}(X)$ in $C(X)$.

\bigskip

\textbf{Theorem 1.10.} \textit{Let the system of linearly independent vectors $%
\mathbf{a}^{1}$ and $\mathbf{a}^{2}$ have a complement to a basis $\{\mathbf{a%
}^{1},...,\mathbf{a}^{n}\}$ in $\mathbb{R}^{n}$ with the property: for any
point $\mathbf{x}^{0}\in X$ and any positive real number $\delta $ there
exist a number $\delta _{0}\in (0,\delta ]$ and a point $\mathbf{x}^{\sigma
} $ in the set}
\begin{equation*}
\sigma =\{\mathbf{x}\in X:\mathbf{a}^{2}\cdot \mathbf{x}^{0}-\delta _{0}\leq
\mathbf{a}^{2}\cdot \mathbf{x}\leq \mathbf{a}^{2}\cdot \mathbf{x}^{0}+\delta
_{0}\},
\end{equation*}
\textit{such that the system}
\begin{equation*}
\left\{
\begin{array}{c}
\mathbf{a}^{2}\cdot \mathbf{x}^{\prime }=\mathbf{a}^{2}\cdot \mathbf{x}%
^{\sigma } \\
\mathbf{a}^{1}\cdot \mathbf{x}^{\prime }=\mathbf{a}^{1}\cdot \mathbf{x} \\
\sum_{i=3}^{n}\left\vert \mathbf{a}^{i}\cdot \mathbf{x}^{\prime }-\mathbf{a}%
^{i}\cdot \mathbf{x}\right\vert <\delta%
\end{array}%
\right. \eqno(1.50)
\end{equation*}%
\textit{has a solution $\mathbf{x}^{\prime }\in \sigma $ for all points $%
\mathbf{x}\in \sigma .$Then the space $\mathcal{R}_{c}(X)$ is proximinal in $%
C(X).$}

\bigskip

\begin{proof}
Introduce the following mappings and sets:
\begin{equation*}
\pi _{i}:X\rightarrow \mathbb{R}\text{, }\pi _{i}(\mathbf{x)=a}^{i}\cdot
\mathbf{x}\text{, }Y_{i}=\pi _{i}(X\mathbf{)}\text{, }i=1,...,n.
\end{equation*}

Since the system of vectors $\{\mathbf{a}^{1},...,\mathbf{a}^{n}\}$ is
linearly independent, the mapping $\pi =(\pi _{1},...\pi _{n})$ is an
injection from $X$ into the Cartesian product $Y_{1} \times ...\times Y_{n} $
. Besides, $\pi $ is linear and continuous. By the open mapping theorem, the
inverse mapping $\pi ^{-1}$ is continuous from $Y=\pi (X)$ onto $X.$ Let $f$
be a continuous function on $X$. Then the composition $f\circ \pi
^{-1}(y_{1},...y_{n})$ will be continuous on $Y,$ where $y_{i}=\pi _{i}(
\mathbf{x),}\ i=1,...,n,$ are the coordinate functions. Consider the
approximation of the function $f\circ \pi ^{-1}$ by elements from
\begin{equation*}
G_{0}=\{g_{1}(y_{1})+g_{2}(y_{2}):\ g_{i}\in C(Y_{i}),\ i=1,2\}
\end{equation*}
over the compact set $Y$. Then one may observe that the function $f$ has an
extremal element in $\mathcal{R}_{c}(X)$ if and only if the function $f\circ
\pi ^{-1}$ has an extremal element in $G_{0}$. Thus the problem of
proximinality of $\mathcal{R}_{c}(X)$ in $C(X)$ is reduced to the problem of
proximinality of $G_{0}$ in $C(Y).$

Let $T,T_{1},...,T_{m+1}$ be metric compact spaces and $T\subset $ $%
T_{1}\times ...\times T_{m+1}.$ For $i=1,...,m,$ let $\varphi _{i}$ be the
continuous mappings from $T$ onto $T_{i}.$ In \cite{35}, the authors obtained
sufficient conditions for proximinality of the set

\begin{equation*}
C_{0}=\{\sum_{i=1}^{n}g_{i}\circ \varphi _{i}:\ g_{i}\in C(T_{i}),\
i=1,...m\}
\end{equation*}%
in the space $C(T)$ of continuous functions on $T.$ Since $Y\subset $ $%
Y_{1}\times Y_{2}\times Z_{3},$ where $Z_{3}=Y_{3}\times ...\times Y_{n},$
we can use this result in our case, for the approximation of the function $%
f\circ \pi ^{-1}$ by elements from $G_{0}$. By this theorem, the set $G_{0}$
is proximinal in $C(Y)$ if for any $y_{2}^{0}\in Y_{2}$ and $\delta >0$
there exists a number $\delta _{0}\in (0,$ $\delta )$ such that the set $%
\sigma (y_{2}^{0},\delta _{0})=[y_{2}^{0}-\delta _{0},y_{2}^{0}+\delta
_{0}]\cap Y_{2}$ has $(2,\delta )$ maximal cross section. The last means
that there exists a point $y_{2}^{\sigma }\in \sigma (y_{2}^{0},\delta _{0})$
with the property: for any point $(y_{1},y_{2},z_{3})\in Y,$ with the second
coordinate $y_{2}$ from the set $\sigma (y_{2}^{0},\delta _{0}),$ there
exists a point $(y_{1}^{\prime },y_{2}^{\sigma },z_{3}^{\prime })\in Y$ such
that $y_{1}=y_{1}^{\prime }$ and $\rho (z_{3},z_{3}^{\prime })<\delta ,$
where $\rho $ is a metrics in $Z_{3}.$ Since these conditions are equivalent
to the conditions of Theorem 1.10, the space $G_{0}$ is proximinal in the
space $C(Y).$ Then by the above conclusion, the space $\mathcal{R}_{c}(X)$
is proximinal in $C(X).$ \end{proof}

Let us give some simple examples of compact sets satisfying the hypothesis
of Theorem 1.10. For the sake of brevity, we restrict ourselves to the case $%
n=3.$

\begin{enumerate}
\item[(a)] Assume $X$ is a closed ball in $\mathbb{R}^{3}$ and $\mathbf{a}^{1}$, $\mathbf{a}^{2}$
are orthogonal directions. Then Theorem 1.10 holds. Note that
in this case, we can take $\delta _{0}=\delta $ and $\mathbf{a}^{3}$ as an orthogonal
vector to both the vectors $\mathbf{a}^{1}$ and $\mathbf{a}^{2}.$

\item[(b)] Let $X$ be the unite cube, $\mathbf{a}^{1}=(1;1;0),\ a^{2}=(1;-1;0).$ Then
Theorem 1.10 also holds. In this case, we can take $\delta _{0}=\delta $ and
$\mathbf{a}^{3}=(0;0;1).$ Note that the unit cube does not satisfy the hypothesis of
the theorem for many directions (take, for example, $\mathbf{a}^{1}=(1;2;0)$ and $\mathbf{a}^{2}=(2;-1;0)$).
\end{enumerate}

In the following example, one can not always chose $\delta _{0}$ as equal to
$\delta $.

\begin{enumerate}
\item[(c)] Let $X=\{(x_{1},x_{2},x_{3}):\ (x_{1},x_{2})\in Q,\ 0\leq
x_{3}\leq 1\},$ where $Q$ is the union of two triangles $A_{1}B_{1}C_{1}$
and $A_{2}B_{2}C_{2}$ with the vertices $A_{1}=(0;0),\ B_{1}=(1;2),\
C_{1}=(2;0),\ A_{2}=(1\frac{1}{2};1),\ B_{2}=(2\frac{1}{2};-1),\ C_{2}=(3%
\frac{1}{2};1).$ Let $\mathbf{a}^{1}=(0;1;0)$ and $\mathbf{a}^{2}=(1;0;0).$ Then it is easy to
see that Theorem 1.10 holds (the vector $\mathbf{a}^{3}$ can be chosen as $(0;0;1)$).
In this case, $\delta _{0}$ can not be always chosen as equal to $\delta $.
Take, for example, $\mathbf{x}^{0}=(1\frac{3}{4};0;0)$ and $\delta =1\frac{3%
}{4}.$ If $\delta _{0}=\delta ,$ then the second equation of the system
(1.50) has not a solution for a point $(1;2;0)$ or a point $(2\frac{1}{2}%
;-1;0).$ But if we take $\delta _{0}$ not more than $\frac{1}{4}$, then for $%
\mathbf{x}^{\sigma }=\mathbf{x}^{0}$ the system has a solution. Note that
the last inequality $\left\vert \mathbf{a}^{3}\cdot \mathbf{x}^{\prime }-%
\mathbf{a}^{3}\cdot \mathbf{x}\right\vert <\delta $ of the system can be
satisfied with the equality $\mathbf{a}^{3}\cdot \mathbf{x}^{\prime }=%
\mathbf{a}^{3}\cdot \mathbf{x}$ if $\mathbf{a}^{3}=(0;0;1).$
\end{enumerate}

It should be remarked that the results of \cite{35} tell nothing about
necessary conditions for proximinality of the spaces considered there. To
fill this gap in our case, we want to give a necessary condition for
proximinality of $\mathcal{R}_{c}(X)$ in $C(X)$. First, let us introduce
some notation. By $\mathcal{R}_{c}^{i},\ i=1,2,$ we will denote the set of
continuous ridge functions $g\left( \mathbf{a}^{i}\cdot \mathbf{x}\right) $
on the given compact set $X\subset \mathbb{R}^{n}.$ Note that $\mathcal{R}%
_{c}=\mathcal{R}_{c}^{1}+\mathcal{R}_{c}^{2}.$ Besides, let $\mathcal{R}%
_{c}^{3}=\mathcal{R}_{c}^{1}\cap \mathcal{R}_{c}^{2}.$ For $i=1,2,3,$ let $%
X_{i}$ be the quotient space obtained by identifying points $y_{1}$ and $%
y_{2}$ in $X$ whenever $f(y_{1})=f(y_{2})$ for each $f$ in $\mathcal{R}%
_{c}^{i}.$ By $\pi _{i}$ denote the natural projection of $X$ onto $X_{i},$ $%
i=1,2,3.$ Note that we have already dealt with the quotient spaces $X_{1}$, $%
X_{2}$ and the projections $\pi _{1},\pi _{2}$ in the previous section.
Recall that the relation on $X$, defined by setting $\ y_{1}\approx y_{2}$
if $y_{1}$ and $y_{2}$ belong to some path, is an equivalence relation and
the equivalence classes are called orbits. By $O(t)$ denote the orbit of $X$
containing $t.$ For $Y\subset X,$ let $var_{Y}\ f$ be the variation of a
function $f$ on the set $Y.$ That is,
\begin{equation*}
\underset{Y}{var}f=\sup\limits_{x,y\in Y}\left\vert f\left( x\right)
-f\left( y\right) \right\vert .
\end{equation*}

The following theorem is valid.

\bigskip

\textbf{Theorem 1.11.} \textit{Suppose that the space $\mathcal{R}_{c}(X)$
is proximinal in $C(X).$Then there exists a positive real number c such that}
\begin{equation*}
\sup_{t\in X}\underset{O\left( t\right) }{var}\mathit{f\leq c}\sup_{t\in X}%
\underset{\pi _{2}^{-1}\left( \pi _{2}\left( t\right) \right) }{var}\mathit{f%
}\eqno(1.51)
\end{equation*}%
\textit{for all $f$ in $\mathcal{R}_{c}^{1}.$}

\bigskip

\begin{proof}
The proof is based on the
following result of Marshall and O'Farrell (see \cite[Proposition 4]{108}): Let $A_{1}\ $and $A_{2}\
$be closed subalgebras of $C(X)\ $that contain the constants. Let $%
(X_{1},\pi _{1}),\ (X_{2},\pi _{2})\ $and $(X_{3},\pi _{3})\ $be the
quotient spaces and projections associated with the algebras $A_{1},$ $%
A_{2}\ $and $A_{3}=A_{1}\cap A_{2}\ $respectively. Then $A_{1}+A_{2}\ $is
closed in $C(X)\ $if and only if there exists a positive real number $c$
such that
\begin{equation*}
\sup\limits_{z\in X_{3}}\underset{\pi _{3}^{-1}\left( z\right) }{var}f\leq
c\sup\limits_{y\in X_{2}}\underset{\pi _{2}^{-1}\left( y\right) }{var}f\eqno%
(1.52)
\end{equation*}%
for all $f\ $in $A_{1}.$

If $\mathcal{R}_{c}(X)$ is proximinal in $C(X),$ then it is necessarily
closed and therefore, by the above proposition, (1.52) holds for the
algebras $A_{1}^{i}=\mathcal{R}_{c}^{i},\ i=1,2,3.$ The right-hand side of
(1.52) is equal to the right-hand side of (1.51). Let $t$ be some point in $%
X $ and $z=\pi _{3}(t).$ Since each function $\ f\in \mathcal{R}_{c}^{3}$ is
constant on the orbit of $t$ (note that $f$ is both of the form ${\
g_{1}\left( \mathbf{a}^{1}{\cdot }\mathbf{x}\right) }$ and of the form ${\
g_{2}\left( \mathbf{a}^{2}{\cdot }\mathbf{x}\right) }$), $O(t)\subset \pi
_{3}^{-1}(z).$ Hence,
\begin{equation*}
\sup_{t\in X}\underset{O\left( t\right) }{var}f\leq c\sup\limits_{z\in X_{3}}%
\underset{\pi _{3}^{-1}\left( z\right) }{var}f\eqno(1.53)
\end{equation*}

From (1.52) and (1.53) we obtain (1.51).
\end{proof}

Note that the inequality (1.52) provides not worse but less practicable
necessary condition for proximinality than the inequality (1.51) does. On
the other hand, there are many cases in which both the inequalities are
equivalent. For example, assume the lengths of irreducible paths of $X$ are
bounded by some positive integer $n_{0}$. In this case, it can be shown that
the inequality (1.52), hence (1.51), holds with the constant $c=\frac{n_{0}}{%
2}$ and moreover $O(t)=\pi _{3}^{-1}(z)$ for all $t\in X$, where $z=\pi
_{3}(t)$ (see the proof of \cite[Theorem 5]{53}). Therefore, the
inequalities (1.51) and (1.52) are equivalent for the considered class of
sets $X.$ The last argument shows that all the compact sets $X\subset $ $%
\mathbb{R}^{n}$ over which $\mathcal{R}_{c}(X)$ is not proximinal in $C(X)$
should be sought in the class of sets having irreducible paths consisting of
sufficiently many points. For example, let $I=[0;1]^{2}$ be the
unit square, $\mathbf{a}^{1}=(1;1)$, $\mathbf{a}^{2}=(1;\frac{1}{2}).$ Consider the path
\begin{equation*}
l_{k}=\{(1;0),(0;1),(\frac{1}{2};0),(0;\frac{1}{2}),(\frac{1}{4};0),...,(0;
\frac{1}{2^{k}})\}.
\end{equation*}

It is clear that $l_{k}$ is an irreducible path with the length $2k+2 $,
where $k$ may be very large. Let $g_{k}$ be a continuous univariate function
on $\mathbb{R}$ satisfying the conditions: $g_{k}(\frac{1}{2^{k-i}} )=i,\
i=0,...,k,$ $g_{k}(t)=0$ if $t<\frac{1}{2^{k}},\ i-1\leq g_{k}(t)\leq i $ if
$t\in (\frac{1}{2^{k-i+1}},\frac{1}{2^{k-i}}),\ i=1,...,k,$ and $g_{k}(t)=k$
if $t>1.$ Then it can be easily verified that
\begin{equation*}
\sup_{t\in X}\underset{\pi _{2}^{-1}\left( \pi _{2}\left( t\right) \right) }{%
var}g_{k}(\mathbf{a}^{1}{\cdot }\mathbf{x})\leq 1.\eqno(1.54)
\end{equation*}

Since $\max_{\mathbf{x}\in I}g_{k}(\mathbf{a}^{1}{\cdot }\mathbf{x} )=k,$ $%
\min_{\mathbf{x}\in I}g_{k}(\mathbf{a}^{1}{\cdot }\mathbf{x})=0$ and $var_{%
\mathbf{x}\in O\left( t_{1}\right) }g_{k}(\mathbf{a}^{1}{\cdot }\mathbf{\ x}%
)=k $ for $t_{1}=(1;0),$ we obtain that
\begin{equation*}
\sup_{t\in X}\underset{O\left( t\right) }{var}g_{k}(\mathbf{a}^{1}{\cdot }%
\mathbf{x})=k.\eqno(1.55)
\end{equation*}

Since $k$ may be very large, from (1.54) and (1.55) it follows that the
inequality (1.51) cannot hold for the function $g_{k}(\mathbf{a}^{1}{\ \cdot
}\mathbf{x})\in \mathcal{R}_{c}^{1}.$ Thus the space $\mathcal{R}_{c}(I)$
with the directions $\mathbf{a}^{1}=(1;1)$ and $\mathbf{a}^{2}=(1;\frac{1}{2})$ is not
proximinal in $C(I)$.

It should be remarked that if a compact set $X\subset $ $\mathbb{R}^{n}$
satisfies the hypothesis of Theorem 1.10, then the length of all irreducible
paths are uniformly bounded (see the proof of Theorem 1.10 and lemma in \cite%
{35}). We have already seen that if the last condition does not hold, then
the proximinality of both $\mathcal{R}_{c}(X)$ in $C(X)$ and $\mathcal{R}%
_{b}(X)$ in $B(X)$ fail for some sets $X$. In addition to the examples given above
and in Section 1.5.2, one can easily construct many other examples of such
sets. All these examples, Theorems 1.9--1.11 and the subsequent
remarks justify the statement of the following conjecture:

\bigskip

\textbf{Conjecture. }\textit{Let $X$ be some subset of $\mathbb{R} ^{n}.$
The space $\mathcal{R}_{b}(X)$ is proximinal in $B(X)$ and the space $%
\mathcal{R}_{c}(X)$ is proximinal in $C(X)$ (in this case, $X$ is considered
to be compact) if and only if the lengths of all irreducible paths of $X$
are uniformly bounded.}

\bigskip

\textbf{Remark 1.2.} Medvedev's result (see \cite[p.58]{76}), which later
came to our attention, in particular, says that the set $R_{c}(X)$ is closed
in $C(X)$ if and only if the lengths of all irreducible paths of $X$ are
uniformly bounded. Thus, in the case of $C(X)$, the necessity of the above
conjecture was proved by Medvedev.

\bigskip

\textbf{Remark 1.3.} Note that there are situations in which a continuous
function (a specific function on a specially constructed set) has an
extremal element in $\mathcal{R}_{b}(X)$, but not in $\mathcal{R}_{c}(X)$
(see \cite[p.73]{76}). One subsection of \cite{76} (see p.68 there) was
devoted to the proximinality of sums of two univariate functions with
continuous and bounded summands in the spaces of continuous and bounded
bivariate functions, respectively. If $X\subset \mathbb{R}^{2}$ and $\mathbf{%
a}^{1},\mathbf{a}^{2} $ be linearly independent directions in $\mathbb{R}^{2}
$, then the linear transformation $y_{1}=$ $\mathbf{a}^{1}\cdot \mathbf{x\,}$%
, $y_{2}=$ $\mathbf{a}^{2}\cdot \mathbf{x}$ reduces the problems of
proximinality of $\mathcal{R}_{b}(X)$ in $B(X)$ and $\mathcal{R}_{c}(X)$ in $%
C(X)$ to the problems considered in that subsection. But in general, when $%
X\subset \mathbb{R}^{n}$ and $n>2$, they cannot be reduced to those in
\cite{76}.

\bigskip

\section{On the approximation by weighted ridge functions}

In this section, we characterize the best $L_{2}$-approximation to a
multivariate function by linear combinations of ridge functions multiplied
by some fixed weight functions. In the special case, when the weight
functions are constants, we obtain explicit formulas for both the best
approximation and approximation error.

\subsection{Problem statement}

Ridge approximation in $L_{2}$ started to be actively studied in the late 90's by
K.I. Oskolkov \cite{114}, V.E. Maiorov \cite{102}, A. Pinkus \cite{118},
V.N. Temlyakov \cite{138}, P. Petrushev \cite{116} and other researchers.

Let $D$ be the unit disk in $\mathbb{R}^{2}$. In \cite{97}, Logan and Shepp
along with other results gave a closed-form expression for the best $L_{2}$-approximation
to a function $f \in L_{2}\left(
D\right)$ from the set $\mathcal{R}\left( \mathbf{a}^{1},...,\mathbf{a}^{r}\right)$.
Their solution requires that the directions $\mathbf{a}^{1},...,\mathbf{a}%
^{r}$ be equally-spaced and involves finite sums of convolutions with
explicit kernels. In the $n$-dimensional case, we obtained an expression of
simpler form for the best $L_{2}$-approximation to square-integrable
multivariate functions over a certain domain, provided that $r=n$ and the
directions $\mathbf{a}^{1},...,\mathbf{a}^{r}$ are linearly independent (see
\cite{52}).

In this section, we consider the approximation by functions from the following more general set
\begin{equation*}
\mathcal{R}\left( \mathbf{a}^{1},...,\mathbf{a}^{r};~w_{1},...,w_{r}\right)
=\left\{ \sum\limits_{i=1}^{r}w_{i}(\mathbf{x})g_{i}\left( \mathbf{a}%
^{i}\cdot \mathbf{x}\right) :g_{i}:\mathbb{R}\rightarrow \mathbb{R}%
,~i=1,...,r\right\} ,
\end{equation*}%
where $w_{1},...,w_{r}$ are fixed multivariate functions. We
characterize the best $L_{2}$-approximation from this set in the case $%
r\leq n.$ Then, in the special case when the weight functions $%
w_{1},...,w_{r}$ are constants, we will prove two theorems on explicit
formulas for the best approximation and the approximation error,
respectively. At present, we do not yet know how to approach
these problems in other possible cases of $r.$

\bigskip

\subsection{Characterization of the best approximation}

Let $X$ be a subset of $\mathbb{R}^{n}$ with a finite Lebesgue measure.
Consider the approximation of a function $f\left( \mathbf{x}\right) =f\left(
x_{1},...,x_{n}\right) $ in $L_{2}\left( X\right) $ by functions from the manifold $%
\mathcal{R}\left( \mathbf{a}^{1},...,\mathbf{a}^{r};~w_{1},...,w_{r}\right) $%
, where $r\leq n.$ We suppose that the functions $w_{i}(\mathbf{x})$ and the
products $w_{i}(\mathbf{x})\cdot g_{i}\left( \mathbf{a}^{i}\cdot \mathbf{x}%
\right) ,~i=1,...,r$, belong to the space $L_{2}\left( X\right) .$ Besides,
we assume that the vectors $\mathbf{a}^{1},...,\mathbf{a}^{r}$ are linearly
independent. We say that a function $g_{w}^{0}=\sum\limits_{i=1}^{r}w_{i}(%
\mathbf{x})g_{i}^{0}\left( \mathbf{a}^{i}\cdot \mathbf{x}\right) $ in $%
\mathcal{R}\left( \mathbf{a}^{1},...,\mathbf{a}^{r};~w_{1},...,w_{r}\right) $
is the best approximation (or extremal) to $f$ if
\begin{equation*}
\left\Vert f-g_{w}^{0}\right\Vert _{L_{2}\left( X\right) }=\inf\limits_{g\in
\mathcal{R}\left( \mathbf{a}^{1},...,\mathbf{a}^{r};~w_{1},...,w_{r}\right)
}\left\Vert f-g\right\Vert _{L_{2}\left( X\right) }.
\end{equation*}

Let the system of vectors $\{\mathbf{a}^{1},...,\mathbf{a}^{r},\mathbf{a%
}^{r+1},...,\mathbf{a}^{n}\}$ be a completion of the system $\{\mathbf{a}%
^{1},...,\mathbf{a}^{r}\}$ to a basis in $\mathbb{R}^{n}.$\ Let $%
J:X\rightarrow \mathbb{R}^{n}$ be the linear transformation given by the
formulas
\begin{equation*}
y_{i}=\mathbf{a}^{i}\cdot \mathbf{x,}\quad \,i=1,...,n.\eqno(1.56)
\end{equation*}%
Since the vectors $\mathbf{a}^{i},$ $i=1,...,n$, are linearly independent,
it is an injection. The Jacobian $\det J$ of this transformation is a
constant different from zero.

Let the formulas
\begin{equation*}
x_{i}=\mathbf{b}^{i}\cdot \mathbf{y},\;\;i=1,...,n,
\end{equation*}%
stand for the solution of linear equations (1.56) with respect to $%
x_{i},\;i=1,...,n.$

Introduce the notation
\begin{equation*}
Y=J\left( X\right)
\end{equation*}%
and
\begin{equation*}
Y_{i}=\left\{ y_{i}\in \mathbb{R}:\;\;y_{i}=\mathbf{a}^{i}\cdot \mathbf{x}%
,\;\;\mathbf{x}\in X\right\} ,\,i=1,...,n.
\end{equation*}
For any function $u\in L_{2}\left( X\right) ,$ put
\begin{equation*}
u^{\ast }=u^{\ast }\left( \mathbf{y}\right) \overset{def}{=}u\left( \mathbf{b%
}^{1}\cdot \mathbf{y},...,\mathbf{b}^{n}\cdot \mathbf{y}\right) .
\end{equation*}%
It is obvious that $u^{\ast }\in L_{2}\left( Y\right) .$ Besides,
\begin{equation*}
\int\limits_{Y}u^{\ast }\left( \mathbf{y}\right) d\mathbf{y}=\left\vert \det
J\right\vert \cdot \int\limits_{X}u\left( \mathbf{x}\right) d\mathbf{x}\eqno%
(1.57)
\end{equation*}%
and
\begin{equation*}
\left\Vert u^{\ast }\right\Vert _{L_{2}\left( Y\right) }=\left\vert \det
J\right\vert ^{1/2}\cdot \left\Vert u\right\Vert _{L_{2}\left( X\right) }.%
\eqno(1.58)
\end{equation*}%

Set
\begin{equation*}
L_{2}^{i}=\{w_{i}^{\ast }(\mathbf{y})g\left( y_{i}\right) \in
L_{2}(Y)\},~i=1,...,r.
\end{equation*}

We need the following auxiliary lemmas.

\bigskip

\textbf{Lemma 1.5.} \textit{Let $f\left( \mathbf{x}\right) \in
L_{2}\left( X\right) $. A function $\sum\limits_{i=1}^{r}w_{i}(\mathbf{x}%
)g_{i}^{0}\left( \mathbf{a}^{i}\cdot \mathbf{x}\right) $ is extremal to the
function $f\left( \mathbf{x}\right) $ if and only if \ $\sum%
\limits_{i=1}^{r}w_{i}^{\ast }(\mathbf{y})g_{i}^{0}\left( y_{i}\right) $ is
extremal from the space $L_{2}^{1}\mathit{\oplus }...\oplus L_{2}^{r}$ to
the function $f^{\ast }\left( \mathbf{y}\right) $.}

\bigskip

Due to (1.58) the proof of this lemma is obvious.

\bigskip

\textbf{Lemma 1.6.} \textit{Let $f\left( \mathbf{x}\right) \in
L_{2}\left( X\right) $. A function $\sum\limits_{i=1}^{r}w_{i}(\mathbf{x}%
)g_{i}^{0}\left( \mathbf{a}^{i}\cdot \mathbf{x}\right) $ is extremal to the
function $f\left( \mathbf{x}\right) $ if and only if}
\begin{equation*}
\int\limits_{X}\left( f\left( \mathbf{x}\right) -\sum\limits_{i=1}^{r}w_{i}(%
\mathbf{x})g_{i}^{0}\left( \mathbf{a}^{i}\cdot \mathbf{x}\right) \right)
w_{j}(\mathbf{x})h\left( \mathbf{a}^{j}\cdot \mathbf{x}\right) d\mathbf{x}%
=0\
\end{equation*}%
\textit{for any ridge function $h\left( \mathbf{a}^{j}\cdot \mathbf{x}%
\right) $ such that $\mathit{w}_{j}\mathit{(x)h}\left( \mathbf{a}^{j}\cdot
\mathbf{x}\right) $$\in L_{2}\left( X\right)$, $j=1,...,r$.}

\bigskip

\textbf{Lemma 1.7.} \textit{The following formula is valid for the error
of approximation to a function $f\left( \mathbf{x}\right) $ in $L_{2}\left(
X\right) $ from $\mathcal{R}\left( \mathbf{a}^{1},...,\mathbf{a}%
^{r};~w_{1},...,w_{r}\right) $:}
\begin{equation*}
E\left( f\right) =\left( \left\Vert f\left( \mathbf{x}\right) \right\Vert
_{L_{2}\left( X\right) }^{2}-\left\Vert \sum\limits_{i=1}^{r}w_{i}(\mathbf{x}%
)g_{i}^{0}\left( \mathbf{a}^{i}\cdot \mathbf{x}\right) \right\Vert
_{L_{2}\left( X\right) }^{2}\right) ^{\frac{1}{2}},
\end{equation*}%
\textit{where $\sum\limits_{i=1}^{r}w_{i}(\mathbf{x})g_{i}^{0}\left( \mathbf{%
a}^{i}\cdot \mathbf{x}\right) $ is the best approximation to $f\left(
\mathbf{x}\right) $.}

\bigskip

Lemmas 1.6 and 1.7 follow from the well-known facts of functional
analysis that the best approximation of an element $x$ in a Hilbert space $H$
from a linear subspace $Z$ of $H$ must be the image of $x$ via the
orthogonal projection onto $Z$ and the sum of squares of norms of orthogonal
vectors is equal to the square of the norm of their sum.

We say that $Y$ is an \textit{$r$-set} if it can be represented as $Y_{1}\times
...\times Y_{r}\times Y_{0},$ where $Y_{0}$ is some set from the space $%
\mathbb{R}^{n-r}.$ In a special case, $Y_{0}$ may be equal to $Y_{r+1}\times
...\times Y_{n},$ but it is not necessary. By $Y^{\left( i\right) },$ we
denote the Cartesian product of the sets $Y_{1},...,Y_{r},Y_{0}$ except for $%
Y_{i},\;i=1,...,r$. That is, $Y^{\left( i\right) }=Y_{1}\times ...\times
Y_{i-1}\times Y_{i+1}\times ...\times Y_{r}\times Y_{0},\,\ i=1,...,r$.

\bigskip

\textbf{Theorem 1.12.} \textit{Let $Y$ be an $r$-set. A function
$\sum\limits_{i=1}^{r}w_{i}(\mathbf{x})g_{i}^{0}\left( \mathbf{a}^{i}\cdot
\mathbf{x}\right) $ is the best approximation to $f(\mathbf{x)}$ if and only
if}
\begin{equation*}
g_{j}^{0}\left( y_{j}\right) =\frac{1}{\int\limits_{Y^{\left( j\right)
}}w_{j}^{\ast 2}(\mathbf{y})d\mathbf{y}^{\left( j\right) }}%
\int\limits_{Y^{\left( j\right) }}\left( f^{\ast }\left( \mathbf{y}\right)
-\sum\limits_{\substack{ i=1  \\ i\neq j}}^{r}w_{i}^{\ast }(\mathbf{y}%
)g_{i}^{0}\left( y_{i}\right) \right) w_{j}^{\ast }(\mathbf{y})d\mathbf{y}%
^{\left( j\right) },\eqno(1.59)
\end{equation*}%
\textit{for $j=1,...,r$.}

\bigskip

\begin{proof} \textit{Necessity.} Let a function $\sum\limits_{i=1}^{r}w_{i}(%
\mathbf{x})g_{i}^{0}\left( \mathbf{a}^{i}\cdot \mathbf{x}\right) $ be
extremal to $f$. Then by Lemma 1.5, the function $\sum%
\limits_{i=1}^{r}w_{i}^{\ast }(\mathbf{y})g_{i}^{0}\left( y_{i}\right) $ in $%
L_{2}^{1}\oplus ...\oplus L_{2}^{r}$ is extremal to $f^{\ast }$. By Lemma
1.6 and equality (1.57),
\begin{equation*}
\int\limits_{Y}f^{\ast }\left( \mathbf{y}\right) w_{j}^{\ast }(\mathbf{y}%
)h\left( y_{j}\right) d\mathbf{y}=\int\limits_{Y}w_{j}^{\ast }(\mathbf{y}%
)h\left( y_{j}\right) \sum\limits_{i=1}^{r}w_{i}^{\ast }(\mathbf{y}%
)g_{i}^{0}\left( y_{i}\right) d\mathbf{y}\eqno(1.60)
\end{equation*}%
for any product $w_{j}^{\ast }(\mathbf{y})h\left( y_{j}\right) $ in $%
L_{2}^{j},\;\;j=1,...,r$. Applying Fubini's theorem to the integrals in
(1.60), we obtain that
\begin{eqnarray*}
&&\int\limits_{Y_{j}}h\left( y_{j}\right) \left[ \int\limits_{Y^{\left(
j\right) }}f^{\ast }\left( \mathbf{y}\right) w_{j}^{\ast }(\mathbf{y})d%
\mathbf{y}^{\left( j\right) }\right] dy_{j} \\
&=&\int\limits_{Y_{j}}h\left( y_{j}\right) \left[ \int\limits_{Y^{\left(
j\right) }}w_{j}^{\ast }(\mathbf{y})\sum\limits_{i=1}^{r}w_{i}^{\ast }(%
\mathbf{y})g_{i}^{0}\left( y_{i}\right) d\mathbf{y}^{\left( j\right) }\right]
dy_{j}.
\end{eqnarray*}%
Since $h\left( y_{j}\right) $ is an arbitrary function such that $%
w_{j}^{\ast }(\mathbf{y})h\left( y_{j}\right) \in L_{2}^{j}$,
\begin{equation*}
\int\limits_{Y^{\left( j\right) }}f^{\ast }\left( \mathbf{y}\right)
w_{j}^{\ast }(\mathbf{y})d\mathbf{y}^{(j)}=\int\limits_{Y^{\left( j\right)
}}w_{j}^{\ast }(\mathbf{y})\sum\limits_{i=1}^{r}w_{i}^{\ast }(\mathbf{y}%
)g_{i}^{0}\left( y_{i}\right) d\mathbf{y}^{\left( j\right) },\;\;j=1,...,r.
\end{equation*}%
Therefore,
\begin{equation*}
\int\limits_{Y^{\left( j\right) }}w_{j}^{\ast 2}(\mathbf{y})g_{j}^{0}\left( {%
y_{j}}\right) d\mathbf{y}^{\left( j\right) }=\int\limits_{Y^{\left( j\right)
}}\left( f^{\ast }\left( \mathbf{y}\right) -\sum\limits_{\substack{ i=1  \\ %
i\neq j}}^{r}w_{i}^{\ast }(\mathbf{y})g_{i}^{0}\left( y_{i}\right) \right)
w_{j}^{\ast }(\mathbf{y})d\mathbf{y}^{\left( j\right) },
\end{equation*}%
for $j=1,...,r.$ Now, since $y_{j}\notin Y^{\left( j\right) }$, we obtain (1.59).

\textit{Sufficiency.} Note that all the equalities in the proof of the necessity can
be obtained in the reverse order. Thus, (1.60) can be obtained from (1.59).
Then by (1.57) and Lemma 1.6, we finally conclude that the function $%
\sum\limits_{i=1}^{r}w_{i}(\mathbf{x})g_{i}^{0}\left( \mathbf{a}^{i}\cdot
\mathbf{x}\right) $ is extremal to $f\left( \mathbf{x}\right) $.
\end{proof}

In the following, $\left\vert Q\right\vert $ will denote the Lebesgue
measure of a measurable set $Q.$ The following corollary is obvious.

\bigskip

\textbf{Corollary 1.6.} \textit{Let $Y$ be an }$r$\textit{-set. A
function $\sum\limits_{i=1}^{r}g_{i}^{0}\left( \mathbf{a}^{i}\cdot \mathbf{x}%
\right) $ in $\mathcal{R}\left( \mathbf{a}^{1},...,\mathbf{a}^{r}\right) $
is the best approximation to $f(\mathbf{x)}$ if and only if}
\begin{equation*}
g_{j}^{0}\left( y_{j}\right) =\frac{1}{\left\vert Y^{\left( j\right)
}\right\vert }\int\limits_{Y^{\left( j\right) }}\left( f^{\ast }\left(
\mathbf{y}\right) -\sum\limits_{\substack{ i=1  \\ i\neq j}}%
^{r}g_{i}^{0}\left( y_{i}\right) \right) d\mathbf{y}^{\left( j\right)
},\;\;j=1,...,r.
\end{equation*}%

\bigskip

In \cite{52}, this corollary was proven for the case $r=n.$

\bigskip

\subsection{Formulas for the best approximation and approximation error}

In this section, we establish explicit formulas for both the
best approximation and approximation error, provided that the weight
functions are constants. In this case, since we vary over $g_{i},$ the set $%
\mathcal{R}\left( \mathbf{a}^{1},...,\mathbf{a}^{r};~w_{1},...,w_{r}\right) $
coincides with $\mathcal{R}\left( \mathbf{a}^{1},...,\mathbf{a}^{r}\right) .$
Thus, without loss of generality, we may assume that $w_{i}(\mathbf{x})=1$
for $i=1,...,r.$

For brevity of the further exposition, introduce the notation
\begin{equation*}
A=\int\limits_{Y}f^{\ast }\left( \mathbf{y}\right) d\mathbf{y}\text{ and \ }%
f_{i}^{\ast }=f_{i}^{\ast }(y_{i})=\int\limits_{Y^{\left( i\right) }}f^{\ast
}\left( \mathbf{y}\right) d\mathbf{y}^{\left( i\right) },~i=1,...,r.
\end{equation*}
The following theorem is a generalization of the main result of \cite%
{52} from the case $r=n$ to the cases $r<n.$

\bigskip

\textbf{Theorem 1.13.} \textit{Let $Y$ be an }$r$\textit{-set. Set the
functions}
\begin{equation*}
g_{1}^{0}\left( y_{1}\right) =\frac{1}{\left\vert Y^{\left( 1\right)
}\right\vert }f_{1}^{\ast }-\left( r-1\right) \frac{A}{\left\vert
Y\right\vert }
\end{equation*}%
\textit{and}
\begin{equation*}
g_{j}^{0}\left( y_{j}\right) =\frac{1}{\left\vert Y^{\left( j\right)
}\right\vert }f_{j}^{\ast },\;j=2,...,r.
\end{equation*}%
\textit{Then the function $\sum\limits_{i=1}^{r}g_{i}^{0}\left( \mathbf{a}%
^{i}\cdot \mathbf{x}\right) $ is the best approximation from $\mathcal{R}%
\left( \mathbf{a}^{1},...,\mathbf{a}^{r}\right) $ to $f\left( \mathbf{x}%
\right) $.}

\bigskip

The proof is simple. It is sufficient to verify that
the functions $g_{j}^{0}\left( y_{j}\right) ,\;j=1,...,r$, satisfy the
conditions of Corollary 1.6. This becomes obvious if note that
\begin{equation*}
\sum\limits_{\underset{i\neq j}{i=1}}^{r}\frac{1}{\left\vert Y^{\left(
j\right) }\right\vert }\frac{1}{\left\vert Y^{\left( i\right) }\right\vert }%
\int\limits_{Y^{\left( j\right) }}\left[ \int\limits_{Y^{\left( i\right)
}}f^{\ast }\left( \mathbf{y}\right) d\mathbf{y}^{\left( i\right) }\right] d%
\mathbf{y}^{\left( j\right) }=\left( r-1\right) \frac{1}{\left\vert
Y\right\vert }\int\limits_{Y}f^{\ast }\left( \mathbf{y}\right) d\mathbf{y}
\end{equation*}%
for $j=1,...,r$.

\bigskip

\textbf{Theorem 1.14.} \textit{Let $Y$ be an $r$-set. Then the
error of approximation to a function $f(x)$ from the set $\mathcal{R}\left(
\mathbf{a}^{1},...,\mathbf{a}^{r}\right) $ can be calculated by the formula}
\begin{equation*}
E(f)=\left\vert \det J\right\vert ^{-1/2}\left( \left\Vert f^{\ast
}\right\Vert _{L_{2}(Y)}^{2}-\sum_{i=1}^{r}\frac{1}{\left\vert Y^{\left(
i\right) }\right\vert ^{2}}\left\Vert f_{i}^{\ast }\right\Vert
_{L_{2}(Y)}^{2}+(r-1)\frac{A^{2}}{\left\vert Y\right\vert }\right) ^{1/2}.
\end{equation*}

\bigskip

\begin{proof} From Eq. (1.58), Lemma 1.7 and Theorem 1.13, it follows that

\begin{equation*}
E(f)=\left\vert \det J\right\vert ^{-1/2}\left( \left\Vert f^{\ast
}\right\Vert _{L_{2}(Y)}^{2}-I\right) ^{1/2},\eqno(1.61)
\end{equation*}%
where

\begin{equation*}
I=\left\Vert \sum_{i=1}^{r}\frac{1}{\left\vert Y^{\left( i\right)
}\right\vert }f_{i}^{\ast }-(r-1)\frac{A}{\left\vert Y\right\vert }%
\right\Vert _{L_{2}(Y)}^{2}.
\end{equation*}%
The integral $I$ can be written as a sum of the following four integrals:

\begin{eqnarray*}
I_{1} &=&\sum_{i=1}^{r}\frac{1}{\left\vert Y^{\left( i\right) }\right\vert
^{2}}\left\Vert f_{i}^{\ast }\right\Vert
_{L_{2}(Y)}^{2},~I_{2}=\sum_{i=1}^{r}\sum\limits_{\substack{ j=1  \\ j\neq i
}}^{r}\frac{1}{\left\vert Y^{\left( i\right) }\right\vert }\frac{1}{%
\left\vert Y^{\left( j\right) }\right\vert }\int\limits_{Y}f_{i}^{\ast
}f_{j}^{\ast }d\mathbf{y,} \\
I_{3} &=&-2(r-1)\frac{1}{\left\vert Y\right\vert }A\sum_{i=1}^{r}\frac{1}{%
\left\vert Y^{\left( i\right) }\right\vert }\int\limits_{Y}f_{i}^{\ast }d%
\mathbf{y,}~I_{4}=(r-1)^{2}\frac{A^{2}}{\left\vert Y\right\vert }.
\end{eqnarray*}%
\qquad \qquad

It is not difficult to verify that

\begin{equation*}
\int\limits_{Y}f_{i}^{\ast }f_{j}^{\ast }d\mathbf{y=}\left\vert Y_{0}\times
\prod\limits_{\substack{ k=1  \\ k\neq i,j}}^{r}Y_{k}\right\vert A^{2},\text{
for }i,j=1,...,r,~i\neq j,\eqno(1.62)
\end{equation*}%
and

\begin{equation*}
\int\limits_{Y}f_{i}^{\ast }d\mathbf{y}=\left\vert Y_{0}\times \prod\limits
_{\substack{ k=1  \\ k\neq i}}^{r}Y_{k}\right\vert A,\text{ for }i=1,...,r.%
\eqno(1.63)
\end{equation*}

Considering (1.62) and (1.63) in the expressions of $I_{2}$ and $I_{3}$
respectively, we obtain that

\begin{equation*}
I_{2}=r(r-1)\frac{A^{2}}{\left\vert Y\right\vert }\text{ and }I_{3}=-2r(r-1)%
\frac{A^{2}}{\left\vert Y\right\vert }.
\end{equation*}%
Therefore,

\begin{equation*}
I=I_{1}+I_{2}+I_{3}+I_{4}=\sum_{i=1}^{r}\frac{1}{\left\vert Y^{\left(
i\right) }\right\vert ^{2}}\left\Vert f_{i}^{\ast }\right\Vert
_{L_{2}(Y)}^{2}-(r-1)\frac{A^{2}}{\left\vert Y\right\vert }.
\end{equation*}%
Now the last equality together with (1.61) complete the proof. \end{proof}

\textbf{Example.} Consider the following set
\begin{equation*}
X=\{\mathbf{x}\in \mathbb{R}^{4}:y_{i}=y_{i}(\mathbf{x})\in \lbrack
0;1],~i=1,...,4\},
\end{equation*}%
where
\begin{equation*}
\left\{
\begin{array}{c}
y_{1}=x_{1}+x_{2}+x_{3}-x_{4} \\
y_{2}=x_{1}+x_{2}-x_{3}+x_{4} \\
y_{3}=x_{1}-x_{2}+x_{3}+x_{4} \\
y_{4}=-x_{1}+x_{2}+x_{3}+x_{4}%
\end{array}%
\right. \eqno(1.64)
\end{equation*}%

Let the function
\begin{equation*}
f=8x_{1}x_{2}x_{3}x_{4}-\sum_{i=1}^{4}x_{i}^{4}+2\sum_{i=1}^{3}%
\sum_{j=i+1}^{4}x_{i}^{2}x_{j}^{2}
\end{equation*}%
be given on $X.$ Consider the approximation of this function by functions from $%
\mathcal{R}\left( \mathbf{a}^{1},\mathbf{a}^{2},\mathbf{a}^{3}\right) ,%
\mathcal{\ }$where $\mathbf{a}^{1}=(1;1;1;-1),~\mathbf{a}^{2}=(1;1;-1;1),~%
\mathbf{a}^{3}=(1;-1;1;1).$ Putting $\mathbf{a}^{4}=(-1;1;1;1),$ we complete
the system of vectors $\mathbf{a}^{1},\mathbf{a}^{2},\mathbf{a}^{3}$ to the
basis $\{\mathbf{a}^{1},\mathbf{a}^{2},\mathbf{a}^{3},\mathbf{a}^{4}\}$ in $%
\mathbb{R}^{4}.$ The linear transformation $J$ defined by (1.64) maps the
set $X$ onto the set $Y=[0;1]^{4}.$ The inverse transformation is given by
the formulas
\begin{equation*}
\left\{
\begin{array}{c}
x_{1}=\frac{1}{4}y_{1}+\frac{1}{4}y_{2}+\frac{1}{4}y_{3}-\frac{1}{4}y_{4} \\
x_{2}=\frac{1}{4}y_{1}+\frac{1}{4}y_{2}-\frac{1}{4}y_{3}+\frac{1}{4}y_{4} \\
x_{3}=\frac{1}{4}y_{1}-\frac{1}{4}y_{2}+\frac{1}{4}y_{3}+\frac{1}{4}y_{4} \\
x_{4}=-\frac{1}{4}y_{1}+\frac{1}{4}y_{2}+\frac{1}{4}y_{3}+\frac{1}{4}y_{4}%
\end{array}%
\right.
\end{equation*}%

It can be easily verified that $f^{\ast }=y_{1}y_{2}y_{3}y_{4}$ and $Y$
is a $3$-set with $Y_{i}=[0;1],$ $i=1,2,3.$ Besides, $Y_{0}=[0;1].$ After
easy calculations we obtain that $A=\allowbreak \frac{1}{16};~$\ $%
f_{i}^{\ast }=\allowbreak \frac{1}{8}y_{i}$ for $i=1,2,3;$ $\det J=-16;$ $%
\left\Vert f^{\ast }\right\Vert _{L_{2}(Y)}^{2}=\frac{1}{81};$ $\left\Vert
f_{i}^{\ast }\right\Vert _{L_{2}(Y)}^{2}=\frac{1}{192},$ $i=1,2,3.$ Now from
Theorems 1.13 and 1.14 it follows that the function $\frac{1}{8}%
\sum_{i=1}^{3}\left( \mathbf{a}^{i}\cdot \mathbf{x}\right) -\allowbreak
\frac{1}{8}$ is the best approximation from $\mathcal{R}\left( \mathbf{a}^{1},%
\mathbf{a}^{2},\mathbf{a}^{3}\right) $ to $f$ and $E(f)=\frac{1}{576}\sqrt{2}%
\sqrt{47}.$

\bigskip

\textbf{Remark 1.4.} Most of the material in this chapter is to be found in
\cite{52,53,54,49,50,47,66,64}.

\newpage

\chapter{The smoothness problem in ridge function representation}

This chapter discusses the following open problem raised in Buhmann and
Pinkus \cite{12}, and Pinkus \cite[p. 14]{117}. Assume we are given a
function $f(\mathbf{x})=f(x_{1},...,x_{n})$ of the form
\begin{equation*}
f(\mathbf{x})=\sum_{i=1}^{k}f_{i}(\mathbf{a}^{i}\cdot \mathbf{x}),\eqno(2.1)
\end{equation*}%
where the $\mathbf{a}^{i},$ $i=1,...,k,$ are pairwise linearly independent
vectors (directions) in $\mathbb{R}^{n}$, $f_{i}$ are arbitrarily behaved
univariate functions and $\mathbf{a}^{i}\cdot \mathbf{x}$ are standard inner
products. Assume, in addition, that $f$ is of a certain smoothness class,
that is, $f\in C^{s}(\mathbb{R}^{n})$, where $s\geq 0$ (with the convention
that $C^{0}(\mathbb{R}^{n})=C(\mathbb{R}^{n})$). Is it true that there will
always exist $g_{i}\in C^{s}(\mathbb{R})$ such that
\begin{equation*}
f(\mathbf{x})=\sum_{i=1}^{k}g_{i}(\mathbf{a}^{i}\cdot \mathbf{x})\text{ ?}%
\eqno(2.2)
\end{equation*}

In this chapter, we solve this problem up to some multivariate polynomial.
In the special case $n=2$, we see that this multivariate polynomial can be
written as a sum of polynomial ridge functions with the given directions $%
\mathbf{a}^{i}$. In addition, we find various conditions on the directions $%
\mathbf{a}^{i}$ guaranteeing a positive solution to the problem. We also
consider the question on constructing $g_{i}$ using the information about
the known functions $f_{i}$.

Most of the material of this chapter may be found in \cite{2,1,A2,A1,120}.

\bigskip

\section{A solution to the problem up to a multivariate polynomial}

In this section, we solve the above problem up to a multivariate polynomial.
That is, we show that if (2.1) holds for $f\in C^{s}(\mathbb{R}^{n})$ and
arbitrarily behaved $f_{i}$, then there exist $g_{i}\in C^{s}(\mathbb{R})$
such that
\begin{equation*}
f(\mathbf{x})=\sum_{i=1}^{k}g_{i}(\mathbf{a}^{i}\cdot \mathbf{x})+P(\mathbf{x%
}),
\end{equation*}%
where $P(\mathbf{x})$ is a polynomial of degree at most $k-1$. In the
special case $n=2$, we see that this multivariate polynomial can be written
as a sum of polynomial ridge functions with the given directions $\mathbf{a}%
^{i}$ and thus (2.2) holds with $g_{i}\in C^{s}(\mathbb{R})$.

\subsection{A brief overview of some results}

We start this subsection with the simple observation that for $k=1$ and $k=2$
the smoothness problem is easily solved. Indeed for $k=1$ by choosing $%
\mathbf{c}\in \mathbb{R}^{n}$ satisfying $\mathbf{a}^{1}\cdot \mathbf{c}=1$,
we have that $f_{1}(t)=f(t\mathbf{c)}$ is in $C^{s}(\mathbb{R})$. The same
argument can be carried out for the case $k=2.$ In this case, since the
vectors $\mathbf{a}^{1}$ and $\mathbf{a}^{2}$ are linearly independent,
there exists a vector $\mathbf{c}\in \mathbb{R}^{n}$ satisfying $\mathbf{a}%
^{1}\cdot \mathbf{c}=1$ and $\mathbf{a}^{2}\cdot \mathbf{c}=0.$ Therefore,
we obtain that the function $f_{1}(t)=f(t\mathbf{c)}-f_{2}(0)$ is in the
class $C^{s}(\mathbb{R})$. Similarly, one can verify that $f_{2}\in C^{s}(%
\mathbb{R})$.

The above cases with one and two ridge functions in (2.1) show that the
functions $f_{i}$ inherit smoothness properties of the given $f$. The
picture is absolutely different if the number of directions $k\geq 3$. For $%
k=3$, there are ultimately smooth functions which decompose into sums of
very badly behaved ridge functions. This phenomena comes from the classical
Cauchy Functional Equation (CFE). This equation,%
\begin{equation*}
h(x+y)=h(x)+h(y),\text{ }h:\mathbb{R\rightarrow R},\eqno(2.3)
\end{equation*}%
looks very simple and has a class of simple solutions $h(x)=cx,$ $c\in
\mathbb{R}$. However, it easily follows from Hamel basis theory that CFE
also has a large class of wild solutions. These solutions are called
\textquotedblleft wild" because they are extremely pathological. They are,
for example, not continuous at a point, not monotone on an interval, not
bounded on any set of positive measure (see, e.g., \cite{A}). Let $h_{1}$ be
any wild solution of the equation (2.3). Then the zero function can be
represented as%
\begin{equation*}
0=h_{1}(x)+h_{1}(y)-h_{1}(x+y).\eqno(2.4)
\end{equation*}%
Note that the functions involved in (2.4) are bivariate ridge functions with
the directions $(1,0)$, $(0,1)$ and $(1,1)$, respectively. This example
shows that for $k\geq 3$ the functions $f_{i}$ in (2.1) may not inherit
smoothness properties of the function $f$, which in the case of (2.4) is the
identically zero function. Thus the above problem arises naturally.

However, it was shown by some authors that, additional conditions on $f_{i}$
or the directions $\mathbf{a}^{i}$ guarantee smoothness of the
representation (2.1). It was first proved by Buhmann and Pinkus \cite{12}
that if in (2.1) $f\in C^{s}(\mathbb{R}^{n})$, $s\geq k-1$ and $f_{i}\in
L_{loc}^{1}(\mathbb{R)}$ for each $i$, then $f_{i}\in C^{s}(\mathbb{R)}$ for
$i=1,...,k.$ Later Pinkus \cite{120} found a strong relationship between CFE
and the problem of smoothness in ridge function representation. He
generalized extensively the previous result of Buhmann and Pinkus \cite{12}.
He showed that the solution is quite simple and natural if the functions $%
f_{i}$ are taken from a certain class $\mathcal{B}$ of real-valued functions
defined on $\mathbb{R}$. $\mathcal{B}$ includes, for example, the set of
continuous functions, the set of bounded functions, the set of Lebesgue
measurable functions (for the precise definition of $\mathcal{B}$ see the
next subsection). The result of Pinkus \cite{120} states that if in (1.1) $%
f\in C^{s}(\mathbb{R}^{n})$ and each $f_{i}\in \mathcal{B}$, then
necessarily $f_{i}\in C^{s}(\mathbb{R)}$ for $i=1,...,k$.

Note that severe restrictions on the directions $\mathbf{a}^{i}$ also
guarantee smoothness of the representation (2.1). For example, in (2.1) the
inclusions $f_{i}\in C^{s}(\mathbb{R})$, $i=1,...,k,$ are automatically
valid if the directions $\mathbf{a}^{i}$ are linearly independent and if
these directions are not linearly independent, then there exists $f\in C^{s}(%
\mathbb{R}^{n})$ of the form (2.1) such that the $f_{i}\notin C^{s}(\mathbb{R%
}),$ $i=1,...,k$ (see \cite{86}). Indeed, if the directions $\mathbf{a}^{i}$
are linearly independent, then for each $i=1,...,k,$ we can choose a vector $%
\mathbf{b}^{i}$ such that $\mathbf{b}^{i}\cdot \mathbf{a}^{i}=1,$ but at the
same time $\mathbf{b}^{i}\cdot \mathbf{a}^{j}=0,$ for all $j=1,...,k,$ $%
j\neq i$. Putting $\mathbf{x}=\mathbf{b}^{i}t$ in (2.1) yields that

\begin{equation*}
f(\mathbf{b}^{i}t)=f_{i}(t)+\sum_{j=1,j\neq i}^{k}f_{j}(0),\text{ }i=1,...,k.
\end{equation*}%
This shows that all the functions $f_{i}$ and $f$ belong to the same
smoothness class. If the directions $\mathbf{a}^{i}$ are not linearly
independent, then there exist numbers $\lambda _{1},...,\lambda _{k}$ such
that $\sum_{i=1}^{k}\left\vert \lambda _{i}\right\vert >0$ and $%
\sum_{i=1}^{k}\lambda _{i}\mathbf{a}^{i}=\mathbf{0}$. Let $h$ be any wild
solution of CFE. Then it is not difficult to verify that

\begin{equation*}
0=\sum_{i=1}^{k}h_{i}(\mathbf{a}^{i}\cdot \mathbf{x}),
\end{equation*}%
where $h_{i}(t)=h(\lambda _{i}t),$ $i=1,...,k.$ Note that in the last
representation, the zero function is an ultimately smooth function, while
all the functions $h_{i}$ are highly nonsmooth.

The above result of Pinkus was a starting point for further research on
continuous and smooth sums of ridge functions. Much work in this direction
was done by Konyagin and Kuleshov \cite{86,K2}, and Kuleshov \cite{K4}. They
mainly analyze the continuity of $f_{i}$, that is, the question of if and
when continuity of $f$ guarantees the continuity of $f_{i}$. There are also
other results concerning different properties, rather than continuity, of $%
f_{i}$. Most results in \cite{86,K2,K4} involve certain subsets (convex open
sets, convex bodies, etc.) of $\mathbb{R}^{n}$ instead of only $\mathbb{R}%
^{n}$ itself.

In \cite{1}, Aliev and Ismailov gave a partial solution to the smoothness
problem. Their solution comprises the cases in which $s\geq 2$ and $k-1$
directions of the given $k$ directions are linearly independent.

Kuleshov \cite{88} generalized Aliev and Ismailov's result \cite[Theorem 2.3]%
{1} to all possible cases of $s$. That is, he proved that if a function $%
f\in C^{s}(\mathbb{R}^{n})$, where $s\geq 0$, is of the form (2.1) and $%
(k-1) $-tuple of the given set of $k$ directions $\mathbf{a}^{i}$ forms a
linearly independent system, then there exist $g_{i}\in C^{s}(\mathbb{R})$, $%
i=1,...,k $, such that (2.2) holds (see \cite[Theorem 3]{88}). In Section
2.2 we give a new constructive proof of Kuleshov's result.

\bigskip

\subsection{A result of A. Pinkus}

In \cite{120}, A. Pinkus considered the smoothness problem in ridge function
representation. For a given function $f$ $:\mathbb{R}^{n}\rightarrow \mathbb{%
R}$, he posed and partially answered the following question. If $f$ belongs
to some smoothness class and (2.1) holds, what can we say about the
smoothness of the functions $f_{i}$? He proved that for a large class of
representing functions $f_{i}$, these $f_{i}$ are smooth. That is, if
apriori we assume that in the representation (2.1) the functions $f_{i}$ is
of a certain class of \textquotedblleft reasonably well behaved functions",
then they have the same degree of smoothness as the function $f.$ As the
mentioned class of \textquotedblleft reasonably well behaved functions" one
may take, e.g., the set of functions that are continuous at a point, the set
of Lebesgue measurable functions, etc. All these classes arise from the
class $\mathcal{B}$ considered by Pinkus \cite{120} and the classical theory
of CFE. In \cite{120}, $\mathcal{B}$ denotes any linear space of real-valued
functions $u$ defined on $\mathbb{R}$, closed under translation, such that
if there is a function $v\in C(\mathbb{R)}$ for which $u-v$ satisfies CFE,
then $u-v$ is necessarily linear, i.e. $u(x)-v(x)=cx,$ for some constant $%
c\in \mathbb{R}$. Such a definition of $\mathcal{B}$ is required in the
proof of the following theorem.

\bigskip

\textbf{Theorem 2.1} (Pinkus \cite{120}). \textit{Assume $f\in C^{s}(\mathbb{%
R}^{n})$ is of the form (2.1). Assume, in addition, that each $f_{i}\in
\mathcal{B}$. Then necessarily $f_{i}\in C^{s}(\mathbb{R)}$ for $i=1,...,k.$}

\bigskip

\begin{proof} We prove this theorem by induction on $k.$ The result is
valid when $k=1$. Indeed, taking any direction $\mathbf{c}$ such that $%
\mathbf{a}^{1}\cdot \mathbf{c}=1$ and putting $x=\mathbf{c}t$ in (2.1), we
obtain that $f_{1}(t)=f(\mathbf{c}t)\in C^{s}(\mathbb{R})$. Assume that the
result is valid for $k-1.$ Let us show that it is valid for $k$.

Chose any vector $\mathbf{e}\in \mathbb{R}^{n}$ satisfying $\mathbf{e\cdot a}%
^{k}=0$ and $\mathbf{e\cdot a}^{i}=b_{i}\neq 0$, for $i=1,...,k-1.$ Clearly,
there exists a vector with this property. The property of $\mathbf{e}$
enables us to write that

\begin{equation*}
f(\mathbf{x+e}t)-f(\mathbf{x)=}\sum_{i=1}^{k-1}f_{i}(\mathbf{a}^{i}\cdot
\mathbf{x}+b_{i}t)-f_{i}(\mathbf{a}^{i}\cdot \mathbf{x}).
\end{equation*}%
Thus%
\begin{equation*}
F(\mathbf{x}):=f(\mathbf{x+e}t)-f(\mathbf{x})=\sum_{i=1}^{k-1}h_{i}(\mathbf{a%
}^{i}\cdot \mathbf{x}),
\end{equation*}%
where%
\begin{equation*}
h_{i}(y)=f_{i}(y+b_{i}t)-f_{i}(y)\text{, }i=1,...,k-1.
\end{equation*}%
Since $f_{i}\in \mathcal{B}$ and $\mathcal{B}$ is translation invariant, $%
h_{i}\in \mathcal{B}$. In addition, since $F\in C^{s}(\mathbb{R}^{n})$, it
follows by our induction assumption that $h_{i}\in C^{s}(\mathbb{R})$. Note
that this inclusion is valid for all $t\in \mathbb{R}$.

In \cite{B1}, de Bruijn proved that if for any $c\in \mathbb{R}$ the
difference $u(y+c)-u(y)$ ($u$ is any real function on $\mathbb{R}$) belongs
to the class $C^{s}(\mathbb{R})$, then $u$ is necessarily of the form $u=v+r$%
, where $v\in $ $C^{s}(\mathbb{R})$ and $r$ satisfies CFE. Thus each
function $f_{i}$ is of the form $f_{i}=v_{i}+r_{i}$, where $v_{i}\in $ $%
C^{s}(\mathbb{R})$ and $r_{i}$ satisfies CFE. By our assumption, each $f_{i}$
is in $\mathcal{B}$, and from the definition of $\mathcal{B}$ it follows
that $r_{i}=f_{i}-v_{i}$ is a linear function. Thus $f_{i}=v_{i}+r_{i}$,
where both $v_{i},r_{i}\in $ $C^{s}(\mathbb{R})$, implying that $f_{i}\in
C^{s}(\mathbb{R})$. This is valid for $i=1,...,k-1$, and hence also for $i=k$.
\end{proof}

\textbf{Remark 2.1. }In de Bruijn \cite{B1,B2}, there are delineated various
classes of real-valued functions $\mathcal{D}$ with the property that if $%
\bigtriangleup _{t}f=f(\cdot +t)-f(\cdot )\in \mathcal{D}$ for all $t\in
\mathbb{R}$, then $f-s\in \mathcal{D}$, for some $s$ satisfying CFE (for
such classes see the next subsection). Some translation invariant classes
among them are $C^{\infty }(\mathbb{R})$ functions; analytic functions;
algebraic polynomials; trigonometric polynomials. Theorem 2.1 can be
suitably restated for any of these classes.

\bigskip

\subsection{Polynomial functions of $k$-th order}

Given $h_{1},...,h_{k}\in \mathbb{R}$, we define inductively the difference
operator $\Delta _{h_{1}...h_{k}}$ as follows
\begin{eqnarray*}
\Delta _{h_{1}}f(x) &:&=f(x+h_{1})-f(x), \\
\Delta _{h_{1}...h_{k}}f &:&=\Delta _{h_{k}}(\Delta _{h_{1}...h_{k-1}}f),%
\text{ }f:\mathbb{R\rightarrow R}.
\end{eqnarray*}%
If $h_{1}=\cdots=h_{k}=h,$ then we write briefly $\Delta _{h}^{k}f$ instead
of $\Delta _{\underset{n\text{ times}}{\underbrace{h...h}}}f$. For various
properties of difference operators see \cite[Section 15.1]{Kuc}.

\bigskip

\textbf{Definition 2.1 }(see \cite{Kuc}). \textit{A function $f:\mathbb{R\rightarrow
R}$ is called a polynomial function of order $k$ ($k\in \mathbb{N}$%
) if for every $x\in \mathbb{R}$ and $h\in \mathbb{R}$ we have}
\begin{equation*}
\Delta _{h}^{k+1}f(x)=0.
\end{equation*}

It can be shown that if $\Delta _{h}^{k+1}f=0$ for any $h\in \mathbb{R}$,
then $\Delta _{h_{1}...h_{k+1}}f=0$ for any $h_{1},...,h_{k+1}\in \mathbb{R}$
(see \cite[Theorem 15.3.3]{Kuc}). A polynomial of degree at most $k$ is a
polynomial function of order $k$ (see \cite[Theorem 15.9.4]{Kuc}). The
polynomial functions generalize ordinary polynomials, and reduce to the
latter under mild regularity assumptions. For example, if a polynomial
function is continuous at one point, or bounded on a set of positive
measure, then it continuous at all points (see \cite{Cies, Kurepa}), and
therefore is a polynomial of degree $k$ (see \cite[Theorem 15.9.4]{Kuc}).

Basic results concerning polynomial functions are due to S. Mazur-W. Orlicz
\cite{Maz}, McKiernan \cite{Mc}, Djokovi\'{c} \cite{Djok}. The following
theorem, which we will use in the sequel, yield implicitly the general
construction of polynomial functions.

\bigskip

\textbf{Theorem 2.2 }(see \cite[Theorems 15.9.1 and 15.9.2]{Kuc}). \textit{A
function $f:\mathbb{R\rightarrow R}$ is a polynomial function of order $k$
if and only if it admits a representation}
\begin{equation*}
f=f_{0}+f_{1}+...+f_{k},
\end{equation*}%
\textit{where $f_{0}$ is a constant and $f_{j}:\mathbb{R\rightarrow R}$, $%
j=1,...,k$, are diagonalizations of $j$-additive symmetric functions $F_{j}:%
\mathbb{R}^{j}\mathbb{\rightarrow R}$, i.e.,}

\begin{equation*}
f_{j}(x)=F_{j}(x,...,x).
\end{equation*}

\bigskip

Note that a function $F_{p}:\mathbb{R}^{p}\mathbb{\rightarrow R}$ is called $%
p$-additive if for every $j,$ $1\leq j\leq p,$ and for every $%
x_{1},...,x_{p},y_{j}\in \mathbb{R}$
\begin{equation*}
F(x_{1},...,x_{j}+y_{j},...,x_{p})=F(x_{1},...,x_{p})+F(x_{1},...,x_{j-1},y_{j},x_{j+1},...,x_{p}),
\end{equation*}%
i.e., $F$ is additive in each of its variables $x_{j}$ (see \cite[p.363]%
{Kuc}). A simple example of a $p$-additive function is given by the product

\begin{equation*}
f_{1}(x_{1})\times\cdots\times f_{p}(x_{p}),
\end{equation*}%
where the univariate functions $f_{j},$ $j=1,...,p$, are additive.

Following de Bruijn, we say that a class $\mathcal{D}$ of real functions has
the \textit{difference property} if any function $f:\mathbb{R\rightarrow R}$
such that $\bigtriangleup _{h}f\in \mathcal{D}$ for all $h\in \mathbb{R}$,
admits a decomposition $f=g+S$, where $g\in \mathcal{D}$ and $S$ satisfies
the Cauchy Functional Equation (2.3). Several classes with the difference
property are investigated in de Bruijn \cite{B1,B2}. Some of these classes
are:

\smallskip

1) $C(\mathbb{R)}$, continuous functions;

2) $C^{s}(\mathbb{R)}$, functions with continuous derivatives up to order $s$%
;

3) $C^{\infty }(\mathbb{R)}$, infinitely differentiable functions;

4) analytic functions;

5) functions which are absolutely continuous on any finite interval;

6) functions having bounded variation over any finite interval;

7) algebraic polynomials;

8) trigonometric polynomials;

9) Riemann integrable functions.

\smallskip

A natural generalization of classes with the difference property are classes
of functions with the difference property of $k$-th order.

\bigskip

\textbf{Definition 2.2 }(see \cite{Gajda}). \textit{A class $\mathcal{F}$ is said to
have the difference property of $k$-th order if any
function $f:\mathbb{R\rightarrow R}$ such that $\bigtriangleup _{h}^{k}f\in
\mathcal{F}$ for all $h\in \mathbb{R}$, admits a decomposition $f=g+H$,
where $g\in \mathcal{F}$ and $H$ is a polynomial function of $k$-th order.}

\bigskip

It is not difficult to see that the class $\mathcal{F}$ has the difference
property of first order if and only if it has the difference property in de
Bruijn's sense. There arises a natural question: which of the above classes
have difference properties of higher orders? Gajda \cite{Gajda} considered
this question in its general form, for functions defined on a locally
compact Abelian group and showed that for any $k\in \mathbb{N}$, continuous
functions have the difference property of $k$-th order (see \cite[Theorem 4]%
{Gajda}). The proof of this result is based on several lemmas, in
particular, on the following lemma, which we will also use in the sequel.

\bigskip

\textbf{Lemma 2.1.} (see \cite[Lemma 5]{Gajda}). \textit{For each $k\in
\mathbb{N}$ the class of all continuous functions defined on $\mathbb{R}$
has the difference property of $k$-th order.}

\bigskip

In fact, Gajda \cite{Gajda} proved this lemma for Banach space valued
functions, but the simplest case with the space $\mathbb{R}$ has all the
difficulties. Unfortunately, the proof of the lemma has an essential gap.
The author of \cite{Gajda} tried to reduce the proof to $\mod1$ periodic
functions, but made a mistake in proving the continuity of the difference $%
\Delta _{h_{1}...h_{k-1}}(f-f^{\ast })$. Here $f^{\ast }:\mathbb{%
R\rightarrow R}$ is a $\mod1$ periodic function defined on the interval $%
[0,1)$ as $f^{\ast }(x)=f(x)$ and extended to the whole $\mathbb{R}$ with
the period $1$. That is, $f^{\ast }(x)=f(x)$ for $x\in \lbrack 0,1)$ and $%
f^{\ast }(x+1)=f^{\ast }(x)$ for $x\in \mathbb{R}$. In the proof, the author
of \cite{Gajda} takes a point $x\in \lbrack m,m+1)$ and writes that

\begin{equation*}
\Delta _{h_{1}...h_{k-1}}(f-f^{\ast })(x)=\Delta
_{h_{1}...h_{k-1}}(f(x)-f(x-m)),
\end{equation*}%
which is not valid. Even though $f^{\ast }(x)=f(x-m)$ for any $x\in \lbrack
m,m+1)$, the differences $\Delta _{h_{1}...h_{k-1}}f^{\ast }(x)$ and $\Delta
_{h_{1}...h_{k-1}}f(x-m)$ are completely different, since the latter may
involve values of $f$ at points outside $[0,1)$, which have no relationship
with the definition of $f^{\ast }$.

In the next section, we give a new proof for Lemma 2.1 (see Theorem 2.3
below). We hope that our proof is free from mathematical errors and thus the
above lemma itself is valid.

\bigskip

\subsection{Some auxiliary results on polynomial functions}

In this section, we do further research on polynomial functions and prove
some auxiliary results.

\bigskip

\textbf{Lemma 2.2.} \textit{If $f:\mathbb{R\rightarrow R}$ is a polynomial
function of order $k$, then for any $p\in $ $\mathbb{N}$ and any fixed $\xi
_{1},...,\xi _{p}\in \mathbb{R}$, the function}%
\begin{equation*}
g(x_{1},...,x_{p})=f(\xi _{1}x_{1}+\cdots +\xi _{p}x_{p}),
\end{equation*}%
\textit{considered on the $p$ dimensional space $\mathbb{Q}^{p}$ of rational
vectors, is an ordinary polynomial of degree at most $k$.}

\bigskip

\begin{proof} By Theorem 2.2,
\begin{equation*}
f=\sum_{m=0}^{k}f_{m},\eqno(2.5)
\end{equation*}%
where $f_{0}$ is a constant and $f_{m}:\mathbb{R\rightarrow R}$, $1,...,m$,
are diagonalizations of $m$-additive symmetric functions $F_{m}:\mathbb{R}%
^{m}\mathbb{\rightarrow R}$, i.e.,

\begin{equation*}
f_{m}(x)=F_{m}(x,...,x).
\end{equation*}%
For a $m$-additive function $F_{m}$ the equality

\begin{equation*}
F_{m}(\xi _{1},...,\xi _{i-1},r\xi _{i},\xi _{i+1},...,\xi _{m})=rF_{m}(\xi
_{1},...,\xi _{m})
\end{equation*}%
holds for all $i=1,...,m$ and any $r\in \mathbb{Q}$, $\xi _{i}\in $ $\mathbb{%
R}$, $i=1,...,m$ (see \cite[Theorem 13.4.1]{Kuc}). Using this, it is not
difficult to verify that for any $(x_{1},...,x_{p})\in \mathbb{Q}^{p}$,
\begin{eqnarray*}
f_{m}(\xi _{1}x_{1}+\cdots +\xi _{p}x_{p}) &=&F_{m}(\xi _{1}x_{1}+\cdots
+\xi _{p}x_{p},...,\xi _{1}x_{1}+\cdots +\xi _{p}x_{p}) \\
&=&\sum_{\substack{ 0\leq s_{i}\leq m,~\overline{i=1,p}  \\ s_{1}+\cdots
+s_{p}=m}}A_{s_{1}...s_{p}}F_{m}(\underset{s_{1}}{\underbrace{\xi
_{1},...,\xi _{1}}},...,\underset{s_{p}}{\underbrace{\xi _{p},...,\xi _{p}}}%
)x_{1}^{s_{1}}...x_{p}^{s_{p}}.
\end{eqnarray*}%
Here $A_{s_{1}...s_{p}}$ are some coefficients, namely $%
A_{s_{1}...s_{p}}=m!/(s_{1}!...s_{p}!).$ Considering the last formula in
(2.5), we conclude that the function $g(x_{1},...,x_{p})$, restricted to $%
\mathbb{Q}^{p}$, is a polynomial of degree at most $k$.
\end{proof}

\textbf{Lemma 2.3.} \textit{Assume $f$ is a polynomial function of order $k$%
. Then there exists a polynomial function $H$ of order $k+1$ such that $%
H(0)=0$ and}

\begin{equation*}
f(x)=H(x+1)-H(x).\eqno(2.6)
\end{equation*}

\bigskip

\begin{proof} Consider the function

\begin{equation*}
H(x):=xf(x)+\sum_{i=1}^{k}(-1)^{i}\frac{x(x+1)...(x+i)}{(i+1)!}\Delta
_{1}^{i}f(x).\eqno(2.7)
\end{equation*}%
Clearly, $H(0)=0.$ We are going to prove that $H$ is a polynomial function
of order $k+1$ and satisfies (2.6).

Let us first show that for any polynomial function $g$ of order $m$ the
function $G_{1}(x)=xg(x)$ is a polynomial function of order $m+1.$ Indeed,
for any $h_{1},...,h_{m+2}\in \mathbb{R}$ we can write that
\begin{equation*}
\Delta _{h_{1}...h_{m+2}}G_{1}(x)=(x+h_{1}+\cdots +h_{m+2})\Delta
_{h_{1}...h_{m+2}}g(x)
\end{equation*}%
\begin{equation*}
+\sum_{i=1}^{m+2}h_{i}\Delta
_{h_{1}...h_{i-1}h_{i+1...}h_{m+2}}g(x).\eqno(2.8)
\end{equation*}%
The last formula is verified directly by using the known product property of
differences, that is, the equality

\begin{equation*}
\Delta _{h}(g_{1}g_{2})=g_{1}\Delta _{h}g_{2}+g_{2}\Delta _{h}g_{1}+\Delta
_{h}g_{1}\Delta _{h}g_{2}.\eqno(2.9)
\end{equation*}%
Now since $g$ is a polynomial function of order $m$, all summands in (2.8)
is equal to zero; hence we obtain that $G_{1}(x)$ is a polynomial function
of order $m+1$. By induction, we can prove that the function $%
G_{p}(x)=x^{p}g(x)$ is a polynomial function of order $m+p.$ Since $\Delta
_{1}^{i}f(x)$ in (2.7) is a polynomial function of order $k-i$, it follows
that all summands in (2.7) are polynomial functions of order $k+1$.
Therefore, $H(x)$ is a polynomial function of order $k+1$.

Now let us prove (2.6). Considering the property (2.9) in (2.7) we can write
that%
\begin{equation*}
\Delta _{1}H(x)=\left[ f(x)+(x+1)\Delta _{1}f(x)\right]
\end{equation*}%
\begin{equation*}
+\sum_{i=1}^{k}(-1)^{i}\left[ \frac{(x+1)...(x+i+1)}{(i+1)!}\Delta
_{1}^{i+1}f(x)+\Delta _{1}\left( \frac{x(x+1)...(x+i)}{(i+1)!}\right) \Delta
_{1}^{i}f(x)\right] .\eqno(2.10)
\end{equation*}

Note that in (2.10)

\begin{equation*}
\Delta _{1}\left( \frac{x(x+1)...(x+i)}{(i+1)!}\right) =\frac{(x+1)...(x+i)}{%
i!}.
\end{equation*}%
Considering this and the assumption $\Delta _{1}^{k+1}f(x)=0$, it follows
from (2.10) that
\begin{equation*}
\Delta _{1}H(x)=f(x),
\end{equation*}%
that is, (2.6) holds.
\end{proof}

\bigskip

The next lemma is due to Gajda \cite{Gajda}.

\bigskip

\textbf{Lemma 2.4 }(see \cite[Corollary 1]{Gajda}). \textit{Let $f:$ $%
\mathbb{R\rightarrow R}$ be a $\mod1$ periodic function such that, for any $%
h_{1},...,h_{k}\in \mathbb{R}$, $\Delta _{h_{1}...h_{k}}f$ is continuous.
Then there exist a continuous function $g:$ $\mathbb{R\rightarrow R}$ and a
polynomial function $H$ of $k$-th order such that $f=g+H$.}

\bigskip

The following theorem generalizes de Bruijn's theorem (see \cite[Theorem 1.1]%
{B1}) on the difference property of continuous functions and shows that
Gajda's above lemma (see Lemma 2.1) is valid. Note that the main result of
\cite{Gajda} also uses this theorem.

\bigskip

\textbf{Theorem 2.3.} \textit{Assume for any $h_{1},...,h_{k}\in \mathbb{R}$%
, the difference $\Delta _{h_{1}...h_{k}}f(x)$ is a continuous function of
the variable $x$. Then there exist a function $g\in C(\mathbb{R})$ and a
polynomial function $H$ of $k$-th order with the property $H(0)=0$ such that}
\begin{equation*}
f=g+H.
\end{equation*}

\bigskip

\begin{proof} We prove this theorem by induction. For $k=1$, the theorem
is the result of de Bruijn: if $f$ is such that, for each $h$, $\Delta
_{h}f(x)$ is a continuous function of $x$, then it can be written in the
form $g+H$, where $g$ is continuous and $H$ is additive (that is, satisfies
the Cauchy Functional Equation). Assume that the theorem is valid for $k-1.$
Let us prove it for $k$. Without loss of generality we may assume that $%
f(0)=f(1)$. Otherwise, we can prove the theorem for $f_{0}(x)=f(x)-\left[
f(1)-f(0)\right] x$ and then automatically obtain its validity for $f$.

Consider the function

\begin{equation*}
F_{1}(x)=f(x+1)-f(x)\text{, }x\in \mathbb{R}.\eqno(2.11)
\end{equation*}%
Since for\ any $h_{1},...,h_{k}\in \mathbb{R}$, $\Delta _{h_{1}...h_{k}}f(x)$
is a continuous function of $x$ and $\Delta _{h_{1}...h_{k-1}}F_{1}=\Delta
_{h_{1}...h_{k-1}1}f$, the difference $\Delta _{h_{1}...h_{k-1}}F_{1}(x)$
will be a continuous function of $x$, as well. By assumption, there exist a
function $g_{1}\in C(\mathbb{R})$ and a polynomial function $H_{1}$ of $%
(k-1) $-th order with the property $H_{1}(0)=0$ such that

\begin{equation*}
F_{1}=g_{1}+H_{1}.\eqno(2.12)
\end{equation*}%
It follows from Lemma 2.3 that there exists a polynomial function $H_{2}$ of
order $k$ such that $H_{2}(0)=0$ and

\begin{equation*}
H_{1}(x)=H_{2}(x+1)-H_{2}(x).\eqno(2.13)
\end{equation*}%
Substituting (2.13) in (2.12) we obtain that

\begin{equation*}
F_{1}(x)=g_{1}(x)+H_{2}(x+1)-H_{2}(x).\eqno(2.14)
\end{equation*}%
It follows from (2.11) and (2.14) that

\begin{equation*}
g_{1}(x)=\left[ f(x+1)-H_{2}(x+1)\right] -\left[ f(x)-H_{2}(x)\right] .\eqno%
(2.15)
\end{equation*}

Consider the function

\begin{equation*}
F_{2}=f-H_{2}.\eqno(2.16)
\end{equation*}%
Since $H_{2}$ is a polynomial function of order $k$ and for any $%
h_{1},...,h_{k}\in \mathbb{R}$ the difference $\Delta _{h_{1}...h_{k}}f(x)$
is a continuous function of $x$, we obtain that $\Delta
_{h_{1}...h_{k}}F_{2}(x)$ is also a continuous function of $x$. In addition,
since $f(0)=f(1)$ and $H_{2}(0)=H_{2}(1)=0$, it follows from (2.16) that $%
F_{2}(0)=F_{2}(1)$. We will use these properties of $F_{2}$ below.

Let us write (2.15) in the form

\begin{equation*}
g_{1}(x)=F_{2}(x+1)-F_{2}(x),\eqno(2.17)
\end{equation*}%
and define the following $\mod1$ periodic function
\begin{eqnarray*}
F^{\ast }(x) &=&F_{2}(x)\text{ for }x\in \lbrack 0,1), \\
F^{\ast }(x+1) &=&F^{\ast }(x)\text{ for }x\in \mathbb{R}.
\end{eqnarray*}

Consider the function

\begin{equation*}
F=F_{2}-F^{\ast }.\eqno(2.18)
\end{equation*}%
Let us show that $F\in C(\mathbb{R)}$. Indeed since $F(x)=0$ for $x\in
\lbrack 0,1)$, $F$ is continuous on $(0,1)$. Consider now the interval $%
[1,2) $. For any $x\in \lbrack 1,2)$ by the definition of $F^{\ast }$ and
(2.17) we can write that

\begin{equation*}
F(x)=F_{2}(x)-F_{2}(x-1)=g_{1}(x-1).\eqno(2.19)
\end{equation*}%
Since $g_{1}\in C(\mathbb{R)}$, it follows from (2.19) that $F$ is
continuous on $(1,2)$. Note that by (2.17) $g_{1}(0)=0$; hence $%
F(1)=g_{1}(0)=0$. Since $F\equiv 0$ on $[0,1)$, $F(1)=0$ and $F\in C(1,2),$
we obtain that $F$ is continuous on $(0,2)$. Consider the interval $[2,3)$.
For any $x\in \lbrack 2,3)$ we can write that

\begin{equation*}
F(x)=F_{2}(x)-F_{2}(x-2)=g_{1}(x-1)+g_{1}(x-2).\eqno(2.20)
\end{equation*}%
Since $g_{1}\in C(\mathbb{R)}$, $F$ is continuous on $(2,3)$. Note that by
(2.19) $\lim_{x\rightarrow 2-}F(x)=g_{1}(1)$ and by (2.20) $F(2)=g_{1}(1).$
We obtain from these arguments that $F$ is continuous on $(0,3)$. In the
same way, we can prove that $F$ is continuous on $(0,m)$ for any $m\in
\mathbb{N}$.

Similar arguments can be used to prove the continuity of $F$ on $(-m,0)$ for
any $m\in \mathbb{N}$. We show it for the first interval $[-1,0)$. For any $%
x\in \lbrack -1,0)$ by the definition of $F^{\ast }$ and (2.17) we can write
that

\begin{equation*}
F(x)=F_{2}(x)-F_{2}(x+1)=-g_{1}(x).
\end{equation*}%
Since $g_{1}\in C(\mathbb{R)}$, it follows that $F$ is continuous on $(-1,0)$%
. Besides, \newline $\lim_{x\rightarrow 0-}F(x)=-g_{1}(0)=0.$ This shows that $F$ is
continuous on $(-1,1)$, since $F\equiv 0$ on $[0,1).$ Combining all the
above arguments we conclude that $F\in C(\mathbb{R)}$.

Since $F\in C(\mathbb{R)}$ and $\Delta _{h_{1}...h_{k}}F_{2}(x)$ is a
continuous function of $x$, we obtain from (2.18) that $\Delta
_{h_{1}...h_{k}}F^{\ast }(x)$ is also a continuous function of $x.$ By Lemma
2.4, there exist a function $g_{2}\in C(\mathbb{R)}$ and a polynomial
function $H_{3}$ of order $k$ such that

\begin{equation*}
F^{\ast }=g_{2}+H_{3}.\eqno(2.21)
\end{equation*}%
It follows from (2.16), (2.18) and (2.21) that

\begin{equation*}
f=F+g_{2}+H_{2}+H_{3}.\eqno(2.22)
\end{equation*}

Introduce the notation
\begin{eqnarray*}
H(x) &=&H_{2}(x)+H_{3}(x)-H_{3}(0), \\
g(x) &=&F(x)+g_{2}(x)+H_{3}(0).
\end{eqnarray*}%
Obviously, $g\in C(\mathbb{R)}$ and $H(0)=0$. It follows from (2.22) and the
above notation that
\begin{equation*}
f=g+H.
\end{equation*}%
This completes the proof of the theorem.
\end{proof}

\bigskip

\subsection{Main results}

We start this subsection with the following lemma.

\bigskip

\textbf{Lemma 2.5.} \textit{Assume we are given pairwise linearly
independent vectors $\mathbf{a}^{i},$ $i=1,...,k,$ and a function $f\in C(%
\mathbb{R}^{n})$ of the form (2.1) with arbitrarily behaved univariate
functions $f_{i}$. Then for any $h_{1},...,h_{k-1}\in \mathbb{R}$, and all
indices $i=1,...,k$, $\Delta _{h_{1}...h_{k-1}}f_{i}\in C(\mathbb{R})$.}

\bigskip

\begin{proof} We prove this lemma for the function $f_{k}.$ It can be
proven for the other functions $f_{i}$ in the same way. Let $%
h_{1},...,h_{k-1}\in \mathbb{R}$ be given. Since the vectors $\mathbf{a}^{i}$
are pairwise linearly independent, for each $j=1,...,k-1,$ there is a vector
$\mathbf{b}^{j}$ such that $\mathbf{b}^{j}\cdot \mathbf{a}^{j}=0$ and $%
\mathbf{b}^{j}\cdot \mathbf{a}^{k}\neq 0$. It is not difficult to see that
for any $\lambda \in \mathbb{R}$, $\Delta _{\lambda \mathbf{b}^{j}}f_{j}(%
\mathbf{a}^{j}\cdot \mathbf{x})=0.$ Therefore, for any $\lambda
_{1},...,\lambda _{k-1}\in \mathbb{R}$, we obtain from (2.1) that

\begin{equation*}
\Delta _{\lambda _{1}\mathbf{b}^{1}...\lambda _{k-1}\mathbf{b}^{k-1}}f(%
\mathbf{x})=\Delta _{\lambda _{1}\mathbf{b}^{1}...\lambda _{k-1}\mathbf{b}%
^{k-1}}f_{k}(\mathbf{a}^{k}\cdot \mathbf{x}).\eqno(2.23)
\end{equation*}%
Note that in multivariate setting the difference operator $\Delta _{\mathbf{h%
}^{1}...\mathbf{h}^{k}}f(\mathbf{x})$ is defined similarly as in the
previous section. If in (2.23) we take

\begin{eqnarray*}
\mathbf{x} &\mathbf{=}&\frac{\mathbf{a}^{k}}{\left\Vert \mathbf{a}%
^{k}\right\Vert ^{2}}t\text{, }t\in \mathbb{R}, \\
\lambda _{j} &=&\frac{h_{j}}{\mathbf{a}^{k}\cdot \mathbf{b}^{j}}\text{, }%
j=1,...,k-1,
\end{eqnarray*}%
we will obtain that $\Delta _{h_{1}...h_{k-1}}f_{k}\in C(\mathbb{R})$.
\end{proof}

The following theorem is valid.

\bigskip

\textbf{Theorem 2.4.} \textit{Assume a function $f\in C(\mathbb{R}^{n})$ is
of the form (2.1). Then there exist continuous functions $g_{i}:\mathbb{%
R\rightarrow R}$, $i=1,...,k$, and a polynomial $P(\mathbf{x})$ of degree at
most $k-1$ such that}

\begin{equation*}
f(\mathbf{x})=\sum_{i=1}^{k}g_{i}(\mathbf{a}^{i}\cdot \mathbf{x})+P(\mathbf{x%
}).\eqno(2.24)
\end{equation*}

\bigskip

\begin{proof} By Lemma 2.5 and Theorem 2.3, for each $i=1,...,k$, there
exists a function $g_{i}\in C(\mathbb{R})$ and a polynomial function $H_{i}$
of $(k-1)$-th order with the property $H_{i}(0)=0$ such that

\begin{equation*}
f_{i}=g_{i}+H_{i}.\eqno(2.25)
\end{equation*}

Consider the function

\begin{equation*}
F(\mathbf{x})=f(\mathbf{x})-\sum_{i=1}^{k}g_{i}(\mathbf{a}^{i}\cdot \mathbf{x%
}).\eqno(2.26)
\end{equation*}%
It follows from (2.1), (2.25) and (2.26) that

\begin{equation*}
F(\mathbf{x})=\sum_{i=1}^{k}H_{i}(\mathbf{a}^{i}\cdot \mathbf{x}).\eqno(2.27)
\end{equation*}

Denote the restrictions of the multivariate functions $H_{i}(\mathbf{a}%
^{i}\cdot \mathbf{x})$ to the space $\mathbb{Q}^{n}$ by $P_{i}(\mathbf{x})$,
respectively. By Lemma 2.2, the functions $P_{i}(\mathbf{x})$ are ordinary
polynomials of degree at most $k-1$. Since the space $\mathbb{Q}^{n}$ is
dense in $\mathbb{R}^{n}$, and the functions $F(\mathbf{x})$, $P_{i}(\mathbf{%
x})$, $i=1,...,k$, are continuous on $\mathbb{R}^{n}$, and the equality

\begin{equation*}
F(\mathbf{x})=\sum_{i=1}^{k}P_{i}(\mathbf{x}),\eqno(2.28)
\end{equation*}%
holds for all $\mathbf{x}\in \mathbb{Q}^{n}$, we obtain that (2.28) holds
also for all $\mathbf{x}\in \mathbb{R}^{n}$. Now (2.24) follows from (2.26)
and (2.28) by putting $P=\sum_{i=1}^{k}P_{i}$.
\end{proof}

Now we generalize Theorem 2.4 from $C(\mathbb{R}^{n})$ to any space $C^{s}(%
\mathbb{R}^{n})$ of $s$-th order continuously differentiable functions.

\bigskip

\textbf{Theorem 2.5.} \textit{Assume $f\in C^{s}(\mathbb{R}^{n})$ is of the
form (2.1). Then there exist functions $g_{i}\in C^{s}(\mathbb{R})$, $%
i=1,...,k$, and a polynomial $P(\mathbf{x})$ of degree at most $k-1$ such
that (2.24) holds.}

\bigskip

The proof is based on Theorems 2.1 and 2.4. On the one hand, it follows from
Theorem 2.4 that the $s$-th order continuously differentiable function $f-P$
can be expressed as $\sum_{i=1}^{k}g_{i}$ with continuous $g_{i}$. On the
other hand, since the class $\mathcal{B}$ in Theorem 2.1, in particular, can
be taken as $C(\mathbb{R}),$ it follows that $g_{i}\in C^{s}(\mathbb{R})$.

\bigskip

Note that Theorem 2.5 solves the problem posed in Buhmann and Pinkus \cite%
{12} and Pinkus \cite[p.14]{117} up to a polynomial. The following theorem
shows that in the two dimensional setting $n=2$ it solves the problem
completely.

\bigskip

\textbf{Theorem 2.6.} \textit{Assume a function $f\in C^{s}(\mathbb{R}^{2})$
is of the form}
\begin{equation*}
f(x,y)=\sum_{i=1}^{k}f_{i}(a_{i}x+b_{i}y),
\end{equation*}%
\textit{where $(a_{i},b_{i})$ are pairwise linearly independent vectors in $%
\mathbb{R}^{2}$ and $f_{i}$ are arbitrary univariate functions. Then there
exist functions $g_{i}\in C^{s}(\mathbb{R})$, $i=1,...,k$, such that}
\begin{equation*}
f(x,y)=\sum_{i=1}^{k}g_{i}(a_{i}x+b_{i}y).\eqno(2.29)
\end{equation*}

\bigskip

The proof of this theorem is not difficult. First we apply Theorem 2.5 and
obtain that
\begin{equation*}
f(x,y)=\sum_{i=1}^{k}\overline{g}_{i}(a_{i}x+b_{i}y)+P(x,y),\eqno(2.30)
\end{equation*}%
where $\overline{g}_{i}\in C^{s}(\mathbb{R})$ and $P(x,y)$ is a bivariate
polynomial of degree at most $k-1$. Then we use the known fact that a
bivariate polynomial $P(x,y)$ of degree $k-1$ is decomposed into a sum of
ridge polynomials with any given $k$ pairwise linearly independent
directions $(a_{i},b_{i}),$ $i=1,...,k$ (see e.g. \cite{97}). That is,
\begin{equation*}
P(x,y)=\sum_{i=1}^{k}p_{i}(a_{i}x+b_{i}y),
\end{equation*}%
where $p_{i}$ are univariate polynomials of degree at most $k-1$.
Considering this in (2.30) gives the desired representation (2.29).

\bigskip

\textbf{Remark 2.2.} Theorem 2.5 can be restated also for the classes $%
C^{\infty }(\mathbb{R})$ of infinitely differentiable functions and $D(%
\mathbb{R})$ of analytic functions. That is, if under the conditions of
Theorem 2.5, we have $f\in C^{\infty }(\mathbb{R}^{n})$ (or $f\in D(\mathbb{R%
}^{n})$), then this function can be represented also in the form (2.24) with
$g_{i}\in C^{\infty }(\mathbb{R})$ (or $g_{i}\in D(\mathbb{R})$). This
follows, similarly to the case $C^{s}(\mathbb{R})$ above, from Theorem 2.4
and Remark 2.1. These arguments are also valid for Theorem 2.6.

\bigskip

\section{A solution to the smoothness problem under certain conditions}

Assume we are given a function $f\in C^{s}(\mathbb{R}^{n})$ of the form
(2.1). In this section, we discuss various conditions on the directions $%
\mathbf{a}^{i}$ guaranteeing the validity of (2.2) with $g_{i}\in C^{s}(%
\mathbb{R})$.

\subsection{Directions with only rational components}

The following theorems, in particular, show that if directions of ridge
functions have only rational coordinates then no polynomial term appears in
Theorems 2.4 and 2.5.

\bigskip

\textbf{Theorem 2.7.} \textit{Assume a function $f\in C(\mathbb{R}^{n})$ is
of the form (2.1) and there is a nonsingular linear transformation $T:$ $%
\mathbb{R}^{n}\rightarrow \mathbb{R}^{n}$ such that $T\mathbf{a}^{i}\in
\mathbb{Q}\mathit{^{n}},$ $i=1,...,k$. Then there exist continuous functions
$g_{i}:\mathbb{R\rightarrow R}$, $i=1,...,k$, such that (2.2) holds.}

\bigskip

\begin{proof} Applying the coordinate change $\mathbf{x\rightarrow y}$,
given by the formula $\mathbf{x}=T\mathbf{y}$, to both sides of (2.1) we
obtain that

\begin{equation*}
\tilde{f}(\mathbf{y})=\sum_{i=1}^{k}f_{i}(\mathbf{b}^{i}\cdot \mathbf{y}),
\end{equation*}%
where $\tilde{f}(\mathbf{y})=f(T\mathbf{y})$ and $\mathbf{b}^{i}=T\mathbf{a}%
^{i},$ $i=1,...,k.$ Let us repeat the proof of Theorem 2.4 for the function $%
\tilde{f}$. Since the vectors $\mathbf{b}^{i}$, $i=1,...,k,$ have rational
coordinates, it is not difficult to see that the restrictions of the
functions $H_{i}$ to $\mathbb{Q}$ are univariate polynomials. Indeed, for
each $\mathbf{b}^{i}$ we can choose a vector $\mathbf{c}^{i}$ with rational
coordinates such that $\mathbf{b}^{i}\cdot \mathbf{c}^{i}=1$. If in the
equality $H_{i}(\mathbf{b}^{i}\cdot \mathbf{x})=P_{i}(\mathbf{x}),$ $\mathbf{%
x}\in \mathbb{Q}^{n}$, we take $\mathbf{x=c}^{i}t$ with $t\in \mathbb{Q}$,
we obtain that $H_{i}(t)=P_{i}(\mathbf{c}^{i}t)$ for all $t\in \mathbb{Q}$.
Now since $P_{i}$ is a multivariate polynomial on $\mathbb{Q}^{n}$, $H_{i}$
is a univariate polynomial on $\mathbb{Q}$. Denote this univariate
polynomial by $L_{i}$. Thus the formula

\begin{equation*}
P_{i}(\mathbf{x})=L_{i}(\mathbf{b}^{i}\cdot \mathbf{x})\eqno(2.31)
\end{equation*}%
holds for each $i=1,...,k$, and all $\mathbf{x}\in \mathbb{Q}^{n}$. Since $%
\mathbb{Q}^{n}$ is dense in $\mathbb{R}^{n}$, we see that (2.31) holds, in
fact, for all $\mathbf{x}\in \mathbb{R}^{n}$. Thus the polynomial $P(\mathbf{%
x})$ in (2.24) can be expressed as $\sum_{i=1}^{k}L_{i}(\mathbf{b}^{i}\cdot
\mathbf{x})$. Considering this in Theorem 2.4, we obtain that

\begin{equation*}
\tilde{f}(\mathbf{y})=\sum_{i=1}^{k}g_{i}(\mathbf{b}^{i}\cdot \mathbf{y}),%
\eqno(2.32)
\end{equation*}%
where $g_{i}$ are continuous functions. Using the inverse transformation $%
\mathbf{y}=T^{-1}\mathbf{x}$ in (2.32) we arrive at (2.2).
\end{proof}

\textbf{Theorem 2.8.} \textit{Assume a function $f\in C^{s}(\mathbb{R}^{n})$
is of the form (2.1) and there is a nonsingular linear transformation $T:$ $%
\mathbb{R}^{n}\rightarrow \mathbb{R}^{n}$ such that $T\mathbf{a}^{i}\in
\mathbb{Q}\mathit{^{n}},$ $i=1,...,k$. Then there exist functions $g_{i}\in
C^{s}(\mathbb{R})$, $i=1,...,k$, such that (2.2) holds.}

\bigskip

The proof of this theorem easily follows from Theorem 2.7 and Theorem 2.5.

\bigskip

\subsection{Linear independence of $k-1$ directions}

We already know that if the given directions $\mathbf{a}^{i}$ form a
linearly independent set, then the smoothness problem has a positive
solution (see Section 2.1.1). What can we say if all $\mathbf{a}^{i}$ are
not linearly independent? In the sequel, we show that if $k-1$ of the
directions $\mathbf{a}^{i}$, $i=1,...,k,$ are linearly independent, then in
(2.1) $f_{i}$ can be replaced with $g_{i}\in C^{s}(\mathbb{R})$. We will
also estimate the modulus of continuity of $g_{i}$ in terms of the modulus
of continuity of a function generated from $f$ under a linear transformation.

Let $F:\mathbb{R}^{n}\rightarrow \mathbb{R}$, $n\geq 1,$ be any function and
$\Omega \subset \mathbb{R}^{n}$. The function
\begin{equation*}
\omega (F;\delta ;\Omega )=\sup \left\{ \left\vert F(\mathbf{x})-F(\mathbf{y}%
)\right\vert :\mathbf{x},\mathbf{y}\in \Omega ,\text{ }\left\vert \mathbf{x}-%
\mathbf{y}\right\vert \leq \delta \right\} ,\text{ }0\leq \delta \leq
diam\Omega ,
\end{equation*}%
is called the modulus of continuity of the function $F(\mathbf{x}%
)=F(x_{1},...,x_{n})$ on the set $\Omega .$ We will also use the notation $%
\omega _{\mathbb{Q}}(F;\delta ;\Omega )$, which stands for the function $%
\omega (F;\delta ;\Omega \cap \mathbb{Q}^{n})$. Here $\mathbb{Q}$ denotes
the set of rational numbers. Note that $\omega _{\mathbb{Q}}(F;\delta
;\Omega )$ makes sense if the set $\Omega \cap \mathbb{Q}^{n}$ is not empty.
Clearly, $\omega _{\mathbb{Q}}(F;\delta ;\Omega )\leq \omega (F;\delta
;\Omega )$. The equality $\omega _{\mathbb{Q}}(F;\delta ;\Omega )=\omega
(F;\delta ;\Omega )$ holds for continuous $F$ and certain sets $\Omega $.
For example, it holds if for any $\mathbf{x},\mathbf{y}\in \Omega $ with $%
\left\vert \mathbf{x}-\mathbf{y}\right\vert \leq \delta $ there exist
sequences $\left\{ \mathbf{x}_{m}\right\} ,\left\{ \mathbf{y}_{m}\right\}
\subset \Omega \cap \mathbb{Q}^{n}$ such that $\mathbf{x}_{m}\rightarrow
\mathbf{x}$, $\mathbf{y}_{m}\rightarrow \mathbf{y}$ and $\left\vert \mathbf{x%
}_{m}-\mathbf{y}_{m}\right\vert \leq \delta ,$ for all $m$. There are many
sets $\Omega $, which satisfy this property.

\bigskip

The following lemma is valid.

\bigskip

\textbf{Lemma 2.6.} \textit{Assume a function $G\in C(\mathbb{R}^{n})$ has
the form}

\begin{equation*}
G(x_{1},...,x_{n})=\sum_{i=1}^{n}g(x_{i})-g(x_{1}+\cdot \cdot \cdot +x_{n}),%
\eqno(2.33)
\end{equation*}%
\textit{where $g$ is an arbitrarily behaved function. Then the following
inequality holds}

\begin{equation*}
\omega _{\mathbb{Q}}(g;\delta ;[-M,M])\leq 2\delta \left\vert
g(1)-g(0)\right\vert +3\omega \left( G;\delta ;[-M,M]^{n}\right) ,\eqno(2.34)
\end{equation*}%
\textit{where $\delta \in \left( 0,\frac{1}{2}\right) \cap \mathbb{Q}$ and $%
M\geq 1$.}

\bigskip

\begin{proof} Consider the function $f(t)=g(t)-g(0)$ and write (2.33) in
the form

\begin{equation*}
F(x_{1},...,x_{n})=\sum_{i=1}^{n}f(x_{i})-f(x_{1}+\cdot \cdot \cdot +x_{n}),%
\eqno(2.35)
\end{equation*}%
where

\begin{equation*}
F(x_{1},...,x_{n})=G(x_{1},...,x_{n})-(n-1)g(0).
\end{equation*}%
Note that the functions $f$ and $g$, as well as the functions $F$ and $G,$
have a common modulus of continuity. Thus we prove the lemma if we prove it
for the pair $\left\langle F,f\right\rangle .$

Since $f(0)=0,$ it follows from (2.35) that

\begin{equation*}
F(x_{1},0,...,0)=F(0,x_{2},0,...,0)=\cdot \cdot \cdot =F(0,0,...,x_{n})=0.%
\eqno(2.36)
\end{equation*}%
For the sake of brevity, introduce the notation $\mathcal{F}(x_{1},x_{2})=$ $%
F(x_{1},x_{2},0,...,0)$. Obviously, for any real number $x,$
\begin{eqnarray*}
\mathcal{F}(x,x) &=&2f(x)-f(2x); \\
\mathcal{F}(x,2x) &=&f(x)+f(2x)-f(3x); \\
&&\cdot \cdot \cdot  \\
\mathcal{F}(x,(k-1)x) &=&f(x)+f((k-1)x)-f(kx).
\end{eqnarray*}

We obtain from the above equalities that
\begin{eqnarray*}
f(2x) &=&2f(x)-\mathcal{F}(x,x), \\
f(3x) &=&3f(x)-\mathcal{F}(x,x)-\mathcal{F}(x,2x), \\
&&\cdot \cdot \cdot \\
f(kx) &=&kf(x)-\mathcal{F}(x,x)-\mathcal{F}(x,2x)-\cdot \cdot \cdot -%
\mathcal{F}(x,(k-1)x).
\end{eqnarray*}%
Thus for any nonnegative integer $k$,

\begin{equation*}
f(x)=\frac{1}{k}f(kx)+\frac{1}{k}\left[ \mathcal{F}(x,x)+\mathcal{F}%
(x,2x)+\cdot \cdot \cdot +\mathcal{F}(x,(k-1)x)\right] .\eqno(2.37)
\end{equation*}

Consider now the simple fraction $\frac{p}{m}\in (0,\frac{1}{2})$ and set $%
m_{0}=\left[ \frac{m}{p}\right] .$ Here $[r]$ denotes the whole number part
of $r$. Clearly, $m_{0}\geq 2$ and the remainder $p_{1}=m-m_{0}p<p.$ Taking $%
x=\frac{p}{m}$ and $k=m_{0}$ in (2.37) gives us the following equality

\begin{equation*}
f\left( \frac{p}{m}\right) =\frac{1}{m_{0}}f\left( 1-\frac{p_{1}}{m}\right)
\end{equation*}%
\begin{equation*}
+\frac{1}{m_{0}}\left[ \mathcal{F}\left( \frac{p}{m},\frac{p}{m}\right) +%
\mathcal{F}\left( \frac{p}{m},\frac{2p}{m}\right) +\cdot \cdot \cdot +%
\mathcal{F}\left( \frac{p}{m},(m_{0}-1)\frac{p}{m}\right) \right] .\eqno%
(2.38)
\end{equation*}%
On the other hand, since

\begin{equation*}
\mathcal{F}\left( \frac{p_{1}}{m},1-\frac{p_{1}}{m}\right) =f\left( \frac{%
p_{1}}{m}\right) +f\left( 1-\frac{p_{1}}{m}\right) -f(1),
\end{equation*}%
it follows from (2.38) that

\begin{equation*}
f\left( \frac{p}{m}\right) =\frac{f(1)}{m_{0}}
\end{equation*}%
\begin{equation*}
+\frac{1}{m_{0}}\left[
\mathcal{F}\left( \frac{p}{m},\frac{p}{m}\right) +\cdot \cdot \cdot +%
\mathcal{F}\left( \frac{p}{m},(m_{0}-1)\frac{p}{m}\right) +\mathcal{F}\left(
\frac{p_{1}}{m},1-\frac{p_{1}}{m}\right) \right]
\end{equation*}
\begin{equation*}
-\frac{1}{m_{0}}f\left( \frac{p_{1}}{m}\right) .\eqno(2.39)
\end{equation*}

Put $m_{1}=\left[ \frac{m}{p_{1}}\right] $, $p_{2}=m-m_{1}p_{1}.$ Clearly, $%
0\leq p_{2}<p_{1}$. Similar to (2.39), we can write that

\begin{equation*}
f\left( \frac{p_{1}}{m}\right) =\frac{f(1)}{m_{1}}
\end{equation*}%
\begin{equation*}
+\frac{1}{m_{1}}\left[
\mathcal{F}\left( \frac{p_{1}}{m},\frac{p_{1}}{m}\right) +\cdot \cdot \cdot +%
\mathcal{F}\left( \frac{p_{1}}{m},(m_{1}-1)\frac{p_{1}}{m}\right) +\mathcal{F%
}\left( \frac{p_{2}}{m},1-\frac{p_{2}}{m}\right) \right]
\end{equation*}
\begin{equation*}
-\frac{1}{m_{1}}f\left( \frac{p_{2}}{m}\right) .\eqno(2.40)
\end{equation*}

Let us make a convention that (2.39) is the $1$-st and (2.40) is the $2$-nd
formula. One can continue this process by defining the chain of pairs $%
(m_{2},p_{3}),$ $(m_{3},p_{4})$ until the pair $(m_{k-1},p_{k})$ with $%
p_{k}=0$ and writing out the corresponding formulas for each pair. For
example, the last $k$-th formula will be of the form
\begin{equation*}
f\left( \frac{p_{k-1}}{m}\right) =\frac{f(1)}{m_{k-1}}
\end{equation*}%
\begin{equation*}
+\frac{1}{m_{k-1}}\left[ \mathcal{F}\left( \frac{p_{k-1}}{m},\frac{p_{k-1}}{m%
}\right) +\cdot \cdot \cdot +\mathcal{F}\left( \frac{p_{k-1}}{m},(m_{k-1}-1)%
\frac{p_{k-1}}{m}\right) +\mathcal{F}\left( \frac{p_{k}}{m},1-\frac{p_{k}}{m}%
\right) \right]
\end{equation*}%
\begin{equation*}
-\frac{1}{m_{k-1}}f\left( \frac{p_{k}}{m}\right) .\eqno(2.41)
\end{equation*}%
Note that in (2.41), $f\left( \frac{p_{k}}{m}\right) =0$ and $\mathcal{F}%
\left( \frac{p_{k}}{m},1-\frac{p_{k}}{m}\right) =0$. Considering now the $k$%
-th formula in the $(k-1)$-th formula, then the obtained formula in the $%
(k-2)$-th formula, and so forth, we will finally arrive at the equality

\begin{equation*}
f\left( \frac{p}{m}\right) =f(1)\left[ \frac{1}{m_{0}}-\frac{1}{m_{0}m_{1}}%
+\cdot \cdot \cdot +\frac{(-1)^{k-1}}{m_{0}m_{1}\cdot \cdot \cdot m_{k-1}}%
\right]
\end{equation*}

\begin{equation*}
+\frac{1}{m_{0}}\left[ \mathcal{F}\left( \frac{p}{m},\frac{p}{m}\right)
+\cdot \cdot \cdot +\mathcal{F}\left( \frac{p}{m},(m_{0}-1)\frac{p}{m}%
\right) +\mathcal{F}\left( \frac{p_{1}}{m},1-\frac{p_{1}}{m}\right) \right]
\end{equation*}

\begin{equation*}
-\frac{1}{m_{0}m_{1}}\left[ \mathcal{F}\left( \frac{p_{1}}{m},\frac{p_{1}}{m}%
\right) +\cdot \cdot \cdot +\mathcal{F}\left( \frac{p_{1}}{m},(m_{1}-1)\frac{%
p_{1}}{m}\right) +\mathcal{F}\left( \frac{p_{2}}{m},1-\frac{p_{2}}{m}\right) %
\right]
\end{equation*}

\begin{equation*}
+\cdot \cdot \cdot +
\end{equation*}%
\begin{equation*}
\frac{(-1)^{k-1}}{m_{0}m_{1}\cdot \cdot \cdot m_{k-1}}\left[ \mathcal{F}%
\left( \frac{p_{k-1}}{m},\frac{p_{k-1}}{m}\right) +\cdot \cdot \cdot +%
\mathcal{F}\left( \frac{p_{k-1}}{m},(m_{k-1}-1)\frac{p_{k-1}}{m}\right) %
\right] .\eqno(2.42)
\end{equation*}%
Taking into account (2.36) and the definition of $\mathcal{F}$, for any
point of the form $\left( \frac{p_{i}}{m},c\right) $, $i=0,1,...,k-1,$ $%
p_{0}=p$, $c\in \lbrack 0,1]$, we can write that
\begin{equation*}
\left\vert \mathcal{F}\left( \frac{p_{i}}{m},c\right) \right\vert
=\left\vert \mathcal{F}\left( \frac{p_{i}}{m},c\right) -\mathcal{F}\left(
0,c\right) \right\vert \leq \omega \left( F;\frac{p_{i}}{m};[0,1]^{n}\right)
\leq \omega \left( F;\frac{p}{m};[0,1]^{n}\right) .
\end{equation*}%
Applying this inequality to each term $\mathcal{F}\left( \frac{p_{i}}{m}%
,\cdot \right) $ in (2.42), we obtain that

\begin{equation*}
\left\vert f\left( \frac{p}{m}\right) \right\vert \leq \left[ \frac{1}{m_{0}}%
-\frac{1}{m_{0}m_{1}}+\cdot \cdot \cdot +\frac{(-1)^{k-1}}{m_{0}m_{1}\cdot
\cdot \cdot m_{k-1}}\right] \left\vert f(1)\right\vert
\end{equation*}

\begin{equation*}
+\left[ 1+\frac{1}{m_{0}}+\cdot \cdot \cdot +\frac{1}{m_{0}\cdot \cdot \cdot
m_{k-2}}\right] \omega \left( F;\frac{p}{m};[0,1]^{n}\right) .\eqno(2.43)
\end{equation*}

Since $m_{0}\leq m_{1}\leq \cdot \cdot \cdot \leq m_{k-1},$ it is not
difficult to see that in (2.43)

\begin{equation*}
\frac{1}{m_{0}}-\frac{1}{m_{0}m_{1}}+\cdot \cdot \cdot +\frac{(-1)^{k-1}}{%
m_{0}m_{1}\cdot \cdot \cdot m_{k-1}}\leq \frac{1}{m_{0}}
\end{equation*}%
and

\begin{equation*}
1+\frac{1}{m_{0}}+\cdot \cdot \cdot +\frac{1}{m_{0}\cdot \cdot \cdot m_{k-2}}%
\leq \frac{m_{0}}{m_{0}-1}.
\end{equation*}%
Considering the above two inequalities in (2.43) we obtain that

\begin{equation*}
\left\vert f\left( \frac{p}{m}\right) \right\vert \leq \frac{\left\vert
f(1)\right\vert }{m_{0}}+\frac{m_{0}}{m_{0}-1}\omega \left( F;\frac{p}{m}%
;[0,1]^{n}\right) .\eqno(2.44)
\end{equation*}%
Since $m_{0}=\left[ \frac{m}{p}\right] \geq 2,$ it follows from (2.44) that

\begin{equation*}
\left\vert f\left( \frac{p}{m}\right) \right\vert \leq \frac{2p\left\vert
f(1)\right\vert }{m}+2\omega \left( F;\frac{p}{m};[0,1]^{n}\right) .\eqno%
(2.45)
\end{equation*}

Let now $\delta \in \left( 0,\frac{1}{2}\right) \cap \mathbb{Q}$ be a
rational increment, $M\geq 1$ and $x,x+\delta $ be two points in $\left[ -M,M%
\right] \cap \mathbb{Q}.$ By (2.45) we can write that

\begin{equation*}
\left\vert f(x+\delta )-f(x)\right\vert \leq \left\vert f(\delta
)\right\vert +\left\vert F(x,\delta ,0,...,0)\right\vert \leq 2\delta
\left\vert f(1)\right\vert +3\omega \left( F;\delta ;[-M,M]^{n}\right) .\eqno%
(2.46)
\end{equation*}%
Now (2.34) follows from (2.46) and the definitions of $f$ and $F$.
\end{proof}

\textbf{Remark 2.3.} The above lemma shows that the restriction of $g$ to
the set of rational numbers $\mathbb{Q}$ is uniformly continuous on any
interval $[-M,M]\cap \mathbb{Q}$.

\bigskip

To prove the main result of this section we need the following lemma.

\bigskip

\textbf{Lemma 2.7.} \textit{Assume a function $G\in C(\mathbb{R}^{n})$ has
the form}

\begin{equation*}
G(x_{1},...,x_{n})=\sum_{i=1}^{n}g(x_{i})-g(x_{1}+\cdot \cdot \cdot +x_{n}),
\end{equation*}%
\textit{where $g$ is an arbitrary function. Then there exists a function $%
F\in C(\mathbb{R})$ such that}

\begin{equation*}
G(x_{1},...,x_{n})=\sum_{i=1}^{n}F(x_{i})-F(x_{1}+\cdot \cdot \cdot +x_{n})%
\eqno(2.47)
\end{equation*}%
\textit{and the following inequality holds}

\begin{equation*}
\omega (F;\delta ;[-M,M])\leq 3\omega \left( G;\delta ;[-M,M]^{n}\right) ,%
\eqno(2.48)
\end{equation*}%
\textit{where $0\leq \delta \leq \frac{1}{2}$ and $M\geq 1$.}

\bigskip

\begin{proof} Consider the function

\begin{equation*}
u(t)=g(t)-\left[ g(1)-g(0)\right] t.
\end{equation*}%
Obviously, $u(1)=u(0)$ and

\begin{equation*}
G(x_{1},...,x_{n})=\sum_{i=1}^{n}u(x_{i})-u(x_{1}+\cdot \cdot \cdot +x_{n}).%
\eqno(2.49)
\end{equation*}%
By Lemma 2.6, the restriction of $u$ to $\mathbb{Q}$ is continuous and
uniformly continuous on every interval $[-M,M]\cap \mathbb{Q}$. Denote this
restriction by $v$.

Let $y$ be any real number and $\{y_{k}\}_{k=1}^{\infty }$ be any sequence
of rational numbers converging to $y$. We can choose $M>0$ so that $y_{k}\in
\lbrack -M,M]$ for any $k\in \mathbb{N}$. It follows from the uniform
continuity of $v$ on $[-M,M]\cap \mathbb{Q}$ that the sequence $%
\{v(y_{k})\}_{k=1}^{\infty }$ is Cauchy. Thus there exits a finite limit $%
\lim_{k\rightarrow \infty }v(y_{k})$. It is not difficult to see that this
limit does not depend on the choice of $\{y_{k}\}_{k=1}^{\infty }$.

Let $F$ denote the following extension of $v$ to the set of real numbers.
\begin{equation*}
F(y)=\left\{
\begin{array}{c}
v(y),\text{ if }y\in \mathbb{Q}\text{;} \\
\lim_{k\rightarrow \infty }v(y_{k}),\text{ if }y\in \mathbb{R}\backslash
\mathbb{Q}\text{ and }\{y_{k}\}\text{ is a sequence in }\mathbb{Q}\text{
tending to }y.%
\end{array}%
\right.
\end{equation*}%
In view of the above arguments, $F$ is well defined on the whole real line.
Let us prove that for this function (2.47) is valid.

Consider an arbitrary point $(x_{1},...,x_{n})\in \mathbb{R}^{n}$ and
sequences of rationale numbers $\{y_{k}^{i}\}_{k=1}^{\infty },i=1,...,n,$
tending to $x_{1},...,x_{n},$ respectively. Taking into account (2.49), we
can write that
\begin{equation*}
G(y_{k}^{1},...,y_{k}^{n})=\sum_{i=1}^{n}v(y_{k}^{i})-v(y_{k}^{1}+\cdot
\cdot \cdot +y_{k}^{n}),\text{ for all }k=1,2,...,\eqno(2.50)
\end{equation*}%
since $v$ is the restriction of $u$ to $\mathbb{Q}$. Tending $k\rightarrow
\infty $ in both sides of (2.50) we obtain (2.47).

Let us now prove that $F\in C(\mathbb{R})$ and (2.48) holds. Since $%
v(1)=v(0) $ we obtain from (2.49) and (2.34) that for $\delta \in \left( 0,%
\frac{1}{2}\right) \cap \mathbb{Q}$, $M\geq 1$ and any numbers $a,b\in
\lbrack -M,M]\cap \mathbb{Q}$, $\left\vert a-b\right\vert \leq \delta ,$ the
following inequality holds

\begin{equation*}
\left\vert v(a)-v(b)\right\vert \leq 3\omega \left( G;\delta
;[-M,M]^{n}\right) .\eqno(2.51)
\end{equation*}%
Consider any real numbers $r_{1}$ and $r_{2}$ satisfying $r_{1},r_{2}\in
\lbrack -M,M]$, $\left\vert r_{1}-r_{2}\right\vert \leq \delta $ and take
sequences $\{a_{k}\}_{k=1}^{\infty }\subset \lbrack -M,M]\cap \mathbb{Q}$, $%
\{b_{k}\}_{k=1}^{\infty }\subset \lbrack -M,M]\cap \mathbb{Q}$ with the
property $\left\vert a_{k}-b_{k}\right\vert \leq \delta ,$ $k=1,2,...,$ and
tending to $r_{1}$ and $r_{2}$, respectively. By (2.51),
\begin{equation*}
\left\vert v(a_{k})-v(b_{k})\right\vert \leq 3\omega \left( G;\delta
;[-M,M]^{n}\right) .
\end{equation*}%
If we take limits\ on both sides of the above inequality, we obtain that

\begin{equation*}
\left\vert F(r_{1})-F(r_{2})\right\vert \leq 3\omega \left( G;\delta
;[-M,M]^{n}\right) ,
\end{equation*}%
which means that $F$ is uniformly continuous on $[-M,M]$ and

\begin{equation*}
\omega \left( F;\delta ;[-M,M]\right) \leq 3\omega \left( G;\delta
;[-M,M]^{n}\right) .
\end{equation*}%
Note that in the last inequality $\delta $ is a rational number from the
interval $\left( 0,\frac{1}{2}\right) .$ It is well known that the modulus
of continuity $\omega (f;\delta ;\Omega )$ of a continuous function $f$ is
continuous from the right for any compact set $\Omega \subset \mathbb{R}^{n}$
and it is continuous from the left for certain compact sets $\Omega $, in
particular for rectangular sets (see \cite{Kol}). It follows immediately
that (2.48) is valid for all $\delta \in \lbrack 0,\frac{1}{2}].$
\end{proof}

The following theorem was first obtained by Kuleshov \cite{88}. Below, we
prove this using completely different ideas. Our proof, which is taken from
\cite{2}, contains a theoretical method for constructing the functions $%
g_{i}\in C^{s}(\mathbb{R})$ in (2.2). Using this method, we will also estimate
the modulus of continuity of $g_{i}$ in terms of the modulus of continuity
of $f$ (see Remark 2.4 below).

\bigskip

\textbf{Theorem 2.9.} \textit{Assume we are given $k$ directions $\mathbf{a}%
^{i}$, $i=1,...,k$, in $\mathbb{R}^{n}\backslash \{\mathbf{0}\}$ and $k-1$
of them are linearly independent.\ Assume that a function $f\in C(\mathbb{R}%
^{n})$ is of the form (2.1). Then $f$ can be represented also in the form
(2.2) with $g_{i}\in C(\mathbb{R})$, $i=1,...,k$.}

\bigskip

\begin{proof} Without loss of generality, we may assume that the first $%
k-1 $ vectors $\mathbf{a}^{1},\mathbf{...},\mathbf{a}^{k-1}$ are linearly
independent. Thus there exist numbers $\lambda _{1},...,\lambda _{k-1}\in
\mathbb{R}$ such that $\mathbf{a}^{k}=\lambda _{1}\mathbf{a}^{1}+\cdot \cdot
\cdot +\lambda _{k-1}\mathbf{a}^{k-1}$. We may also assume that the first $p$
numbers $\lambda _{1},...,\lambda _{p}$, $1\leq p\leq k-1$, are nonzero and
the remaining $\lambda _{j}$s are zero. Indeed, if necessary, we can
rearrange the vectors $\mathbf{a}^{1},\mathbf{...},\mathbf{a}^{k-1}$ so that
this assumption holds. Complete the system $\{\mathbf{a}^{1},...,\mathbf{a}%
^{k-1}\}$ to a basis $\{\mathbf{a}^{1},...,\mathbf{a}^{k-1},\mathbf{b}%
^{k},...,\mathbf{b}^{n}\}$ and consider the linear transformation $\mathbf{y}%
=A\mathbf{x,}$ where $\mathbf{x}=(x_{1},...,x_{n})^{T},$ $\mathbf{y}%
=(y_{1},...,y_{n})^{T}$ and $A$ is the matrix, rows of which are formed by
the coordinates of the vectors $\mathbf{a}^{1},...,\mathbf{a}^{k-1},\mathbf{b%
}^{k},...,\mathbf{b}^{n}.$ Using this transformation, we can write (2.1) in
the form
\begin{equation*}
f(A^{-1}\mathbf{y})=f_{1}(y_{1})+\cdot \cdot \cdot
+f_{k-1}(y_{k-1})+f_{k}(\lambda _{1}y_{1}+\cdot \cdot \cdot +\lambda
_{p}y_{p}).\eqno(2.52)
\end{equation*}

For the brevity of exposition in the sequel, we put $l=k-1$ and use the
notation
\begin{equation*}
w=f_{l+1},\text{ }\Phi (y_{1},...,y_{l})=f(A^{-1}\mathbf{y})\text{ and }%
Y_{j}=(y_{1},...,y_{j-1},y_{j+1},...,y_{l}), j=1,...,l.
\end{equation*}%
Using this notation, we can write (2.52) in the form
\begin{equation*}
\Phi (y_{1},...,y_{l})=f_{1}(y_{1})+\cdot \cdot \cdot
+f_{l}(y_{l})+w(\lambda _{1}y_{1}+\cdot \cdot \cdot +\lambda _{p}y_{p}).\eqno%
(2.53)
\end{equation*}%
In (2.53), taking sequentially $Y_{1}=0,$ $Y_{2}=0,$..., $Y_{l}=0$ we obtain
that
\begin{equation*}
f_{j}(y_{j})=\Phi (y_{1},...,y_{l})|_{Y_{j}=\mathbf{0}}-w(\lambda
_{j}y_{j})-\sum_{\substack{ i=1  \\ i\neq j}}^{l}f_{i}(0),\text{ }j=1,...,l.%
\eqno(2.54)
\end{equation*}%
Substituting (2.54) in (2.53), we obtain the equality

\begin{equation*}
\left.
\begin{array}{c}
w(\lambda _{1}y_{1}+\cdot \cdot \cdot +\lambda
_{p}y_{p})-\sum_{j=1}^{p}w(\lambda _{j}y_{j})-(l-p)w(0)= \\
=\Phi (y_{1},...,y_{l})-\sum_{j=1}^{l}\Phi (y_{1},...,y_{l})|_{Y_{j}=\mathbf{%
0}}+(l-1)\sum_{j=1}^{l}f_{j}(0).%
\end{array}%
\right. \eqno(2.55)
\end{equation*}%
We see that the right hand side of (2.55) depends only on the variables $%
y_{1},y_{2},...,y_{p}.$ Denote the right hand side of (2.55) by $%
H(y_{1},...,y_{p}).$ That is, set%
\begin{equation*}
H(y_{1},...,y_{p})\overset{def}{=}\Phi (y_{1},...,y_{l})-\sum_{j=1}^{l}\Phi
(y_{1},...,y_{l})|_{Y_{j}=\mathbf{0}}+(l-1)\sum_{j=1}^{l}f_{j}(0).\eqno(2.56)
\end{equation*}%
We will use the following identity, which follows from (2.55) and (2.56)%
\begin{equation*}
H(y_{1},...,y_{p})=w(\lambda _{1}y_{1}+\cdot \cdot \cdot +\lambda
_{p}y_{p})-\sum_{j=1}^{p}w(\lambda _{j}y_{j})-(l-p)w(0).\eqno(2.57)
\end{equation*}

It follows from (2.56) and the continuity of $f$ that the function $H$ is
continuous on $\mathbb{R}^{p}$. Then, defining the function

\begin{equation*}
G(y_{1},...,y_{p})=H(\frac{y_{1}}{\lambda _{1}},...,\frac{y_{p}}{\lambda _{p}%
})+(l-p)w(0)\eqno(2.58)
\end{equation*}%
and applying Lemma 2.7, we obtain that there exists a function $F\in C(%
\mathbb{R})$ such that

\begin{equation*}
G(y_{1},...,y_{p})=F(y_{1}+\cdot \cdot \cdot +y_{p})-\sum_{j=1}^{p}F(y_{p}).%
\eqno(2.59)
\end{equation*}%
It follows from the formulas (2.57)-(2.59) that

\begin{equation*}
w(\lambda _{1}y_{1}+\cdot \cdot \cdot +\lambda
_{p}y_{p})-\sum_{j=1}^{p}w(\lambda _{j}y_{j})=F(\lambda _{1}y_{1}+\cdot
\cdot \cdot +\lambda _{p}y_{p})-\sum_{j=1}^{p}F(\lambda _{j}y_{j}).\eqno%
(2.60)
\end{equation*}

Let us introduce the following functions
\begin{equation*}
\left.
\begin{array}{c}
g_{j}(y_{j})=\Phi (y_{1},...,y_{l})|_{Y_{j}=\mathbf{0}}-F(\lambda
_{j}y_{j})-\sum_{\substack{ i=1  \\ i\neq j}}^{l}f_{i}(0),\text{ }j=1,...,p,
\\
g_{j}(y_{j})=\Phi (y_{1},...,y_{l})|_{Y_{j}=\mathbf{0}}-\sum_{\substack{ i=1
\\ i\neq j}}^{l}f_{i}(0)-w(0),\text{ }j=p+1,...,l.%
\end{array}%
\right. \eqno(2.61)
\end{equation*}

Note that $g_{j}\in C(\mathbb{R}),$ $j=1,...,l.$ Considering (2.55), (2.60)
and (2.61) it is not difficult to verify that
\begin{equation*}
\Phi (y_{1},...,y_{l})=g_{1}(y_{1})+\cdot \cdot \cdot
+g_{l}(y_{l})+F(\lambda _{1}y_{1}+\cdot \cdot \cdot +\lambda _{p}y_{p}).\eqno%
(2.62)
\end{equation*}%
In (2.62), denoting $F=g_{k}$, recalling the definition of $\Phi $ and going
back to the variable $\mathbf{x}=(x_{1},...,x_{n})$ by using again the
linear transformation $\mathbf{y}=A\mathbf{x}$, we finally obtain (2.2).
\end{proof}

\textbf{Remark 2.4.} Using Theorem 2.9 and Lemma 2.7, one can estimate the
modulus of continuity of the functions $g_{i}$ in representation (2.2) in
terms of the modulus of continuity of $\Phi $. To show how one can do this,
assume $C>0,$ $M\geq 1,$ $0\leq \delta \leq \frac{1}{2\max \{\left\vert
\lambda _{j}\right\vert \}}$ and introduce the following sets
\begin{equation*}
\mathcal{M}_{j}=\left\{ (x_{1},...,x_{l})\in \mathbb{R}^{l}:x_{j}\in \lbrack
-M/\left\vert \lambda _{j}\right\vert ,M/\left\vert \lambda _{j}\right\vert ]%
\text{, }x_{i}=0\text{ for }i\neq j\right\},
\end{equation*}%
\begin{equation*}
j=1,...,p;
\end{equation*}%
\begin{equation*}
\mathcal{C}_{j}=\left\{ (x_{1},...,x_{l})\in \mathbb{R}^{l}:x_{j}\in \lbrack
-C,C]\text{, }x_{i}=0\text{ for }i\neq j\right\} ,\text{ }j=p+1,...,l.
\end{equation*}%
It can be easily obtained from (2.61) that
\begin{equation*}
\omega (g_{j};\delta ;[-M/\left\vert \lambda _{j}\right\vert ,M/\left\vert
\lambda _{j}\right\vert ])\leq \omega \left( \Phi ;\delta ;\mathcal{M}%
_{j}\right) +\omega (F;\delta _{1};[-M,M]),\text{ }j=1,...,p,\eqno(2.63)
\end{equation*}%
\begin{equation*}
\omega (g_{j};\delta ;[-C,C])\leq \omega \left( \Phi ;\delta ;\mathcal{C}%
_{j}\right) ,\text{ }j=p+1,...,l,\eqno(2.64)
\end{equation*}%
where $\delta _{1}=\delta \cdot \max \{\left\vert \lambda _{j}\right\vert
\}. $ To estimate $\omega (F;\delta _{1};[-M,M])$ in (2.63), we refer to
Lemma 2.7. Applying Lemma 2.7 to the function $G$ in (2.58) we obtain that
in addition to (2.59) the following inequality holds.
\begin{equation*}
\omega (F;\delta _{1};[-M,M])\leq 3\omega \left( G;\delta
_{1};[-M,M]^{p}\right) .\eqno(2.65)
\end{equation*}%
Note that here $0\leq \delta _{1}\leq \frac{1}{2}$ as in Lemma 2.7. It
follows from (2.58) and (2.65) that
\begin{equation*}
\omega (F;\delta _{1};[-M,M])\leq 3\omega \left( H;\delta _{2};\mathcal{K}%
\right) ,\eqno(2.66)
\end{equation*}%
where $\mathcal{K}=[-M/\left\vert \lambda _{1}\right\vert ,M/\left\vert
\lambda _{1}\right\vert ]\times \cdot \cdot \cdot \times \lbrack
-M/\left\vert \lambda _{p}\right\vert ,M/\left\vert \lambda _{p}\right\vert
] $ and \newline
$\delta _{2}=\delta _{1}/\min \{\left\vert \lambda _{j}\right\vert \} $.
Further, (2.66) and (2.56) together yield that
\begin{equation*}
\omega (F;\delta _{1};[-M,M])\leq (3l+3)\omega \left( \Phi ;\delta _{2};%
\mathcal{S}\right) ,\eqno(2.67)
\end{equation*}%
where $\mathcal{S}=\left\{ (x_{1},...,x_{l})\in \mathbb{R}%
^{l}:(x_{1},...,x_{p})\in \mathcal{K}\text{, }x_{i}=0\text{ for }i>p\right\}
$. Now it follows from (2.63) and (2.67) that
\begin{equation*}
\omega (g_{j};\delta ;[-M/\left\vert \lambda _{j}\right\vert ,M/\left\vert
\lambda _{j}\right\vert ])\leq \omega \left( \Phi ;\delta ;\mathcal{M}%
_{j}\right) +(3l+3)\omega \left( \Phi ;\delta _{2};\mathcal{S}\right) ,\text{
}j=1,...,p.\eqno(2.68)
\end{equation*}%
Formulas (2.64), (2.67) and (2.68) provide us with upper estimates for the
modulus of continuity of the functions $g_{j},$ $j=1,...,k,$ in terms of the
modulus of continuity of $\Phi $. Recall that in these estimates $l=k-1$, $%
F=g_{k}$ and $\lambda _{j}$ are coefficients\textit{\ }in the expression $%
\mathbf{a}^{k}=\lambda _{1}\mathbf{a}^{1}+\cdot \cdot \cdot +\lambda _{p}%
\mathbf{a}^{p}$.

\bigskip

Theorems 2.1 and 2.9 together give the following result.

\bigskip

\textbf{Theorem 2.10. }\textit{Assume we are given $k$ directions $\mathbf{a}%
^{i}$, $i=1,...,k$, in $\mathbb{R}^{n}\backslash \{\mathbf{0}\}$ and $k-1$
of them are linearly independent.\ Assume that a function $f\in C^{s}(%
\mathbb{R}^{n})$ is of the form (2.1). Then $f$ can be represented also in
the form (2.2), where the functions $g_{i}\in C^{s}(\mathbb{R})$, $i=1,...,k$%
.}

\bigskip

Indeed, on the one hand, it follows from Theorem 2.9 that $f$ can be
expressed as (2.2) with continuous $g_{i}$. On the other hand, since the
class $\mathcal{B}$ in Theorem 2.1, in particular, can be taken as $C(%
\mathbb{R}),$ it follows that $g_{i}\in C^{s}(\mathbb{R})$.

\bigskip

\textbf{Remark 2.5.} In addition to the above $C^{s}(\mathbb{R})$, Theorems
2.9 can be restated also for some other subclasses of the space of
continuous functions. These are $C^{\infty }(\mathbb{R})$ functions;
analytic functions; algebraic polynomials; trigonometric polynomials. More
precisely, assume $\mathcal{H}(\mathbb{R})$ is any of these subclasses and $%
\mathcal{H}(\mathbb{R}^{n})$ is the $n$-variable analog of the $\mathcal{H}(%
\mathbb{R})$. If under the conditions of Theorem 2.9, we have $f\in \mathcal{%
H}(\mathbb{R}^{n})$, then this function can be represented in the form (2.2)
with $g_{i}\in \mathcal{H}(\mathbb{R}).$ This follows, similarly to the case
$C^{s}(\mathbb{R})$ above, from Theorem 2.9 and Remark 2.1.

\bigskip

\section{A constructive analysis of the smoothness problem}

Note that Theorems 2.4-2.10 are generally existence results. They tell about
existence of smooth ridge functions $g_{i}$ in the corresponding
representation formula (2.2) or (2.24). They are uninformative if we want to
construct explicitly these functions.

In this section, we give two theorems which do not only address the
smoothness problem, but also are useful in constructing the mentioned $g_{i}$%
.

\subsection{Bivariate case}

We start with the constructive analysis of the smoothness problem for
bivariate functions. We show that if a bivariate function of a certain
smoothness class is represented by a sum of finitely many, arbitrarily
behaved ridge functions, then, under suitable conditions, it also can be
represented by a sum of ridge functions of the same smoothness class and
these ridge functions can be constructed explicitely.

\bigskip

\textbf{Theorem 2.11.} \textit{Assume $(a_{i},b_{i})$, $i=1,...,n$ \ are
pairwise linearly independent vectors in $\mathbb{R}^{2}$. Assume that a
function $f\in C^{s}(\mathbb{R}^{2})$ has the form}
\begin{equation*}
f(x,y)=\sum_{i=1}^{n}f_{i}(a_{i}x+b_{i}y),
\end{equation*}%
\textit{where $f_{i}$ are arbitrary univariate functions and $s\geq n-2.$
Then $f$ can be represented also in the form}
\begin{equation*}
f(x,y)=\sum_{i=1}^{n}g_{i}(a_{i}x+b_{i}y),\eqno(2.69)
\end{equation*}%
\textit{where the functions $g_{i}\in C^{s}(\mathbb{R})$, $i=1,...,n$. In
(2.69), the functions $g_{i}$, $i=1,...,n,$ can be constructed by the
formulas}
\begin{eqnarray*}
g_{p} &=&\varphi _{p,n-p-1},\text{ }p=1,...,n-2; \\
g_{n-1} &=&h_{1,n-1};\text{ }g_{n}=h_{2,n-1}.
\end{eqnarray*}%
\textit{Here all the involved functions $\varphi _{p,n-p-1}$\textit{, }$%
h_{1,n-1}$ and $h_{2,n-1}$ can be found inductively as follows}%
\begin{eqnarray*}
h_{1,1}(t) &=&\frac{\partial ^{n-2}}{\partial l_{1}\cdot \cdot \cdot
\partial l_{n-2}}f^{\ast }(t,0),~ \\
h_{2,1}(t) &=&\frac{\partial ^{n-2}}{\partial l_{1}\cdot \cdot \cdot
\partial l_{n-2}}f^{\ast }(0,t)-\frac{\partial ^{n-2}}{\partial l_{1}\cdot
\cdot \cdot \partial l_{n-2}}f^{\ast }(0,0); \\
h_{1,k+1}(t) &=&\frac{1}{e_{1}\cdot l_{k}}\int_{0}^{t}h_{1,k}(z)dz,\text{ }%
k=1,...,n-2; \\
h_{2,k+1}(t) &=&\frac{1}{e_{2}\cdot l_{k}}\int_{0}^{t}h_{2,k}(z)dz,\text{ }%
k=1,...,n-2;
\end{eqnarray*}%
\textit{and}
\begin{equation*}
\varphi _{p,1}(t)=\frac{\partial ^{n-p-2}f^{\ast }}{\partial l_{p+1}\cdot
\cdot \cdot \partial l_{n-2}}\left( \frac{\widetilde{a}_{p}t}{\widetilde{a}%
_{p}^{2}+\widetilde{b}_{p}^{2}},\frac{\widetilde{b}_{p}t}{\widetilde{a}%
_{p}^{2}+\widetilde{b}_{p}^{2}}\right) -h_{1,p+1}\left( \frac{\widetilde{a}%
_{p}t}{\widetilde{a}_{p}^{2}+\widetilde{b}_{p}^{2}}\right)
\end{equation*}%
\begin{equation*}
-h_{2,p+1}\left( \frac{\widetilde{b}_{p}t}{\widetilde{a}_{p}^{2}+\widetilde{b%
}_{p}^{2}}\right)-\sum_{j=1}^{p-1}\varphi _{j,p-j+1}\left( \frac{\widetilde{a%
}_{j}\widetilde{a}_{p}+\widetilde{b}_{j}\widetilde{b}_{p}}{\widetilde{a}%
_{p}^{2}+\widetilde{b}_{p}^{2}}t\right),
\end{equation*}
\begin{equation*}
p=1,...,n-2\left( \text{for }p=n-2\text{, }\frac{\partial ^{n-p-2}f^{\ast }}{%
\partial l_{p+1}\cdot \cdot \cdot \partial l_{n-2}}:=f^{\ast }\right) ;
\end{equation*}

\begin{equation*}
\varphi _{p,k+1}(t)=\frac{1}{(\widetilde{a}_{p},\widetilde{b}_{p})\cdot
l_{k+p}}\int_{0}^{t}\varphi _{p,k}(z)dz,\text{ }p=1,...,n-3,\text{ }%
k=1,...,n-p-2.
\end{equation*}

\textit{In the above formulas}
\begin{equation*}
\widetilde{a}_{p}=\frac{a_{p}b_{n}-a_{n}b_{p}}{a_{n-1}b_{n}-a_{n}b_{n-1}};~%
\widetilde{b}_{p}=\frac{a_{n-1}b_{p}-a_{p}b_{n-1}}{a_{n-1}b_{n}-a_{n}b_{n-1}}%
,~p=1,...,n-2,
\end{equation*}

\begin{equation*}
l_{p}=\left( \frac{\widetilde{b}_{p}}{\sqrt{\widetilde{a}_{p}^{2}+\widetilde{%
b}_{p}^{2}}},\frac{-\widetilde{a}_{p}}{\sqrt{\widetilde{a}_{p}^{2}+%
\widetilde{b}_{p}^{2}}}\right) ,\text{ }p=1,...,n-2.
\end{equation*}

\begin{equation*}
f^{\ast }(x,y)=f\left( \frac{b_{n}x-b_{n-1}y}{a_{n-1}b_{n}-a_{n}b_{n-1}},%
\frac{a_{n}x-a_{n-1}y}{a_{n}b_{n-1}-a_{n-1}b_{n}}\right) .
\end{equation*}

\bigskip

\begin{proof} Since the vectors $(a_{n-1},b_{n-1})$ and $(a_{n},b_{n})$
are linearly independent, there is a nonsingular linear transformation $%
S:(x,y)\rightarrow (x^{^{\prime }},y^{^{\prime }})$ such that $%
S:(a_{n-1},b_{n-1})\rightarrow (1,0)$ and $S:(a_{n},b_{n})\rightarrow (0,1).$
Thus, without loss of generality we may assume that the vectors $%
(a_{n-1},b_{n-1})$ and $(a_{n},b_{n})$ coincide with the coordinate vectors $%
e_{1}=(1,0)$ and $e_{2}=(0,1)$ respectively. Therefore, to prove the first
part of the theorem it is enough to show that if a function $f\in C^{s}(%
\mathbb{R}^{2})$ is expressed in the form

\begin{equation*}
f(x,y)=\sum_{i=1}^{n-2}f_{i}(a_{i}x+b_{i}y)+f_{n-1}(x)+f_{n}(y),
\end{equation*}%
with arbitrary $f_{i}$, then there exist functions $g_{i}$ $\in C^{s}(%
\mathbb{R})$, $i=1,...,n$, such that $f$ is also expressed in the form

\begin{equation*}
f(x,y)=\sum_{i=1}^{n-2}g_{i}(a_{i}x+b_{i}y)+g_{n-1}(x)+g_{n}(y).\eqno(2.70)
\end{equation*}

By $\Delta _{l}^{(\delta )}F$ we denote the increment of a function $F$ in a
direction $l=(l^{\prime },l^{\prime \prime }).$ That is,

\begin{equation*}
\Delta _{l}^{(\delta )}F(x,y)=F(x+l^{\prime }\delta ,y+l^{\prime \prime
}\delta )-F(x,y).
\end{equation*}%
We also use the notation $\frac{\partial F}{\partial l}$ which denotes the
derivative of $F$ in the direction $l$.

It is easy to check that the increment of a ridge function $g(ax+by)$ in a
direction perpendicular to $(a,b)$ is zero. Let $l_{1},...,l_{n-2}$ be unit
vectors perpendicular to the vectors $(a_{1},b_{1}),...,(a_{n-2},b_{n-2})$
correspondingly. Then for any set of numbers $\delta _{1},...,\delta
_{n-2}\in \mathbb{R}$ we have

\begin{equation*}
\Delta _{l_{1}}^{(\delta _{1})}\cdot \cdot \cdot \Delta _{l_{n-2}}^{(\delta
_{n-2})}f(x,y)=\Delta _{l_{1}}^{(\delta _{1})}\cdot \cdot \cdot \Delta
_{l_{n-2}}^{(\delta _{n-2})}\left[ f_{n-1}(x)+f_{n}(y)\right] .\eqno(2.71)
\end{equation*}

Denote the left hand side of (2.71) by $S(x,y).$ That is, set%
\begin{equation*}
S(x,y)\overset{def}{=}\Delta _{l_{1}}^{(\delta _{1})}\cdot \cdot \cdot
\Delta _{l_{n-2}}^{(\delta _{n-2})}f(x,y).
\end{equation*}%
Then from (2.71) it follows that for any real numbers $\delta _{n-1}$and $%
\delta _{n}$,

\begin{equation*}
\Delta _{e_{1}}^{(\delta _{n-1})}\Delta _{e_{2}}^{(\delta _{n})}S(x,y)=0,
\end{equation*}%
or in expanded form,

\begin{equation*}
S(x+\delta _{n-1},y+\delta _{n})-S(x,y+\delta _{n})-S(x+\delta
_{n-1},y)+S(x,y)=0.
\end{equation*}%
Putting in the last equality $\delta _{n-1}=-x,$ $\delta _{n}=-y$, we obtain
that

\begin{equation*}
S(x,y)=S(x,0)+S(0,y)-S(0,0).
\end{equation*}%
This means that

\begin{equation*}
\Delta _{l_{1}}^{(\delta _{1})}\cdot \cdot \cdot \Delta _{l_{n-2}}^{(\delta
_{n-2})}f(x,y)
\end{equation*}
\begin{equation*}
=\Delta _{l_{1}}^{(\delta _{1})}\cdot \cdot \cdot \Delta
_{l_{n-2}}^{(\delta _{n-2})}f(x,0)+\Delta _{l_{1}}^{(\delta _{1})}\cdot
\cdot \cdot \Delta _{l_{n-2}}^{(\delta _{n-2})}f(0,y)-\Delta
_{l_{1}}^{(\delta _{1})}\cdot \cdot \cdot \Delta _{l_{n-2}}^{(\delta
_{n-2})}f(0,0).
\end{equation*}

By the hypothesis of the theorem, the derivative $\frac{\partial ^{n-2}}{%
\partial l_{1}\cdot \cdot \cdot \partial l_{n-2}}f(x,y)$ exists at any point
$(x,y)\in $ $\mathbb{R}^{2}$. Thus, it follows from the above formula that

\begin{equation*}
\frac{\partial ^{n-2}f}{\partial l_{1}\cdot \cdot \cdot \partial l_{n-2}}%
(x,y)=h_{1,1}(x)+h_{2,1}(y),\eqno(2.72)
\end{equation*}%
where $h_{1,1}(x)=\frac{\partial ^{n-2}}{\partial l_{1}\cdot \cdot \cdot
\partial l_{n-2}}f(x,0)$ and $h_{2,1}(y)=\frac{\partial ^{n-2}}{\partial
l_{1}\cdot \cdot \cdot \partial l_{n-2}}f(0,y)-\frac{\partial ^{n-2}}{%
\partial l_{1}\cdot \cdot \cdot \partial l_{n-2}}f(0,0)$. Note that $h_{1,1}$
and $h_{2,1}$ belong to the class $C^{s-n+2}(\mathbb{R}).$

By $h_{1,2}$ and $h_{2,2}$ denote the antiderivatives of $h_{1,1}$ and $%
h_{2,1}$ satisfying the condition $h_{1,2}(0)=h_{2,2}(0)=0$ and multiplied
by the numbers $1/(e_{1}\cdot l_{1})$ and $1/(e_{2}\cdot l_{1})$
correspondingly. That is,
\begin{eqnarray*}
h_{1,2}(x) &=&\frac{1}{e_{1}\cdot l_{1}}\int_{0}^{x}h_{1,1}(z)dz; \\
h_{2,2}(y) &=&\frac{1}{e_{2}\cdot l_{1}}\int_{0}^{y}h_{2,1}(z)dz.
\end{eqnarray*}%
Here $e\cdot l$ denotes the scalar product between vectors $e$ and $l$.
Obviously, the function

\begin{equation*}
F_{1}(x,y)=h_{1,2}(x)+h_{2,2}(y)
\end{equation*}%
obeys the equality

\begin{equation*}
\frac{\partial F_{1}}{\partial l_{1}}(x,y)=h_{1,1}(x)+h_{2,1}(y).\eqno(2.73)
\end{equation*}%
From (2.72) and (2.73) we obtain that

\begin{equation*}
\frac{\partial }{\partial l_{1}}\left[ \frac{\partial ^{n-3}f}{\partial
l_{2}\cdot \cdot \cdot \partial l_{n-2}}-F_{1}\right] =0.
\end{equation*}%
Hence, for some ridge function $\varphi _{1,1}(a_{1}x+b_{1}y),$
\begin{equation*}
\frac{\partial ^{n-3}f}{\partial l_{2}\cdot \cdot \cdot \partial l_{n-2}}%
(x,y)=h_{1,2}(x)+h_{2,2}(y)+\varphi _{1,1}(a_{1}x+b_{1}y).\eqno(2.74)
\end{equation*}%
Here all the functions $h_{2,1},h_{2,2}(y),\varphi _{1,1}\in C^{s-n+3}(%
\mathbb{R}).$

Set the following functions
\begin{eqnarray*}
h_{1,3}(x) &=&\frac{1}{e_{1}\cdot l_{2}}\int_{0}^{x}h_{1,2}(z)dz; \\
h_{2,3}(y) &=&\frac{1}{e_{2}\cdot l_{2}}\int_{0}^{y}h_{2,2}(z)dz; \\
\varphi _{1,2}(t) &=&\frac{1}{(a_{1},b_{1})\cdot l_{2}}\int_{0}^{t}\varphi
_{1,1}(z)dz.
\end{eqnarray*}%
Note that the function
\begin{equation*}
F_{2}(x,y)=h_{1,3}(x)+h_{2,3}(y)+\varphi _{1,2}(a_{1}x+b_{1}y)
\end{equation*}%
obeys the equality

\begin{equation*}
\frac{\partial F_{2}}{\partial l_{2}}(x,y)=h_{1,2}(x)+h_{2,2}(y)+\varphi
_{1,1}(a_{1}x+b_{1}y).\eqno(2.75)
\end{equation*}%
From (2.74) and (2.75) it follows that

\begin{equation*}
\frac{\partial }{\partial l_{2}}\left[ \frac{\partial ^{n-4}f}{\partial
l_{3}\cdot \cdot \cdot \partial l_{n-2}}-F_{2}\right] =0.
\end{equation*}%
The last equality means that for some ridge function $\varphi
_{2,1}(a_{2}x+b_{2}y),$
\begin{equation*}
\frac{\partial ^{n-4}f}{\partial l_{3}\cdot \cdot \cdot \partial l_{n-2}}%
(x,y)=h_{1,3}(x)+h_{2,3}(y)+\varphi _{1,2}(a_{1}x+b_{1}y)+\varphi
_{2,1}(a_{2}x+b_{2}y).\eqno(2.76)
\end{equation*}%
Here all the functions $h_{1,3},$ $h_{2,3},$ $\varphi _{1,2},$ $\varphi
_{2,1}\in C^{s-n+4}(\mathbb{R}).$

Note that in the left hand sides of (2.72), (2.74) and (2.76) we have the
mixed directional derivatives of $f$ and the order of these derivatives is
decreased by one in each consecutive step. Continuing the above process,
until it reaches the function $f$, we obtain the desired representation
(2.70).

The formulas for $g_{i}$ are obtained in the process of the above proof.
These formulas involve certain functions which can be found inductively as
described in the proof. The validity of the formulas for the functions $%
h_{1,k}$ and $h_{2,k}$, $k=1,...,n-1,$ is obvious. The formulas for $\varphi
_{p,1}$ and $\varphi _{p,k+1}$ can be obtained from (2.74), (2.76) and the
subsequent (assumed but not written) equations if we put $x=\widetilde{a}%
_{p}t/(\widetilde{a}_{p}^{2}+\widetilde{b}_{p}^{2})$ and $y=\widetilde{b}%
_{p}t/(\widetilde{a}_{p}^{2}+\widetilde{b}_{p}^{2})$. Note that $(\widetilde{%
a}_{p},\widetilde{b}_{p}),$ $p=1,...,n-2,$ are the images of vectors $%
(a_{p},b_{p})$ under the linear transformation $S$ which takes the vectors $%
(a_{n-1},b_{n-1})$ and $(a_{n},b_{n})$ to the coordinate vectors $e_{1}=(1,0)$ and
$e_{2}=(0,1),$ respectively. Besides, note that for $p=1,...,n-2,$ the
vectors $l_{p}$ are perpendicular to the vectors $(\widetilde{a}_{p},%
\widetilde{b}_{p})$, respectively and $f^{\ast }$ is the function generated
from $f$ by the above liner transformation.
\end{proof}

Theorem 2.11 can be applied to some higher order partial differential
equations in two variables, e.g., to the following homogeneous equation

\begin{equation*}
\prod\limits_{i=1}^{r}\left( \alpha _{i}\frac{\partial }{\partial x}+\beta
_{i}\frac{\partial }{\partial y}\right) u(x,y)=0,\eqno(2.77)
\end{equation*}%
where $(\alpha _{i},\beta _{i}),~i=1,...,r,$ are pairwise linearly
independent vectors in $\mathbb{R}^{2}$. Clearly, the general solution to
this equation are all functions of the form

\begin{equation*}
u(x,y)=\sum\limits_{i=1}^{r}v_{i}(\beta _{i}x-\alpha _{i}y),\eqno(2.78)
\end{equation*}%
where $v_{i}\in C^{r}(\mathbb{R})$, $i=1,...,r$. Based on Theorem 2.11, for
the general solution, one can demand only smoothness of the sum $u$ and
dispense with smoothness of the summands $v_{i}$. More precisely, the
following corollary is valid.

\bigskip

\textbf{Corollary 2.1.} \textit{Assume a function $u\in C^{r}(\mathbb{R}%
^{2}) $ is of the form (2.78) with arbitrarily behaved $v_{i}$. Then $u$ is
a solution to Equation (2.77).}

\bigskip

\textbf{Remark 2.6.} If in Theorem 2.11 $s\geq n-1,$ then the functions $%
g_{i}$, $i=1,...,n,$ can be constructed (up to polynomials) by the method
discussed in Buhmann and Pinkus \cite{12}. This method is based on the fact
that for a direction $\mathbf{c}=(c_{1},...,c_{m})$ orthogonal to a given
direction $\mathbf{a}\in \mathbb{R}^{m}\backslash \{\mathbf{0}\},$ the
operator
\begin{equation*}
D_{\mathbf{c}}=\sum_{k=1}^{m}c_{k}\frac{\partial }{\partial x_{k}}
\end{equation*}%
acts on $m$-variable ridge functions $g(\mathbf{a}\cdot \mathbf{x})$ as
follows
\begin{equation*}
D_{\mathbf{c}}g(\mathbf{a}\cdot \mathbf{x})=\left( \mathbf{c}\cdot \mathbf{a}%
\right) g^{\prime }(\mathbf{a}\cdot \mathbf{x}).
\end{equation*}%
Thus, if in our case for fixed $r\in \{1,...,n\},$ vectors $l_{k},$ $k\in
\{1,...,n\}$, $k\neq r$, are perpendicular to the vectors $(a_{k},b_{k})$,
then
\begin{equation*}
\prod\limits_{\substack{ k=1  \\ k\neq r}}^{n}D_{l_{k}}f(x,y)=\prod\limits
_{\substack{ k=1  \\ k\neq r}}^{n}D_{l_{k}}\sum_{i=1}^{n}g_{i}(a_{i}x+b_{i}y)
\end{equation*}%
\begin{equation*}
=\sum_{i=1}^{n}\left( \prod\limits _{\substack{ k=1  \\ k\neq r}}^{n}\left(
(a_{i},b_{i})\cdot l_{k}\right) \right)
g_{i}^{(n-1)}(a_{i}x+b_{i}y)=\prod\limits_{\substack{ k=1  \\ k\neq r}}%
^{n}\left( (a_{r},b_{r})\cdot l_{k}\right) g_{r}^{(n-1)}(a_{r}x+b_{r}y).
\end{equation*}%
Now $g_{r}$ can be easily constructed from the above formula (up to a
polynomial of degree at most $n-2$). Note that this method is not feasible
if in Theorem 2.11 the function $f$ is of the class $C^{n-2}(\mathbb{R}^{2})$%
.

\bigskip

\subsection{Multivariate case}

In this subsection, we generalize ideas from the previous subsection to
prove constructively that if a multivariate function of a certain smoothness
class is represented by a sum of $k$ arbitrarily behaved ridge functions,
then, under suitable conditions, it can be represented by a sum of ridge
functions of the same smoothness class and some polynomial of a certain
degree. The appearance of a polynomial term is mainly related to the fact
that in $\mathbb{R}^{n}$ ($n\geq 3)$ there are many directions orthogonal to
a given direction. Such a result was proved nonconstructively in Section
2.1.5 (see Theorem 2.5), but here under a mild hypothesis on the degree of
smoothness, we give a new proof for this theorem, which will provide us with
a recipe for constructing the functions $g_{i}$ in (2.24).

The following theorem is valid.

\bigskip

\textbf{Theorem 2.12.} \textit{Assume $f\in C^{s}(\mathbb{R}^{n})$ is of the
form (2.1). Let $s\geq k-p+1,$ where $p$ is the number of vectors $\mathbf{a}%
^{i}$ forming a maximal linearly independent system. Then there exist
functions $g_{i}\in C^{s}(\mathbb{R})$ and a polynomial $P(\mathbf{x})$ of
total degree at most $k-p+1$ such that (2.24) holds and $g_{i}$ can be
constructed algorithmically.}

\bigskip

\begin{proof} We start the proof by choosing a maximal linearly
independent system in $\{\mathbf{a}^{1},....,\mathbf{a}^{k}\}$. The case
when the system $\{\mathbf{a}^{1},....,\mathbf{a}^{k}\}$ itself is linearly
independent is obvious (see Section 2.1.1). Thus we omit this special case
here. Without loss of generality we may assume that the first $p$ vectors $%
\mathbf{a}^{1},....,\mathbf{a}^{p}$, $p<k$, are linearly independent. Thus,
the vectors $\mathbf{a}^{j},$ $j=p+1,...,k,$ can be expressed as linear
combinations $\lambda _{1}^{j}\mathbf{a}^{1}+\cdot \cdot \cdot +\lambda
_{p}^{j}\mathbf{a}^{p}$, where $\lambda _{1}^{j},...,\lambda _{p}^{j}$ are
real numbers. In addition, we can always apply a nonsingular linear
transformation $S$ of the coordinates such that $S:\mathbf{a}^{i}\rightarrow
\mathbf{e}_{i},$ $i=1,...,p,$ where $\mathbf{e}_{i}$ denotes the $i$-th unit
vector. This reduces the initial representation (2.1) to the following
simpler form

\begin{equation*}
f(\mathbf{x})=f_{1}(x_{1})+\cdot \cdot \cdot
+f_{p}(x_{p})+\sum_{i=1}^{m}f_{p+i}(\mathbf{a}^{i}\cdot \mathbf{x}).\eqno%
(2.79)
\end{equation*}%
Note that we keep the notation of (2.1), but here $\mathbf{x}%
=(x_{1},...,x_{p}),$ $\mathbf{a}^{i}=(\lambda _{1}^{i}\mathbf{,...,}\lambda
_{p}^{i})\in \mathbb{R}^{p}$ and $m=k-p.$ Obviously, we prove Theorem 2.12
if we prove it for the representation (2.79). Thus, in the sequel, we prove
that if $f\in C^{s}(\mathbb{R}^{n})$ is of the form (2.79) and $s\geq m+1,$
then there exist functions $g_{i}\in C^{s}(\mathbb{R})$ and a polynomial $P(%
\mathbf{x})$ of total degree at most $m+1$ such that

\begin{equation*}
f(\mathbf{x})=g_{1}(x_{1})+\cdot \cdot \cdot
+g_{p}(x_{p})+\sum_{i=1}^{m}g_{p+i}(\mathbf{a}^{i}\cdot \mathbf{x})+P(%
\mathbf{x}).
\end{equation*}

In the process of the proof, we also see how these $g_{i}$ are constructed.

For each $i=1,...,m,$ let $\{\mathbf{e}_{1}^{(i)},...,\mathbf{e}%
_{p-1}^{(i)}\}$ denote an orthonormal basis in the hyperplane perpendicular
to $\mathbf{a}^{i}.$ By $\Delta _{\mathbf{e}}^{(\delta )}F$ we denote the
increment of a function $F$ in a direction $\mathbf{e}$ of length $\delta .$
That is,

\begin{equation*}
\Delta _{\mathbf{e}}^{(\delta )}F(\mathbf{x})=F(\mathbf{x}+\delta \mathbf{e}%
)-F(\mathbf{x}).
\end{equation*}%
We also use the notation $\frac{\partial F}{\partial \mathbf{e}}$ to denote
the derivative of $F$ in a direction $\mathbf{e}$.

It is easy to check that the increment of a ridge function $g(\mathbf{a\cdot
x})$ in any direction perpendicular to $\mathbf{a}$ is zero. For example,

\begin{equation*}
\Delta _{\mathbf{e}_{j}^{(i)}}^{(\delta )}g(\mathbf{a}^{i}\mathbf{\cdot x)}%
=0,
\end{equation*}%
for all $i=1,...,m,$ $j=1,...,p-1.$ Therefore, for any indices $%
i_{1},...,i_{m}\in \{1,...,p-1\}$, $q\in \{1,...,p\}$ and numbers $\delta
_{1},...,\delta _{m},\delta \in \mathbb{R}$ we have the formula

\begin{equation*}
\Delta _{\mathbf{e}_{i_{1}}^{(1)}}^{(\delta _{1})}\Delta _{\mathbf{e}%
_{i_{2}}^{(2)}}^{(\delta _{2})}\cdot \cdot \cdot \Delta _{\mathbf{e}%
_{i_{m}}^{(m)}}^{(\delta _{m})}\Delta _{\mathbf{e}_{q}}^{(\delta )}f(\mathbf{%
x})=\Delta _{\mathbf{e}_{i_{1}}^{(1)}}^{(\delta _{1})}\Delta _{\mathbf{e}%
_{i_{2}}^{(2)}}^{(\delta _{2})}\cdot \cdot \cdot \Delta _{\mathbf{e}%
_{i_{m}}^{(m)}}^{(\delta _{m})}\Delta _{\mathbf{e}_{q}}^{(\delta
)}f_{q}(x_{q}),
\end{equation*}%
where $\mathbf{e}_{q}$ denotes the $q$-th unit vector. This means that for
each $q=1,...,p,$ the mixed directional derivative

\begin{equation*}
\frac{\partial ^{m+1}f}{\partial \mathbf{e}_{i_{1}}^{(1)}\cdot \cdot \cdot
\partial \mathbf{e}_{i_{m}}^{(m)}\partial x_{q}}(\mathbf{x})
\end{equation*}%
depends only on the variable $x_{q}.$ Denote this derivative by $%
h_{i_{1},...,i_{m}}^{0,q}(x_{q})$:

\begin{equation*}
h_{i_{1},...,i_{m}}^{0,q}(x_{q})=\frac{\partial ^{m+1}f}{\partial \mathbf{e}%
_{i_{1}}^{(1)}\cdot \cdot \cdot \partial \mathbf{e}_{i_{m}}^{(m)}\partial
x_{q}}(\mathbf{x}),\text{ }q=1,...,p.\eqno(2.80)
\end{equation*}%
Since $f\in C^{s}(\mathbb{R}^{p}),$ we obtain that $%
h_{i_{1},...,i_{m}}^{0,q}\in C^{s-m-1}(\mathbb{R}).$ It follows from (2.80)
that

\begin{equation*}
d\left( \frac{\partial ^{m}f}{\partial \mathbf{e}_{i_{1}}^{(1)}\cdot \cdot
\cdot \partial \mathbf{e}_{i_{m}}^{(m)}}\right)
=h_{i_{1},...,i_{m}}^{0,1}(x_{1})dx_{1}+\cdot \cdot \cdot
+h_{i_{1},...,i_{m}}^{0,p}(x_{p})dx_{p}.\eqno(2.81)
\end{equation*}%
We conclude from (2.81) that

\begin{equation*}
\frac{\partial ^{m}f}{\partial \mathbf{e}_{i_{1}}^{(1)}\cdot \cdot \cdot
\partial \mathbf{e}_{i_{m}}^{(m)}}(\mathbf{x}%
)=h_{i_{1},...,i_{m}}^{1,1}(x_{1})+\cdot \cdot \cdot
+h_{i_{1},...,i_{m}}^{1,p}(x_{p})+c_{i_{1},...,i_{m}},\eqno(2.82)
\end{equation*}%
where the functions $h_{i_{1},...,i_{m}}^{1,q}(x_{q})$, $q=1,...,p,$ are
antiderivatives of $h_{i_{1},...,i_{m}}^{0,q}(x_{q})$ satisfying the
condition $h_{i_{1},...,i_{m}}^{1,q}(0)=0$ and $c_{i_{1},...,i_{m}}$ is a
constant. Note that $h_{i_{1},...,i_{m}}^{1,q}\in C^{s-m}(\mathbb{R}),$ $%
q=1,...,p.$ Obviously, for any pair $k,t\in \{1,...,p-1\}$,

\begin{equation*}
\frac{\partial ^{m+1}f}{\partial \mathbf{e}_{i_{1}}^{(1)}\cdot \cdot \cdot
\partial \mathbf{e}_{i_{m-1}}^{(m-1)}\partial \mathbf{e}_{k}^{(m)}\partial
\mathbf{e}_{t}^{(m)}}=\frac{\partial ^{m+1}f}{\partial \mathbf{e}%
_{i_{1}}^{(1)}\cdot \cdot \cdot \partial \mathbf{e}_{i_{m-1}}^{(m-1)}%
\partial \mathbf{e}_{t}^{(m)}\partial \mathbf{e}_{k}^{(m)}}\eqno(2.83)
\end{equation*}%
It follows from (2.82) and (2.83) that

\begin{equation*}
(\mathbf{e}_{q}\cdot \mathbf{e}_{k}^{(m)})\left(
h_{i_{1},...,i_{m-1},t}^{1,q}\right) ^{^{\prime }}(x_{q})=(\mathbf{e}%
_{q}\cdot \mathbf{e}_{t}^{(m)})\left( h_{i_{1},...,i_{m-1},k}^{1,q}\right)
^{^{\prime }}(x_{q})+c,
\end{equation*}%
where $c$ is a constant depending on the parameters $i_{1},...,i_{m-1},k,t$
and $q.$ Recall that by construction, $h_{i_{1},...,i_{m}}^{1,q}(0)=0.$ Hence

\begin{equation*}
(\mathbf{e}_{q}\cdot \mathbf{e}%
_{k}^{(m)})h_{i_{1},...,i_{m-1},t}^{1,q}(x_{q})=(\mathbf{e}_{q}\cdot \mathbf{%
e}_{t}^{(m)})h_{i_{1},...,i_{m-1},k}^{1,q}(x_{q})+cx_{q}.\eqno(2.84)
\end{equation*}

Since for each $q=1,...,p,$ the vectors $\mathbf{e}_{q}$ and $\mathbf{a}^{m}$
are linearly independent, there exists an index $i_{m}(q)\in \{1,...,p-1\}$
such that the vector $\mathbf{e}_{i_{m}(q)}^{(m)}$ is not orthogonal to $%
\mathbf{e}_{q}.$ That is, $\mathbf{e}_{q}\cdot \mathbf{e}_{i_{m}(q)}^{(m)}$ $%
\neq 0.$ For each $q=1,...,p,$ fix the index $i_{m}(q)$ and define the
following functions

\begin{equation*}
h_{i_{1},...,i_{m-1}}^{2,q}(x_{q})=\frac{1}{\mathbf{e}_{q}\cdot \mathbf{e}%
_{i_{m}(q)}^{(m)}}\int_{0}^{x_{q}}h_{i_{1},...,i_{m-1},i_{m}(q)}^{1,q}(z)dz.%
\eqno(2.85)
\end{equation*}%
and

\begin{equation*}
F_{i_{1},...,i_{m-1}}(\mathbf{x})=h_{i_{1},...,i_{m-1}}^{2,1}(x_{1})+\cdot
\cdot \cdot +h_{i_{1},...,i_{m-1}}^{2,p}(x_{p}).
\end{equation*}

\bigskip It is easy to obtain from (2.84) and (2.85) that for any $i_{m}\in
\{1,...,p-1\},$%
\begin{equation*}
\frac{\partial F_{i_{1},...,i_{m-1}}}{\partial \mathbf{e}_{i_{m}}^{(m)}}(%
\mathbf{x})=h_{i_{1},...,i_{m}}^{1,1}(x_{1})+\cdot \cdot \cdot
+h_{i_{1},...,i_{m}}^{1,p}(x_{p})+P_{i_{1},...,i_{m}}^{(1)},\eqno(2.86)
\end{equation*}%
where $P_{i_{1},...,i_{m}}^{(1)}$ is a polynomial of total degree not
greater than $1$. It follows from (2.82) and (2.86) that

\begin{equation*}
\frac{\partial }{\partial \mathbf{e}_{i_{m}}^{(m)}}\left[ \frac{\partial
^{m-1}f}{\partial \mathbf{e}_{i_{1}}^{(1)}\cdot \cdot \cdot \partial \mathbf{%
e}_{i_{m-1}}^{(m-1)}}-F_{i_{1},...,i_{m-1}}\right] (\mathbf{x}%
)=c_{i_{1},...,i_{m}}-P_{i_{1},...,i_{m}}^{(1)}(\mathbf{x}).\eqno(2.87)
\end{equation*}%
Note that the last equality is valid for all vectors $\mathbf{e}%
_{i_{m}}^{(m)},$ which form a basis in the hyperplane orthogonal to $\mathbf{%
a}^{m}$. Thus from (2.87) we conclude that the following expansion is valid

\begin{equation*}
\left.
\begin{array}{c}
\frac{\partial ^{m-1}f}{\partial \mathbf{e}_{i_{1}}^{(1)}\cdot \cdot \cdot
\partial \mathbf{e}_{i_{m-1}}^{(m-1)}}(\mathbf{x}%
)=h_{i_{1},...,i_{m-1}}^{2,1}(x_{1})+\cdot \cdot \cdot
+h_{i_{1},...,i_{m-1}}^{2,p}(x_{p}) \\
+\varphi _{i_{1},...,i_{m-1}}^{2,1}(\mathbf{a}^{m}\cdot \mathbf{x}%
)+P_{i_{1},...,i_{m-1}}^{(2)}(\mathbf{x}).%
\end{array}%
\right. \eqno(2.88)
\end{equation*}%
Here all the functions $%
h_{i_{1},...,i_{m-1}}^{2,1},...,h_{i_{1},...,i_{m-1}}^{2,p}(x_{p}),\varphi
_{i_{1},...,i_{m-1}}^{2,1}\in C^{s-m+1}(\mathbb{R})$ and $%
P_{i_{1},...,i_{m-1}}^{(2)}$ is a polynomial of total degree not greater
than $2.$

Since for each $q=1,...,p,$ the vector $\mathbf{e}_{q}$ is not collinear to $%
\mathbf{a}^{m-1},$ there is an index $i_{m-1}(q)\in \{1,...,p-1\}$ such that
$\mathbf{e}_{i_{m-1}(q)}^{(m-1)}$ is not orthogonal to $\mathbf{e}_{q}$.
Similarly, since $\mathbf{a}^{m-1}$ is not collinear to $\mathbf{a}^{m},$
there is an index $i_{m-1}(m)\in \{1,...,p-1\}$ such that $\mathbf{e}%
_{i_{m-1}(m)}^{(m-1)}$ is not orthogonal to $\mathbf{a}^{m}$. Fix the
indices $i_{m-1}(q)$, $i_{m-1}(m)$ and consider the following functions

\begin{equation*}
h_{i_{1},...,i_{m-2}}^{3,q}(x_{q})=\frac{1}{\mathbf{e}_{q}\cdot \mathbf{e}%
_{i_{m-1}(q)}^{(m-1)}}%
\int_{0}^{x_{q}}h_{i_{1},...,i_{m-2},i_{m-1}(q)}^{2,q}(z)dz,\text{ }%
q=1,...,p,\eqno(2.89)
\end{equation*}

\begin{equation*}
\varphi _{i_{1},...,i_{m-2}}^{3,1}(t)=\frac{1}{\mathbf{a}^{m}\cdot \mathbf{e}%
_{i_{m-1}(m)}^{(m-1)}}\int_{0}^{t}\varphi
_{i_{1},...,i_{m-2},i_{m-1}(m)}^{2,1}(z)dz,\eqno(2.90)
\end{equation*}%
and

\begin{equation*}
F_{i_{1},...,i_{m-2}}(\mathbf{x})=h_{i_{1},...,i_{m-2}}^{3,1}(x_{1})+\cdot
\cdot \cdot +h_{i_{1},...,i_{m-2}}^{3,p}(x_{p})+\varphi
_{i_{1},...,i_{m-2}}^{3,1}(\mathbf{a}^{m}\cdot \mathbf{x}).\eqno(2.91)
\end{equation*}

Similar to (2.84), one can easily verify that for any pair $k,t\in
\{1,...,p-1\}$ and for all $q=1,...,p,$ the following equalities are valid.
\begin{equation*}
\left.
\begin{array}{c}
(\mathbf{e}_{q}\cdot \mathbf{e}%
_{k}^{(m-1)})h_{i_{1},...,i_{m-2},t}^{2,q}(x_{q})=(\mathbf{e}_{q}\cdot
\mathbf{e}_{t}^{(m-1)})h_{i_{1},...,i_{m-2},k}^{2,q}(x_{q})+H_{q}(x_{q}), \\
(\mathbf{a}^{m}\cdot \mathbf{e}_{k}^{(m-1)})\varphi
_{i_{1},...,i_{m-2},t}^{2,1}(\mathbf{a}^{m}\cdot \mathbf{x})=(\mathbf{a}%
^{m}\cdot \mathbf{e}_{t}^{(m-1)})\varphi _{i_{1},...,i_{m-2},k}^{2,1}(%
\mathbf{a}^{m}\cdot \mathbf{x})+\Phi (\mathbf{x}),%
\end{array}%
\right. \eqno(2.92)
\end{equation*}%
where $H_{q}$ and $\Phi $ are univariate and $n$-variable polynomials of
degree not greater than $2.$ Indeed, applying the Schwarz formula
\begin{equation*}
\frac{\partial ^{m+1}f}{\partial \mathbf{e}_{i_{1}}^{(1)}\cdot \cdot \cdot
\partial \mathbf{e}_{i_{m-2}}^{(m-2)}\partial \mathbf{e}_{k}^{(m-1)}\partial
\mathbf{e}_{t}^{(m-1)}\partial \mathbf{e}_{i_{m}(q)}^{(m)}}=\frac{\partial
^{m+1}f}{\partial \mathbf{e}_{i_{1}}^{(1)}\cdot \cdot \cdot \partial \mathbf{%
e}_{i_{m-2}}^{(m-2)}\partial \mathbf{e}_{t}^{(m-1)}\partial \mathbf{e}%
_{k}^{(m-1)}\partial \mathbf{e}_{i_{m}(q)}^{(m)}}
\end{equation*}%
on the symmetry of derivatives, it follows from (2.82) that for any $k,t\in
\{1,...,p-1\}$
\begin{equation*}
(\mathbf{e}_{q}\cdot \mathbf{e}_{k}^{(m-1)})\left(
h_{i_{1},...,i_{m-2},t,i_{m}(q)}^{1,q}\right) ^{^{\prime }}(x_{q})=(\mathbf{e%
}_{q}\cdot \mathbf{e}_{t}^{(m-1)})\left(
h_{i_{1},...,i_{m-2},k,i_{m}(q)}^{1,q}\right) ^{^{\prime }}(x_{q})+d,
\end{equation*}%
where $d$ is a constant depending on the parameters $i_{1},...,i_{m-2},k,t$
and $i_{m}(q).$ Since, by construction, $h_{i_{1},...,i_{m}}^{1,q}(0)=0,$ we
obtain that

\begin{equation*}
(\mathbf{e}_{q}\cdot \mathbf{e}%
_{k}^{(m-1)})h_{i_{1},...,i_{m-2},t,i_{m}(q)}^{1,q}(x_{q})=(\mathbf{e}%
_{q}\cdot \mathbf{e}%
_{t}^{(m-1)})h_{i_{1},...,i_{m-2},k,i_{m}(q)}^{1,q}(x_{q})+dx_{q}.
\end{equation*}%
The last equality together with (2.85) yield that

\begin{equation*}
(\mathbf{e}_{q}\cdot \mathbf{e}_{k}^{(m-1)})\left(
h_{i_{1},...,i_{m-2},t}^{2,q}\right) ^{^{\prime }}(x_{q})=(\mathbf{e}%
_{q}\cdot \mathbf{e}_{t}^{(m-1)})\left( h_{i_{1},...,i_{m-2},k}^{2,q}\right)
^{^{\prime }}(x_{q})+dx_{q}.
\end{equation*}%
Therefore, the first equality in (2.92) holds. Considering this and applying
the corresponding Schwarz formula to (2.88) we obtain the second equality in
(2.92).

Taking into account the definitions (2.89), (2.90) and the relations (2.92),
we obtain from (2.91) that for any $i_{m-1}\in \{1,...,p-1\},$%
\begin{equation*}
\frac{\partial F_{i_{1},...,i_{m-2}}}{\partial \mathbf{e}_{i_{m-1}}^{(m-1)}}(%
\mathbf{x})
\end{equation*}
\begin{equation*}
=h_{i_{1},...,i_{m-1}}^{2,1}(x_{1})+\cdot \cdot \cdot
+h_{i_{1},...,i_{m-1}}^{2,p}(x_{p})+\varphi _{i_{1},...,i_{m-1}}^{2,1}(%
\mathbf{a}^{m}\cdot \mathbf{x})+\widetilde{P}_{i_{1},...,i_{m-1}}^{(2)}(%
\mathbf{x}),\eqno(2.93)
\end{equation*}%
where $\widetilde{P}_{i_{1},...,i_{m-1}}^{(2)}$ is a polynomial of degree
not greater than $2.$ It follows from (2.88) and (2.93) that
\begin{equation*}
\frac{\partial }{\partial \mathbf{e}_{i_{m-1}}^{(m-1)}}\left[ \frac{\partial
^{m-2}f}{\partial \mathbf{e}_{i_{1}}^{(1)}\cdot \cdot \cdot \partial \mathbf{%
e}_{i_{m-2}}^{(m-2)}}-F_{i_{1},...,i_{m-2}}\right] (\mathbf{x}%
)=P_{i_{1},...,i_{m-1}}^{(2)}(\mathbf{x})-\widetilde{P}%
_{i_{1},...,i_{m-1}}^{(2)}(\mathbf{x}).\eqno(2.94)
\end{equation*}%
Note that the last equality is valid for all vectors $\mathbf{e}%
_{i_{m-1}}^{(m-1)},$ which form a basis in the hyperplane orthogonal to $%
\mathbf{a}^{m-1}$. Considering this, from (2.94) we derive the following
representation
\begin{equation*}
\left.
\begin{array}{c}
\frac{\partial ^{m-2}f}{\partial \mathbf{e}_{i_{1}}^{(1)}\cdot \cdot \cdot
\partial \mathbf{e}_{i_{m-2}}^{(m-2)}}(\mathbf{x}%
)=h_{i_{1},...,i_{m-2}}^{3,1}(x_{1})+\cdot \cdot \cdot
+h_{i_{1},...,i_{m-2}}^{3,p}(x_{p}) \\
+\varphi _{i_{1},...,i_{m-2}}^{3,1}(\mathbf{a}^{m}\cdot \mathbf{x})+\varphi
_{i_{1},...,i_{m-2}}^{3,2}(\mathbf{a}^{m-1}\cdot \mathbf{x}%
)+P_{i_{1},...,i_{m-2}}^{(3)}(\mathbf{x}).%
\end{array}%
\right. \eqno(2.95)
\end{equation*}%
Here all the functions \newline $%
h_{i_{1},...,i_{m-2}}^{3,1},...,h_{i_{1},...,i_{m-2}}^{3,p}(x_{p}),$ $%
\varphi _{i_{1},...,i_{m-2}}^{3,1},$ $\varphi _{i_{1},...,i_{m-2}}^{3,2}\in
C^{s-m+2}(\mathbb{R})$ and $P_{i_{1},...,i_{m-2}}^{(3)}$ is a polynomial of
total degree not greater than $3.$

Note that in the left hand sides of (2.82), (2.88) and (2.95) we have the
mixed directional derivatives of $f$ and the order of these derivatives is
decreased by one at each consecutive step. Continuing the above process,
until it reaches the function $f$, we obtain the desired representation.
Note that the above proof gives a recipe for constructing the smooth ridge
functions $g_{i}$. Writing out explicit recurrent formulas for $g_{i}$, as
in Theorem 2.11, is technically cumbersome here and hence is avoided.
\end{proof}

\textbf{Remark 2.7.} Note that using Theorem 2.12, the degree of polynomial $%
P(\mathbf{x})$ in Theorem 2.5 can be reduced. Indeed, it follows from (2.27)
and (2.28) that the the above polynomial $P(\mathbf{x})$ is of the form
(2.1). On the other hand, by Theorem 2.12 there exist functions $g_{i}^{\ast
}\in C^{s}(\mathbb{R})$, $i=1,...,k$, and a polynomial $G(\mathbf{x})$ of
degree at most $k-p+1$ such that

\begin{equation*}
P(\mathbf{x})=\sum_{i=1}^{k}g_{i}^{\ast }(\mathbf{a}^{i}\cdot \mathbf{x})+G(%
\mathbf{x}).
\end{equation*}%
Now considering this in (2.24) we see that our assertion is true.

\bigskip

At the end of this chapter, we want to draw the reader's attention to the
following uniqueness question. Assume we are given pairwise linearly
independent vectors $\mathbf{a}^{i},$ $i=1,...,k,$ in $\mathbb{R}^{n}$ and a
function $f:\mathbb{R}^{n}\rightarrow \mathbb{R}$ of the form (2.1). How
many different ways can $f$ be written as a sum of ridge functions with the
directions $\mathbf{a}^{i}$? Clearly, representation (2.1) is not
unique, since we can always add some constants $c_{i}$ to $f_{i}$ without
changing the resulting sum in (2.1) provided that $\sum_{i=1}^{k}c_{i}=0$. It turns out
that under minimal requirements representation (2.1) is unique up to
polynomials of degree at most $k-2$. More precisely, if, in addition to (2.1), $%
f$ also has the form (2.2) and $%
f_{i},g_{i}\in \mathcal{B}$, $i=1,...,k$, then the functions $f_{i}-g_{i}$ are univariate
polynomials of degree at most $k-2$. This result is due to Pinkus
\cite[Theorem 3.1]{117}. It follows immediately from
this result that in Theorems 2.6--2.10 the functions $g_{i}$ is unique up to
a univariate polynomial. This is also valid for $g_{i}$ in Theorems 2.4 and
2.5, but in this case for the proof we must apply a slightly different result of Pinkus
\cite[Corollary 3.2]{117}: Assume a multivariate
polynomial $f$ of degree $m$ is of the form (2.1) and $f_{i}\in \mathcal{B}$
for $i=1,...,k$. Then $f_{i}$ are univariate polynomials of degree at most $%
l=\max \left\{ m,k-2\right\} $.

A different uniqueness problem, in a more general setting, will be analyzed
in Chapter 4. In that problem we will look for sets $Q\subset \mathbb{R}^{n}$
for which representation (2.1), considered on $Q$, is unique.

\newpage

\chapter{Approximation of multivariate functions by sums of univariate
functions}

It is clear that in the special case, when directions of ridge functions
coincide with the coordinate directions, the problem of approximation by
linear combinations of these functions turn into the problem of
approximation by sums of univariate functions. This is also the simplest
case in ridge function approximation. The simplicity of the approximation
guarantees its practicability in application areas, where complicated
multivariate functions are main obstacles. In mathematics, this type of
approximation has arisen, for example, in connection with the classical
functional equations \cite{11}, the numerical solution of certain PDE
boundary value problems \cite{9}, dimension theory \cite{133,132}, etc. In
this chapter, we obtain some results concerning the problem of best
approximation by sums of univariate functions.

Most of the material of this chapter is taken from \cite{59,55,56,48}.

\bigskip

\section{Characterization of some bivariate function classes by formulas for
the error of approximation}

This section is devoted to calculation formulas for the error of
approximation of bivariate functions by sums of univariate functions.
Certain classes of bivariate functions depending on some numerical parameter
are constructed and characterized in terms of the approximation error
calculation formulas.

\subsection{Exposition of the problem}

The approximation problem considered here is to approximate a continuous and
real-valued function of two variables by sums of two continuous functions of
one variable. To make the problem precise, let $Q$ be a compact set in the $%
xOy$ plane. Consider the approximation of a continuous function $f \in C(Q)$
by functions from the manifold $D=\left\{ \varphi (x)+\psi (y)\right\} ,$
where $\varphi (x),\psi (y)$ are defined and continuous on the projections
of $Q$ into the coordinate axes $x$ and $y$, respectively. The approximation
error is defined as the distance from $f$ to $D:$
\begin{equation*}
E(f)=dist(f,D)=\inf\limits_{D}\left\Vert f-\varphi -\psi \right\Vert _{C(Q)}=
\end{equation*}%
\begin{equation*}
=\inf\limits_{D}\max\limits_{(x,y)\in Q}\left\vert f(x,y)-\varphi (x)-\psi
(y)\right\vert.
\end{equation*}%
A function $\varphi _{0}(x)+\psi _{0}(y)$ from $D$, if it exists, is called
an extremal element or a best approximating sum if
\begin{equation*}
E(f)=\left\Vert f-\varphi _{0}-\psi _{0}\right\Vert _{C(Q)}.
\end{equation*}%
To show that $E(f)$ depends also on $Q$, in some cases to avoid confusion,
we will write $E(f,Q)$ instead of $E(f)$.

In this section we deal with calculation formulas for $E(f)$. In 1951
Diliberto and Straus published a paper \cite{26}, in which along with other
results they established a formula for $E(f,R)$, where $R$ here and
throughout this section is a rectangle with sides parallel to the coordinate
axes, containing supremum over all closed lightning bolts. Later the same
formula was established by other authors differently, in cases of both
rectangle (see \cite{113}) and more general sets (see \cite{79,107}). Although the formula was valid for all continuous functions, it was not
easily calculable. Some authors started to seek easily calculable formulas
for the approximation error for some subsets of continuous functions. Rivlin
and Sibner \cite{121} proved a result, which allow one to find the exact
value of $E(f,R)$ for a function $f(x,y)$ having the continuous and
nonnegative derivative $\frac{\partial ^{2}f}{\partial x\partial y}$. This
result in a more general case (for functions of $n$ variables) was proved by
Flatto \cite{30}. Babaev \cite{6} generalized Rivlin and Sibner's result (as
well as Flatto's result, see \cite{7}). More precisely, he considered the
class $M(R)$ of continuous functions $f(x,y) $ with the property
\begin{equation*}
\Delta
_{h_{1},h_{1}}f=f(x,y)+f(x+h_{1},y+h_{2})-f(x,y+h_{2})-f(x+h_{1},y)\geq 0
\end{equation*}%
for each rectangle $\left[ x,x+h_{1}\right] \times \left[ y,y+h_{2}\right]
\subset R$, and proved that if $f(x,y)$ belongs to $M(R)$, where $R=\left[
a_{1},b_{1}\right] \times \left[ a_{2},b_{2}\right] $, then
\begin{equation*}
E(f,R)=\frac{1}{4}\left[
f(a_{1},a_{2})+f(b_{1},b_{2})-f(a_{1},b_{2})-f(b_{1},a_{2})\right] .
\end{equation*}%
As seen from this formula, to calculate $E(f)$ it is sufficient to find only
values of $f(x,y)$ at the vertices of $R$. One can see that the formula also
gives a sufficient condition for membership in the class $M(R)$, i.e. if
\begin{equation*}
E(f,S)=\frac{1}{4}\left[
f(x_{1},y_{1})+f(x_{2},y_{2})-f(x_{1},y_{2})-f(x_{2},y_{1})\right] ,
\end{equation*}%
for a given $f$ and for each $S=\left[ x_{1},x_{2}\right] \times \left[
y_{1},y_{2}\right] \subset R$, then the function $f(x,y)$ is from $M(R)$.

Our purpose is to construct new classes of continuous functions, which will
depend on a numerical parameter, and characterize each class in terms of the
approximation error calculation formulas. The mentioned parameter will show
which points of $R$ the calculation formula involves. We will also construct
a best approximating sum $\varphi _{0}+\psi _{0}$ to a function from
constructed classes.

\bigskip

\subsection{Definition of the main classes}

Let throughout this section $R=\left[ a_{1},b_{1}\right] \times \left[
a_{2},b_{2}\right] $ be a rectangle and $c\in (a_{1},b_{1}]$. Denote $R_{1}=%
\left[ a_{1},c\right] \times \left[ a_{2},b_{2}\right] $ and $R_{2}=\left[
c,b_{1}\right] \times \left[ a_{2},b_{2}\right] $. It is clear that $%
R=R_{1}\cup R_{2}$ and if $c=b_{1}$, then $R=R_{1}$.

We associate each rectangle $S=\left[ x_{1},x_{2}\right] \times \left[
y_{1},y_{2}\right] $ lying in $R$ with the following functional:
\begin{equation*}
L(f,S)=\frac{1}{4}\left[
f(x_{1},y_{1})+f(x_{2},y_{2})-f(x_{1},y_{2})-f(x_{2},y_{1})\right] .
\end{equation*}

\bigskip

\textbf{Definition 3.1.} \textit{We say that a continuous function $f(x,y)$
belongs to the class $V_{c}(R)$ if }

\textit{1) $L(f,S)\geq 0$, for each $S\subset R_{1}$; }

\textit{2) $L(f,S)\leq 0$, for each $S\subset R_{2}$; }

\textit{3) $L(f,S)\geq 0$, for each $S=\left[ a_{1},b_{1}\right] \times %
\left[ y_{1},y_{2}\right] ,~\ S\subset R$.}

\bigskip

It can be shown that for any $c\in (a_{1},b_{1}]$ the class $V_{c}(R)$ is
not empty. Indeed, one can easily verify that the function
\begin{equation*}
v _{c}(x,y)=\left\{
\begin{array}{c}
w(x,y)-w(c,y),\;\ (x,y)\in R_{1} \\
w(c,y)-w(x,y),\;\ (x,y)\in R_{2}%
\end{array}%
\right.
\end{equation*}
where $w(x,y)=\left( \frac{x-a_{1}}{b_{1}-a_{1}}\right) ^{\frac{1}{n}}\cdot
y $ and $n\geq \log _{2}\frac{b_{1}-a_{1}}{c-a_{1}}$, satisfies conditions
1)-3) and therefore belongs to $V_{c}(R)$. The class $V_{c}(R)$ has the
following obvious properties:\newline
a) For given functions $f_{1},f_{2}\in V_{c}(R)$ and numbers $\alpha
_{1},\alpha _{2}\geq 0$, $\alpha _{1}f_{1}+\alpha _{2}f_{2}\in V_{c}(R)$. $%
V_{c}(R)$ is a closed subset of the space of continuous functions.\newline
b) $V_{b_{1}}(R)=M(R)$.\newline
c) If $f$ is a common element of $V_{c_{1}}(R)$ and $V_{c_{2}}(R)$, $%
a_{1}<c_{1}<c_{2}\leq b_{1}$ then $f(x,y)=\varphi (x)+\psi (y)$ on the
rectangle $\left[ c_{1},c_{2}\right] \times \left[ a_{2},b_{2}\right] $.

The properties a) and b) are clear. The property c) also becomes clear if
note that according to the definition of the classes $V_{c_{1}}(R)$ and $%
V_{c_{2}}(R)$, for each rectangle
\begin{equation*}
S\subset \left[ c_{1},c_{2}\right] \times \left[ a_{2},b_{2}\right]
\end{equation*}
we have
\begin{equation*}
L(f,S)\leq 0\;\; \mbox{and}\;\; L(f,S)\geq 0,
\end{equation*}
respectively. Hence
\begin{equation*}
L(f,S)=0\;\; \mbox{for each}\;\; S\subset \left[ c_{1},c_{2}\right] \times %
\left[ a_{2},b_{2}\right].
\end{equation*}
Thus it is not difficult to understand that $f$ is of the form $\varphi
(x)+\psi (y)$ on the rectangle $\left[ c_{1},c_{2}\right] \times \left[
a_{2},b_{2}\right] $.

\bigskip

\textbf{Lemma 3.1.} \textit{Assume a function $f(x,y)$ has the continuous
derivative $\frac{\partial ^{2}f}{\partial x\partial y}$ on the rectangle $R$
and satisfies the following conditions}

\textit{1) $\frac{\partial ^{2}f}{\partial x\partial y}\geq 0$, for all $%
(x,y)\in R_{1}$; }

\textit{2) $\frac{\partial ^{2}f}{\partial x\partial y}\leq 0$, for all $%
(x,y)\in R_{2}$; }

\textit{3) $\frac{df(a_{1},y)}{dy}\leq \frac{df(b_{1},y)}{dy}$, for all $%
y\in \left[ a_{2},b_{2}\right]$.}

\textit{Then $f(x,y)$ belongs to $V_{c}(R)$.}

\bigskip

The proof of this lemma is very simple and can be obtained by integrating
both sides of inequalities in conditions 1)-3) through sets $\left[
x_{1},x_{2}\right] \times \left[ y_{1},y_{2}\right] \subset R_{1}$, $\left[
x_{1},x_{2}\right] \times \left[ y_{1},y_{2}\right] \subset R_{2}$ and $%
\left[ y_{1},y_{2}\right] \subset \left[ a_{2},b_{2}\right]$, respectively.

\bigskip

\textbf{Example 3.1.} Consider the function $f(x,y)=y\sin \pi x$ on the unit
square $K=\left[ 0,1\right] \times \left[ 0,1\right] $ and rectangles $K_{1}=%
\left[ 0,\frac{1}{2}\right] \times \left[ 0,1\right] ,K_{2}=\left[ \frac{1}{2%
},1\right] \times \left[ 0,1\right] $. It is not difficult to verify that
this function satisfies all conditions of the lemma and therefore belongs to
$V_{\frac{1}{2}}(K)$.

\bigskip

\subsection{Construction of an extremal element}

The following theorem is valid.

\bigskip

\textbf{Theorem 3.1.}\ \textit{The approximation error of a function $f(x,y)$
from the class $V_{c}(R)$ can be calculated by the formula
\begin{equation*}
E(f,R)=L(f,R_{1})=\frac{1}{4}\left[
f(a_{1},a_{2})+f(c,b_{2})-f(a_{1},b_{2})-f(c,a_{2})\right] .
\end{equation*}%
Let $y_{0}$ be any solution from $\left[ a_{2},b_{2}\right] $ of the
equation
\begin{equation*}
L(f,Y)=\frac{1}{2}L(f,R_{1}),\;\;Y=\left[ a_{1},c\right] \times \left[
a_{2},y\right] .
\end{equation*}%
Then the function $\varphi _{0}(x)+\psi _{0}(y)$, where
\begin{equation*}
\varphi _{0}(x)=f(x,y_{0}),
\end{equation*}%
\begin{equation*}
\psi _{0}(y)=\frac{1}{2}\left[ f(a_{1},y)+f(c,y)-f(a_{1},y_{0})-f(c,y_{0})%
\right]
\end{equation*}%
is a best approximating sum from the manifold $D$ to $f$.}

\bigskip

To prove this theorem we need the following lemma.

\bigskip

\textbf{Lemma 3.2.}\ \textit{Let $f(x,y)$ be a function from $V_{c}(R)$ and $%
X=\left[ a_{1},x\right] \times \left[ y_{1},y_{2}\right] $ be a rectangle
with fixed $y_{1},y_{2}\in \left[ a_{2},b_{2}\right] $. Then the function $%
h(x)=L(f,X)$ has the properties: }

\textit{1) $h(x)\geq 0$, for any $x\in \left[ a_{1},b_{1}\right] $; }

\textit{2) $\max\limits_{\left[ a_{1},b_{1}\right] }h(x)=h(c)$ and $%
\min\limits_{\left[ a_{1},b_{1}\right] }h(x)=h(a_{1})=0$.}

\bigskip

\textbf{Proof.} If $X\subset R_{1}$, then the validity of $h(x)\geq 0$
follows from the definition of $V_{c}(R)$. If $X$ is from $R$ but not lying
in $R_{1}$, then by denoting $X^{\prime }=\left[ x,b_{1}\right] \times \left[
y_{1},y_{2}\right] ,S=X\cup X^{\prime }$ and using the obvious equality
\begin{equation*}
L(f,S)=L(f,X)+L(f,X^{\prime })
\end{equation*}
we deduce from the definition of $V_{c}(R)$ that $h(x)\geq 0$.

To prove the second part of the lemma, it is enough to show that $h(x)$
increases on the interval $\left[ a_{1},c\right] $ and decreases on the
interval $\left[ c,b_{1}\right] $. Indeed, if $a_{1}\leq x_{1}\leq x_{2}\leq
c$, then
\begin{equation*}
h(x_{2})=L(f,X_{2})=L(f,X_{1})+L(f,X_{12}),\eqno(3.1)
\end{equation*}%
where $X_{1}=\left[ a_{1},x_{1}\right] \times \left[ y_{1},y_{2}\right] ,$ $%
X_{2}=\left[ a_{1},x_{2}\right] \times \left[ y_{1},y_{2}\right] ,$ $X_{12}=%
\left[ x_{1},x_{2}\right] \times \left[ y_{1},y_{2}\right] $. Taking into
consideration that $L(f,X_{1})=h(x_{1})$ and $X_{12}$ lies in $R_{1}$ we
obtain from (3.1) that $h(x_{2})\geq h(x_{1})$. If $c\leq x_{1}\leq
x_{2}\leq b_{1}$, then $X_{12}$ lies in $R_{2}$ and we obtain from (3.1)
that $h(x_{2})\leq h(x_{1})$.

\bigskip

\textbf{Proof of Theorem 3.1.} It is obvious that $L(f,R_{1})=L(f-\varphi
-\psi ,R_{1})$ for each sum $\varphi (x)+\psi (y)$. Hence
\begin{equation*}
L(f,R_{1})\leq \left\Vert f-\varphi -\psi \right\Vert _{C(R_{1})}\leq
\left\Vert f-\varphi -\psi \right\Vert _{C(R)}.
\end{equation*}%
Since a sum $\varphi (x)+\psi (y)$ is arbitrary, $L(f,R_{1})\leq E(f,R)$. To
complete the proof it is sufficient to construct a sum $\varphi _{0}(x)+\psi
_{0}(y)$ for which the equality
\begin{equation*}
\left\Vert f-\varphi _{0}-\psi _{0}\right\Vert _{C(R)}=L(f,R_{1})\eqno(3.2)
\end{equation*}%
holds.

Consider the function
\begin{equation*}
g(x,y)=f(x,y)-f(x,a_{2})-f(a_{1},y)+f(a_{1},a_{2}).
\end{equation*}%
This function has the following obvious properties

1) $g(x,a_{2})=g(a_{1},y)=0$;

2) $L(f,R_{1})=L(g,R_{1})=\frac{1}{4}g(c,b_{2})$;

3) $E(f,R)=E(g,R)$;

4) The function of one variable $g(c,y)$ increases on the interval $\left[
a_{2},b_{2}\right] $.

The last property of $g$ allows us to write that
\begin{equation*}
0=g(c,a_{2})\leq \frac{1}{2}g(c,b_{2})\leq g(c,b_{2}).
\end{equation*}%
Since $g(x,y)$ is continuous, there exists at least one solution $y=y_{0}$
of the equation
\begin{equation*}
g(c,y)=\frac{1}{2}g(c,b_{2})
\end{equation*}%
or, in other notation, of the equation
\begin{equation*}
L(f,Y)=\frac{1}{2}L(f,R_{1}),\;\;\text{where}\;\;Y=\left[ a_{1},c\right]
\times \left[ a_{2},y\right] ,
\end{equation*}%
Introduce the functions
\begin{equation*}
\varphi _{1}(x)=g(x,y_{0}),
\end{equation*}%
\begin{equation*}
\psi _{1}(y)=\frac{1}{2}\left( g(c,y)-g(c,y_{0})\right) ,
\end{equation*}%
\begin{equation*}
G(x,y)=g(x,y)-\varphi _{1}(x)-\psi _{1}(y).
\end{equation*}%
Calculate the norm of $G(x,y)$ on $R$. Consider the rectangles $R^{\prime }=%
\left[ a_{1},b_{1}\right] \times \left[ y_{0},b_{2}\right] $ and $R^{\prime
\prime }=\left[ a_{1},b_{1}\right] \times \left[ a_{2},y_{0}\right] $. It is
clear that
\begin{equation*}
\left\Vert G\right\Vert _{C(R)}=\max \left\{ \left\Vert G\right\Vert
_{C(R^{\prime })},\left\Vert G\right\Vert _{C(R^{\prime \prime })}\right\} .
\end{equation*}%
First calculate the norm $\left\Vert G\right\Vert _{C(R^{\prime })}$:
\begin{equation*}
\left\Vert G\right\Vert _{C(R^{\prime })}=\max\limits_{(x,y)\in R^{\prime
}}\left\vert G(x,y)\right\vert =\max\limits_{y\in \left[ y_{0},b_{2}\right]
}\max\limits_{x\in \left[ a_{1},b_{1}\right] }\left\vert G(x,y)\right\vert .%
\eqno(3.3)
\end{equation*}%
For a fixed point $y$ (we keep it fixed until (3.6)) from the interval $%
\left[ y_{0},b_{2}\right] $ we can write that
\begin{equation*}
\max\limits_{x\in \left[ a_{1},b_{1}\right] }G(x,y)=\max\limits_{x\in \left[
a_{1},b_{1}\right] }\left( g(x,y)-g(x,y_{0})\right) -\psi _{1}(y)\eqno(3.4)
\end{equation*}%
and
\begin{equation*}
\min\limits_{x\in \left[ a_{1},b_{1}\right] }G(x,y)=\min\limits_{x\in \left[
a_{1},b_{1}\right] }\left( g(x,y)-g(x,y_{0})\right) -\psi _{1}(y).\eqno(3.5)
\end{equation*}%
By Lemma 3.2, the function
\begin{equation*}
h_{1}(x)=4L(f,X)=g(x,y)-g(x,y_{0}),\;\;\mbox{where}\;\;X=\left[ a_{1},x%
\right] \times \left[ y_{0},y\right] ,
\end{equation*}%
reaches its maximum on $x=c$ and minimum on $x=a_{1}$:
\begin{equation*}
\max\limits_{x\in \left[ a_{1},b_{1}\right] }h_{1}(x)=g(c,y)-g(c,y_{0})
\end{equation*}%
\begin{equation*}
\min\limits_{x\in \left[ a_{1},b_{1}\right]
}h_{1}(x)=g(a_{1},y)-g(a_{1},y_{0})=0.
\end{equation*}%
Considering these facts in (3.4) and (3.5) we obtain that
\begin{equation*}
\max\limits_{x\in \left[ a_{1},b_{1}\right] }G(x,y)=g(c,y)-g(c,y_{0})-\psi
_{1}(y)=\frac{1}{2}\left( g(c,y)-g(c,y_{0})\right) ,
\end{equation*}%
\begin{equation*}
\min\limits_{x\in \left[ a_{1},b_{1}\right] }G(x,y)=-\psi _{1}(y)=-\frac{1}{2%
}\left( g(c,y)-g(c,y_{0})\right) .
\end{equation*}%
Consequently,
\begin{equation*}
\max\limits_{x\in \left[ a_{1},b_{1}\right] }\left\vert G(x,y)\right\vert =%
\frac{1}{2}\left( g(c,y)-g(c,y_{0})\right) .\eqno(3.6)
\end{equation*}%
Taking (3.6) and the $4$-th property of $g$ into account in (3.3) yields
\begin{equation*}
\left\Vert G\right\Vert _{C(R^{\prime })}=\frac{1}{2}\left(
g(c,b_{2})-g(c,y_{0})\right) =\frac{1}{4}g(c,b_{2}).
\end{equation*}%
Similarly it can be shown that
\begin{equation*}
\left\Vert G\right\Vert _{C(R^{\prime \prime })}=\frac{1}{4}g(c,b_{2}).
\end{equation*}%
Hence
\begin{equation*}
\left\Vert G\right\Vert _{C(R)}=\frac{1}{4}g(c,b_{2})=L(f,R_{1}).
\end{equation*}%
But by the definition of $G$,
\begin{equation*}
G(x,y)=g(x,y)-\varphi _{1}(x)-\psi _{1}(y)=f(x,y)-\varphi _{0}(x)-\psi
_{0}(y),
\end{equation*}%
where
\begin{equation*}
\varphi _{0}(x)=\varphi
_{1}(x)+f(x,a_{2})-f(a_{1},a_{2})+f(a_{1},y_{0})=f(x,y_{0}),
\end{equation*}%
\begin{equation*}
\psi _{0}(y)=\psi _{1}(y)+f(a_{1},y)-f(a_{1},y_{0})=
\end{equation*}%
\begin{equation*}
=\frac{1}{2}\left( f(a_{1},y)+f(c,y)-f(a_{1},y_{0})-f(c,y_{0})\right) .
\end{equation*}%
Therefore,
\begin{equation*}
\left\Vert f-\varphi _{0}-\psi _{0}\right\Vert _{C(R)}=L(f,R_{1}).
\end{equation*}%
We proved (3.2) and hence Theorem 3.1. Note that the function $\varphi
_{0}(x)+\psi _{0}(y)$ is a best approximating sum from the manifold ${D}$ to
$f$.

\bigskip

\textbf{Remark 3.1.} In the special case $c=b_{1}$, Theorem 3.1 turns into
Babaev's result from \cite{6}.

\bigskip

\textbf{Corollary 3.1.} \textit{Let a function $f(x,y)$ have the continuous
derivative $\frac{\partial ^{2}f}{\partial x\partial y}$ on the rectangle $R$
and satisfy the following conditions }

\textit{1) $\frac{\partial ^{2}f}{\partial x\partial y}\geq 0$, for all $%
(x,y)\in R_{1}$; }

\textit{2) $\frac{\partial ^{2}f}{\partial x\partial y}\leq 0$, for all $%
(x,y)\in R_{2}$; }

\textit{3) $\frac{df(a_{1},y)}{dy}\leq \frac{df(b_{1},y)}{dy}$, for all $%
y\in \left[ a_{2},b_{2}\right] $. }

\textit{Then}
\begin{equation*}
E(f,R)=L(f,R_{1})=\frac{1}{4}\left[
f(a_{1},a_{2})+f(c,b_{2})-f(a_{1},b_{2})-f(c,a_{2})\right] .
\end{equation*}

The proof of this corollary can be obtained directly from Lemma 3.1 and
Theorem 3.1.

\bigskip

\textbf{Remark 3.2.} Rivlin and Sibner \cite{121} proved Corollary 3.1 in
the special case $c=b_{1}$.

\bigskip

\textbf{Example 3.2.} As we know (see Example 3.1) the function $f=y\sin \pi
x$ belongs to $V_{\frac{1}{2}}(K)$, where $K=\left[ 0,1\right] \times \left[
0,1\right] $. By Theorem 3.1, $E(f,K)=\frac{1}{4}$ and the function $\frac{1%
}{2}\sin \pi x+\frac{1}{2}y-\frac{1}{4}$ is a best approximating sum.

\bigskip

The following theorem shows that in some cases the approximation error
formula in Theorem 3.1 is valid for more general sets than rectangles with
sides parallel to the coordinate axes.

\bigskip

\textbf{Theorem 3.2.} \textit{Let $f(x,y)$ be a function from $V_{c}(R)$ and
$Q\subset R$ be a compact set which contains all vertices of $R_{1}$ (points
$(a_{1},a_{2}),(a_{1},b_{2}),(c,a_{2}),(c,b_{2})$). Then}
\begin{equation*}
E(f,Q)=L(f,R_{1})=\frac{1}{4}\left[
f(a_{1},a_{2})+f(c,b_{2})-f(a_{1},b_{2})-f(c,a_{2})\right] .
\end{equation*}

\textbf{Proof.} Since $Q\subset R,$ $E(f,Q)\leq E(f,R)$. On the other hand
by Theorem 3.1, $E(f,R)=L(f,R_{1})$. Hence $E(f,Q)\leq L(f,R_{1})$. It can
be shown, as it has been shown in the proof of Theorem 3.1, that $%
L(f,R_{1})\leq E(f,Q)$. But then automatically $E(f,Q)=L(f,R_{1})$.

\bigskip

\textbf{Example 3.3.} Calculate the approximation error of the function $%
f(x,y)=-(x-2)^{2n}y^{m}$ ($n$ and $m$ are positive integers) on the domain
\begin{equation*}
Q=\left\{ (x,y):0\leq x\leq 2,0\leq y\leq (x-1)^{2}+1\right\} .
\end{equation*}%
It can be easily verified that $f \in V_{2}(R)$, where $R=\left[ 0,4\right]
\times \left[ 0,2\right] $. Besides, $Q$ contains all vertices of $R_{1}=%
\left[ 0,2\right] \times \left[ 0,2\right] $. Consequently, by Theorem 3.2, $%
E(f,Q)=L(f,R_{1})=2^{2(n-1)+m}$.

\bigskip

\subsection{Characterization of $V_{c}(R)$}

The following theorem characterizes the class $V_{c}(R)$ in terms of the
approximation error calculation formulas.

\bigskip

\textbf{Theorem 3.3.} \textit{The following conditions are necessary and
sufficient for a continuous function $f(x,y)$ belong to $V_{c}(R):$ }

\textit{1) $E(f,S)=L(f,S)$, for each rectangle $S=\left[ x_{1},x_{2}\right]
\times \left[ y_{1},y_{2}\right] ,S\subset R_{1}$; }

\textit{2) $E(f,S)=-L(f,S)$, for each rectangle $S=\left[ x_{1},x_{2}\right]
\times \left[ y_{1},y_{2}\right] ,S\subset R_{2}$; }

\textit{3) $E(f,S)=L(f,S_{1})$, for each rectangle $S=\left[ a_{1},b_{1}%
\right] \times \left[ y_{1},y_{2}\right] ,S\subset R$ and $S_{1}=\left[
a_{1},c\right] \times \left[ y_{1},y_{2}\right] $.}

\bigskip

\textbf{Proof.} The necessity easily follows from the definition of $V_{c}(R)
$, Babaev's above-mentioned result (see Section 3.1.1) and Theorem 3.1. The
sufficiency is clear if pay attention to the fact that $E(f,S)\geq 0 $.

\bigskip

\subsection{Classes $V_{c}^{-}(R),U(R)$ and $U_{c}^{-}(R)$}

By $V_{c}^{-}(R)$ we denote the class of functions $f(x,y)$ such that $-f
\in V_{c}(R)$. It is clear that $E(f,R)=-L(f,R_{1})$ for each $f\in
V_{c}^{-}(R)$.

We define $U_{c}(R),a_{1}\leq c< b_{1}$, as a class of continuous functions $%
f(x,y)$ with the properties

1) $L(f,S)\leq 0$, for each rectangle $S=\left[ x_{1},x_{2}\right] \times %
\left[ y_{1},y_{2}\right] ,\;S\subset R_{1};$

2) $L(f,S)\geq 0$, for each rectangle $S=\left[ x_{1},x_{2}\right] \times %
\left[ y_{1},y_{2}\right] ,\;S\subset R_{2};$

3) $L(f,S)\geq 0$, for each rectangle $S=\left[ a_{1},b_{1}\right] \times %
\left[ y_{1},y_{2}\right] ,\;S\subset R.$

Using the same techniques in the proof of Theorem 3.1 it can be shown that
the following theorem is valid:

\bigskip

\textbf{Theorem 3.4.} \textit{The approximation error of a function $f(x,y)$
from the class $U_{c}(R)$ can be calculated by the formula
\begin{equation*}
E(f,R)=L(f,R_{2})=\frac{1}{4}\left[
f(c,a_{2})+f(b_{1},b_{2})-f(c,b_{2})-f(b_{1},a_{2})\right] .
\end{equation*}%
Let $y_{0}$ be any solution from $\left[ a_{2},b_{2}\right] $ of the
equation
\begin{equation*}
L(f,Y)=\frac{1}{2}L(f,R_{2}),\qquad Y=\left[ c,b_{1}\right] \times \left[
a_{2},y\right] .
\end{equation*}%
Then the function $\varphi _{0}(x)+\psi _{0}(y)$, where
\begin{equation*}
\varphi _{0}(x)=f(x,y_{0}),\quad \psi _{0}(y)=\frac{1}{2}\left[
f(c,y)+f(b_{1},y)-f(c,y_{0})-f(b_{1},y_{0})\right] ,
\end{equation*}%
is a best approximating sum from the manifold $D$ to $f$. }

\bigskip

By $U_{c}^{-}(R)$ denote the class of functions $f(x,y)$ such that $-f \in
U_{c}(R)$. It is clear that $E(f,R)=-L(f,R_{2})$ for each $f\in U_{c}^{-}(R)$%
.

\bigskip

\textbf{Remark 3.3.} The correspondingly modified versions of Theorems 2.2,
2.3 and Corollary 3.1 are valid for the classes $V_{c}^{-}(R),U_{c}(R)$ and $%
U_{c}^{-}(R)$.

\bigskip

\textbf{Example 3.4.} Consider the function $f(x,y)=\left( x-\frac{1}{2}%
\right) ^{2}y$ on the unit square $K=\left[ 0,1\right] \times \left[ 0,1%
\right] $. It can be easily verified that $f\in U_{\frac{1}{2}}(K)$. Hence,
by Theorem 3.4, $E(f,K)=\frac{1}{16}$ and the function $\frac{1}{2}\left( x-%
\frac{1}{2}\right) ^{2}+\frac{1}{8}y-\frac{1}{16}$ is a best approximating
function.

\bigskip

\section{Approximation by sums of univariate functions on certain domains}

The purpose of this section is to develop a method for obtaining explicit
formulas for the error of approximation of bivariate functions by sums of
univariate functions. It should be remarked that formulas of this type were
known only for functions defined on a rectangle with sides parallel to the
coordinate axes. Our method, based on a maximization process over closed
bolts, allows the consideration of functions defined on hexagons, octagons
and stairlike polygons with sides parallel to the coordinate axes.

\subsection{Problem statement}

Let $Q$ be a compact set in $\mathbb{R}^2$. Consider the approximation of a
continuous function $f \in C(Q)$ by functions from the set $D=\left\{
\varphi (x)+\psi (y)\right\} ,$ where $\varphi (x),\psi (y)$ are defined and
continuous on the projections of $Q$ into the coordinate axes $x$ and $y$,
respectively. The approximation error is defined as follows
\begin{equation*}
E(f,Q)=\inf\limits_{\varphi+\psi \in D}\left\Vert f-\varphi -\psi
\right\Vert _{C(Q)}.
\end{equation*}%
Our purpose is to develop a method for obtaining explicit formulas providing
precise and easy computation of $E(f,Q)$ for polygons $Q$ with sides
parallel to the coordinate axes. This method will be based on the herein
developed \textit{closed bolts maximization process} and can be used in
alternative proofs of the known results from \cite{6}, \cite{59} and \cite%
{121}. First, we show efficiency of the method in the example of a hexagon
with sides parallel to the coordinate axes. Then we formulate an analogous
theorem for staircase polygons and two theorems for octagons, which can be
proved in a similar way, and touch some aspects of the question about the
case of an arbitrary polygon with sides parallel to the coordinate axes. The
condition posed on sides of polygons (being parallel to the coordinate axes)
is essential for our method. This has several reasons, which get clear
through the proof of Theorem 3.5. Here we are able to explain one of these
reasons: by \cite[Theorem 3]{34}, a continuous function $f(x,y)$ defined
on a polygon with sides parallel to the coordinate axes has an extremal
element, the existence of which is required in our method. Now let $K$ be a
rectangle (not speaking about polygons) with sides not parallel to the
coordinate axes. Does any function $f\in C(K)$ have an extremal element? No one knows (see \cite{34}).

In the sequel, all the considered polygons are supposed to have sides
parallel to the coordinate axes.

\bigskip

\subsection{The maximization process}

Let $H$ be a closed hexagon. It is clear that $H$ can be uniquely
represented in the form
\begin{equation*}
H=R_{1}\cup R_{2},\eqno(3.7)
\end{equation*}%
where $R_{1},R_{2}$ are rectangles and there does not exist any rectangle $R$
such that $R_{1}\subset R\subset H$ or $R_{2}\subset R\subset H$.

We associate each closed bolt $p=\left\{ p_{1},p_{2},\cdots p_{2n}\right\} $
with the following functional
\begin{equation*}
l(f,p)=\frac{1}{2n}\sum\limits_{k=1}^{2n}(-1)^{k-1}f(p_{k}).
\end{equation*}

Denote by $M(H)$ the class of bivariate continuous functions $f$ on $H$
satisfying the condition
\begin{equation*}
f(x_{1},y_{1})+f(x_{2},y_{2})-f(x_{1},y_{2})-f(x_{2},y_{1})\geq 0
\end{equation*}
for any rectangle $\left[ x_{1},x_{2}\right] \times \left[ y_{1},y_{2}\right]
\subset H.$\newline

\textbf{Theorem 3.5.}\ \textit{Let $H$ be a hexagon and (3.7) be its
representation. Let $f \in M(H)$. Then
\begin{equation*}
E(f,H)=\max \left\{ \left\vert l(f,h)\right\vert ,\left\vert
l(f,r_{1})\right\vert ,\left\vert l(f,r_{2})\right\vert \right\} ,\eqno(3.8)
\end{equation*}%
where $h,r_{1},r_{2}$ are closed bolts formed by vertices of the polygons $%
H,R_{1},R_{2}$ respectively. }

\bigskip

\begin{proof} Without loss of generality, we may assume that the rectangles $%
R_{1}$ and $R_{2}$ are of the following form
\begin{equation*}
R_{1}=\left[ a_{1},a_{2}\right] \times \left[ b_{1},b_{3}\right] ,~\ R_{2}=%
\left[ a_{1},a_{3}\right] \times \left[ b_{1},b_{2}\right] ,~\
a_{1}<a_{2}<a_{3},\;b_{1}<b_{2}<b_{3}.
\end{equation*}%
Introduce the notation
\begin{equation*}
\begin{array}{c}
f_{11}=f\left( a_{1},b_{1}\right) ,~\ f_{12}=-f\left( a_{1},b_{2}\right)
,\;f_{13}=-f\left( a_{1},b_{3}\right) ; \\
f_{21}=-f\left( a_{2},b_{1}\right) ,\;f_{22}=-f\left( a_{2},b_{2}\right) ,~\
f_{23}=f\left( a_{2},b_{3}\right) ; \\
f_{31}=-f\left( a_{3},b_{1}\right) ,\;f_{32}=f\left( a_{3},b_{2}\right) .%
\end{array}%
\eqno(3.9)
\end{equation*}%
It is clear that
\begin{equation*}
\begin{array}{c}
\left\vert l(f,r_{1})\right\vert =\dfrac{1}{4}\left(
f_{11}+f_{13}+f_{23}+f_{21}\right) , \\
\left\vert l(f,r_{2})\right\vert =\dfrac{1}{4}\left(
f_{11}+f_{12}+f_{32}+f_{31}\right) , \\
\left\vert l(f,h)\right\vert =\dfrac{1}{6}\left(
f_{11}+f_{13}+f_{23}+f_{22}+f_{32}+f_{31}\right) .%
\end{array}%
\eqno(3.10)
\end{equation*}

Let $p=\left\{ p_{1},p_{2},\cdots p_{2n}\right\} $ be any closed bolt. We
group the points $p_{1},p_{2},\cdots p_{2n}$ by putting
\begin{equation*}
p_{+}=\left\{ p_{1},p_{3},\cdots p_{2n-1}\right\} ,\;p_{-}=\left\{
p_{2},p_{4},\cdots p_{2n}\right\}.
\end{equation*}

First, assume that $l(f,p)\geq 0$. We apply the following algorithm, which we call \textit{the maximization process over closed bolts}, to $p$.

\textbf{Step 1.} Consider sequentially the units $p_ip_{i+1}$ $\left(i=\overline{%
1,2n}, p_{2n+1}=p_1\right)$ with the vertices $p_{i}\left(
x_{i},y_{i}\right) ,~\ p_{i+1}\left( x_{i+1},y_{i+1}\right) $ having equal
abscissae: $x_{i}=x_{i+1}$. Four cases are possible.

1) $p_{i}\in p_{+}$ and $y_{i+1}>y_{i}$. In this case, replace the unit $%
p_{i}p_{i+1}$ by a new unit $q_{i}q_{i+1}$ with the vertices $%
q_{i}=(a_{1},y_{i}),\;\ q_{i+1}=(a_{1},y_{i+1})$.

2) $p_{i}\in p_{+}$ and $y_{i+1}<y_{i}$. In this case, replace the unit $%
p_{i}p_{i+1}$ by a new unit $q_{i}q_{i+1}$ with the vertices $q_{i}=\left(
a_{2},y_{i}\right) ,\;q_{i+1}=(a_{2},y_{i+1})$ if $b_{2}<y_{i}\leq b_{3}$ or
with the vertices $q_{i}=\left( a_{3},y_{i}\right),
q_{i+1}=(a_{3},y_{i+1})$ if $b_{1}\leq y_{i}\leq b_{2}$.

3) $p_{i}\in p_{-}$ and $y_{i+1}<y_{i}$. In this case, replace $p_{i}p_{i+1}$
by a new unit $q_{i}q_{i+1}$ with the vertices $%
q_{i}=(a_{1},y_{i}),~q_{i+1}=(a_{1},y_{i+1})$.

4) $p_{i}\in p_{-}$ and $y_{i+1}>y_{i}$. In this case, replace $p_{i}p_{i+1}$
by a new unit $q_{i}q_{i+1}$ with the vertices $%
q_{i}=(a_{2},y_{i}),~q_{i+1}=(a_{2},y_{i+1})$\ if\ $b_{2}<y_{i+1}\leq b_{3}$
or with the vertices $q_{i}=(a_{3},y_{i}),~q_{i+1}=(a_{3},y_{i+1})$ if $%
b_{1}\leq y_{i+1}\leq b_{2}$.

Since $f\in M(H)$, it is not difficult to verify that
\begin{equation*}
\begin{array}{c}
f(p_{i})-f(p_{i+1})\leq f(q_{i})-f(q_{i+1})\ \ \mbox{for cases 1)
and 2)}, \\
-f(p_{i})+f(p_{i+1})\leq -f(q_{i})+f(q_{i+1})\ \ \mbox{for cases 3) and 4)}%
\end{array}%
\eqno(3.11)
\end{equation*}

It is clear that after Step 1 the bolt $p$ will be replaced by the ordered set $%
q=\left\{ q_{1},q_{2},\cdots ,q_{2n}\right\} $. We do not say a bolt but an
ordered set because of a possibility of coincidence of some successive
points $q_{i}, q_{i+1}$ (this, for example, may happen if the 1-st case
takes place for the units $p_{i-1}p_{i}$ and $p_{i+1}p_{i+2}$). Let us exclude
simultaneously successive and coincident points from $q$. Then we obtain
some closed bolt, which we denote by $q^{\prime }=\left\{ q_{1}^{\prime
},q_{2}^{\prime },\cdots ,q_{2m}^{\prime }\right\}$. It is not difficult to
understand that all points of the bolt $q^{\prime}$ are located on straight
lines $x=a_{1},~x=a_{2},~x=a_{3}$.

From inequalities (3.11) and the fact that $2m\leq 2n,$ we deduce that
\begin{equation*}
l(f,p)\leq l(f,q^{\prime }).\eqno(3.12)
\end{equation*}

\textbf{Step 2.} Consider sequentially units $q_{i}^{\prime }q_{i+1}^{\prime
}\;\left( i=\overline{1,2m},q_{2m+1}^{\prime }=q_{1}^{\prime }\right) $ with
the vertices $q_{i}^{\prime }=\left( x_{i}^{\prime },y_{i}^{\prime }\right)
,~\ q_{i+1}^{\prime }\left( x_{i+1}^{\prime },y_{i+1}^{\prime }\right) $
having equal ordinates: $y_{i}^{\prime }=y_{i+1}^{\prime }$. The following
four cases are possible.

1) $q_{i}^{\prime }\in q_{+}^{\prime }$ and $x_{i+1}^{\prime }>x_{i}^{\prime
}$. In this case, replace the unit $q_{i}^{\prime }q_{i+1}^{\prime }$ by a
new unit $p_{i}^{\prime }p_{i+1}^{\prime }$ with the vertices $p_{i}^{\prime
}=\left( x_{i}^{\prime },b_{1}\right) ,~\ p_{i+1}^{\prime }=\left(
x_{i+1}^{\prime },b_{1}\right) $.

2) $q_{i}^{\prime }\in q_{+}^{\prime }$ and $x_{i+1}^{\prime }<x_{i}^{\prime
}$. In this case, replace the unit $q_{i}^{\prime }q_{i+1}^{\prime }$ by a
new unit $p_{i}^{\prime }p_{i+1}^{\prime }$ with the vertices $p_{i}^{\prime
}=\left( x_{i}^{\prime },b_{2}\right) ,$ $\ p_{i+1}^{\prime }=\left(
x_{i+1}^{\prime },b_{2}\right) $ if $x_{i}^{\prime }=a_{3}$ and with the
vertices $p_{i}^{\prime }=\left( x_{i}^{\prime },b_{3}\right) ,$ $%
p_{i+1}^{\prime }=\left( x_{i+1}^{\prime },b_{3}\right) $ if $x_{i}^{\prime
}=a_{2}$.

3) $q_{i}^{\prime }\in q_{-}^{\prime }$ and $x_{i+1}^{\prime }<x_{i}^{\prime
}$. In this case, replace $q_{i}^{\prime }q_{i+1}^{\prime }$ by a new unit $%
p_{i}^{\prime }p_{i+1}^{\prime }$ with the vertices $p_{i}^{\prime }=\left(
x_{i}^{\prime },b_{1}\right) ,~\ p_{i+1}^{\prime }=\left( x_{i+1}^{\prime
},b_{1}\right)$.

4) $q_{i}^{\prime }\in q_{-}^{\prime }$ and $x_{i+1}^{\prime }>x_{i}^{\prime
}$. In this case, replace $q_{i}^{\prime }q_{i+1}^{\prime }$ by a new unit $%
p_{i}^{\prime }p_{i+1}^{\prime }$ with the vertices $p_{i}^{\prime }=\left(
x_{i}^{\prime },b_{2}\right) ,~\ p_{i+1}^{\prime }=\left( x_{i+1}^{\prime
},b_{2}\right)$ if $x_{i+1}^{\prime }=a_{3}$ and with the vertices $%
p_{i}^{\prime }=\left( x_{i}^{\prime },b_{3}\right) ,~\ p_{i+1}^{\prime
}=\left( x_{i+1}^{\prime },b_{3}\right) $ if $x_{i+1}^{\prime }=a_{2}$.

It is easy to see that after Step 2 the bolt $q^{\prime }$ will be replaced by
the bolt $p^{\prime }=\left\{ p_{1}^{\prime },p_{2}^{\prime },\cdots
p_{2m}^{\prime }\right\} $ and
\begin{equation*}
l(f,q^{\prime })\leq l(f,p^{\prime }).\eqno(3.13)
\end{equation*}

From (3.12) and (3.13) we obtain that
\begin{equation*}
l(f,p)\leq l(f,p^{\prime }).\eqno(3.14)
\end{equation*}

It is clear that each point of the set $p_{+}^{\prime }$ coincides with one
of the points $\left( a_{1},b_{1}\right) ,~\left( a_{2},b_{3}\right) ,$ $%
\left( a_{3},b_{2}\right) $ and each point of the set $p_{-}^{\prime }$
coincides with one of the points $\left( a_{1},b_{2}\right) ,~\left(
a_{1},b_{3}\right) ,$ $~\left( a_{2},b_{1}\right) ,~\left(
a_{2},b_{2}\right) ,~\left( a_{3},b_{1}\right) .$ Denote by $m_{ij}$ the
number of points of the bolt $p^{\prime }$ coinciding with the point $\left(
a_{i},b_{j}\right) ,~i,j=\overline{1,3},~i+j\neq 6$. By (3.9), we can write that
\begin{equation*}
l(f,p^{\prime })=\frac{1}{2m}\sum\limits_{\substack{ i,j=\overline{1,3}  \\ %
i+j\leq 5}}m_{ij}f_{ij}.\eqno(3.15)
\end{equation*}

On the straight line $x=a_{i}~\ $or $\ y=b_{i},~i=\overline{1,3}$, the
number of points of the set $p_{+}^{\prime }$ is equal to the number of
points of the set $p_{-}^{\prime }$. Hence
\begin{equation*}
m_{11}=m_{12}+m_{13}=m_{21}+m_{31};\ m_{23}=m_{22}+m_{21}=m_{13};\
m_{32}=m_{31}=m_{12}+m_{22}.
\end{equation*}%
From these equalities we deduce that
\begin{equation*}
m_{11}=m_{12}+m_{21}+m_{22};\ m_{13}=m_{21}+m_{22};\ m_{23}=m_{21}+m_{22};\
m_{31}=m_{12}+m_{22}.\eqno(3.16)
\end{equation*}%
Consequently,
\begin{equation*}
2m=\sum\limits_{\substack{ i,j=\overline{1,3}  \\ i+j\leq 5}}%
m_{ij}=4m_{12}+4m_{21}+6m_{22}.\eqno(3.17)
\end{equation*}%
Considering (3.16) and (3.17) in (3.15) and taking (3.10) into account, we
obtain that
\begin{equation*}
l(f,p^{\prime })=\dfrac{4m_{12}\left\vert l(f,r_{2})\right\vert
+4m_{21}\left\vert l(f,r_{1})\right\vert +6m_{22}\left\vert
l(f,h)\right\vert }{4m_{12}+4m_{21}+6m_{22}}
\end{equation*}%
\begin{equation*}
\leq \max \left\{ \left\vert l(f,r_{1})\right\vert ,\left\vert
l(f,r_{2})\right\vert ,\left\vert l(f,h)\right\vert \right\} .
\end{equation*}%
Therefore, due to (3.14),
\begin{equation*}
l(f,p)\leq \max \left\{ \left\vert l(f,r_{1})\right\vert ,\left\vert
l(f,r_{2})\right\vert ,\left\vert l(f,h)\right\vert \right\} .\eqno(3.18)
\end{equation*}

Note that in the beginning of the proof the bolt $p$ has been chosen so that
$l(f,p)\geq 0$. Let now $p=\left\{ p_{1},p_{2},\cdots p_{2n}\right\} $ be
any closed bolt such that $l(f,p)\leq 0$. Since $l(f,p^{\prime
\prime })=$ $-l(f,p)\geq 0$ for the bolt $p^{\prime \prime }=\left\{
p_{2},p_{3},\cdots ,p_{2n},p_{1}\right\} $,we obtain from (3.18) that
\begin{equation*}
-l(f,p)\leq \max \left\{ \left\vert l(f,r_{1})\right\vert ,\left\vert
l(f,r_{2})\right\vert ,\left\vert l(f,h)\right\vert \right\} .\eqno(3.19)
\end{equation*}%
From (3.18) and (3.19) we deduce on the strength of arbitrariness of $p$
that
\begin{equation*}
\sup\limits_{p\subset H}\left\{ \left\vert l(f,p)\right\vert \right\} =\max
\left\{ \left\vert l(f,r_{1})\right\vert ,\left\vert l(f,r_{2})\right\vert
,\left\vert l(f,h)\right\vert \right\} ,\eqno(3.20)
\end{equation*}%
where the $sup$ is taken over all closed bolts of the hexagon $H$.

The hexagon $H$ satisfies the conditions of Theorem 1.10 on the
existence of a best approximation. By \cite[Theorem 2]{79} (see Section 3.3), we obtain that
\begin{equation*}
E(f,H)=\sup\limits_{p\subset H}\left\{ \left\vert l(f,p)\right\vert \right\}
.\eqno(3.21)
\end{equation*}%
From (3.20) and (3.21) we finally conclude that
\begin{equation*}
E(f,H)=\max \left\{ \left\vert l(f,r_{1})\right\vert ,\left\vert
l(f,r_{2})\right\vert ,\left\vert l(f,h)\right\vert \right\} .
\end{equation*}%

\end{proof}

\textbf{Corollary 3.2.}\ \textit{Let a function $f(x,y)$ have the continuous
nonnegative derivative $\dfrac{\partial ^{2}f}{\partial x\partial y}$ on $H$%
. Then the formula (3.8) is valid.}

\bigskip

The proof is very simple and can be obtained by integrating the inequality
\linebreak $\dfrac{\partial ^{2}f}{\partial x\partial y}\!\geq \!0$ over an
arbitrary rectangle $\left[ x_{1},x_{2}\right] \times \left[ y_{1},y_{2}%
\right] \subset H$ and applying Theorem 3.5.

The method used in the proof of Theorem 3.5 can be generalized to obtain
similar results for stairlike polygons. For example, let $S$ be a closed
polygon of the following form

\begin{equation*}
S=\bigcup\limits_{i=1}^{N-1}P_{i},
\end{equation*}%
where $N\geq 2,$ $P_{i}=\left[ a_{i},a_{i+1}\right] \times \left[
b_{1},b_{N+1-i}\right] ,$ $i=\overline{1,N-1},$ $a_{1}<a_{2}<\dots <a_{N},$ $%
b_{1}<b_{2}<\dots <b_{N}$. Such polygons will be called \textit{stairlike
polygons} (see \cite{56}).

A closed $2m$-gon $F$ with sides parallel to the coordinate axes is called a
maximal $2m$-gon of the polygon $S$ if $F\subset S$ and there is no another $%
2m$-gon $F^{\prime }$ such that $F\subset F^{\prime }\subset S$. Clearly, if
$F$ is a maximal $2m$-gon of the polygon $S$, then $m\leq N.$ A closed bolt
formed by the vertices of a maximal polygon $F$ is called a maximal bolt of $%
S$. By $S^{B}$ denote the set of all maximal bolts of the stairlike polygon $%
S.$

\bigskip

\textbf{Theorem 3.6.} \textit{Let $S$ be a stairlike polygon. The
approximation error of a function $f \in M(S)$ can be computed by the formula%
}
\begin{equation*}
E\left( f,S\right) =\max \left\{ \left\vert r(f,h)\right\vert ,\;h\in
S^{B}\right\} .
\end{equation*}

\bigskip

For the proof of this theorem see \cite{56}.

\bigskip

\subsection{$E$-bolts}

The main idea in the proof of Theorem 3.5 can be successfully used in
obtaining formulas of type (3.8) for functions $f(x,y)$ defined on another
simple polygons. The following two theorems include cases of some octagons
and can be proved in a similar way.

\bigskip

\textbf{Theorem 3.7.}\ \textit{Let $a_{1}<a_{2}<a_{3}<a_{4},$ $%
b_{1}<b_{2}<b_{3}$ and $Q$ be an octagon of the following form
\begin{equation*}
Q=\bigcup\limits_{i=1}^{4}R_{i},\ \ \ where
\end{equation*}%
$R_{1}=\left[ a_{1},a_{2}\right] \times \left[ b_{1},b_{2}\right] ,R_{2}=%
\left[ a_{2},a_{3}\right] \times \left[ b_{1},b_{2}\right] ,R_{3}=\left[
a_{3},a_{4}\right] \times \left[ b_{1},b_{2}\right] ,R_{4}=\left[ a_{2},a_{3}%
\right] \times \left[ b_{2},b_{3}\right] $.} \textit{Let $f\in M(Q)$. Then
the following formula holds
\begin{equation*}
E(f,Q)=\max \left\{ \left| l(f,q)\right| ,\left| l(f,r_{123} )\right|
,\left| l(f,r_{124} )\right| ,\left| l(f,r_{234} )\right| ,\left| l(f,r_{24}
)\right| \right\},
\end{equation*}
where $q,$ $r_{123},$ $r_{124},$ $r_{234},$ $r_{24} $ are closed bolts
formed by the vertices of the polygons $Q,$ $R_{1} \cup R_{2}\cup R_{3}
,R_{1}\cup R_{2}\cup R_{4} ,R_{2}\cup R_{3}\cup R_{4} $ and $R_{2}\cup R_{4}
$, respectively.}

\bigskip

\textbf{Theorem 3.8.}\ \textit{Let $a_{1}<a_{2}<a_{3}<a_{4},\
b_{1}<b_{2}<b_{3}$ and $Q$ be an octagon of the following form
\begin{equation*}
Q=\bigcup_{i=1}^{3}R_{i},
\end{equation*}%
where $R_{1}=\left[ a_{1},a_{4}\right] \times \left[ b_{1},b_{2}\right]
,R_{2}=\left[ a_{1},a_{2}\right] \times \left[ b_{2},b_{3}\right] ,R_{3}=%
\left[ a_{3},a_{4}\right] \times \left[ b_{2},b_{3}\right] $.} \textit{Let $%
f\in M(Q)$. Then
\begin{equation*}
E(f,Q)=\max \left\{ \left| l(f,r)\right| ,\left| l(f,r_{12} )\right| ,\left|
l(f,r_{13} )\right| \right\},
\end{equation*}
where $r,r_{12} ,r_{13} $ are closed bolts formed by the vertices of the
polygons $R=\left[ a_{1} ,a_{4} \right] \times \left[ b_{1} ,b_{3} \right],$
$R_{1}\cup R_{2} ,R_{1}\cup R_{3}$, respectively.}

\bigskip

Although the closed bolts maximization process can be applied to bolts of an
arbitrary polygon, some combinatorial difficulties arise when grouping
values at points of maximized bolts (bolts obtained after the maximization
process, see (3.15)-(3.18)). While we do not know a complete answer to this
problem, we can describe points of a polygon $F$ with which points of
maximized bolts coincide and state a conjecture concerning the approximation
error.

Let $F=A_{1}A_{2}...A_{2n}$ be any polygon with sides parallel to the
coordinate axes. The vertices $A_{1},$ $A_{2},$ $...,$ $A_{2n}$ in the given
order form a closed bolt, which we denote by $r_{F}$. By $\left[ r_{F}\right]
$ denote the length of $r_{F}$. In our case, $\left[ r_{F}\right] =2n$.

\bigskip

\textbf{Definition 3.2.}\ \textit{Let $F$ and $S$ be polygons with sides
parallel to the coordinate axes. We say that the closed bolt $r_{F}$ is an $%
e $-bolt (extended bolt) of $S$ if\ $r_{F}\subset S$\ and there does not
exist any polygon $F^{^{\prime }}$ such that $F\subset F^{^{\prime }},\ \
r_{F^{^{\prime }}}\subset S,\ \ \left[ r_{F^{^{\prime }}}\right] \leq \left[
r_{F}\right] .$}

\bigskip

For example, in Theorem 3.8 the octagon $Q$ has $3$ $e$-bolts. They are $%
r,r_{12}$ and $r_{13}$. In Theorem 3.7, the octagon $Q$ has $5$ $e$-bolts,
which are $q,r_{123},r_{124},r_{234}$ and $r_{24}$ . The polygon $%
S_{2n}=\bigcup\limits_{i=1}^{n-1}R_{i}$, where $R_{i}=\left[ a_{i},a_{i+1}%
\right] \times \left[ b_{1},b_{n+1-i}\right] ,i=\overline{1,n-1}%
,a_{1}<a_{2}<...<a_{n},b_{1}<b_{2}<...<b_{n}$ has exactly $2^{n-1}-1$ $e$%
-bolts. It is not difficult to observe that the set of points of a closed
bolt obtained after the maximization process is a subset of the set of
points of all $e$-bolts. This condition and Theorems 2.5-2.8 justify the
statement of the following conjecture:

\textit{Let $S$ be any polygon with sides parallel to the coordinate axes
and $f \in M(S)$. Then
\begin{equation*}
E(f,S)=\max_{h\in S^{E}}\left\{ \left\vert l(f,h)\right\vert \right\},
\end{equation*}%
where $S^{E}$ is a set of all $e$-bolts of the polygon $S$.}

\bigskip

\subsection{Error estimates}

Theorem 3.5 allows us to consider classes wider than $M(H)$ and establish
sharp estimates for the approximation error.

\bigskip

\textbf{Theorem 3.9.}\ \textit{Let $H$ be a hexagon and (3.7) be its
representation. The following sharp estimates are valid for a function $%
f(x,y)$ having the continuous derivative $\dfrac{\partial ^{2}f}{\partial
x\partial y}$ on $H$:
\begin{equation*}
A\leq E(f,H)\leq BC+\frac{3}{2}\left( B\left\vert l(g,h)\right\vert
-\left\vert l(f,h)\right\vert \right) ,\eqno(3.22)
\end{equation*}%
where
\begin{equation*}
B=\max_{(x,y)\in H}\left\vert \frac{\partial ^{2}f(x,y)}{\partial x\partial y%
}\right\vert ,\ \ \ g=g(x,y)=x\cdot y,
\end{equation*}%
\begin{equation*}
A=\max \left\{ \left\vert l(f,h)\right\vert ,\ \left\vert
l(f,r_{1})\right\vert ,\left\vert l(f,r_{2})\right\vert \right\} ,\ C=\max
\left\{ \left\vert l(g,h)\right\vert ,\left\vert l(g,r_{1})\right\vert ,\
\left\vert l(g,r_{2})\right\vert \right\},
\end{equation*}%
where $h,r_{1},r_{2}$ are closed bolts formed by vertices of the polygons $%
H,R_{1}$ and $R_{2}$, respectively.}

\bigskip

\textbf{Remark 3.4.} Inequalities similar to (3.22) were established in
Babaev \cite{8} for the approximation of a function $f(x)=f(x_{1},...,x_{n})$%
, defined on a parallelepiped with sides parallel to the coordinate axes, by
sums $\sum\limits_{i=1}^{n}\varphi _{i}(x\backslash x_{i})$. For the
approximation of bivariate functions, Babaev's result contains only
rectangular case.

\bigskip

\textbf{Remark 3.5.} Estimates (3.22) are easily calculable in contrast to
those established in \cite{5} for continuous functions defined on certain
domains, which are different from polygons.

\bigskip

To prove Theorem 3.9 we need the following lemmas.

\bigskip

\textbf{Lemma 3.3.}\ \textit{Let $X$ be a normed space, $F$ be a subspace of
$X$. The following inequality is valid for an element $x=x_{1}+x_{2}$ from $X
$:
\begin{equation*}
\left\vert E(x_{1})-E(x_{2})\right\vert \leq E(x)\leq E(x_{1})+E(x_{2}),
\end{equation*}%
where
\begin{equation*}
E(x)=E(x,F)=\inf_{y\in F}\left\Vert x-y\right\Vert .
\end{equation*}%
} \bigskip

\textbf{Lemma 3.4.}\ \textit{If $f\in M(H)$, then
\begin{equation*}
\left\vert l(f,r_{i})\right\vert \leq \frac{3}{2}\left\vert
l(f,h)\right\vert ,i=1,2.
\end{equation*}%
}

Lemma 3.3 is obvious. To prove Lemma 3.4, note that for any $f\in M(H)$
\begin{equation*}
6\left\vert l(f,h)\right\vert =4\left\vert l(f,r_{i})\right\vert
+4\left\vert l(f,r_{3})\right\vert ,\ \ i=1,2,
\end{equation*}%
where $r_{3}$ is a closed bolt formed by the vertices of the rectangle $%
R_{3}=H\backslash R_{i}.$

Now let us prove Theorem 3.9.

\begin{proof} It is not difficult to verify that if $\frac{\partial ^{2}u}{%
\partial x\partial y}\geq 0$ on $H$ for some $u(x,y),$ $\frac{\partial
^{2}u(x,y)}{\partial x\partial y}\in C(H)$, then $u\in M(H)$ (see the proof
of Corollary 3.2). Set $f_{1}=f+Bg$. Since $\frac{\partial ^{2}f_{1}}{%
\partial x\partial y}\geq 0$ on $H$, $f_{1}\in M(H)$. By Lemma 3.4,
\begin{equation*}
\left\vert l(f_{1},r_{i})\right\vert \leq \frac{3}{2}\left\vert
l(f_{1},h)\right\vert ,i=1,2.\eqno(3.23)
\end{equation*}

Theorem 3.5 implies that
\begin{equation*}
E(f_{1},H)=\max \left\{ \left\vert l(f_{1},h)\right\vert ,\left\vert
l(f_{1},r_{1})\right\vert ,\left\vert l\left( f_{1},r_{2}\right) \right\vert
\right\} .\eqno(3.24)
\end{equation*}
We deduce from (3.23) and (3.24) that
\begin{equation*}
E(f_{1},H)\leq \frac{3}{2}\left\vert l(f_{1},h)\right\vert .
\end{equation*}

First, let the closed bolt $h$ start at the point $(a_{1},b_{1})$. Then it
is clear that
\begin{equation*}
E(f_{1},H)\leq \frac{3}{2}l(f_{1},h).\eqno(3.25)
\end{equation*}
By Lemma 3.3,
\begin{equation*}
E(f,H)-E(Bg,H)\leq E(f_{1},H).\eqno(3.26)
\end{equation*}
Inequalities (3.25) and (3.26) yield
\begin{equation*}
E(f,H)\leq BE(g,H)+\frac{3}{2}l(f_{1},h).\eqno(3.27)
\end{equation*}

Since the functional $l(f,h)$ is linear,
\begin{equation*}
l(f_{1},h)=l(f,h)+Bl(g,h).
\end{equation*}%
Considering this expression of $l(f_{1},h)$ in (3.27), we obtain that
\begin{equation*}
E(f,H)\leq BE(g,H)+\frac{3}{2}Bl(g,h)+\frac{3}{2}l(f,h).\eqno(3.28)
\end{equation*}

Now consider the function $f_{2}=Bg-f$. Obviously, $\frac{\partial ^{2}f_{2}}{\partial x\partial y}%
\geq 0$ on $H$. It can be shown, in the same way as (3.28) has been
obtained, that
\begin{equation*}
E(f,H)\leq BE(g,H)+\frac{3}{2}Bl(g,h)-\frac{3}{2}l(f,h).\eqno(3.29)
\end{equation*}
From (3.28) and (3.29) it follows that
\begin{equation*}
E(f,H)\leq BE(g,H)+\frac{3}{2}Bl(g,h)-\frac{3}{2}\left\vert
l(f,h)\right\vert .\eqno(3.30)
\end{equation*}%
Since $g\in M(H)$ and $h$ starts at the point $(a_{1},b_{1}),$ we have $%
l(g,h)\geq 0$.

Let now $h$ start at a point such that $l(u,h)\leq 0$ for any $u\in M(H)$.
Then in a similar way as above we can prove that
\begin{equation*}
E(f,H)\leq BE(g,H)-\frac{3}{2}Bl(g,h)-\frac{3}{2}\left\vert
l(f,h)\right\vert ,\eqno(3.31)
\end{equation*}%
where $l(g,h)\leq 0$. From (3.30), (3.31) and the fact that $E(g,H)=C$ (in
view of Theorem 3.5), it follows that
\begin{equation*}
E(f,H)\leq BC+\frac{3}{2}\left( B\left\vert l(g,h)\right\vert -\left\vert
l(f,h)\right\vert \right) .
\end{equation*}%
The upper bound in (3.22) has been established. Note that it is attained
by $f=g=xy$.

The proof of the lower bound in (3.22) is simple. One of the obvious properties
of the functional $l(f,p)$ is that $\left\vert l(f,p)\right\vert \leq E(f,H)$
for any continuous function $f$ on $H$ and a closed bolt $p$. Hence,
\begin{equation*}
A=\max \left\{ \left\vert l(f,h)\right\vert ,\left\vert
l(f,r_{1})\right\vert ,\left\vert l(f,r_{2})\right\vert \right\} \leq E(f,H).
\end{equation*}

Note that by Theorem 3.5 the lower bound in (3.22) is attained by an arbitrary
function from $M(H)$. \end{proof}

\textbf{Remark 3.6.} Using Theorems 2.7 and 2.8 one can obtain sharp
estimates of type (3.22) for bivariate functions defined on the
corresponding simple polygons with sides parallel to the coordinate axes.

\bigskip

\section{On the theorem of M. Golomb}

Let $X_{1},...,X_{n}$ be compact spaces and $X=X_{1}\times \cdots \times
X_{n}.$ Consider the approximation of a function $f\in C(X)$ by sums $%
g_{1}(x_{1})+\cdots +g_{n}(x_{n}),$ where $g_{i}\in C(X_{i}),$ $i=1,...,n.$
In \cite{37}, M.Golomb obtained a formula for the error of this
approximation in terms of measures constructed on special points of $X$,
called ``projection cycles". However, his proof had a gap, which was pointed
out later by Marshall and O'Farrell \cite{107}. But the question if the
formula was correct, remained open. The purpose of this section is to prove
that Golomb's formula is valid, and moreover it holds in a stronger form.

\subsection{History of Golomb's formula}

Let $X_{i},i=1,...,n,$ be compact Hausdorff spaces. Consider the
approximation to a continuous function $f$, defined on $X=X_{1}\times \cdots
\times X_{n}$, from the manifold
\begin{equation*}
M=\left\{ \sum_{i=1}^{n}g_{i}(x_{i}):g_{i}\in C(X_{i}),~~i=1,...,n\right\} .
\end{equation*}%
The approximation error is defined as the distance from $f$ to $M$:
\begin{equation*}
E(f)\overset{def}{=}dist(f,M)=\underset{g\in M}{\inf }\left\Vert
f-g\right\Vert _{C(X)}.
\end{equation*}

The well-known duality relation says that

\begin{equation*}
E(f)=\underset{\left\Vert \mu \right\Vert \leq 1}{\underset{\mu \in M^{\bot }%
}{\sup }}\left\vert \int\limits_{X}fd\mu \right\vert ,\eqno(3.32)
\end{equation*}%
where $M^{\bot }$ is the space of regular Borel measures annihilating all
functions in $M$ and $\left\Vert \mu \right\Vert $ stands for the total
variation of a measure $\mu $. It should be noted that the $\sup $ in (3.32)
is attained by some measure $\mu ^{\ast }$ with total variation $\left\Vert
\mu ^{\ast }\right\Vert =1.$ We are interested in the problem: is it
possible to replace in (3.32) the class $M^{\bot }$ by some subclass of it
consisting of measures of simple structure? For the case $n=2,$ this problem
was first considered by Diliberto and Straus \cite{26}. They showed that the
measures generated by closed bolts are sufficient for the equality (3.32).

In case of general topological spaces, a lightning bolt is defined similarly
to the case $\mathbb{R}^2$. Let $X=X_{1}\times X_{2}$ and $\pi _{i}$ be the
projections of $X$ onto $X_{i},$ $i=1,2.$ A lightning bolt (or, simply, a
bolt) is a finite ordered set $\{a_{1},...,a_{k}\}$ contained in $X$, such
that $a_{i}\neq a_{i+1}$, for $i=1,2,...,k-1$, and either $\pi
_{1}(a_{1})=\pi _{1}(a_{2}),$ $\pi _{2}(a_{2})=\pi _{2}(a_{3})$, $\pi
_{1}(a_{3})=\pi _{1}(a_{4}),...,$ or $\pi _{2}(a_{1})=\pi _{2}(a_{2}),$ $\pi
_{1}(a_{2})=\pi _{1}(a_{3})$, $\pi _{2}(a_{3})=\pi _{2}(a_{4}),...$ A bolt $%
\{a_{1},...,a_{k}\}$ is said to be closed if $k$ is an even number and the
set $\{a_{2},...,a_{k},a_{1}\}$ is also a bolt.

Let $l=\{a_{1},...,a_{2k}\}$ be a closed bolt. Consider a measure $\mu _{l}$
having atoms $\pm \frac{1}{2k}$ with alternating signs at the vertices of $l$%
. That is,
\begin{equation*}
\mu _{l}=\frac{1}{2k}\sum_{i=1}^{2k}(-1)^{i-1}\delta _{a_{i}}\text{ \ or \ }%
\mu _{l}=\frac{1}{2k}\sum_{i=1}^{2k}(-1)^{i}\delta _{a_{i}},
\end{equation*}%
where $\delta _{a_{i}}$ is a point mass at $a_{i}.$ It is clear that $\mu
_{l}\in M^{\bot }$ and $\left\Vert \mu _{l}\right\Vert \leq 1$. $\left\Vert
\mu _{l}\right\Vert =1$ if and only if the set of vertices of the bolt $l$
having even indices does not intersect with that having odd indices. The
following duality relation was first established by Diliberto and Straus
\cite{26}
\begin{equation*}
E(f)=\underset{l\subset X}{\sup }\left\vert \int\limits_{X}fd\mu
_{l}\right\vert ,\eqno(3.33)
\end{equation*}%
where $X=X_{1}\times X_{2}$ and the $\sup $ is taken over all closed bolts
of $X$. In fact, Diliberto and Straus obtained the formula (3.33) for the
case when $X$ is a rectangle in $\mathbb{R}^{2}$ with sides parallel to the
coordinate axis. The same result was independently proved by Smolyak (see
\cite{113}). Yet another proof of (3.33), in the case when $X$ is a
Cartesian product of two compact Hausdorff spaces, was given by Light and
Cheney \cite{93}. For $X$'s other than a rectangle in $\mathbb{R}^{2}$, the
theorem under some additional assumptions appeared in the works \cite%
{62,79,107}. But we shall not discuss these works here.

Golomb's paper \cite{37} made a start to a systematic study of approximation
of multivariate functions by various compositions, including sums of
univariate functions. Golomb generalized the notion of a closed bolt to the $%
n$-dimensional case and obtained the analogue of formula (3.33) for the
error of approximation from the manifold $M$. The objects introduced in \cite%
{37} were called \textit{projection cycles} and they are defined as sets of
the form
\begin{equation*}
p=\{b_{1},...,b_{k};~c_{1},...,c_{k}\}\subset X,\eqno(3.34)
\end{equation*}%
with the property that $b_{i}\neq c_{j}$, $i,j=1,...,k$ and for all $\nu
=1,...,n,$ the group of the $\nu $-th coordinates of $c_{1},...,c_{k}$ is a
permutation of that of the $\nu $-th coordinates of $b_{1},...,b_{k}.$ Some
points in the $b$-part $\left( b_{1},...,b_{k}\right) $ or $c$-part $\left(
c_{1},...,c_{k}\right) $ of $p$ may coincide. The measure associated with $p$
is
\begin{equation*}
\mu _{p}=\frac{1}{2k}\left( \sum_{i=1}^{k}\delta
_{b_{i}}-\sum_{i=1}^{k}\delta _{c_{i}}\right).
\end{equation*}

It is clear that $\mu _{p}\in M^{\bot }$ and $\left\Vert \mu _{p}\right\Vert
=1.$ Besides, if $n=2,$ then a projection cycle is the union of closed bolts
after some suitable permutation of its points. Golomb's result states that
\begin{equation*}
E(f)=\underset{p\subset X}{\sup }\left\vert \int\limits_{X}fd\mu
_{p}\right\vert ,\eqno(3.35)
\end{equation*}%
where $X=X_{1}\times \cdots \times X_{n}$ and the $\sup $ is taken over all
projection cycles of $X$. It can be proved that in the case $n=2,$ the
formulas (3.33) and (3.35) are equivalent. Unfortunately, the proof of
(3.35) had a gap, which was pointed out many years later by Marshall and
O'Farrell \cite{107}. But the question if the formula (3.35) was correct,
remained unsolved (see also the monograph by Khavinson \cite{76}%
). Note that Golomb's result was used and cited in the literature, for
example, in works \cite{75,126}.

In the following subsection, we will construct families of normalized
measures (that is, measures with the total variation equal to $1$) on
projection cycles. Each measure $\mu _{p}$ defined above will be a member of
some family. We will also consider minimal projection cycles and measures
constructed on them. By properties of these measures, we show that Golomb's
formula (3.35) is valid in a stronger form.

\bigskip

\subsection{Measures supported on projection cycles}

Let us give an equivalent definition of a projection cycle. This will be
useful in constructing of certain measures having simple structure and
capability of approximating arbitrary measures in $M^{\bot }$.

In the sequel, $\chi _{a}$ will denote the characteristic function of a
single point set $\{a\}\subset \mathbb{R}$.

\bigskip

\textbf{Definition 3.3.} \textit{Let $X=X_{1}\times \cdots \times X_{n}$ and
$\pi _{i}$ be the projections of $X$ onto the sets $X_{i},$ $i=1,...,n.$ We
say that a set $p=\{x_{1},...,x_{m}\}\subset X$ is a projection cycle if
there exists a vector $\lambda =(\lambda _{1},...,\lambda _{m})$ with
nonzero real coordinates such that}
\begin{equation*}
\sum_{j=1}^{m}\lambda _{j}\chi _{\pi _{i}(x_{j})}=0,\text{ \ }i=1,...,n.\eqno%
(3.36)
\end{equation*}

\bigskip

Let us give some explanatory remarks concerning Definition 3.3. Fix the
subscript $i.$ Let the set $\{\pi _{i}(x_{j})$, $j=1,...,m\}$ have $s_{i}$
different values, which we denote by $\gamma _{1}^{i},\gamma
_{2}^{i},...,\gamma _{s_{i}}^{i}.$ Then (3.36) implies that
\begin{equation*}
\sum_{j}\lambda _{j}=0,
\end{equation*}%
where the sum is taken over all $j$ such that $\pi _{i}(x_{j})=\gamma
_{k}^{i},$ $k=1,...,s_{i}.$ Thus for fixed $i$, we have $s_{i}$ homogeneous
linear equations in $\lambda _{1},...,\lambda _{m}.$ The coefficients of
these equations are the integers $0$ and $1.$ By varying $i$, we obtain $%
s=\sum_{i=1}^{n}s_{i}$ such equations. Hence (3.36), in its expanded form,
stands for the system of these equations. One can observe that if this
system has a solution $(\lambda _{1},...,\lambda _{m})$ with nonzero real
components $\lambda _{i},$ then it also has a solution $(n_{1},...,n_{m})$
with nonzero integer components $n_{i},$ $i=1,...,m.$ This means that in
Definition 3.3, we can replace the vector $\lambda $ by the vector $%
n=(n_{1},...,n_{m})\,$, where $n_{i}\in \mathbb{Z}\backslash \{0\},$ $%
i=1,...,m.$ Thus, Definition 3.3 is equivalent to the following definition.

\bigskip

\textbf{Definition 3.4.} \textit{A set $p=\{x_{1},...,x_{m}\}\subset X$ is
called a projection cycle if there exist nonzero integers $n_{1},...,n_{m}$
such that}

\begin{equation*}
\sum_{j=1}^{m}n_{j}\chi _{\pi _{i}(x_{j})}=0,\text{ \ }i=1,...,n.\eqno(3.37)
\end{equation*}

\bigskip

\textbf{Lemma 3.5.} \textit{Definition 3.4 is equivalent to Golomb's
definition of a projection cycle.}

\bigskip

\begin{proof} Let $p=\{x_{1},...,x_{m}\}$ be a projection cycle with respect
to Definition 3.4. By $b$ and $c$ denote the set of all points $x_{i}$ such
that the integers $n_{i}$ associated with them in (3.37) are positive and
negative correspondingly. Write out each point $x_{i}$ $n_{i}$ times if $%
n_{i}>0$ and $-n_{i}$ times if $n_{i}<0.$ Then the set $\{b;c\}$ is a
projection cycle with respect to Golomb's definition. The
inverse is also true. Let a set $p_{1}=\{b_{1},...,b_{k};~c_{1},...,c_{k}\}$
be a projection cycle with respect to Golomb's definition. Here, some points
$b_{i}$ or $c_{i}$ may be repeated. Let $p=\{x_{1},...,x_{m}\}$ stand for
the set $p_{1}$, but with no repetition of its points. Let $n_{i}$ show how
many times $x_{i}$ appear in $p_{1}.$ We take $n_{i}$ positive if $x_{i}$
appears in the $b$-part of $p_{1}$ and negative if it appears in the $c$%
-part of $p_{1}.$ Clearly, the set $\{x_{1},...,x_{m}\}$ is a projection
cycle with respect to Definition 3.4, since the integers $n_{i},$ $%
i=1,...,m, $ satisfy (3.37). \end{proof}

In the sequel, we will use Definition 3.3. A pair $\left\langle p,\lambda
\right\rangle ,$ where $p$ is a projection cycle in $X$ and $\lambda $ is a
vector associated with $p$ by (3.36), will be called a ``projection
cycle-vector pair" of $X.$ To each such pair $\left\langle p,\lambda
\right\rangle $ with $p=\{x_{1},...,x_{m}\}$ and $\lambda =(\lambda
_{1},...,\lambda _{m})$, we correspond the measure

\begin{equation*}
\mu _{p,\lambda }=\frac{1}{\sum_{j=1}^{m}\left\vert \lambda _{j}\right\vert }%
\sum_{j=1}^{m}\lambda _{j}\delta _{x_{j}}.\eqno(3.38)
\end{equation*}

Clearly, $\mu _{p,\lambda }\in M^{\bot }$ and $\left\Vert \mu _{p,\lambda
}\right\Vert =1$. We will also deal with measures supported on some certain
subsets of projection cycles called \textit{minimal projection cycles}. A
projection cycle is said to be minimal if it does not contain any projection
cycle as its proper subset. For example, the set $p=%
\{(0,0,0),~(0,0,1),~(0,1,0),~(1,0,0),~(1,1,1)\}$ is a minimal projection
cycle in $\mathbb{R}^{3},$ since the vector $\lambda =(2,-1,-1,-1,1)$
satisfies Eq. (3.36) and there is no such vector for any other subset of $p$%
. Adding one point $(0,1,1)$ from the right to $p$, we will also have a
projection cycle, but not minimal. Note that in this case, $\lambda $ can be
taken as $(3,-1,-1,-2,2,-1).$

\bigskip

\textbf{Remark 3.7.} A minimal projection cycle under the name of a \textit{%
loop} was introduced and used in the works of Klopotowski, Nadkarni, Rao
\cite{81,80}.

\bigskip

To prove our main result we need some auxiliary facts.

\bigskip

\textbf{Lemma 3.6.} (1)\textit{\ The vector $\lambda =(\lambda
_{1},...,\lambda _{m})$ associated with a minimal projection cycle $%
p=(x_{1},...,x_{m})$ is unique up to multiplication by a constant.}

(2)\textit{\ If in (1), $\sum_{j=1}^{m}\left\vert \lambda _{j}\right\vert
=1, $ then all the numbers $\lambda _{j}$, $j=1,...,m,$ are rational.}

\bigskip

\begin{proof} Let $\lambda ^{1}=(\lambda _{1}^{1},...,\lambda _{m}^{1})$ and $%
\lambda ^{2}=(\lambda _{1}^{2},...,\lambda _{m}^{2})$ be any two vectors
associated with $p.$ That is,
\begin{equation*}
\sum_{j=1}^{m}\lambda _{j}^{1}\chi _{\pi _{i}(x_{j})}=0\text{ and }%
\sum_{j=1}^{m}\lambda _{j}^{2}\chi _{\pi _{i}(x_{j})}=0,\text{ \ }i=1,...,n.
\end{equation*}%
After multiplying the second equality by $c=\frac{\lambda _{1}^{1}}{\lambda
_{1}^{2}}$ and subtracting from the first, we obtain that
\begin{equation*}
\sum_{j=2}^{m}(\lambda _{j}^{1}-c\lambda _{j}^{2})\chi _{\pi _{i}(x_{j})}=0%
\text{, \ }i=1,...,n.
\end{equation*}%
Now since the cycle $p$ is minimal, $\lambda _{j}^{1}=c\lambda _{j}^{2},$
for all $j=1,...,m.$

The second part of the lemma is a consequence of the first part. Indeed, let
$n=(n_{1},...,n_{m})$ be a vector with the nonzero integer coordinates
associated with $p.$ Then the vector $\lambda ^{^{\prime }}=(\lambda
_{1}^{^{\prime }},...,\lambda _{m}^{^{\prime }}),$ where $\lambda
_{j}^{^{\prime }}=\frac{n_{j}}{\sum_{j=1}^{m}\left\vert n_{j}\right\vert },$
$j=1,...,m,$ is also associated with $p.$ All coordinates of $\lambda
^{^{\prime }}$ are rational and therefore by the first part of the lemma, it
is the unique vector satisfying $\sum_{j=1}^{m}\left\vert \lambda
_{j}^{^{\prime }}\right\vert =1.$ \end{proof}

By this lemma, a minimal projection cycle $p$ uniquely (up to a sign)
defines the measure
\begin{equation*}
~\mu _{p}=\sum_{j=1}^{m}\lambda _{j}\delta _{x_{j}},\text{ \ }%
\sum_{j=1}^{m}\left\vert \lambda _{j}\right\vert =1.
\end{equation*}

\bigskip

\textbf{Lemma 3.7}. \textit{Let $\mu $ be a normalized orthogonal measure on
a projection cycle $l\subset X$. Then it is a convex combination of
normalized orthogonal measures on minimal projection cycles of $l$. That is,}
\begin{equation*}
\mu =\sum_{i=1}^{s}t_{i}\mu _{l_{i}},\text{ }\sum_{i=1}^{s}t_{i}=1,~t_{i}>0,
\end{equation*}
\textit{where $l_{i},$ $i=1,...,s,$ are minimal projection cycles in $l.$}

\bigskip

This lemma follows from the result of Navada (see \cite[Theorem 2]{112}):
Let $S\subset X_{1}\times \cdots \times X_{n}$ be a finite set. Then any
extreme point of the convex set of measures $\mu $ on $S$, $\mu \in M^{\bot }
$, $\left\Vert \mu \right\Vert \leq 1$, has its support on a minimal
projection cycle contained in $S$.

\bigskip

\textbf{Remark 3.8.} In the case $n=2$, Lemma 3.7 was proved by Medvedev
(see \cite[p.77]{76}).

\bigskip

\textbf{Lemma 3.8} (see \cite[p.73]{76}). \textit{Let $X=X_{1}\times \cdots
\times X_{n}$ and $\pi _{i}$ be the projections of $X$ onto the sets $X_{i},$
$i=1,...,n.$ In order that a measure $\mu \in C(X)^{\ast }$ be orthogonal to
the subspace $M$, it is necessary and sufficient that}
\begin{equation*}
\mu \circ \pi _{i}^{-1}=0,\text{ }i=1,...,n.
\end{equation*}

\bigskip

\textbf{Lemma 3.9} (see \cite[p.75]{76}). \textit{Let $\mu \in M^{\bot }$
and $\left\Vert \mu \right\Vert =1.$ Then there exist a net of measures $%
\{\mu _{\alpha }\}\subset M^{\bot }$ weak$^{\text{*}}$ converging in $%
C(X)^{\ast }$ to $\mu $ and satisfying the following properties:}

1) $\left\Vert \mu _{\alpha }\right\Vert =1;$

2) \textit{The closed support of each $\mu _{\alpha }$ is a finite set.}

\bigskip

Our main result is the following theorem.

\bigskip

\textbf{Theorem 3.10.} \textit{The error of approximation from the manifold $%
M$ obeys the equality}
\begin{equation*}
E(f)=\underset{l\subset X}{\sup }\left\vert \int\limits_{X}fd\mu
_{l}\right\vert ,
\end{equation*}%
\textit{where the $\sup $ is taken over all minimal projection cycles of $X.$%
}

\bigskip

\begin{proof} Let $\overset{\sim }{\mu }$ be a measure with finite support $%
\{x_{1},...,x_{m}\}$ and orthogonal to the space $M.$ Put $\lambda _{j}=%
\overset{\sim }{\mu }(x_{j}),$ $j=1,...m.$ By Lemma 3.8, $\overset{\sim }{%
\mu }(\pi _{i}^{-1}(\pi _{i}(x_{j})))=0,$ for all $i=1,...,n,$ $j=1,...,m.$
Fix the indices $i$ and $j.$ Then we have the equation $\sum_{k}\lambda
_{k}=0,$ where the sum is taken over all indices $k$ such that $\pi
_{i}(x_{k})=\pi _{i}(x_{j}).$ Varying $i$ and $j,$ we obtain a system of
such equations, which concisely can be written as
\begin{equation*}
\sum_{k=1}^{m}\lambda _{k}\chi _{\pi _{i}(x_{k})}=0,\text{ \ }i=1,...,n.
\end{equation*}%
This means that the finite support of $\overset{\sim }{\mu }$ forms a
projection cycle. Therefore, a net of measures approximating the given
measure $\mu $ in Lemma 3.9 are all of the form (3.38).

Let now $\mu _{p,\lambda }$ be any measure of the form (3.38). Since $\mu
_{p,\lambda }\in M^{\bot }$ and $\left\Vert \mu _{p,\lambda }\right\Vert =1,$
we can write

\begin{equation*}
\left\vert \int\limits_{X}fd\mu _{p,\lambda }\right\vert =\left\vert
\int\limits_{X}(f-g)d\mu _{p,\lambda }\right\vert \leq \left\Vert
f-g\right\Vert ,\eqno(3.39)
\end{equation*}%
where $g$ is an arbitrary function in $M$. It follows from (3.39) that

\begin{equation*}
\underset{\left\langle p,\lambda \right\rangle }{\sup }\left\vert
\int\limits_{X}fd\mu _{p,\lambda }\right\vert \leq E(f),\eqno(3.40)
\end{equation*}%
where the $\sup $ is taken over all projection cycle-vector pairs of $X.$

Consider the general duality relation (3.32). Let $\mu _{0}$ be a measure
attaining the supremum in (3.32) and $\left\{ \mu _{p,\lambda }\right\} $ be
a net of measures of the form (3.38) approximating $\mu _{0}$ in the weak$^{%
\text{*}}$ topology of $C(X)^{\ast }.$ We already know that this is
possible. For any $\varepsilon >0,$ there exists a measure $\mu
_{p_{0},\lambda _{0}}$ in $\left\{ \mu _{p,\lambda }\right\} $ such that
\begin{equation*}
\left\vert \int\limits_{X}fd\mu _{0}-\int\limits_{X}fd\mu _{p_{0},\lambda
_{0}}\right\vert <\varepsilon .
\end{equation*}%
From the last inequality we obtain that

\begin{equation*}
\left\vert \int\limits_{X}fd\mu _{p_{0},\lambda _{0}}\right\vert >\left\vert
\int\limits_{X}fd\mu _{0}\right\vert -\varepsilon =E(f)-\varepsilon .
\end{equation*}%
Hence,
\begin{equation*}
\underset{\left\langle p,\lambda \right\rangle }{\sup }\left\vert
\int\limits_{X}fd\mu _{p,\lambda }\right\vert \geq E(f).\eqno(3.41)
\end{equation*}%
From (3.40) and (3.41) it follows that
\begin{equation*}
\underset{\left\langle p,\lambda \right\rangle }{\sup }\left\vert
\int\limits_{X}fd\mu _{p,\lambda }\right\vert =E(f).\eqno(3.42)
\end{equation*}

By Lemma 3.7,
\begin{equation*}
\mu _{p,\lambda }=\sum_{i=1}^{s}t_{i}\mu _{l_{i}},
\end{equation*}%
where $l_{i}$, $i=1,...,s,$ are minimal projection cycles in $p$ and $%
\sum_{i=1}^{s}t_{i}=1,~t_{i}>0.$ Let $k$ be an index in the set $\{1,...,s\}
$ such that

\begin{equation*}
\left\vert \int\limits_{X}fd\mu _{l_{k}}\right\vert =\max \left\{ \left\vert
\int\limits_{X}fd\mu _{l_{i}}\right\vert ,\text{ }i=1,...,s\right\} .
\end{equation*}%
Then

\begin{equation*}
\left\vert \int\limits_{X}fd\mu _{p,\lambda }\right\vert \leq \left\vert
\int\limits_{X}fd\mu _{l_{k}}\right\vert .\eqno(3.43)
\end{equation*}%
Now since

\begin{equation*}
\left\vert \int\limits_{X}fd\mu _{l}\right\vert \leq E(f),
\end{equation*}%
for any minimal cycle $l,$ from (3.42) and (3.43) we obtain the assertion of
the theorem. \end{proof}

\textbf{Remark 3.9.} Theorem 3.10 not only proves Golomb's formula, but also
improves it. Indeed, based on Lemma 3.5, one can easily observe that the
formula (3.35) is equivalent to the formula

\begin{equation*}
E(f)=\underset{\left\langle p,\lambda \right\rangle }{\sup }\left\vert
\int\limits_{X}fd\mu _{p,\lambda }\right\vert ,
\end{equation*}%
where the $\sup $ is taken over all projection cycle-vector pairs $%
\left\langle p,\lambda \right\rangle $ of $X$ provided that all the numbers $%
\lambda _{i}\diagup \sum_{j=1}^{m}\left\vert \lambda _{j}\right\vert $, $%
i=1,...,m,$ are rational. But by Lemma 3.6, minimal projection cycles enjoy
this property.

\newpage

\chapter{Generalized ridge functions and linear superpositions}

A ridge function $g(\mathbf{a}\cdot \mathbf{x})$ with a direction $\mathbf{a}%
\in \mathbb{R}^{d}\backslash \{\mathbf{0}\}$ admits a natural generalization
to a multivariate function of the form $g(\alpha _{1}(x_{1})+\cdot \cdot
\cdot +\alpha _{d}(x_{d}))$, where $\alpha _{i}(x_{i})$, $i=\overline{1,d},$
are real, presumably well behaved, fixed univariate functions. We know from
Chapter 1 that finitely many directions $\mathbf{a}^{j}$ are not enough for
sums $\sum g_{j}\left( \mathbf{a}^{j}\cdot \mathbf{x}\right) $ to
approximate multivariate functions. However, we will see in this chapter
that sums of the form $\sum g_{j}(\alpha _{1}^{j}(x_{1})+\cdot \cdot \cdot
+\alpha _{d}^{j}(x_{d}))$ with finitely many $\alpha _{i}^{j}(x_{i})$ is
capable not only approximating multivariate functions but also precisely
representing them. First we study the problem of representation of a
function $f:X\rightarrow \mathbb{R}$, where $X$ is any set, as a linear
superposition $\sum_{j}g_{j}(h_{j}(x))$ with arbitrary but fixed functions $%
h_{j}:X\rightarrow {{\mathbb{R}}}$. Then we apply the obtained result and
the famous Kolmogorov superposition theorem to prove representability of an
arbitrarily behaved multivariate function in the form of a generalized ridge
function $\sum g_{j}(\alpha _{1}^{j}(x_{1})+\cdot \cdot \cdot +\alpha
_{d}^{j}(x_{d}))$. We also study the uniqueness of representation
of functions by linear superpositions.

The material of this chapter is taken from \cite{49,Ism}.

\bigskip

\section{Representation theorems}

In this section, we study some problems of representation of real functions by
linear superpositions and linear combinations of generalized ridge functions.

\subsection{Problem statement and historical notes}

Let $X$ be any set and $h_{i}:X\rightarrow {{\mathbb{R}}},~i=1,...,r,$ be
arbitrarily fixed functions. Consider the set
\begin{equation*}
\mathcal{B}(X)=\mathcal{B}(h_{1},...,h_{r};X)=\left\{
\sum\limits_{i=1}^{r}g_{i}(h_{i}(x)),~x\in X,~g_{i}:\mathbb{R}\rightarrow
\mathbb{R},~i=1,...,r\right\} \eqno(4.1)
\end{equation*}%
Members of this set will be called linear superpositions with respect to the
functions $h_{1},...,h_{r}$ (see \cite{141}).
For a detailed study of linear superpositions and their approximation-theoretic
properties we refer the reader to the monograph by Khavinson \cite{76}.
Note that sums of generalized ridge functions $\sum g_{j}(\alpha _{1}^{j}(x_{1})+\cdot \cdot \cdot
+\alpha _{d}^{j}(x_{d}))$ with fixed $\alpha _{i}^{j}(x_{i})$
are a special case of linear superpositions. In Section 1.2,
we considered linear superpositions defined on a subset of the $d$%
-dimensional Euclidean space, while here $X$ is a set of arbitrary nature.
As in Section 1.2, we are interested in the question: what conditions on $X$
guarantee that each function on $X$ will be in the set $\mathcal{B}(X)$? The
simplest case $X\subset \mathbb{R}^{d},~r=d$ and $h_{i}$ are the coordinate
functions was solved in \cite{81}. See also \cite[p.57]{76} for the case $%
r=2.$

By $\mathcal{B}_{c}(X)$ and $\mathcal{B}_{b}(X)$ denote the right hand side
of (4.1) with continuous and bounded $g_{i}:\mathbb{R}\rightarrow \mathbb{R}%
,~i=1,...,r,$ respectively. Our starting point is the well-known
superposition theorem of Kolmogorov \cite{83}. It states that for the unit
cube $\mathbb{I}^{d},~\mathbb{I}=[0,1],~d\geq 2,$ there exists $2d+1$
functions $\{s_{q}\}_{q=1}^{2d+1}\subset C(\mathbb{I}^{d})$ of the form
\begin{equation*}
s_{q}(x_{1},...,x_{d})=\sum_{p=1}^{d}\varphi _{pq}(x_{p}),~\varphi _{pq}\in
C(\mathbb{I}),~p=1,...,d,~q=1,...,2d+1\eqno(4.2)
\end{equation*}%
such that each function $f\in C(\mathbb{I}^{d})$ admits the representation
\begin{equation*}
f(x)=\sum_{q=1}^{2d+1}g_{q}(s_{q}(x)),~x=(x_{1},...,x_{d})\in \mathbb{I}%
^{d},~g_{q}\in C({{\mathbb{R)}}}.\eqno(4.3)
\end{equation*}

Note that the functions $g_{q}(s_{q}(x))$, involved in the right hand side of
(4.3), are generalized ridge functions. In our notation, (4.3) means that $%
\mathcal{B}_{c}(s_{1},...,s_{2d+1};\mathbb{I}^{d})=C(\mathbb{I}^{d}).$ This
surprising and deep result, which solved (negatively) Hilbert's 13-th
problem, was improved and generalized in several directions. It was first
observed by Lorentz \cite{98} that the functions $g_{q}$ can be replaced by
a single continuous function $g.$ Sprecher \cite{128} showed that the
theorem can be proven with constant multiples of a single function $\varphi $
and translations. Specifically, $\varphi _{pq}$ in (4.2) can be chosen as $%
\lambda ^{p}\varphi (x_{p}+\varepsilon q),$ where $\varepsilon $ and $%
\lambda $ are some positive constants. Fridman \cite{31} succeeded in
showing that the functions $\varphi _{pq}$ can be constructed to belong to
the class $Lip(1).$ Vitushkin and Henkin \cite{141} showed that $\varphi
_{pq}$ cannot be taken to be continuously differentiable.

Ostrand \cite{115} extended the Kolmogorov theorem to general compact metric
spaces. In particular, he proved that for each compact $d$-dimensional
metric space $X$ there exist continuous real functions $\{\alpha
_{i}\}_{i=1}^{2d+1}\subset C(X)$ such that $\mathcal{B}_{c}(\alpha
_{1},...,\alpha _{2d+1};X)=C(X).$ Sternfeld \cite{130} showed that the
number $2d+1$ cannot be reduced for any $d$-dimensional space $X.$ Thus the
number of terms in the Kolmogorov superposition theorem is the best possible.

Some papers of Sternfeld were devoted to the representation of continuous
and bounded functions by linear superpositions. Let $C(X)$ and $B(X)$ denote
the space of continuous and bounded functions on some set $X$ respectively
(in the first case, $X$ is supposed to be a compact metric space). Let $%
F=\{h\}$ be a family of functions on $X.$ $F$ is called a uniformly
separating family (\textit{u.s.f.}) if there exists a number $0<\lambda \leq
1$ such that for each pair $\{x_{j}\}_{j=1}^{m}$, $\{z_{j}\}_{j=1}^{m}$ of
disjoint finite sequences in $X$, there exists some $h\in F$ so that if from
the two sequences $\{h(x_{j})\}_{j=1}^{m}$and $\{h(z_{j})\}_{j=1}^{m}$ in $%
h(X)$ we remove a maximal number of pairs of points $h(x_{j_{1}})$ and $%
h(z_{j_{2}})$ with $h(x_{j_{1}})=h(z_{j_{2}}),$ there remains at least $%
\lambda m$ points in each sequence (or , equivalently, at most $(1-\lambda )m
$ pairs can be removed). Sternfeld \cite{132} proved that for a finite
family $F=\{h_{1},...,h_{r}\}$ of functions on $X$, being a \textit{u.s.f.}
is equivalent to the equality $\mathcal{B}_{b}(h_{1},...,h_{r};X)=B(X),$ and
that in the case where $X$ is a compact metric space and the elements of $F$
are continuous functions on $X$, the equality $\mathcal{B}%
_{c}(h_{1},...,h_{r};X)=C(X)$ implies that $F$ is a \textit{u.s.f.} Thus, in
particular, Sternfeld obtained that the formula (4.3) is valid for all
bounded functions, where $g_{q}$ are bounded functions depending on $f$ (see
also \cite[p.21]{76}).

Let $X$ be a compact metric space. The family $F=\{h\}\subset C(X)$ is said
to be a measure separating family (\textit{m.s.f.}) if there exists a number
$0<\lambda \leq 1$ such that for any measure $\mu $ in $\ C(X)^{\ast },$ the
inequality $\left\Vert \mu \circ h^{-1}\right\Vert \geq \lambda \left\Vert
\mu \right\Vert $ holds for some $h\in F.$ Sternfeld \cite{131} proved that $%
\mathcal{B}_{c}(h_{1},...,h_{r};X)=C(X)$ if and only if the family $%
\{h_{1},...,h_{r}\}$ is a \textit{m.s.f.} In \cite{132}, it was also shown
that if $r=2,$ then the properties \textit{u.s.f.} and \textit{m.s.f.} are
equivalent. Therefore, the equality $\mathcal{B}_{b}(h_{1},h_{2};X)=B(X)$ is
equivalent to $\mathcal{B}_{c}(h_{1},h_{2};X)=C(X).$ But for $r\,>2$, these
two properties are no longer equivalent. That is, $\mathcal{B}%
_{b}(h_{1},...,h_{r};X)=B(X)$ does not always imply $\mathcal{B}%
_{c}(h_{1},...,h_{r};X)=C(X)$ (see \cite{131}).

Our purpose is to consider the above mentioned problem of representation by
linear superpositions without involving any topology (that of continuity or
boundedness). We start with characterization of those sets $X$ for which $%
\mathcal{B}(h_{1},...,h_{r};X)=T(X),$ where $T(X)$ is the space of all
functions on $X.$ As in Section 1.2, this will be done in terms of cycles.
We claim that nonexistence of cycles in $X$ is equivalent to the equality $%
\mathcal{B}(X)=T(X)$ for an arbitrary set $X$. In particular, we show that $%
\mathcal{B}_{c}(X)=C(X)$ always implies $\mathcal{B}(X)=T(X).$ This
implication will enable us to obtain some new results, namely extensions of
the previously known theorems from continuous to discontinuous multivariate
functions. For example, we will prove that the formula (4.3) is valid for
all discontinuous multivariate functions $f$ defined on the unite cube $%
\mathbb{I}^{d},$ where $g_{q}$ are univariate functions depending on $f.$

\bigskip

\subsection{Extension of Kolmogorov's superposition theorem}

In this subsection, we show that if some representation by linear
superpositions holds for continuous functions, then it holds for all
functions. This will lead us to natural extensions of some known
superposition theorems (such as Kolmogorov's superposition theorem,
Ostrand's superposition theorem, etc) from continuous to discontinuous
functions.

In the sequel, by $\chi _{A}$ we will denote the characteristic function of
a set $\ A\subset \mathbb{R}.$ That is,
\begin{equation*}
\chi _{A}(y)=\left\{
\begin{array}{c}
1,~if~y\in A \\
0,~if~y\notin A.%
\end{array}%
\right.
\end{equation*}

The following definition is a generalized version of Definition 1.1 from
Section 1.2, where in connection with ridge functions only subsets of $%
\mathbb{R}^{d}$ were considered.

\bigskip

\textbf{Definition 4.1.} \textit{Given an arbitrary set $X$ and functions $%
h_{i}:X\rightarrow \mathbb{R},~i=1,...,r$. A set of points $%
\{x_{1},...,x_{n}\}\subset X$ is called to be a cycle with respect to the
functions $h_{1},...,h_{r}$ (or, concisely, a cycle if there is no
confusion), if there exists a vector $\lambda =(\lambda _{1},...,\lambda
_{n})$ with the nonzero real coordinates $\lambda _{i},~i=1,...,n,$ such
that }
\begin{equation*}
\sum_{j=1}^{n}\lambda _{j}\chi _{h_{i}(x_{j})}=0,~i=1,...,r.\eqno(4.4)
\end{equation*}

\textit{A cycle $p=\{x_{1},...,x_{n}\}$ is said to be minimal if $p$ does
not contain any cycle as its proper subset.}

\bigskip

Note that in this definition the vector $\lambda =(\lambda _{1},\ldots
,\lambda _{n})$ can be chosen so that it has only integer components.
Indeed, let for $i=1,...,r,$ the set $\{h_{i}(x_{j}),~j=1,...,n\}$ have $%
k_{i}$ different values. Then it is not difficult to see that Eq. (4.4)
stands for a system of $\sum_{i=1}^{r}k_{i}$ homogeneous linear equations in
unknowns $\lambda _{1},...,\lambda _{n}.$ This system can be written in the
matrix form $(\lambda _{1},\ldots ,\lambda _{n})\times C=0,$ where $C$ is an
$n$ by $\sum_{i=1}^{r}k_{i}$ matrix. The basic property of this matrix is
that all of its entries are 0's and 1's and no row or column of $C$ is
identically zero. Since Eq. (4.4) has a nontrivial solution $(\lambda
_{1}^{^{\prime }},\ldots ,\lambda _{n}^{^{\prime }})\in \mathbf{R}^{n}$ and
all entries of $C$ are integers, by applying the Gauss elimination method we
can see that there always exists a nontrivial solution $(\lambda _{1},\ldots
,\lambda _{n})$ with the integer components $\lambda _{i}$, $i=1,...,n$.

For a number of simple examples, see Section 1.2.

Let $T(X)$ denote the set of all functions on $X.$ With each pair $%
\left\langle p,\lambda \right\rangle ,$ where $p=\{x_{1},...,x_{n}\}$ is a
cycle in $X$ and $\lambda =(\lambda _{1},...,\lambda _{n})$ is a vector
known from Definition 4.1, we associate the functional

\begin{equation*}
G_{p,\lambda }:T(X)\rightarrow \mathbb{R},~~G_{p,\lambda
}(f)=\sum_{j=1}^{n}\lambda _{j}f(x_{j}).
\end{equation*}%
In the following, such pairs $\left\langle p,\lambda \right\rangle $ will be
called \textit{cycle-vector pairs} of $X.$ It is clear that the functional $%
G_{p,\lambda }$ is linear. Besides, $G_{p,\lambda }(g)=0$ for all functions $%
g\in \mathcal{B}(h_{1},...,h_{r};X).$ Indeed, assume that (4.4) holds. Given
$i\leq r$, let $z=h_{i}(x_{j})$ for some $j$. Hence, $%
\sum_{j~(h_{i}(x_{j})=z)}\lambda _{j}=0$ and $\sum_{j~(h_{i}(x_{j})=z)}%
\lambda _{j}g_{i}(h_{i}(x_{j}))=0$. A summation yields $G_{p,\lambda
}(g_{i}\circ h_{i})=0$. Since $G_{p,\lambda }$ is linear, we obtain that $%
G_{p,\lambda }(\sum_{i=1}^{r}g_{i}\circ h_{i})=0$.

A minimal cycle $p=\{x_{1},...,x_{n}\}$ has the following obvious properties:

\begin{description}
\item[(a)] \textit{The vector $\lambda $ associated with $p$ by Eq. (4.4) is
unique up to multiplication by a constant;}

\item[(b)] \textit{If in (4.4), $\sum_{j=1}^{n}\left\vert \lambda
_{j}\right\vert =1,$ then all the numbers $\lambda _{j},~j=1,...,n,$ are
rational.}
\end{description}

Thus, a minimal cycle $p$ uniquely (up to a sign) defines the functional

\begin{equation*}
~G_{p}(f)=\sum_{j=1}^{n}\lambda _{j}f(x_{j}),\text{ \ }\sum_{j=1}^{n}\left%
\vert \lambda _{j}\right\vert =1.
\end{equation*}

\bigskip

\textbf{Proposition 4.1.} \textit{1) Let $X$ have cycles. A function $%
f:X\rightarrow \mathbb{R}$ belongs to the space $\mathcal{B}%
(h_{1},...,h_{r};X)$ if and only if $G_{p}(f)=0$ for any minimal cycle $%
p\subset X$ with respect to the functions $h_{1},...,h_{r}$.}

\textit{2) Let $X$ has no cycles. Then $\mathcal{B}(h_{1},...,h_{r};X)=T(X).$%
}

\bigskip

\textbf{Proposition 4.2.} \textit{$\mathcal{B}(h_{1},...,h_{r};X)=T(X)$ if
and only if $X$ has no cycles.}

\bigskip

These propositions are proved by the same way as Theorems 1.1 and 1.2. We
use these propositions to obtain our main result (see Theorem 4.1 below).

The condition whether $X$ have cycles or not, depends both on $X$ and the
functions $h_{1},...,h_{r}$. In the following, we see that if $%
h_{1},...,h_{r}$ are ``nice" functions (smooth functions with the simple
structure. For example, ridge functions) and $X\subset \mathbb{R}^{d}$ is a
``rich" set (for example, the set with interior points), then $X$ has always
cycles. Thus the representability by linear combinations of univariate
functions with the fixed ``nice" multivariate functions requires at least
that $X$ should not possess interior points. The picture is quite different
when the functions $h_{1},...,h_{r}$ are not ``nice". Even in the case when
they are continuous, we will see that many sets of $\mathbb{R}^{d}$ (the
unite cube, any compact subset of that, or even the whole space $\mathbb{R}%
^{d}$ itself) may have no cycles. If disregard the continuity, there exists
even one function $h$ such that every multivariate function is representable
as $g\circ h$ over any subset of $\mathbb{R}^{d}$. First, let us introduce
the following definition.

\bigskip

\textbf{Definition 4.2.} \textit{Let $X$ be a set and $h_{i}:X\rightarrow
\mathbb{R}, $ $i=1,...,r,$ be arbitrarily fixed functions. A class $A(X)$ of
functions on $X$ will be called a ``permissible function class" if for any
minimal cycle $p\subset X$ with respect to the functions $h_{1},...,h_{r}$
(if it exists), there is a function $f_{0}$ in $A(X)$ such that $%
G_{p}(f_{0})\neq 0. $}

\bigskip

Clearly, $C(X)$ and $B(X)$ are both permissible function classes (in case of
$C(X),$ $X$ is considered to be a normal topological space).

\bigskip

\textbf{Theorem 4.1.} \textit{Let $A(X)$ be a permissible function class. If
$A(X) \subset \mathcal{B}(h_{1},...,h_{r};X)$, then $\mathcal{B}%
(h_{1},...,h_{r};X)=T(X).$}

\bigskip

The proof is simple and based on Propositions 4.1 and 4.2. Assume for a
moment that $X$ admits a cycle $p$. By Proposition 4.1, the functional $G_{p}
$ annihilates all members of the set $B(h_{1},...,h_{r};X).$ By Definition
4.2 of permissible function classes, $A(X)\ $contains a function $f_{0}$
such that $G_{p}(f_{0})\neq 0.$ Therefore, $f_{0}\notin B(h_{1},...,h_{r};X)$%
. We see that the embedding $A(X) \subset B(h_{1},...,h_{r};X)$ is
impossible if $X$ has a cycle. Thus $X$ has no cycles. Then by Proposition
4.2, $\mathcal{B}(h_{1},...,h_{r};X)=T(X).$

In the ``if part" of Theorem 4.1, instead of $\mathcal{B}(h_{1},...,h_{r};X)$
and $A(X)$ one can take $\mathcal{B}_{c}(h_{1},...,h_{r};X)$ and $C(X)$ (or $%
\mathcal{B}_{b}(h_{1},...,h_{r};X)$ and $B(X)$) respectively. That is, the
following corollaries are valid.

\bigskip

\textbf{Corollary 4.1.} \textit{Let $X$ be a set and $h_{i}:X\rightarrow
\mathbb{R},$ $i=1,...,r,$ be arbitrarily fixed bounded functions\textit{. If
$\mathcal{B}_{b}(h_{1},...,h_{r};X)=B(X)$, then $\mathcal{B}%
(h_{1},...,h_{r};X)=T(X).$}}

\bigskip \textbf{Corollary 4.2.} \textit{Let $X$ be a normal topological
space and $h_{i}:X\rightarrow \mathbb{R},$ $i=1,...,r,$ be arbitrarily fixed
continuous functions\textit{. If $\mathcal{B}_{c}(h_{1},...,h_{r};X)=C(X)$,
then $\mathcal{B}(h_{1},...,h_{r};X)=T(X).$}}

\bigskip

The main advantage of Theorem 4.1 is that we need not check directly if the
set $X$ has no cycles, which in many cases may turn out to be very tedious
task. Using this theorem, we can extend free-of-charge\ the existing
superposition theorems from the classes $B(X)$ or $C(X)$ (or some other
permissible function classes) to all functions defined on $X.$ For example,
this theorem allows us to extend the Kolmogorov superposition theorem from
continuous to all multivariate functions.

\bigskip

\textbf{Theorem 4.2.} \textit{Let $d\geq 2$, $\mathbb{I}=[-1;1]$, and $%
~\varphi _{pq}, ~p=1,...,d, ~q=1,...,2d+1$, be the universal continuous
functions in (4.2). Then each multivariate function $f:\mathbb{I}%
^{d}\rightarrow \mathbb{R}$ can be represented in the form}
\begin{equation*}
f(x)=\sum_{q=1}^{2d+1}g_{q}(\sum_{p=1}^{d}\varphi
_{pq}(x_{p})),~x=(x_{1},...,x_{d})\in \mathbb{I}^{d}.
\end{equation*}%
\textit{where $g_{q}$ are univariate functions depending on $f.$}

\bigskip

It should be remarked that Sternfeld \cite{132}, in particular, obtained
that the formula (4.3) is valid for functions $f\in B(\mathbb{I}^{d})$
provided that $g_{q}$ are bounded functions depending on $f$ (see \cite[%
Chapter 1]{76} for more detailed information and interesting discussions).

Let $X$ be a compact metric space and $h_{i}\in C(X)$, $i=1,...,r.$ The
result of Sternfeld (see Section 4.1) and Corollary 4.1 give us the
implications
\begin{equation*}
\mathcal{B}_{c}(h_{1},...,h_{r};X)=C(X)\Rightarrow \mathcal{B}%
_{b}(h_{1},...,h_{r};X)=B(X)
\end{equation*}
\begin{equation*}
\Rightarrow \mathcal{B}(h_{1},...,h_{r};X)=T(X).
\end{equation*}

The first implication is invertible when $r=2$ (see \cite{132}). We want to
show that the second is not invertible even in the case $r=2.$ The following
interesting example is due to Khavinson \cite[p.67]{76}.

Let $X\subset \mathbb{R}^{2}$ consist of a broken line whose sides are
parallel to the coordinate axis and whose vertices are
\begin{equation*}
(0;0),(1;0),(1;1),(1+\frac{1}{2^{2}};1),(1+\frac{1}{2^{2}};1+\frac{1}{2^{2}}
),(1+\frac{1}{2^{2}}+\frac{1}{3^{2}};1+\frac{1}{2^{2}}),...
\end{equation*}

We add to this line the limit point of the vertices $(\frac{\pi ^{2}}{6},%
\frac{\pi ^{2}}{6})$. Let $r=2$ and $h_{1},h_{2}$ be the coordinate
functions. Then the set $X$ has no cycles with respect to $h_{1}$ and $%
h_{2}. $ By Proposition 4.1, every function $f$ on $X$ is of the form $%
g_{1}(x_{1})+g_{2}(x_{2})$, $(x_{1},x_{2})\in X$. Now construct a function $%
f_{0}$ on $X$ as follows. On the link joining $(0;0)$ to $(1;0)$ $%
f_{0}(x_{1},x_{2})$ continuously increases from $0$ to $1$; on the link from
$(1;0)$ to $(1;1)$ it continuously decreases from $1$ to $0$; on the link
from $(1;1)$ to $(1+\frac{1}{2^{2}};1)$ it increases from $0$ to $\frac{1}{2}
$; on the link from $(1+\frac{1}{2^{2}};1)$ to $(1+\frac{1}{2^{2}};1+\frac{1%
}{2^{2}})$ it decreases from $\frac{1}{2}$ to $0$; on the next link it
increases from $0$ to $\frac{1}{3}$, etc. At the point $(\frac{\pi ^{2}}{6},%
\frac{\pi ^{2}}{6})$ set the value of $f_{0}$ equal to $0.$ Obviously, $%
f_{0} $ is a continuous functions and by the above argument, $%
f_{0}(x_{1},x_{2})=g_{1}(x_{1})+g_{2}(x_{2}).$ But $g_{1}$ and $g_{2}$
cannot be chosen as continuous functions, since they get unbounded as $x_{1}$
and $x_{2}$ tends to $\frac{\pi ^{2}}{6}$. Thus, $\mathcal{B}%
(h_{1},h_{2};X)=T(X)$, but at the same time $\mathcal{B}_{c}(h_{1},h_{2};X)%
\neq C(X)$ (or, equivalently, $\mathcal{B}_{b}(h_{1},h_{2};X)\neq B(X)$).

\bigskip

\subsection{Some other superposition theorems}

We have seen in the previous subsection that the unit cube in $\mathbb{R}%
^{d} $ has no cycles with respect to some $2d+1$ continuous functions
(namely, the Kolmogorov functions $s_{q}$ (4.2)). From the result of Ostrand
\cite{115} (see Section 4.1) and Corollary 4.2 it follows that compact sets $%
X$ of finite dimension also lack cycles with respect to a certain family of
finitely many continuous functions on $X$. Namely, the following
generalization of Ostrand's theorem is valid.

\bigskip

\textbf{Theorem 4.3.} \textit{For $p=1,2,...,m$ let $X_{p}$ be a compact
metric space of finite dimension $d_{p}$ and let $n=\sum_{p=1}^{n}d_{p}.$
There exist continuous functions $\alpha _{pq}:X_{p}\rightarrow \lbrack
0,1], $ $p=1,...,m,$ $q=1,...,2n+1,$ such that every real function $f$
defined on $\Pi _{p=1}^{m}X_{p}$ is representable in the form}

\begin{equation*}
f(x_{1},...,x_{m})=\sum_{q=1}^{2n+1}g_{q}(\sum_{p=1}^{m}\alpha _{pq}(x_{p})).%
\eqno(4.5)
\end{equation*}%
\textit{where $g_{q}$ are real functions depending on $f$. If $f$ is
continuous, then the functions $g_{q}$ can be chosen continuous.}

\bigskip

Note that Ostrand proved ``if $f$ is continuous..." part of Theorem 4.3,
while we prove the validity of (4.5) for discontinuous $f$.

One may ask if there exists a finite family of functions $\{h_{i}:\mathbb{R}%
^{d}\rightarrow \mathbb{R}\}_{i=1}^{n}$ such that any subset of $\mathbb{R}%
^{d}$ does not admit cycles with respect to this family? The answer is
positive. This follows from the result of Demko \cite{23}: there exist $2d+1$
continuous functions $\varphi _{1},...,\varphi _{2d+1}$ defined on $\mathbb{R%
}^{d}$ such that every bounded continuous function on $\mathbb{R}^{d}$ is
expressible in the form $\sum_{i=1}^{2d+1}g\circ \varphi _{i}$ for some $%
g\in C(\mathbb{R})$. This theorem together with Corollary 4.1 yield that
every function on $\mathbb{R}^{d}$ is expressible in the form $%
\sum_{i=1}^{2d+1}g_{i}\circ \varphi _{i}$ for some $g_{i}:\mathbb{R}%
\rightarrow \mathbb{R},~i=1,...,2d+1$. We do not yet know if $g_{i}$ here
can be replaced by a single univariate function. We also don't know if the
number $2d+1$ can be reduced so that the whole space of $\mathbb{R}^{d}$ (or
any $d$-dimensional compact subset of that, or at least the unit cube $%
\mathbb{I}^{d}$) has no cycles with respect to some continuous functions $%
\varphi _{1},...,\varphi _{k}:\mathbb{R}^{d}\rightarrow \mathbb{R}$, where $%
k<2d+1$. One of the basic results of Sternfeld \cite{130} says that the
dimension of a compact metric space $X$ equals $d$ if and only if there
exist functions $\varphi _{1},...,\varphi _{2d+1}\in C(X)$ such that $%
\mathcal{B}_{c}(\varphi _{1},...,\varphi _{2d+1};X)=C(X)$ and for any fmily $%
\{\psi _{i}\}_{i=1}^{k}\subset C(X),$ $k<2d+1$, we have $\mathcal{B}%
_{c}(\psi _{1},...,\psi _{k};X)\neq C(X).$ In particular, from this result
it follows that the number of terms in the Kolmogorov superposition theorem
cannot be reduced. But since the equalities $\mathcal{B}_{c}(X)=C(X)$ and $%
\mathcal{B}(X)=T(X)$ are not equivalent, the above question on the
nonexistence of cycles in $\mathbb{R}^{d}$ with respect to less than $2d+1$
continuous functions is far from trivial.

If disregard the continuity, one can construct even one function $\varphi :%
\mathbb{R}^{d}\rightarrow \mathbb{R}$ such that the whole space $\mathbb{R}%
^{d}$ will not possess cycles with respect to $\varphi $ and therefore,
every function $f:\mathbb{R}^{d}\rightarrow \mathbb{R}$ will admit the
representation $f=g\circ \varphi $ with some univariate $g$ depending on $f$%
. Our argument easily follows from Corollary 4.2 and the result of Sprecher
\cite{127}: for any natural number $d$, $d\geq 2$, there exist functions $%
h_{p}:\mathbb{I}\rightarrow \mathbb{R}$, $p=1,...,d,$ such that every
function $f\in C(\mathbb{I}^{d})$ can be represented in the form

\begin{equation*}
f(x_{1},...,x_{d})=g\left( \sum_{p=1}^{d}h_{p}(x_{p})\right) ,\eqno(4.6)
\end{equation*}%
where $g$ is a univariate (generally discontinuous) function depending on $f$%
.

Note that the function involved in the right hand side of (4.6) is a generalized ridge function. Thus, the result of Sprecher together with our result means
that every multivariate function $f$ is representable as a generalized ridge
function $g\left( \cdot \right) $ and if $f$ is continuous, then $g$ can be
chosen continuous as well.

\bigskip

\textbf{Remark 4.1.} Concerning ordinary ridge functions $g(\mathbf{a}\cdot
\mathbf{x})$, representation of every multivariate function by linear
combinations of such functions may not be possible over many sets in $%
\mathbb{R}^{d}$. For example, this is not possible for sets having interior
points. More precisely, assume we are given finitely many nonzero directions
$\mathbf{a}^{1},...,\mathbf{a}^{r}$ in $\mathbb{R}^{d}$. Then $\mathcal{R}%
\left( \mathbf{a}^{1},...,\mathbf{a}^{r};X\right) \neq T(X)$ for any set $%
X\subset \mathbb{R}^{d}$ with a nonempty interior. Indeed, let $\mathbf{y}$
be a point in the interior of $X$. Consider vectors $\mathbf{b}^{i}$, $%
i=1,...,r,$ with sufficiently small coordinates such that $\mathbf{a}%
^{i}\cdot \mathbf{b}^{i}=0$, $i=1,...,r$. Note that the vectors $\mathbf{b}%
^{i}$, $i=1,...,r,$ can be chosen pairwise linearly independent. With each
vector $\mathbf{\varepsilon }=(\varepsilon _{1},...,\varepsilon _{r})$, $%
\varepsilon _{i}\in \{0,1\}$, $i=1,...,r,$ we associate the point
\begin{equation*}
\mathbf{x}_{\mathbf{\varepsilon }}=\mathbf{y+}\sum_{i=1}^{r}\varepsilon _{i}%
\mathbf{b}^{i}.
\end{equation*}%
Since the coordinates of $\mathbf{b}^{i}$ are sufficiently small, we may
assume that all the points $\mathbf{x}_{\mathbf{\varepsilon }}$ are in the
interior of $X$. We correspond each point $\mathbf{x}_{\mathbf{\varepsilon }%
} $ to the number $(-1)^{\left\vert \mathbf{\varepsilon }\right\vert }$,
where $\left\vert \mathbf{\varepsilon }\right\vert =\varepsilon _{1}+\cdots
+\varepsilon _{r}.$ One may easily verify that the pair $\left\langle \{%
\mathbf{x}_{\mathbf{\varepsilon }}\},\{(-1)^{\left\vert \mathbf{\varepsilon }%
\right\vert }\}\right\rangle $ is a cycle-vector pair of $X$. Therefore, by
Proposition 4.2, $\mathcal{R}\left( \mathbf{a}^{1},...,\mathbf{a}%
^{r};X\right) \neq T(X).$

Note that the above method of construction of the set $\{\mathbf{x}_{\mathbf{%
\varepsilon }}\}$ is due to Lin and Pinkus \cite{95}.

\bigskip

\textbf{Remark 4.2.} A different generalization of ridge functions was
considered in Lin and Pinkus \cite{95}. This generalization involves
multivariate functions of the form $g(A\mathbf{x})$, where $\mathbf{x}\in
\mathbb{R}^{d}$ is the variable, $A$ is a fixed $d\times n$ matrix, $1\leq
n<d$, and $g$ is a real-valued function defined on $\mathbb{R}^{n}$. For $%
n=1,$ this reduces to a ridge function.

\bigskip

\section{Uniqueness theorems}

Let $Q$ be a set such that every function on $Q$ can be represented by
linear superpositions. This representation is generally not unique. But for
some sets it may be unique provided that initial values of the representing
functions are prescribed at some point of $Q$. In this section, we are going
to study properties of such sets. All the obtained results are valid, in
particular, for linear combinations of generalized ridge functions.

\subsection{Formulation of the problem}

Assume $X$ is an arbitrary set, $h_{i}:X\rightarrow \mathbb{R}$, $i=1,\ldots
,r$, are fixed functions and $\mathcal{B}(X)$ is the set defined in (4.1).
Let $T(X)$ denote the set of all real functions on $X$. {Obviously, }$%
\mathcal{B}(X)${\ is a linear subspace of $T(X)$. For a set }$Q\subset X$,
let $T(Q)$ and $\mathcal{B}(Q)$ denote the restrictions of $T(X)$ and $%
\mathcal{B}(X)$ to $Q$, respectively. Sets $Q$ with the property $\mathcal{B}%
(Q)=T(Q)$ will be called \textit{representation sets}. Recall that
Proposition 4.2 gives a complete characterization of such sets. For a
representation set $Q$, we will also use the notation $Q\in RS.$ Here, $RS$
stands for the set of all representation sets in $X$.

Let $Q\in RS.$ Clearly for a function $f$ defined on $Q$ the representation
\begin{equation*}
f(x)=\sum_{i=1}^{r}g_{i}(h_{i}(x)),~x\in Q\eqno(4.7)
\end{equation*}%
is not unique. We are interested in the uniqueness of such representation
under some reasonable restrictions on the functions $g_{i}\circ h_{i}$.
These restrictions may be various, but in this section, we require that the
values of $g_{i}\circ h_{i}$ are prescribed at some point $x_{0}\in Q$. That
is, we require that
\begin{equation*}
g_{i}(h_{i}(x_{0}))=a_{i},~i=1,...,r-1,\eqno(4.8)
\end{equation*}%
where $a_{i}$ are arbitrarily fixed real numbers. Is representation (4.7)
subject to initial conditions (4.8) always unique? Obviously, not. We are
going to identify those representation sets $Q$ for which representation
(4.7) subject to conditions (4.8) is unique for all functions $%
f:Q\rightarrow \mathbb{R}$. In the sequel, such sets $Q$ will be called
\textit{unicity sets}.

\bigskip

\subsection{Complete representation sets}

From Proposition 4.2 it is easy to obtain the following set-theoretic
properties of representation sets:

\bigskip

(1) $Q\in RS$ $\Longleftrightarrow $ $A\in RS$ for every finite set $%
A\subset Q$;

(2) The union of any linearly ordered (under inclusion) system of
representation sets is also a representation set

(3) For any representation set $Q$ there is a maximal representation set,
that is, a set $M\in RS$ such that $Q\subset M$ and for any $P\supset M$, $%
P\in RS$ we have $P=M$.

(4) If $M\subset X$ is a maximal representation set, then $h_{i}(M)=h_{i}(X)$%
, $i=1,...,r$.

\bigskip

Properties (1) and (2) are obvious, since any cycle is a finite set. The
(3)-rd property follows from (2) and Zorn's lemma. To prove property (4)
note that if $x_{0}\in X$ and $h_{i}(x_{0})\notin h_{i}(M)$ for some $i$,
one can construct the representation set $M\cup \{x_{0}\}$, which is bigger
than $M$. But this is impossible, since $M$ is maximal.

\bigskip

\textbf{Definition 4.3.} \textit{A set $Q\subset X$ is called a complete
representation set if $Q$ itself is a representation set and there is no
other representation set $P$ such that $Q\subset P$ and $h_{i}(P)=h_{i}(Q)$,
$i=1,...,r$.}

\bigskip

The set of all complete representation sets of $X$ will be denoted by $CRS$.
Obviously, every representation set is contained in a complete
representation set. That is, if $A\in RS$, then there exists $B\in CRS$ such
that $h_{i}(B)=h_{i}(A),$ $i=1,...,r.$ It turns out that for the functions $%
h_{1},...,h_{r}$, complete representation sets entirely characterize unicity
sets. To prove this fact we need some auxiliary lemmas.

\bigskip

\textbf{Lemma 4.1.} \textit{Let $Q\subset X$ be a representation set and for
some point $x_{0}\in Q$ the zero function representation}
\begin{equation*}
0=\sum_{i=1}^{r}g_{i}(h_{i}(x)),\text{ \ }x\in Q,
\end{equation*}%
\textit{is unique, provided that $g_{i}(h_{i}(x_{0}))=0,$ $i=1,...,r-1$.
That is, all the functions $g_{i}\equiv 0$ on the sets $h_{i}(Q)$, $%
i=1,...,r.$ Then $Q\in CRS.$}

\bigskip

\begin{proof} Assume that $Q\notin CRS$. Then there exists a point $p\in X$ such
that $p\notin Q$, $h_{i}(p)\in h_{i}(Q)$, for all $i=1,...,r,$ and $%
Q^{^{\prime }}=Q\cup \{p\}$ is also a representation set. Consider a
function $f_{0}:Q^{^{\prime }}\rightarrow \mathbb{R}$ such that $f_{0}(q)=0$%
, for any $q\in Q$ and $f_{0}(p)=1.$ Since $Q^{^{\prime }}\in RS$,

\begin{equation*}
f_{0}(x)=\sum_{i=1}^{r}s_{i}(h_{i}(x)),\text{ \ }x\in Q^{^{\prime }}.
\end{equation*}%
Then

\begin{equation*}
f_{0}(x)=\sum_{i=1}^{r}g_{i}(h_{i}(x)),\text{ \ }x\in Q^{^{\prime }},\eqno%
(4.9)
\end{equation*}%
where

\begin{equation*}
g_{i}(h_{i}(x))=s_{i}(h_{i}(x))-s_{i}(h_{i}(x_{0})),\text{ }i=1,...,r-1
\end{equation*}%
and

\begin{equation*}
g_{r}(h_{r}(x))=s_{r}(h_{r}(x))+\sum_{i=1}^{r-1}s_{i}(h_{i}(x_{0})).
\end{equation*}%
\qquad

A restriction of representation (4.9) to the set $Q$ gives the equality

\begin{equation*}
\sum_{i=1}^{r}g_{i}(h_{i}(x))=0,\text{ for all }x\in Q.\eqno(4.10)
\end{equation*}%
Note that $g_{i}(h_{i}(x_{0}))=0,$ $i=1,...,r-1.$ It follows from the
hypothesis of the lemma that representation (4.10) is unique. Hence, $%
g_{i}(h_{i}(x))=0,$ for all $x\in Q$ and $i=1,...,r.$ But from (4.9) it
follows that

\begin{equation*}
\sum_{i=1}^{r}g_{i}(h_{i}(p))=f_{0}(p)=1.
\end{equation*}%
Since $h_{i}(p)\in h_{i}(Q)$ for all $i=1,...,r,$ the above relation
contradicts that the functions $g_{i}$ are identically zero on the sets $%
h_{i}(Q)$, $i=1,...,r.$ This means that our assumption is not true and $Q\in
CRS.$
\end{proof}

The following lemma is a strengthened version of Lemma 4.1.

\bigskip

\textbf{Lemma 4.2.} \textit{Let $Q\in RS$ and for some point $x_{0}\in Q$,
numbers $c_{1},c_{2},...,c_{r-1}\in \mathbb{R}$ and a function $v\in T(Q)$
the representation}
\begin{equation*}
v(x)=\sum_{i=1}^{r}v_{i}(h_{i}(x))
\end{equation*}%
\textit{is unique under the initial conditions $v_{i}(h_{i}(x_{0}))=c_{i},$ $%
i=1,...,r-1$. Then for any numbers $b_{1},b_{2}...,b_{r-1}\in \mathbb{R}$
and an arbitrary function $f\in T(Q)$ the representation}
\begin{equation*}
f(x)=\sum_{i=1}^{r}f_{i}(h_{i}(x))
\end{equation*}%
\textit{is also unique, provided that $f_{i}(h_{i}(x_{0}))=b_{i},$ $%
i=1,...,r-1$. Besides, $Q\in CRS.$}

\bigskip

\begin{proof} Assume the contrary. Assume that there is a function $f\in T(Q)$
having two different representations subject to the same initial conditions.
That is,
\begin{equation*}
f(x)=\sum_{i=1}^{r}f_{i}(h_{i}(x))=\sum_{i=1}^{r}f_{i}^{^{\prime }}(h_{i}(x))
\end{equation*}%
with $f_{i}(h_{i}(x_{0}))=f_{i}^{^{\prime }}(h_{i}(x_{0}))=b_{i},$ $%
i=1,...,r-1$ and $f_{i}\neq f_{i}^{^{\prime }}$ for some indice $i\in
\{1,...,r\}.$ In this case, the function $v(x)$ will possess the following
two different representations
\begin{equation*}
v(x)=\sum_{i=1}^{r}v_{i}(h_{i}(x))=\sum_{i=1}^{r}\left[
v_{i}(h_{i}(x))+f_{i}(h_{i}(x))-f_{i}^{^{\prime }}(h_{i}(x))\right] .
\end{equation*}%
both satisfying the initial conditions. The obtained contradiction and above
Lemma 4.1 complete the proof.
\end{proof}

In the sequel, we will assume that for any points $t_{i}\in h_{i}(X),$ $%
i=1,...,r,$ the system of equations $h_{i}(x)=t_{i}$, $i=1,...,r,$ has at
least one solution.

\bigskip

\textbf{Lemma 4.3.} \textit{Let $Q\in CRS.$ Then for any point $x_{0}\in Q$
the representation}
\begin{equation*}
0=\sum_{i=1}^{r}g_{i}(h_{i}(x)),\text{ }x\in Q,\eqno(4.11)
\end{equation*}%
\textit{subject to the conditions}
\begin{equation*}
g_{i}(h_{i}(x_{0}))=0,\text{ }i=1,...,r-1,\eqno(4.12)
\end{equation*}%
\textit{is unique. That is, $g_{i}\equiv 0$ on the sets $h_{i}(Q)$, $%
i=1,...,r.$}

\bigskip

\begin{proof} Assume the contrary. Assume that representation (4.11) subject to
(4.12) is not unique, or in other words, not all of $g_{i}$ are identically
zero. Without loss of generality, we may suppose that $g_{r}(h_{r}(y))\neq 0,
$ for some $y\in Q.$ Let $\xi \in X$ be a solution of the system of
equations $h_{i}(x)=h_{i}(x_{0}),$ $i=1,...,r-1,$ and $h_{r}(x)=h_{r}(y)$.
Therefore, $g_{i}(h_{i}(\xi ))=0,$ $i=1,...,r-1,$ and $g_{r}(h_{r}(\xi
))\neq 0.$ Obviously, $\xi \notin Q.$ Otherwise, we may have $%
g_{r}(h_{r}(\xi ))=0.$

We are going to prove that $Q^{\prime }=Q\cup \{\xi \}$ is a representation
set. For this purpose, consider an arbitrary function $f:Q^{\prime
}\rightarrow \mathbb{R}$. The restriction of $f$ to the set $Q$ admits a
decomposition

\begin{equation*}
f(x)=\sum_{i=1}^{r}t_{i}(h_{i}(x)),\text{ }x\in Q.
\end{equation*}

One is allowed to fix the values $t_{i}(h_{i}(x_{0}))=0,$ $i=1,...,r-1.$
Note that then $t_{i}(h_{i}(\xi ))=0,$ $i=1,...,r-1.$ Consider now the
functions

\begin{equation*}
v_{i}(h_{i}(x))=t_{i}(h_{i}(x))+\frac{f(\xi )-t_{r}(h_{r}(\xi ))}{%
g_{r}(h_{r}(\xi ))}g_{i}(h_{i}(x)),\text{ }x\in Q^{\prime },\text{ }%
i=1,...,r.
\end{equation*}

It can be easily verified that

\begin{equation*}
f(x)=\sum_{i=1}^{r}v_{i}(h_{i}(x)),\text{ }x\in Q^{\prime }.
\end{equation*}%
Since $f$ is arbitrary, we obtain that $Q^{\prime }\in RS,$ where $Q^{\prime
}\supset Q$ and $h_{i}(Q^{\prime })=h_{i}(Q),$ $i=1,...,r.$ But this
contradicts the hypothesis of the lemma that $Q\in CRS$.
\end{proof}

The following theorem is valid.

\bigskip

\textbf{Theorem 4.4.} \textit{$Q\in CRS$ if and only if for any $x_{0}\in Q,$
any $f\in T(Q)$ and any $a_{1},...,a_{r-1}\in \mathbb{R}$ the representation}
\begin{equation*}
f(x)=\sum_{i=1}^{r}g_{i}(h_{i}(x)),\text{ }x\in Q,
\end{equation*}%
\textit{subject to the conditions $g_{i}(h_{i}(x_{0}))=a_{i},$ $i=1,...,r-1,$
is unique. Equivalently, a set $Q\in CRS$ if and only if it is a unicity set.%
}

\bigskip

Theorem 4.4 is an obvious consequence of Lemmas 4.2 and 4.3.

\bigskip

\textbf{Remark 4.3.} In Theorem 4.4, all the words "any" can be replaced
with the word "some".

\bigskip

\textbf{Remark 4.4. }For the case $X=X_{1}\times \cdot \cdot \cdot \times
X_{n}$, the possibility and uniqueness of the representation by sums $%
\sum_{i=1}^{n}u_{i}(x_{i})$, \thinspace $u_{i}:X_{i}\rightarrow \mathbb{R}$,
$i=1,...,n$, were investigated in \cite{81,80}.

\bigskip

\textbf{Examples.} Let $r=2,$ $X=\mathbb{R}^{2},$ $%
h_{1}(x_{1},x_{2})=x_{1}+x_{2},$ $h_{2}(x_{1},x_{2})=x_{1}-x_{2},$ $Q$ be
the graph of the function $x_{2}=\arcsin (\sin x_{1}).$ The set $Q$ has no
cycles with respect to the functions $h_{1}$ and $h_{2}.$ Therefore, by
Proposition 4.2, $Q\in RS.$ By adding a point $p\notin Q$, we obtain the set
$Q\cup \{p\},$ which contains a cycle and hence is not a representation set.
Thus, $Q\in CRS$ and hence $Q$ is a unicity set.

Let now $r=2,$ $X=\mathbb{R}^{2},$ $h_{1}(x_{1},x_{2})=x_{1},$ $%
h_{2}(x_{1},x_{2})=x_{2},$ and $Q$ be the graph of the function $x_{2}=x_{1}.
$ Clearly, $Q\in RS$ and $Q\notin CRS.$ By the definition of complete
representation sets, there is a set $P\supset Q$ such that $P\in RS$ and for
any $T\supset P$, $T$ is not a representation set. There are many sets $P$
with this property. One of them can be obtained by adding to $Q$ any
straight line $l$ parallel to one of the coordinate axes. Indeed, if $%
y\notin Q\cup l,$ then the set $Q_{1}=Q\cup l\cup \{y\}$ contains a
four-point cycle (with one vertex as $y$, two vertices lying on $l$ and one
vertex lying on $Q$). This means that $Q_{1}\notin RS$ and hence $Q\cup l\in
CRS.$

\bigskip

The following corollary can be easily obtained from Theorem 4.4 and Lemma
4.2.

\bigskip

\textbf{Corollary 4.3.} \textit{$Q\in CRS$ if and only if $Q\in RS$ and in
the representation}
\begin{equation*}
0=\sum_{i=1}^{r}g_{i}(h_{i}(x)),\text{ }x\in Q,
\end{equation*}%
\textit{all the functions $g_{i},$ $i=1,...,r,$ are constants.}

\bigskip

We have seen that complete representation sets enjoy the unicity property.
Let us study some other properties of these sets. The following properties
are valid.

\bigskip

(a) If $Q_{1},Q_{2}\in CRS,$ $Q_{1}\cap Q_{2}\neq \emptyset $ and $Q_{1}\cup
Q_{2}\in RS$, then $Q_{1}\cup Q_{2}\in CRS.$

(b) Let $\{Q_{\alpha }\},$ $\alpha \in \Phi ,$ be a family of complete
representation sets such that $\cap _{\alpha \in \Phi }Q_{\alpha }\neq
\emptyset $ and $\cup _{\alpha \in \Phi }Q_{\alpha }\in RS.$ Then $\cup
_{\alpha \in \Phi }Q_{\alpha }\in CRS.$

\bigskip

The above two properties follow from Corollary 4.3. Note that (b) is a
generalization of (a). The following property is a consequence of (b) and
property (2) of representation sets.

\bigskip

(c) Let $\{Q_{\alpha }\},$ $\alpha \in \Phi ,$ be a totally ordered (under
inclusion) family of complete representation sets. Then $\cup _{\alpha \in
\Phi }Q_{\alpha }\in CRS.$

\bigskip

We know that every representation set $A$ is contained in a complete
representation set $Q$ such that $h_{i}(A)=h_{i}(Q),$ $i=1,...,r.$ What can
we say about the set $Q\backslash A$? Clearly, $Q\backslash A\in RS.$ But
can we chose $Q$ so that $Q\backslash A\in CRS$? The following theorem
answers this question.

\bigskip

\textbf{Theorem 4.5.} \textit{Let $A\in RS$ and $A\notin CRS.$ Then there
exists a set $B\in CRS$ such that $A\subset B,$ $h_{i}(A)=h_{i}(B),$ $%
i=1,...,r,$ and $B\backslash A\in CRS.$}

\bigskip

\begin{proof} Since the representation set $A$ is not complete, there exists a
point $p\notin A$ such that $h_{i}(p)\in h_{i}(A),$ $i=1,...,r,$ and $%
A^{\prime }=A\cup \{p\}\in RS$. By $\mathcal{M}$ denote the collection of
sets $M$ such that

1) $A\subset M$ and $M\in RS$;

2) $h_{i}(M)=h_{i}(A)$ for all $i=1,...,r$;

3) $M\backslash A\in CRS.$

Obviously, $\mathcal{M}$ is not empty. It contains the above set $A^{\prime
} $. Consider the partial order on $\mathcal{M}$ defined by inclusion. Let $%
\{M_{\beta }\},$ $\beta \in \Gamma $, be any chain in $\mathcal{M}$. The set
$\cup _{\beta \in \Gamma }M_{\beta }$ is an upper bound for this chain. To
see this, let us check that $\cup _{\beta \in \Gamma }M_{\beta }$ belongs to
$\mathcal{M}$. That is, all the above conditions 1)-3) are satisfied. Indeed,

1) $A\subset \cup _{\beta \in \Gamma }M_{\beta }$ and $\cup _{\beta \in
\Gamma }M_{\beta }\in RS.$ This follows from property (2) of representation
sets;

2) $h_{i}(\cup _{\beta \in \Gamma }M_{\beta })=\cup _{\beta \in \Gamma
}h_{i}(M_{\beta })=\cup _{\beta \in \Gamma }h_{i}(A)=h_{i}(A),$ $i=1,...,r$;

3) $\cup _{\beta \in \Gamma }M_{\beta }\backslash A\in CRS$. This follows
from property (c) of complete representation sets and the facts that $%
M_{\beta }\backslash A\in CRS$ for any $\beta \in \Gamma $ and the system $%
\{M_{\beta }\backslash A\}$, $\beta \in \Gamma $, is totally ordered under
inclusion.

Thus we see that any chain in $\mathcal{M}$ has an upper bound. By Zorn's
lemma, there are maximal sets in $\mathcal{M}$. Assume $B$ is one of such
sets.

Let us now prove that $B\in CRS$. Assume on the contrary that $B\notin CRS$.
Then by Lemma 4.2, for any point $x_{0}\in B$ the representation

\begin{equation*}
0=\sum_{i=1}^{r}g_{i}(h_{i}(x)),\text{ }x\in B,\eqno(4.13)
\end{equation*}%
subject to the conditions $g_{i}(h_{i}(x_{0}))=0,$ $i=1,...,r-1,$ is not
unique. That is, there is a point $y\in B$ such that for some index $i,$ $%
g_{i}(h_{i}(y))\neq 0.$ Without loss of generality we may assume that $%
g_{r}(h_{r}(y))\neq 0$. Clearly, $y$ cannot belong to $B\backslash A$, since
$B\backslash A\in CRS$ and over complete representation sets, the zero
function has a trivial representation provided that conditions (4.12) hold.
Thus, $y\in A$. Let $\xi \in X$ be a point such that $h_{i}(\xi
)=h_{i}(x_{0}),$ $i=1,...,r-1$, and $h_{r}(\xi )=h_{r}(y).$ The point $\xi
\notin B,$ otherwise from (4.13) we would obtain that $%
g_{r}(h_{r}(y))=g_{r}(h_{r}(\xi ))=0$. Following the techniques in the proof
of Lemma 4.3, it can be shown that $B_{1}=B\cup \{\xi \}\in RS$.

Now we prove that $B_{1}\backslash A\in CRS$. Consider the representation

\begin{equation*}
0=\sum_{i=1}^{r}g_{i}^{\prime }(h_{i}(x)),\text{ }x\in B_{1}\backslash A,%
\eqno(4.14)
\end{equation*}%
subject to the conditions $g_{i}^{\prime }(h_{i}(x_{0}))=0,$ $i=1,...,r-1,$
where \ $x_{0}$ is some point in $B\backslash A.$ Such representation holds
uniquely on $B\backslash A,$ since $B\backslash A\in CRS$. That is, all the
functions $g_{i}^{\prime }$ are identically zero on $h_{i}(B\backslash A),$ $%
i=1,...,r$. On the other hand, since $g_{i}^{\prime }(h_{i}(\xi
))=g_{i}^{\prime }(h_{i}(x_{0}))=0$, for all $i=1,...,r-1$, we obtain that $%
g_{r}^{\prime }(h_{r}(\xi ))=0.$ This means that representation (4.14)
subject to the conditions $g_{i}^{\prime }(h_{i}(x_{0}))=0,$ $i=1,...,r-1,$
is unique on $B_{1}\backslash A.$ That is, all the functions $g_{i}^{\prime }
$ in (4.14) are zero functions on $h_{i}(B_{1}\backslash A),$ $i=1,...,r.$
Hence by Lemma 4.1, $B_{1}\backslash A\in CRS$. Thus, $B_{1}\in \mathcal{M}$%
. But the set $B$ was chosen as a maximal set in $\mathcal{M}$. We see that
our assumption $B\notin CRS$ leads to the contradiction that there is a set $%
B_{1}\in \mathcal{M}$ bigger than the maximal set $B$. Thus, in fact, $B\in
CRS$.
\end{proof}

\bigskip

\subsection{$C$-orbits and $C$-trips}

Let $A$ be a representation set. The relation on $A$ defined by setting $%
x\sim y$ if there is a finite complete representation subset of $A$
containing both $x$ and $y$, is an equivalence relation. Indeed, it is
reflexive and symmetric. It is transitive by property (a) of complete
representation sets. The equivalence classes we call $C$\textit{-orbits}. In
the case $r=2$, $C$-orbits turn into classical orbits considered by Marshall
and O'Farrell \cite{108,107}, which have a very nice
geometric interpretation in terms of paths (see Section 1.3). A classical
orbit consists of all possible traces of an arbitrary point in it traveling
alternatively in the level sets of $h_{1}$ and $h_{2}.$ In the general
setting, one partial case of $C$-orbits were introduced by Klopotowski,
Nadkarni, Rao \cite{80} under the name of \textit{%
related components}. The case considered in \cite{80} requires that
$A\subset X=X_{1}\times \cdot \cdot \cdot \times
X_{n}$ and $h_{i}$ be the canonical projections of $X$ onto $X_{i},$ $%
i=1,...,r,$ respectively.

Finite complete representation sets containing $x$ and $y$ will be called $C$%
\textit{-trips} connecting $x$ and $y$. A $C$-trip of the smallest
cardinality connecting $x$ and $y$ will be called a \textit{minimal }$C$%
\textit{-trip}.

\bigskip

\textbf{Theorem 4.6.} \textit{Let $A$ be a representation set and $x$ and $y$
be any two points of some $C$-orbit in $A$. Then there is only one minimal $C
$-trip connecting them.}

\bigskip

\begin{proof} Assume that $L_{1}$ and $L_{2}$ are two minimal $C$-trips connecting $%
x$ and $y.$ By the definition, $L_{1}$ and $L_{2}$ are complete
representation sets. Note that $L_{1}\cup L_{2}$ is also complete. Let us
prove that the set $L_{1}\cap L_{2}$ is complete. Clearly, $L_{1}\cap
L_{2}\in RS.$ Let $x_{0}\in L_{1}\cap L_{2}$. In particular, $x_{0}$ can be
one of the points $x$ and $y$. Consider the representation

\begin{equation*}
0=\sum_{i=1}^{r}g_{i}(h_{i}(x)),\text{ }x\in L_{1}\cap L_{2},\eqno(4.15)
\end{equation*}%
subject to $g_{i}(h_{i}(x_{0}))=0,$ $i=1,...,r-1$. On the strength of Lemma
4.1, it is enough to prove that this representation is unique. For $i=1,...,r
$, let $g_{i}^{\prime }$ be any extension of $g_{i}$ from the set $%
h_{i}(L_{1}\cap L_{2})$ to the set $h_{i}(L_{1})$. Construct the function

\begin{equation*}
f^{\prime }(x)=\sum_{i=1}^{r}g_{i}^{\prime }(h_{i}(x)),\text{ }x\in L_{1}.%
\eqno(4.16)
\end{equation*}%
Since $f^{\prime }(x)=0$ on $L_{1}\cap L_{2}$, the following function is
well defined

\begin{equation*}
f(x)=\left\{
\begin{array}{c}
f^{\prime }(x),\text{ }x\in L_{1}, \\
0,\text{ }x\in L_{2}.%
\end{array}%
\right.
\end{equation*}%
Since $L_{1}\cup L_{2}\in CRS$, the representation

\begin{equation*}
f(x)=\sum_{i=1}^{r}w_{i}(h_{i}(x)),\text{ }x\in L_{1}\cup L_{2}.\eqno(4.17)
\end{equation*}%
subject to

\begin{equation*}
w_{i}(h_{i}(x_{0}))=0,\text{ }i=1,...,r-1.\eqno(4.18)
\end{equation*}%
is unique. Besides, since $L_{1}\in CRS$ and $g_{i}^{\prime
}(h_{i}(x_{0}))=g_{i}(h_{i}(x_{0}))=0,$ $i=1,...,r-1$, representation (4.16)
is unique. This means that for each function $g_{i}$, there is only one
extension $g^{\prime }$. Note that

\begin{equation*}
f(x)=f^{\prime }(x)=\sum_{i=1}^{r}w_{i}(h_{i}(x)),\text{ }x\in L_{1}.
\end{equation*}%
Now from the uniqueness of representation (4.16) we obtain that

\begin{equation*}
w_{i}(h_{i}(x))=g_{i}^{\prime }(h_{i}(x)),\text{ }i=1,...,r,\text{ }x\in
L_{1}.\eqno(4.19)
\end{equation*}

A restriction of formula (4.17) to the set $L_{2}$ gives

\begin{equation*}
0=\sum_{i=1}^{r}w_{i}(h_{i}(x)),\text{ }x\in L_{2}.\eqno(4.20)
\end{equation*}%
Since $L_{2}\in CRS$, representation (4.20) subject to conditions (4.18) is
unique, whence

\begin{equation*}
w_{i}(h_{i}(x))=0,\text{ }i=1,...,r\text{, \ }x\in L_{2}.\eqno(4.21)
\end{equation*}%
From (4.19) and (4.21) it follows that

\begin{equation*}
g_{i}(h_{i}(x))=g_{i}^{\prime }(h_{i}(x))=0,\text{ }i=1,...,r,\text{ }x\in
L_{1}\cap L_{2}.
\end{equation*}%
Thus, we see that representation (4.15) subject to the conditions $%
g_{i}(h_{i}(x_{0}))=0,$ $i=1,...,r-1$ is unique on the intersection $%
L_{1}\cap L_{2}.$ Therefore by Lemma 4.1, $L_{1}\cap L_{2}\in CRS.$

Let the cardinalities of $L_{1}$ and $L_{2}$ be equal to $n.$ Since $x,y\in
L_{1}\cap L_{2}$ and $L_{1}\cap L_{2}\in CRS$, we obtain from the definition
of minimal $C$-trips that the cardinality of $L_{1}\cap L_{2}$ is also $n.$
Hence, $L_{1}\cap L_{2}=L_{1}=L_{2}.$
\end{proof}

Let $Q$ be a representation set. That is, each function $f:Q\rightarrow
\mathbb{R}$ enjoys representation (4.7). Can we construct the functions $%
g_{i},$ $i=1,...,r,$ for a given $f$? There is a procedure for constructing
one certain collection of $g_{i}$, provided that $Q$ consists of a single $C$%
-orbit, that is, any two points of $Q$ can be connected by a $C$-trip. To
describe this procedure, take a point $x_{0}\in Q$ and fix it. We are going
to find $g_{i}$ from (4.7) and conditions (4.8). Let $y$ be any point in $Q$%
. To find the values of $g_{i}$ at the points $h_{i}(y),$ $i=1,...,r,$
connect $x_{0}$ and $y$ by a minimal $C$-trip $S=\{x_{1},...,x_{n}\},$ where
$x_{1}=x_{0}$ and $x_{n}=y.$ Since $S$ is a complete representation set,
equation (4.7) subject to (4.8) has a unique solution on $S$. That is, we
can find $g_{i}(h_{i}(y)),$ $i=1,...,r,$ by solving the system of linear
equations

\begin{equation*}
\sum_{i=1}^{r}g_{i}(h_{i}(x_{j}))=f(x_{j}),\text{ }j=1,...,n\text{.}
\end{equation*}

We see that each minimal $C$-trip containing $x_{0}$ generates a system of
linear equations, which is uniquely solvable. Since any point in $Q$ can be
connected with $x_{0}$ by such a trip, we can find $g_{i}(t)$ at each point $%
t\in h_{i}(Q),$ $i=1,...,r.$

The above procedure can still be effective for some particular
representation sets $Q$ consisting of many $C$-orbits. Let $\{C_{\alpha }\},$
$\alpha \in \Lambda ,$ denote the set of all $C$-orbits of $Q$. Fix some
points $x_{\alpha }\in C_{\alpha },$ $\alpha \in \Lambda $, one in each
orbit. Let $y_{\alpha }$ be any points of $C_{\alpha },$ $\alpha \in \Lambda
,$ respectively. We can apply the above procedure of finding the values of $%
g_{i}$ at each $y_{\alpha }$ if $h_{i}(y_{\alpha })\neq h_{i}(y_{\beta })$
for all $i$ and $\alpha \neq \beta $. For $h_{i}(y_{\alpha })=h_{i}(y_{\beta
}),$ one cannot guarantee that after solving the corresponding systems of
linear equations (associated with $y_{\alpha }$ and $y_{\beta }$), the
solutions $g_{i}(h_{i}(y_{\alpha })$ and $g_{i}(h_{i}(y_{\beta }))$ will be
equal. That is, for the case $h_{i}(y_{\alpha })=h_{i}(y_{\beta })$, the
constructed functions $g_{i}$ may not be well defined.

\bigskip

\textbf{Remark 4.5.} All the results in this section are valid, in
particular, for linear combinations of generalized ridge functions.

\newpage

\chapter{Applications to neural networks}

Neural networks have increasingly been used in many areas of applied
sciences. Most of the applications employ neural networks to approximate
complex nonlinear functional dependencies on a high dimensional data set.
The theoretical justification for such applications is that any continuous
function can be approximated within an arbitrary precision by carefully
selecting parameters in the network. The most commonly used model of neural
networks is the \textit{multilayer feedforward perceptron} (MLP) model. This model
consists of a finite number of successive layers. The first and the last
layers are called the input and the output layers, respectively. The
intermediate layers are called hidden layers. MLP models are usually
classified not by their number of layers, but by their number of hidden
layers. In this chapter, we study approximation properties of the single and
two hidden layer feedforward perceptron models. Our analysis is based on
ridge functions and the Kolmogorov superposition theorem.

The material of this chapter may be found in \cite{39,46,43,65}.

\bigskip

\section{Single hidden layer neural networks}

In this section, we consider single hidden layer neural networks with a set
of weights consisting of a finite number of directions or straight lines. For certain activation functions,
we characterize compact sets $X$ in the $d$-dimensional space such that the corresponding neural
network can approximate any continuous function on $X$.

\subsection{Problem statement}

Approximation capabilities of neural networks have
been investigated in a great deal of works over the last 30 years (see,
e.g., \cite%
{Alm,An,B,15,17,Ch,19,20,21,24,GI,39,40,41,67,68,71,70,89,104,110,119,123,135}).
In this section, we are interested in questions of density of a single
hidden layer perceptron model. A typical density result shows that this
model can approximate an arbitrary function in a given class with any degree
of accuracy.

\textit{A single hidden layer perceptron model} with $r$ units in the hidden layer
and input $\mathbf{x}=(x_{1},...,x_{d})$ evaluates a function of the form
\begin{equation*}
\sum_{i=1}^{r}c_{i}\sigma (\mathbf{w}^{i}\mathbf{\cdot x}-\theta _{i}),\eqno%
(5.1)
\end{equation*}%
where the \textit{weights} $\mathbf{w}^{i}$ are vectors in $\mathbb{R}^{d}$, the
\textit{thresholds} $\theta _{i}$ and the coefficients $c_{i}$ are real numbers and \
the \textit{activation function} $\sigma $ is a univariate function, which is
considered to be continuous here. Note that in Eq (5.1) each function $%
\sigma (\mathbf{w}^{i}\mathbf{\cdot x}-\theta _{i})$ is a ridge function with the direction $\mathbf{w}^{i}$.
For various activation functions $\sigma $, it has been proved in a number
of papers that one can approximate arbitrarily well a given continuous
function by functions of the form (5.1) ($r$ is not fixed!) over any compact
subset of $\mathbb{R}^{d}$. In other words, the set
\begin{equation*}
\mathcal{M}(\sigma )=span\text{\ }\{\sigma (\mathbf{w\cdot x}-\theta ):\
\theta \in \mathbb{R}\text{, }\mathbf{w\in }\mathbb{R}^{d}\}
\end{equation*}%
is dense in the space $C(\mathbb{R}^{d})$ in the topology of uniform
convergence on compact sets (see, e.g., \cite{17,21,41,67,68}). The most
general result of this type belongs to Leshno, Lin, Pinkus and Schocken \cite%
{89}. They proved that a necessary and sufficient condition for a
continuous activation function to have the density property is that it not
be a polynomial. This result shows the efficacy of the single hidden layer
perceptron model within all possible choices of the activation function $%
\sigma $, provided that $\sigma $ is continuous. In fact, density of the set
$\mathcal{M}(\sigma )$ also holds for some reasonable sets of weights and
thresholds. (see\cite{119}).

Some authors showed that a single hidden layer perceptron with a suitably
restricted set of weights can also have the density property
(or, in neural network terminology, the \textit{universal approximation property}).
For example, White and Stinchcombe \cite{135} proved that a
single layer network with a polygonal, polynomial spline or analytic
activation function and a bounded set of weights has the density property. Ito \cite{68}
investigated this property of networks using monotone sigmoidal functions
(tending to $0$ at minus infinity and $1$ at infinity), with only weights
located on the unit sphere. We see that weights required for the density property are
not necessary to be of an arbitrarily large magnitude. But what if they are
too restricted. How can one learn approximation properties of networks with
an arbitrarily restricted set of weights? This problem is too difficult to be
solved completely in this general formulation. But there are some cases that deserve a
special attention. The most interesting case is, of course, neural networks
with weights varying on a finite set of directions or lines. To the best of
our knowledge, approximation capabilities of such networks have not been
studied yet. More precisely, let $W$ be a set of weights consisting of a
finite number of vectors (or straight lines) in $\mathbb{R}^{d}$. It is
clear that if $w$ varies only in $W$, the set $\mathcal{M}(\sigma )$ can not
be dense in $C(\mathbb{R}^{d})$ in the topology of uniform convergence on compacta (compact sets).
In this case, one may want to determine boundaries of efficacy of the model. Over
which compact sets $X\subset \mathbb{R}^{d}$ does the model preserve its
general propensity to approximate arbitrarily well every continuous
multivariate function? In Section 5.1.2, we will consider this problem and
give both sufficient and necessary conditions for well approximation
(approximation with arbitrary precision) by
networks with weights from a finite set of directions or lines. For a set $W$
of weights consisting of two vectors, we show that there is a geometrically
explicit solution to the problem. In Section 5.1.3, we discuss some
aspects of the exact representation by neural networks with weights varying on finitely many straight lines.

\bigskip

\subsection{Density results}

In this subsection we give a sufficient and also a necessary conditions for
approximation by neural networks with finitely many weights and with weights
varying on a finite set of straight lines (through the origin).

Let $X$ be a compact subset of $\mathbb{R}^{d}$. Consider the following set
functions
\begin{equation*}
\tau _{i}(Z)=\{\mathbf{x}\in Z:~|p_{i}^{-1}(p_{i}(\mathbf{x}))\bigcap Z|\geq
2\},\quad Z\subset X,~i=1,\ldots ,k,
\end{equation*}%
where $p_{i}(\mathbf{x})=\mathbf{a}^{i}\cdot \mathbf{x}$, $|Y|$ denotes the
cardinality of a considered set $Y$. Define $\tau (Z)$ to be $%
\bigcap_{i=1}^{k}\tau _{i}(Z)$ and define $\tau ^{2}(Z)=\tau (\tau (Z))$, $%
\tau ^{3}(Z)=\tau (\tau ^{2}(Z))$ and so on inductively. These functions
first appeared in the work \cite{132} by Sternfeld, where he investigated
problems of representation by linear superpositions. Clearly, $\tau
(Z)\supseteq \tau ^{2}(Z)\supseteq \tau ^{3}(Z)\supseteq ...$It is possible
that for some $n$, $\tau ^{n}(Z)=\emptyset .$ In this case, one can see that
$Z$ does not contain a cycle. In general, if some set $Z\subset X$ forms a
cycle, then $\tau ^{n}(Z)=Z.$ But the reverse is not true. Indeed, let $%
Z=X=\{(0,0,\frac{1}{2}),(0,0,1),(0,1,0),(1,0,1),(1,1,0),(\frac{1}{2},\frac{1%
}{2},0),(\frac{1}{2},\frac{1}{2},\frac{1}{2})\}$, $\mathbf{a}^{i},i=1,2,3,$
are the coordinate directions in $\mathbb{R}^{3}$. It is not difficult to
verify that $X$ does not possess cycles with respect to these directions and
at the same time $\tau (X)=X$ (and so $\tau ^{n}(X)=X$ for every $n)$.

Consider the linear combinations of ridge functions with fixed directions $%
\mathbf{a}^{1},...,\mathbf{a}^{k}$
\begin{equation*}
\mathcal{R}\left( \mathbf{a}^{1},...,\mathbf{a}^{k}\right) =\left\{
\sum\limits_{i=1}^{k}g_{i}\left( \mathbf{a}^{i}\cdot \mathbf{x}\right)
:g_{i}\in C(\mathbb{R)},~i=1,...,k\right\} .\eqno(5.2)
\end{equation*}

Let $K$ be a family of functions defined on $\mathbb{R}^{d}$ and $X$ be a
subset of $\mathbb{R}^{d}.$ By $K_{X}$ we will denote the restriction of
this family to $X.$ Thus $\mathcal{R}_{X}\left( \mathbf{a}^{1},...,\mathbf{a}%
^{k}\right) $ stands for the set of sums of ridge functions in (5.2) defined
on $X$.

The following theorem is a particular case of the known general result of
Sproston and Strauss \cite{129} established for the sum of subalgebras of $%
C(X)$.

\bigskip

\textbf{Theorem 5.1. }\textit{Let $X$ be a compact subset of $\mathbb{R}^{d}$%
. If $\cap _{n=1,2,...}\tau ^{n}(X)=\emptyset $, then the set $\mathcal{R}%
_{X}\left( \mathbf{a}^{1},...,\mathbf{a}^{k}\right) $ is dense in $C(X)$.}

\bigskip

In our analysis, we need the following lemma.

\bigskip

\textbf{Lemma 5.1. } \textit{If $\mathcal{R}_{X}\left( \mathbf{a}^{1},...,%
\mathbf{a}^{k}\right) $ is dense in $C(X),$ then the set $X$ does not
contain a cycle with respect to the directions $\mathbf{a}^{1},...,\mathbf{a}%
^{k}$.}

\bigskip

\begin{proof}
Suppose the contrary. Suppose that the set $X$ contains cycles. Each cycle $%
l=(x_{1},\ldots ,x_{n})$ and the associated vector $\lambda =(\lambda
_{1},\ldots ,\lambda _{n})$ generate the functional
\begin{equation*}
G_{l,\lambda }(f)=\sum_{j=1}^{n}\lambda _{j}f(x_{j}),\quad f\in C(X).
\end{equation*}

Clearly, $G_{l,\lambda }$ is linear and continuous with the norm $%
\sum_{j=1}^{n}|\lambda _{j}|.$It is not difficult to verify that $%
G_{l,\lambda }(g)=0$ for all functions $g\in \mathcal{R}\left( \mathbf{a}%
^{1},...,\mathbf{a}^{k}\right) .$ Let $f_{0}$ be a continuous function such
that $f_{0}(x_{j})=1$ if $\lambda _{j}>0$ and $f_{0}(x_{j})=-1$ if $\lambda
_{j}<0$, $j=1,\ldots ,n$. For this function, $G_{l,\lambda }(f_{0})\neq 0$.
Thus, we have constructed a nonzero linear functional which belongs to the
annihilator of the manifold $\mathcal{R}_{X}\left( \mathbf{a}^{1},...,%
\mathbf{a}^{k}\right) $. This means that $\mathcal{R}_{X}\left( \mathbf{a}%
^{1},...,\mathbf{a}^{k}\right) $ is not dense in $C(X)$. The obtained
contradiction proves the lemma.
\end{proof}

Now we are ready to step forward from ridge function approximation to neural
networks. Let $\sigma \in C(\mathbb{R)}$ be a continuous activation
function. For a subset $W\subset \mathbb{R}^{d},$ let $\mathcal{M}(\sigma ;W,%
\mathbb{R})$ stand for the set of neural networks with weights from $W.$
That is,
\begin{equation*}
\mathcal{M}(\sigma ;W,\mathbb{R})=span\{\sigma (\mathbf{w}\cdot \mathbf{x}%
-\theta ):~\mathbf{w}\in W,~\theta \in \mathbb{R}\}.
\end{equation*}

\bigskip

\textbf{Theorem 5.2.} \textit{Let $\sigma \in C(\mathbb{R})\cap L_{p}(%
\mathbb{R)}$, where $1\leq p<\infty $, or $\sigma $ be a continuous,
bounded, nonconstant function, which has a limit at infinity (or minus
infinity). Let $W=\{\mathbf{a}^{1},...,\mathbf{a}^{k}\} \subset \mathbb{R}^{d}$
be the given set of weights and $X$ be a
compact subset of $\mathbb{R}^{d}$. The following assertions are valid:}

\textit{(1) if $\cap _{n=1,2,...}\tau ^{n}(X)=\emptyset $, then the set $%
\mathcal{M}_{X}(\sigma ;W,\mathbb{R})$ is dense in the space of all
continuous functions on $X$.}

\textit{(2) if $\mathcal{M}_{X}(\sigma ;W,\mathbb{R})$ is dense in $C(X)$,
then the set $X$ does not contain cycles.}

\bigskip

\begin{proof}
Part (1). Let $X$ be a compact subset of $\mathbb{R}^{d}$ for which \textit{%
\ }$\cap _{n=1,2,...}\tau ^{n}(X)=\emptyset $. By Theorem 5.1, the set $%
\mathcal{R}_{X}\left( \mathbf{a}^{1},...,\mathbf{a}^{k}\right) $ is dense in
$C(X)$. This means that for any positive real number $\varepsilon $ there
exist continuous univariate functions $g_{i},$ $i=1,...,k$ such that

\begin{equation*}
\left\vert f(\mathbf{x})-\sum_{i=1}^{k}{g_{i}\left( \mathbf{a}^{i}{\cdot }%
\mathbf{x}\right) }\right\vert <\frac{\varepsilon }{k+1}\eqno(5.3)
\end{equation*}%
for all $\mathbf{x}\in X$. Since $X$ is compact, the sets $Y_{i}=\{\mathbf{a}%
^{i}{\cdot }\mathbf{x:\ x}\in X\},\ i=1,2,...,k$ are also compacts. In 1947,
Schwartz \cite{124} proved that continuous and $p$-th degree Lebesgue integrable
univariate functions or continuous, bounded, nonconstant functions having a
limit at infinity (or minus infinity) are not mean-periodic. Note that a
function $f\in C(\mathbb{R}^{d})$ is called mean periodic if the set $span$\
$\{f(\mathbf{x}-\mathbf{b}):\ \mathbf{b}\in \mathbb{R}^{d}\}$ is not dense
in $C(\mathbb{R}^{d})$ in the topology of uniform convergence on compacta
(see \cite{124}). Thus, Schwartz proved that the set
\begin{equation*}
span\text{\ }\{\sigma (y-\theta ):\ \theta \in \mathbb{R}\}
\end{equation*}%
is dense in $C(\mathbb{R)}$ in the topology of uniform convergence. We
learned about this result from Pinkus \cite[p.162]{119}. This density result
means that for the given $\varepsilon $ there exist numbers $c_{ij},\theta
_{ij}\in \mathbb{R}$, $i=1,2,...,k$, $j=1,...,m_{i}$ such that%
\begin{equation*}
\left\vert g_{i}(y)-\sum_{j=1}^{m_{i}}c_{ij}\sigma (y-\theta
_{ij})\right\vert \,<\frac{\varepsilon }{k+1}\eqno(5.4)
\end{equation*}%
for all $y\in Y_{i},\ i=1,2,...,k.$ From (5.3) and (5.4) we obtain that

\begin{equation*}
\left\Vert f(\mathbf{x})-\sum_{i=1}^{k}\sum_{j=1}^{m_{i}}c_{ij}\sigma (%
\mathbf{a}^{i}{\cdot }\mathbf{x}-\theta _{ij})\right\Vert
_{C(X)}<\varepsilon .\eqno(5.5)
\end{equation*}%
Hence $\overline{\mathcal{M}_{X}(\sigma ;W,\mathbb{R})}=C(X).$

\bigskip

Part (2). Let $X$ be a compact subset of $\mathbb{R}^{d}$ and the set $%
\mathcal{M}_{X}(\sigma ;W,\mathbb{R})$ be dense in $C(X).$ Then for an
arbitrary positive real number $\varepsilon $, inequality (5.5) holds with
some coefficients $c_{ij},\theta _{ij},\ i=1,2,\ j=1,...,m_{i}.$ Since for
each $i=1,2,...,k$, the function $\sum_{j=1}^{m_{i}}c_{ij}\sigma (\mathbf{a}%
^{i}{\cdot }\mathbf{x}-\theta _{ij})$ is a function of the form $g_{i}(%
\mathbf{a}^{i}{\cdot }\mathbf{x}),$ the subspace $\mathcal{R}_{X}\left(
\mathbf{a}^{1},...,\mathbf{a}^{k}\right) $ is dense in $C(X)$. Then by Lemma
5.1, the set $X$ contains no cycles.
\end{proof}

The above theorem still holds if the set of weights $W=\{\mathbf{a}^{1},...,\mathbf{a}^{k}\}$
is replaced by the set $W_{1}=\{t_{1}\mathbf{a}^{1},...,t_{k}\mathbf{a}^{k}:\ t_{1},...,t_{k}\in \mathbb{R}\}$.
In fact, for $W_{1}$, the above restrictions on the activation function $\sigma $ may be weakened.

\bigskip

\textbf{Theorem 5.3.} \textit{Assume $\sigma \in C(\mathbb{R})$ is not a polynomial. Let
$W_{1}=\{t_{1}\mathbf{a}^{1},...,t_{k}\mathbf{a}^{k}:\ t_{1},...,t_{k}\in \mathbb{R}\}$ be the
given set of weights and $X$ be a compact subset of $\mathbb{R}^{d}$.
The following assertions are valid:}

\textit{(1) if $\cap _{n=1,2,...}\tau ^{n}(X)=\emptyset $, then the set $%
\mathcal{M}_{X}(\sigma ;W_{1},\mathbb{R})$ is dense in the space of all
continuous functions on $X$.}

\textit{(2) if $\mathcal{M}_{X}(\sigma ;W_{1},\mathbb{R})$ is dense in $C(X)$%
, then the set $X$ does not contain cycles.}

\bigskip

The proof of this theorem is similar to that of Theorem 5.2 and based on the
following result of Leshno, Lin, Pinkus and Schocken \cite{89}: if $\sigma $
is not a polynomial, then the set
\begin{equation*}
span\text{\ }\{\sigma (ty-\theta ):\ t,\theta \in \mathbb{R}\}
\end{equation*}%
is dense in $C(\mathbb{R)}$ in the topology of uniform convergence on compacta.

The above example with the set
\begin{equation*}
\{(0,0,\frac{1}{2}),(0,0,1),(0,1,0),(1,0,1),(1,1,0),(\frac{1}{2},\frac{1}{2}%
,0),(\frac{1}{2},\frac{1}{2},\frac{1}{2})\}
\end{equation*}%
shows that the sufficient condition in part (1) of Theorem 5.2 is not
necessary. The necessary condition in part (2), in general, is not
sufficient. But it is not easily seen. Here, is the nontrivial example
showing that nonexistence of cycles is not sufficient for the density $%
\overline{\mathcal{M}_{X}(\sigma ;W,\mathbb{R})}=C(X).$ For the sake of
simplicity, we restrict ourselves to $\mathbb{R}^{2}.$ Let $\mathbf{a}%
^{1}=(1;1),$ $\mathbf{a}^{2}=(1;-1)$ and the set of weights $W=\{\mathbf{a}%
^{1},\mathbf{a}^{2}\}.$ The set $X$ can be constructed as follows. Let $%
X_{1} $ be the union of the four line segments $[(-3;0),(-1;0)],$ $%
[(-1;2),(1;2)],$ $[(1;0),(3;0)]$ and $[(-1;-2),(1;-2)].$ Rotate one segment
in $X_{1}$ $90^{\circ }$ about its center and remove the middle one-third
from each line segment. The obtained set denote by $X_{2}$. By the same way,
one can construct $X_{3},X_{4},$ and so on. It is clear that the set $X_{i}$
has $2^{i+1}$ line segments. Let $X$ be a limit of the sets $X_{i}$, $%
i=1,2,...$. Note that there are no cycles.

By $S_{i}$, $i=\overline{1,4},$ denote the closed discs with the unit radius
and centered at the points $(-2;0),$ $(0;2),$ $(2;0)$ and $(0;-2)$
respectively. Consider a continuous function $f_{0}$ such that $f_{0}(%
\mathbf{x})=1$ for $\mathbf{x}\in (S_{1}\cup S_{3})\cap X$, $f_{0}(\mathbf{x}%
)=-1$ for $\mathbf{x}\in (S_{2}\cup S_{4})\cap X$, and $-1<f_{0}(\mathbf{x}%
)<1$ elsewhere on $\mathbb{R}^{2}$. Let $p=(\mathbf{y}^{1},\mathbf{y}%
^{2},...)$ be any infinite path in $X.$ Note that the points $\mathbf{y}%
^{i}, $ $i=1,2,...,$ are alternatively in the sets $(S_{1}\cup S_{3})\cap X$
and $(S_{2}\cup S_{4})\cap X$. Obviously,
\begin{equation*}
E(f_{0},X)\overset{def}{=}\inf_{g\in \mathcal{R}_{X}\left( \mathbf{a}^{1},%
\mathbf{a}^{2}\right) }\left\Vert f_{0}-g\right\Vert _{C(X)}\leq \left\Vert
f_{0}\right\Vert _{C(X)}=1.\eqno(5.6)
\end{equation*}

For each positive integer $k=1,2,...$, set $p_{k}=(\mathbf{y}^{1},...,%
\mathbf{y}^{k})$ and consider the path functionals
\begin{equation*}
G_{p_{k}}(f)=\frac{1}{k}\sum_{i=1}^{k}(-1)^{i-1}f(\mathbf{y}^{i}).
\end{equation*}

$G_{p_{k}}$ is a continuous linear functional obeying the following obvious
properties:

\begin{enumerate}
\item[(1)] $\left\Vert G_{p_{k}}\right\Vert =G_{p_{k}}(f_{0})=1;$

\item[(2)] $G_{p_{k}}(g_{1}+g_{2})\leq \frac{2}{k}(\left\Vert
g_{1}\right\Vert +\left\Vert g_{2}\right\Vert )$ for ridge functions $g_{1}={%
g_{1}\left( \mathbf{a}^{1}{\cdot }\mathbf{x}\right) }$ and $g_{2}={%
g_{2}\left( \mathbf{a}^{2}{\cdot }\mathbf{x}\right) .}$
\end{enumerate}

By property (1), the sequence $\{G_{p_{k}}\}_{k=1}^{\infty }$ has a weak$^{%
\text{*}}$ cluster point. This point will be denoted by $G.$ By property
(2), $G\in \mathcal{R} _{X}\left( \mathbf{a}^{1},\mathbf{a}^{2}\right)
^{\bot}.$ Therefore,
\begin{equation*}
1=G(f_{0})=G(f_{0}-g)\leq \left\Vert f_{0}-g\right\Vert _{C(X)}\text{ \ for
any }g\in \mathcal{R}_{X}\left( \mathbf{a}^{1},\mathbf{a}^{2}\right) .
\end{equation*}

Taking $\inf$ over $g$ in the right-hand side of the last inequality, we
obtain that $1\leq E(f_{0},X).$ Now it follows from (5.6) that $%
E(f_{0},X)=1. $ Recall that $\mathcal{M}_{X}(\sigma ;W,\mathbb{R})\subset
\mathcal{R}_{X}\left( \mathbf{a}^{1},\mathbf{a}^{2}\right) .$ Thus%
\begin{equation*}
\inf_{h\in \mathcal{M}_{X}(\sigma ;W,\mathbb{R})}\left\Vert f-h\right\Vert
_{C(X)}\geq 1.
\end{equation*}

The last inequality finally shows that $\overline{\mathcal{M}_{X}(\sigma ;W,%
\mathbb{R})}\neq C(X).$

\bigskip

For neural networks with weights consisting of only two vectors (or
directions) the problem of density becomes more clear. In this case, under
some minor restrictions on $X,$ the necessary condition in part (2) of
Theorem 5.2 (nonexistence of cycles) is also sufficient for the density of $%
\mathcal{M}_{X}(\sigma ;W,\mathbb{R})$ in $C(X)$. These restrictions are
imposed on the following equivalent classes of $X$ induced by paths. The
relation $\mathbf{x}\thicksim \mathbf{y}$ when $\mathbf{x}$ and $\mathbf{y}$
belong to some path in a given compact set $X\subset \mathbb{R}^{d}$ defines
an equivalence relation. Recall that the equivalence classes are called
orbits (see Section 1.3.4).

\bigskip

\textbf{Theorem 5.4.} \textit{Let $\sigma \in C(\mathbb{R})\cap L_{p}(%
\mathbb{R)}$, where $1\leq p<\infty $, or $\sigma $ be a continuous,
bounded, nonconstant function, which has a limit at infinity (or minus
infinity). Let $W=\{\mathbf{a}^{1},\mathbf{a}^{2}\} \subset \mathbb{R}^{d}$
be the given set of weights and $X$ be a compact subset of $\mathbb{R}%
^{d}$ with all its orbits closed. Then $\mathcal{M}_{X}(\sigma ;W,\mathbb{R}%
) $ is dense in the space of all continuous functions on $X$ if and only
if $X$ contains no closed paths with respect to the directions
$\mathbf{a}^{1}$ and $\mathbf{a}^{2}$.}

\bigskip

\begin{proof} \textit{Sufficiency.} Let $X$ be a compact subset of $\mathbb{R}%
^{d}$ with all its orbits closed. Besides, let $X$ contain no closed paths.
By Theorem 1.6 (see Section 1.3.4), the set $\mathcal{R}_{X}\left( \mathbf{a}%
^{1},\mathbf{a}^{2}\right) $ is dense in $C(X)$. This means that for any
positive real number $\varepsilon $ there exist continuous univariate
functions $g_{1}$ and $g_{2}$ such that
\begin{equation*}
\left\vert f(\mathbf{x})-{g_{1}\left( \mathbf{a}^{1}{\cdot }\mathbf{x}%
\right) -g_{2}\left( \mathbf{a}^{2}{\cdot }\mathbf{x}\right) }\right\vert <%
\frac{\varepsilon }{3}\eqno(5.7)
\end{equation*}%
for all $\mathbf{x}\in X$. Since $X$ is compact, the sets $Y_{i}=\{\mathbf{a}%
^{i}{\cdot }\mathbf{x:\ x}\in X\},\ i=1,2,$ are also compacts. As mentioned
above, Schwartz \cite{124} proved that continuous and $p$-th degree Lebesgue
integrable univariate functions or continuous, bounded, nonconstant
functions having a limit at infinity (or minus infinity) are not
mean-periodic. Thus, the set
\begin{equation*}
span\text{\ }\{\sigma (y-\theta ):\ \theta \in \mathbb{R}\}
\end{equation*}%
is dense in $C(\mathbb{R)}$ in the topology of uniform convergence. This
density result means that for the given $\varepsilon $ there exist numbers $%
c_{ij},\theta _{ij}\in \mathbb{R}$, $i=1,2,$ $j=1,\dots ,m_{i}$ such that%
\begin{equation*}
\left\vert g_{i}(y)-\sum_{j=1}^{m_{i}}c_{ij}\sigma (y-\theta
_{ij})\right\vert \,<\frac{\varepsilon }{3}\eqno(5.8)
\end{equation*}%
for all $y\in Y_{i},\ i=1,2.$ From (5.7) and (5.8) we obtain that
\begin{equation*}
\left\Vert f(\mathbf{x})-\sum_{i=1}^{2}\sum_{j=1}^{m_{i}}c_{ij}\sigma (%
\mathbf{a}^{i}{\cdot }\mathbf{x}-\theta _{ij})\right\Vert
_{C(X)}<\varepsilon .\eqno(5.9)
\end{equation*}%
Hence $\overline{\mathcal{M}_{X}(\sigma ;W,\mathbb{R})}=C(X).$

\bigskip

\textit{Necessity.} Let $X$ be a compact subset of $\mathbb{R}^{n}$ with all
its orbits closed and the set $\mathcal{M}_{X}(\sigma ;W,\mathbb{R})$ be
dense in $C(X).$ Then for an arbitrary positive real number $\varepsilon $,
inequality (5.9) holds with some coefficients $c_{ij},\theta _{ij},\ i=1,2,\
j=1,\dots ,m_{i}.$ Since for $i=1,2,$ $\sum_{j=1}^{m_{i}}c_{ij}\sigma (%
\mathbf{a}^{i}{\cdot }\mathbf{x}-\theta _{ij})$ is a function of the form $%
g_{i}(\mathbf{a}^{i}{\cdot }\mathbf{x}),$ the subspace $\mathcal{R}%
_{X}\left( \mathbf{a}^{1},\mathbf{a}^{2}\right) $ is dense in $C(X)$. Then
by Theorem 1.6, the set $X$ contains no closed paths.
\end{proof}

\bigskip

\textbf{Remark 5.1.} It can be shown that the necessity of the theorem is
valid without any restriction on orbits of $X$. Indeed if $X$ contains a
closed path, then it contains a closed path $p=(\mathbf{x}^{1},\dots ,%
\mathbf{x}^{2m})$ with different points. The functional $G_{p}=%
\sum_{i=1}^{2m}(-1)^{i-1}f(\mathbf{x}^{i})$ belongs to the annihilator of
the subspace $\mathcal{R}_{X}\left( \mathbf{a}^{1},\mathbf{a}^{2}\right) .$
There exist nontrivial continuous functions $f_{0}$ on $X$ such that $%
G_{p}(f_{0})\neq 0$ (take, for example, any continuous function $f_{0}$
taking values $+1$ at $\{\mathbf{x}^{1},\mathbf{x}^{3},\dots ,\mathbf{x}%
^{2m-1}\}$, $-1$ at $\{\mathbf{x}^{2},\mathbf{x}^{4},\dots ,\mathbf{x}%
^{2m}\} $ and $-1<f_{0}(\mathbf{x})<1$ elsewhere). This shows that the
subspace $\mathcal{R}_{X}\left( \mathbf{a}^{1},\mathbf{a}^{2}\right) $ is
not dense in $C(X)$. But in this case, the set $\mathcal{M}_{X}(\sigma ;W,%
\mathbb{R})$ cannot be dense in $C(X)$. The obtained contradiction means
that our assumption is not true and $X$ contains no closed paths.

\bigskip

Theorem 5.4 remains valid if the set of weights $\ W=\{\mathbf{a}^{1},%
\mathbf{a}^{2}\}$ is replaced by the set $W_{1}=\{t_{1}\mathbf{a}^{1},t_{2}%
\mathbf{a}^{2}:\ t_{1},t_{2}\in \mathbb{R\}}$. In fact, for the set $W_{1}$,
the required conditions on $\sigma $ may be weakened. As in Theorem 5.3, the
activation function $\sigma $ can be taken only non-polynomial.

\bigskip

\textbf{Theorem 5.5.} \textit{Assume $\sigma \in C(\mathbb{R})$
is not a polynomial. Let $\mathbf{a}^{1}$ and $\mathbf{a}^{2}$ be
fixed vectors and $W_{1}=\{t_{1}\mathbf{a}^{1},t_{2}%
\mathbf{a}^{2}:\ t_{1},t_{2}\in \mathbb{R\}}$
be the set of weights. Let $X$ be a compact subset of $\mathbb{R%
}^{d}$ with all its orbits closed. Then $\mathcal{M}_{X}(\sigma ;W_{1},%
\mathbb{R})$ is dense in the space of all continuous functions on $X$ if
and only if $X$ contains no closed paths with respect to the directions
$\mathbf{a}^{1}$ and $\mathbf{a}^{2}$.}

\bigskip

The proof is analogous to that of Theorem 5.4 and based on the above
mentioned result of Leshno, Lin, Pinkus and Schocken \cite{89}.

\bigskip

\textbf{Examples:}

\begin{description}
\item[(a)] Let $\mathbf{a}^{1}$ and $\mathbf{a}^{2}$ be two noncollinear
vectors in $\mathbb{R}^{2}.$ Let $B=B_{1}...B_{k}$ be a broken line with the
sides $B_{i}B_{i+1},\ i=1,...,k-1,$ alternatively perpendicular to $\mathbf{a%
}^{1}$ and $\mathbf{a}^{2}$. Besides, let $B$ does not contain vertices of
any parallelogram with sides perpendicular to these vectors. Then the set $%
\mathcal{M}_{B}(\sigma ;\{\mathbf{a}^{1},\mathbf{a}^{2}\},\mathbb{R})$ is
dense in $C(B).$

\item[(b)] Let $\mathbf{a}^{1}$ and $\mathbf{a}^{2}$ be two noncollinear
vectors in $\mathbb{R}^{2}.$ If $X$ is the union of two parallel line
segments, not perpendicular to any of the vectors $\mathbf{a}^{1}$ and $%
\mathbf{a}^{2}$, then the set $\mathcal{M}_{X}(\sigma ;\{\mathbf{a}^{1},
\mathbf{a}^{2}\},\mathbb{R})$ is dense in $C(X).$

\item[(c)] Let now $\mathbf{a}^{1}$ and $\mathbf{a}^{2}$ be two collinear
vectors in $\mathbb{R}^{2}.$ Note that in this case any path consisting of
two points is automatically closed. Thus the set $\mathcal{M}_{X}(\sigma ;\{%
\mathbf{a}^{1},\mathbf{a}^{2}\},\mathbb{R})$ is dense in $C(X)$ if and only
if $X$ contains no path different from a singleton. A simple example is a
line segment not perpendicular to the given direction.

\item[(d)] Let $X$ be a compact set with an interior point. Then Theorem 5.4
fails, since any such set contains vertices of some parallelogram with sides
perpendicular to the given directions $\mathbf{a}^{1}$ and $\mathbf{a}^{2}$,
that is a closed path.
\end{description}

\bigskip

\subsection{A necessary condition for the representation by neural networks}

In this subsection we give a necessary condition for the representation of
functions by neural networks with weights from a finitely many straight
lines. Before formulating our result, we introduce new objects, namely \textit{semicycles}
with respect to directions $\mathbf{a}^{1},...,%
\mathbf{a}^{k}\in \mathbb{R}^{d}\backslash \{\mathbf{0}\}$.

\bigskip

\textbf{Definition 5.1.} \textit{A set of points $l=(\mathbf{x}^{1},\ldots ,%
\mathbf{x}^{n})\subset \mathbb{R}^{d}$ is called a semicycle with respect to
the directions $\mathbf{a}^{1},...,\mathbf{a}^{k}$ if there
exists a vector $\lambda =(\lambda _{1},\ldots ,\lambda _{n})\in \mathbf{Z}%
^{n}\setminus \{\mathbf{0}\}$ such that for any $i=1,\ldots ,k,$ we have}
\begin{equation*}
\sum_{j=1}^{n}\lambda _{j}\delta _{\mathbf{a}^{i}\cdot \mathbf{x}%
^{j}}=\sum_{s=1}^{r_{i}}\lambda _{i_{s}}\delta _{\mathbf{a}^{i}\cdot \mathbf{%
x}^{i_{s}}},\quad where~r_{i}\leq k.\eqno(5.10)
\end{equation*}%
Here $\delta _{a}$ is the characteristic function of the single point set $%
\{a\}$. Note that for $i=1,\ldots ,k$, the set $\{\lambda
_{i_{s}},~s=1,...,r_{i}\}$ is a subset of the set $\{\lambda
_{j},~j=1,...,n\}$. Thus, Eq. (5.10) means that for each $i$, we actually
have at most $k$ terms in the sum $\sum_{j=1}^{n}\lambda _{j}\delta _{%
\mathbf{a}^{i}\cdot \mathbf{x}^{j}}$.

Recall that if in (5.10) for any $i=1,\ldots ,k$, we have
\begin{equation*}
\sum_{j=1}^{n}\lambda _{j}\delta _{\mathbf{a}^{i}\cdot \mathbf{x}^{j}}=0,
\end{equation*}%
then the set $l=(\mathbf{x}^{1},\ldots ,\mathbf{x}^{n})$ is a cycle with
respect to the directions $\mathbf{a}^{1},...,\mathbf{a}^{k}$
(see Section 1.2). Thus a cycle is a special case of a semicycle.

Let us give a simple example of a semicycle. Assume $k=2$ and $\mathbf{a}^{1}\cdot \mathbf{x}^{1}=\mathbf{a}%
^{1}\cdot \mathbf{x}^{2}$, $\mathbf{a}^{2}\cdot \mathbf{x}^{2}=\mathbf{a}%
^{2}\cdot \mathbf{x}^{3}$, $\mathbf{a}^{1}\cdot \mathbf{x}^{3}=\mathbf{a}%
^{1}\cdot \mathbf{x}^{4}$,..., $\mathbf{a}^{2}\cdot \mathbf{x}^{n-1}=\mathbf{%
a}^{2}\cdot \mathbf{x}^{n}$. Then it is not difficult to see that for a
vector $\lambda =(\lambda _{1},\ldots ,\lambda _{n})$ with the components $\lambda _{j}=(-1)^{j},$
the following equalities hold:
\begin{eqnarray*}
\sum_{j=1}^{n}\lambda _{j}\delta _{\mathbf{a}^{1}\cdot \mathbf{x}^{j}}
&=&\lambda _{n}\delta _{\mathbf{a}^{1}\cdot \mathbf{x}^{n}}, \\
\sum_{j=1}^{n}\lambda _{j}\delta _{\mathbf{a}^{2}\cdot \mathbf{x}^{j}}
&=&\lambda _{1}\delta _{\mathbf{a}^{2}\cdot \mathbf{x}^{1}}.
\end{eqnarray*}%
Thus, by Definition 5.1, the set $l=\{\mathbf{x}^{1},\ldots ,\mathbf{x}%
^{n}\}$ is a semicycle with respect to the directions $\mathbf{a}^{1}$
and $\mathbf{a}^{2}$. Note that this set, in the given order of its points, forms
a path with respect to the directions $\mathbf{a}^{1}$ and $\mathbf{a}^{2}$ (see
Section 1.3). It is not difficult to see that any path with respect to
$\mathbf{a}^{1}$ and $\mathbf{a}^{2}$ is a semicycle with respect
to these directions. But semicycles may also involve some union of paths.

Note that one can construct many semicycles by adding not more than $k$
arbitrary points to a cycle with respect to the directions $\mathbf{a}^{1},%
\mathbf{a}^{2},...,\mathbf{a}^{k}$.

A cycle (or semicycle) $l$ is called a \textit{$q$-cycle} (\textit{$q$%
-semicycle}) if the vector $\lambda $ associated with $l$ can be chosen so
that $\left\vert \lambda _{i}\right\vert \leq q,$ $i=1,...,n,$ and $q$ is
the minimal number with this property.

The semicycle considered above is a $1$-semicycle. If in that example, $%
\mathbf{a}^{2}\cdot \mathbf{x}^{n-1}=\mathbf{a}^{2}\cdot \mathbf{x}^{1}$,
then the set $\{x_{1},x_{2},...,x_{n-1}\}$ is a $1$-cycle. Let us give a
simple example of a $2$-cycle with respect to the directions $\mathbf{a}%
^{1}=(1,0)$ and $\mathbf{a}^{2}=(0,1)$. Consider the union
\begin{equation*}
\{0,1\}^{2}\cup \{0,2\}^{2}=\{(0,0),(1,1),(2,2),(0,1),(1,0),(0,2),(2,0)\}.
\end{equation*}

It is easy to see that this set is a $2$-cycle with the associated vector $%
(2,1,1,-1,-1,-1,-1).$ Similarly, one can construct a $q$-cycle or $q$%
-semicycle for any positive integer $q$.

\bigskip

\textbf{Theorem 5.6.} \textit{Assume $W=\{t_{1}\mathbf{a}^{1},...,t_{k}\mathbf{a}^{k}:\ t_{1},...,t_{k}\in \mathbb{R}\}$
is the given set of weights. If $\mathcal{M}_{X}(\sigma ;W,\mathbb{R})=C(X)$, then $X$ contains no cycles
and the lengths (number of points) of all $q$-semicycles in $X$ are bounded
by some positive integer.}

\bigskip

\begin{proof} Let $\mathcal{M}_{X}(\sigma ;W,\mathbb{R})=C(X).$\textit{\ Then}
$\mathcal{R}_{1}+\mathcal{R}_{2}+...+\mathcal{R}_{k}=C\left( X\right) $,
where
\begin{equation*}
\mathcal{R}_{i}=\{g_{i}(\mathbf{a}^{i}\cdot \mathbf{x):~}g_{i}\in C(\mathbb{%
R)}\},~i=1,2,...,k.
\end{equation*}

Consider the linear space%
\begin{equation*}
\mathcal{U}=\prod_{i=1}^{k}\mathcal{R}_{i}=\{(g_{1},\ldots ,g_{k}):~g_{i}\in
\mathcal{R}_{i},~i=1,\ldots ,k\}
\end{equation*}%
endowed with the norm
\begin{equation*}
\Vert (g_{1},\ldots ,g_{k})\Vert =\Vert g_{1}\Vert +\cdots +\Vert g_{k}\Vert
.
\end{equation*}

By $\mathcal{U}^{\ast}$ denote the dual space of $\mathcal{U}$. Each functional $%
F\in \mathcal{U}^{\ast }$ can be written as
\begin{equation*}
F=F_{1}+\cdots +F_{k},
\end{equation*}%
where the functionals $F_{i}\in \mathcal{R}_{i}^{\ast }$ and
\begin{equation*}
F_{i}(g_{i})=F[(0,\ldots ,g_{i},\ldots ,0)],\quad i=1,\ldots ,k.
\end{equation*}%
We see that the functional $F$ determines the collection $%
(F_{1},\ldots ,F_{k})$. Conversely, every collection $(F_{1},\ldots ,F_{k})$
of continuous linear functionals $F_{i}\in \mathcal{R}_{i}^{\ast }$, $%
i=1,\ldots ,k$, determines the functional $F_{1}+\cdots +F_{k},$ on $%
\mathcal{U}$. Considering this, in what follows, elements of $\mathcal{U}%
^{\ast }$ will be denoted by $(F_{1},\ldots ,F_{k})$.

It is not difficult to verify that
\begin{equation*}
\Vert (F_{1},\ldots ,F_{k})\Vert =\max \{\Vert F_{1}\Vert ,\ldots ,\Vert
F_{k}\Vert \}.\eqno(5.11)
\end{equation*}

Let $l=(\mathbf{x}^{1},\ldots ,\mathbf{x}^{n})$ be any $q$-semicycle (with
respect to the directions $\mathbf{a}^{1}$,...,$\mathbf{a}%
^{k}$) in $X$ and $\lambda =(\lambda _{1},\ldots ,\lambda _{n})$ be a vector
associated with it. Consider the following functional

\begin{equation*}
G_{l,\lambda }(f)=\sum_{j=1}^{n}\lambda _{j}f(\mathbf{x}^{j}),\quad f\in
C(X).
\end{equation*}

Since $l$ satisfies (5.10), for each function $g_{i}\in \mathcal{R}_{i}$, $%
i=1,\ldots ,k$, we have
\begin{equation*}
G_{l,\lambda }(g_{i})=\sum_{j=1}^{n}\lambda _{j}g_{i}(\mathbf{a}^{i}\cdot
\mathbf{x}^{j})=\sum_{s=1}^{r_{i}}\lambda _{i_{s}}g_{i}(\mathbf{a}^{i}\cdot
\mathbf{x}^{i_{s}}),\eqno(5.12)
\end{equation*}%
where $r_{i}\leq k$. That is, for each set $\mathcal{R}_{i}$, $G_{l,\lambda
} $ can be reduced to a functional defined with the help of not more than $k$
points of the semicycle $l$.

Consider the operator
\begin{equation*}
A:\mathcal{U}\rightarrow C(X),\quad A[(g_{1},\ldots ,g_{k})]=g_{1}+\cdots
+g_{k}.
\end{equation*}%
Clearly, $A$ is a linear continuous operator with the norm $\Vert A\Vert =1$%
. Besides, since $\mathcal{R}_{1}+\mathcal{R}_{2}+...+\mathcal{R}_{k}=C(X)$,
$A$ is a surjection. Consider also the conjugate operator
\begin{equation*}
A^{\ast }:C(X)^{\ast }\rightarrow \mathcal{U}^{\ast },~A^{\ast
}[H]=(F_{1},\ldots ,F_{k}),
\end{equation*}%
where $F_{i}(g_{i})=H(g_{i})$, for any $g_{i}\in \mathcal{R}_{i}$, $%
i=1,\ldots ,k$. Set $A^{\ast }[G_{l,\lambda }]=(G_{1},\ldots ,G_{k})$. From
(5.12) it follows that
\begin{equation*}
|G_{i}(g_{i})|=|G_{l,\lambda }(g_{i})|\leq \Vert g_{i}\Vert
\sum_{s=1}^{r_{i}}|\lambda _{i_{s}}|\leq kq\Vert g_{i}\Vert ,\quad
i=1,\ldots ,k,
\end{equation*}%
Therefore,
\begin{equation*}
\Vert G_{i}\Vert \leq kq,\quad i=1,\ldots ,k.
\end{equation*}%
From (5.11) we obtain that
\begin{equation*}
\Vert A^{\ast }[G_{l,\lambda }]\Vert =\Vert (G_{1},\ldots ,G_{k})\Vert \leq
kq.\eqno(5.13)
\end{equation*}%
Since $A$ is a surjection, there exists a positive real number $\delta $
such that
\begin{equation*}
\Vert A^{\ast }[H]\Vert >\delta \Vert H\Vert
\end{equation*}%
for any functional $H\in C(X)^{\ast }$(see \cite[p.100]{122}). Taking into account that $%
\Vert G_{l,\lambda }\Vert =\sum_{j=1}^{n}|\lambda _{j}|$, for the functional
$G_{l,\lambda }$ we have
\begin{equation*}
\Vert A^{\ast }[G_{l,\lambda }]\Vert >\delta \sum_{j=1}^{n}|\lambda _{j}|.%
\eqno(5.14)
\end{equation*}%
It follows from (5.13) and (5.14) that
\begin{equation*}
\delta <\frac{kq}{\sum_{j=1}^{n}|\lambda _{j}|}.
\end{equation*}%
The last inequality shows that $n$ (the length of the arbitrarily chosen $q$%
-semicycle $l$) cannot be as great as possible, otherwise $\delta =0$. This
simply means that there must be some positive integer bounding the lengths
of all $q$-semicycles in $X$.

It remains to show that there are no cycles in $X$. Indeed, if $l=(\mathbf{x}%
^{1},\ldots ,\mathbf{x}^{n})$ is a cycle in $X$ and $\lambda =(\lambda
_{1},\ldots ,\lambda _{n})$ is a vector associated with it, then the above
functional $G_{l,\lambda }$ annihilates all functions from $\mathcal{R}_{1}+%
\mathcal{R}_{2}+...+\mathcal{R}_{k}$. On the other hand, $G_{l,\lambda
}(f)=\sum_{j=1}^{n}|\lambda _{j}|\neq 0$ for a continuous function $f$ on $X$
satisfying the conditions $f(\mathbf{x}^{j})=1$ if $\lambda _{j}>0$ and $f(%
\mathbf{x}^{j})=-1$ if $\lambda _{j}<0$, $j=1,\ldots ,n$. This implies that $%
\mathcal{R}_{1}+\mathcal{R}_{2}+...+\mathcal{R}_{k}\neq C\left( X\right) $.
Since $\mathcal{M}_{X}(\sigma ;W,\mathbb{R})\subseteq \mathcal{R}_{1}+%
\mathcal{R}_{2}+...+\mathcal{R}_{k}$, we obtain that $\mathcal{M}_{X}(\sigma
;W,\mathbb{R})\neq C\left( X\right) $ on the contrary to our assumption.
\end{proof}

\bigskip

\textbf{Remark 5.2.} Assume $\mathcal{M}_{X}(\sigma ;W,\mathbb{R})$ is dense
in $C(X).$ Is it necessarily closed? Theorem 5.6 may describe cases when it
is not. For example, let $\mathbf{a}^{1}=(1;-1),\ \mathbf{a}^{2}=(1;1),$ $%
W=\{\mathbf{a}^{1},\mathbf{a}^{2}\}$ and $\sigma $ be any continuous,
bounded and nonconstant function, which has a limit at infinity. Consider
the set
\begin{eqnarray*}
X &=&\{(2;\frac{2}{3}),(\frac{2}{3};\frac{2}{3}),(0;0),(1;1),(1+\frac{1}{2}%
;1-\frac{1}{2}),(1+\frac{1}{2}+\frac{1}{4};1-\frac{1}{2}+\frac{1}{4}), \\
&&(1+\frac{1}{2}+\frac{1}{4}+\frac{1}{8};1-\frac{1}{2}+\frac{1}{4}-\frac{1}{8%
}),...\}.
\end{eqnarray*}

It is clear that $X$ is a compact set with all its orbits closed. (In fact,
there is only one orbit, which coincides with $X$). Hence, by Theorem 5.4, $%
\overline{\mathcal{M}_{X}(\sigma ;W,\mathbb{R})}=C(X).$ But by Theorem 5.6, $%
\mathcal{M}_{X}(\sigma ;W,\mathbb{R})\neq C(X).$ Therefore, the set $%
\mathcal{M}_{X}(\sigma ;W,\mathbb{R})$ is not closed in $C(X).$

\bigskip

\section{Two hidden layer neural networks}

A single hidden layer perceptron is able to approximate a given data with
any degree of accuracy. But in applications it is necessary to define how
many neurons one should take in a hidden layer. The more the number of
neurons, the more the probability of the network to give precise results.
Unfortunately, practicality decreases with the increase of the number of
neurons in the hidden layer. In other words, single hidden layer perceptrons
are not always effective if the number of neurons in the hidden layer is
prescribed. In this section, we show that this phenomenon is no longer true
for perceptrons with two hidden layers. We prove that a two hidden layer
neural network with $d$ inputs, $d$ neurons in the first hidden layer, $2d+2$
neurons in the second hidden layer and with a specifically constructed
sigmoidal and infinitely differentiable activation function can approximate
any continuous multivariate function with arbitrary accuracy.

\subsection{Relation of the Kolmogorov superposition theorem to two hidden
layer neural networks}

Note that if $r$ is fixed in (5.1), then the set

\begin{equation*}
\mathcal{M}_{r}(\sigma )=\left\{ \sum_{i=1}^{r}c_{i}\sigma (\mathbf{w}^{i}%
\mathbf{\cdot x}-\theta _{i}):~c_{i},\theta _{i}\in \mathbb{R},\mathbf{w\in }%
\mathbb{R}^{d}\right\}
\end{equation*}%
is no longer dense in in the space $C(\mathbb{R}^{d})$ (in the topology of
uniform convergence on compact sets) for any activation function $\sigma $.
The set $\mathcal{M}_{r}(\sigma )$ will not be dense even if we variate
over all univariate continuous functions $\sigma $ (see \cite[Theorem 5.1]%
{95}). In the following, we will see that this property of single hidden
layer neural networks does not carry over to networks with more than one
hidden layer.

A two hidden layer network is defined by iteration of the single hidden layer neural network
model. The output of two hidden layer perceptron with $r$ units in the first
layer, $s$ units in the second layer and the input $x=(x_{1},...,x_{d})$ is
\begin{equation*}
\sum_{i=1}^{s}d_{i}\sigma \left( \sum_{j=1}^{r}c_{ij}\sigma (\mathbf{w}%
^{ij}\cdot \mathbf{x-}\theta _{ij})-\gamma _{i}\right) .
\end{equation*}%
Here $d_{i},c_{ij},\theta _{ij},\gamma _{i}$ are real numbers, $\mathbf{w}%
^{ij}$ are vectors of $\mathbb{R}^{d}$ and $\sigma $ is a fixed univariate
function.

In many applications, it is convenient to take the activation function $%
\sigma $ as a \textit{sigmoidal function} which is defined as
\begin{equation*}
\lim_{t\rightarrow -\infty }\sigma (t)=0\quad \text{ and }\quad
\lim_{t\rightarrow +\infty }\sigma (t)=1.
\end{equation*}%
The literature on neural networks abounds with the use of such functions and
their superpositions. The following are typical examples of sigmoidal
functions:
\begin{align*}
\sigma (t)& =\frac{1}{1+e^{-t}} & & \text{(the squashing function),} \\
\sigma (t)& =%
\begin{cases}
0, & t\leq -1, \\
\dfrac{t+1}{2}, & -1\leq t\leq 1, \\
1, & t\geq 1%
\end{cases}
& & \text{(the piecewise linear function),} \\
\sigma (t)& =\frac{1}{\pi }\arctan t+\frac{1}{2} & & \text{(the arctan
sigmoid function),} \\
\sigma (t)& =\frac{1}{\sqrt{2\pi }}\int\limits_{-\infty }^{t}e^{-x^{2}/2}dx
& & \text{(the Gaussian function).}
\end{align*}

In this section, we prove that there exists a two hidden layer neural
network model with $d$ units in the first layer and $2d+2$ units in the
second layer such that it has the ability to approximate any $d$-variable
continuous function with arbitrary accuracy. As an activation function for
this model we take a specific sigmoidal function. The idea behind the proof
of this result is very much connected to the Kolmogorov superposition
theorem (see Section 4.1). This theorem has been much discussed in neural
network literature (see, e.g., \cite{119}). In our opinion, the most
remarkable application of the Kolmogorov superposition theorem to neural
networks was given by Maiorov and Pinkus \cite{99}. They showed that there
exists a sigmoidal, strictly increasing, analytic activation function, for
which a fixed number of units in both hidden layers are sufficient to
approximate arbitrarily well any continuous multivariate function. Namely,
the authors of \cite{99} proved the following theorem.

\bigskip

\textbf{Theorem 5.7 (Maiorov and Pinkus \cite{99}).} \textit{There exists an
activation function $\sigma $ which is analytic, strictly increasing and
sigmoidal and has the following property: For any $f\in C[0,1]^{d}$ and $%
\varepsilon >0,$ there exist constants $d_{i},$ $c_{ij},$ $\theta _{ij},$ $%
\gamma _{i}$, and vectors $\mathbf{w}^{ij}\in \mathbb{R}^{d}$ for which}

\begin{equation*}
\left\vert f(\mathbf{x})-\sum_{i=1}^{6d+3}d_{i}\sigma \left(
\sum_{j=1}^{3d}c_{ij}\sigma (\mathbf{w}^{ij}\cdot \mathbf{x-}\theta
_{ij})-\gamma _{i}\right) \right\vert <\varepsilon \eqno(5.15)
\end{equation*}%
\textit{for all $\mathbf{x}=(x_{1},...,x_{d})\in \lbrack 0,1]^{d}.$}

\bigskip

This theorem is based on the following version of the Kolmogorov
superposition theorem given by Lorentz \cite{98} and Sprecher \cite{128}.

\bigskip

\textbf{Theorem 5.8 (Kolmogorov's superposition theorem).} \textit{For the
unit cube $\mathbb{I}^{d},~\mathbb{I}=[0,1],~d\geq 2,$ there exists
constants $\lambda _{q}>0,$ $q=1,...,d,$ $\sum_{q=1}^{d}\lambda _{q}=1,$ and
nondecreasing continuous functions $\phi _{p}:[0,1]\rightarrow \lbrack 0,1],$
$p=1,...,2d+1,$ such that every continuous function $f:\mathbb{I}%
^{d}\rightarrow \mathbb{R}$ admits the representation}
\begin{equation*}
f(x_{1},...x_{d})=\sum_{p=1}^{2d+1}g\left( \sum_{q=1}^{d}\lambda _{q}\phi
_{p}(x_{q})\right) \eqno(5.16)
\end{equation*}%
\textit{for some $g\in C[0,1]$ depending on $f.$}

\bigskip

In the next subsection, using the general ideas developed in \cite{99}, we
show that the bounds of units in hidden layers in (5.15) may be chosen even
equal to the bounds in the Kolmogorov superposition theorem. More precisely,
these bounds can be taken as $2d+2$ and $d$ instead of $6d+3$ and $3d$. To
attain this purpose, we change the ``analyticity" of $\sigma $ to ``infinite
differentiability". In addition, near infinity we assume that $\sigma $ is ``$%
\lambda $-strictly increasing" instead of being ``strictly increasing".

\bigskip

\subsection{The main result}

\smallskip

We begin this subsection with a definition of a \textit{$\lambda$-monotone
function}. Let $\lambda $ be any nonnegative number. A real function $f$
defined on $(a,b)$ is called \textit{$\lambda $-increasing} (\textit{$\lambda $-decreasing})
if there exists an increasing (decreasing) function $u:(a,b)\rightarrow
\mathbb{R}$ such that $\left\vert f(x)-u(x)\right\vert \leq \lambda ,$ for
all $x\in (a,b)$. If $u$ is strictly increasing (or strictly decreasing),
then the above function $f$ is called a \textit{$\lambda $-strictly increasing} (or \textit{$%
\lambda $-strictly decreasing}) function. Clearly, $0$-monotonicity coincides
with the usual concept of monotonicity and a $\lambda _{1}$-monotone
function is $\lambda _{2}$-monotone if $\lambda _{1}\leq \lambda _{2}$. It
is also clear from the definition that a $\lambda $-monotone function
behaves like a usual monotone function as $\lambda $ gets very small.

Our purpose is to prove the following theorem.

\bigskip

\textbf{Theorem 5.9.} \textit{For any positive numbers $\alpha $ and $%
\lambda $,\textit{\ }there exists a $C^{\infty }(\mathbb{R}),$
sigmoidal activation function $\sigma :$ $\mathbb{R\rightarrow R}$ which is
strictly increasing on $(-\infty ,\alpha )$, $\lambda$-strictly
increasing on $[\alpha ,+\infty )$, and satisfies the following property:
For any $f\in C[0,1]^{d}$ and $\varepsilon >0,$ there exist constants $%
d_{p}, $ $c_{pq},$ $\theta _{pq},$ $\gamma _{p}$, and vectors $\mathbf{w}%
^{pq}\in \mathbb{R}^{d}$ for which}

\begin{equation*}
\left\vert f(\mathbf{x})-\sum_{p=1}^{2d+2}d_{p}\sigma \left(
\sum_{q=1}^{d}c_{pq}\sigma (\mathbf{w}^{pq}\cdot \mathbf{x-}\theta
_{pq})-\gamma _{p}\right) \right\vert <\varepsilon \eqno(5.17)
\end{equation*}%
\textit{for all $\mathbf{x}=(x_{1},...,x_{d})\in \lbrack 0,1]^{d}.$}

\bigskip

\begin{proof} Let $\alpha $ be any positive number. Divide the interval $%
[\alpha ,+\infty )$ into the segments $[\alpha ,2\alpha ],$ $[2\alpha
,3\alpha ],...$. Let $h(t)$ be any strictly increasing, infinitely
differentiable function on $[\alpha ,+\infty )$ with the properties

\bigskip

1) $0<h(t)<1$ for all $t\in \lbrack \alpha ,+\infty )$;

2) $1-h(\alpha )\leq \lambda ;$

3) $h(t)\rightarrow 1,$ as $t\rightarrow +\infty .$

\bigskip

The existence of a strictly increasing smooth function satisfying these
properties is easy to verify. Note that from conditions (1)-(3) it follows
that any function $f(t)$ satisfying the inequality $h(t)<f(t)<1$ for all $%
t\in \lbrack \alpha ,+\infty ),$ is $\lambda $-strictly increasing and $%
f(t)\rightarrow 1,$ as $t\rightarrow +\infty .$

We are going to construct $\sigma $ obeying the required properties in
stages. Let $\{u_{n}(t)\}_{n=1}^{\infty }$ be the sequence of all
polynomials with rational coefficients defined on $[0,1].$ First, we define $%
\sigma $ on the closed intervals $[(2m-1)\alpha ,2m\alpha ],$ $m=1,2,...$,
as the function
\begin{equation*}
\sigma (t)=a_{m}+b_{m}u_{m}(\frac{t}{\alpha }-2m+1),\text{ }t\in \lbrack
(2m-1)\alpha ,2m\alpha ],\eqno(5.18)
\end{equation*}%
or equivalently,
\begin{equation*}
\sigma (\alpha t+(2m-1)\alpha )=a_{m}+b_{m}u_{m}(t),\text{ }t\in \lbrack
0,1],\eqno(5.19)
\end{equation*}%
where $a_{m}$ and $b_{m}\neq 0$ are appropriately chosen constants. These
constants are determined from the condition
\begin{equation*}
h(t)<\sigma (t)<1,\eqno(5.20)
\end{equation*}%
\bigskip for all $t\in \lbrack (2m-1)\alpha ,2m\alpha ].$ There is a simple
procedure for determining a suitable pair of $a_{m}$ and $b_{m}$. Indeed,
let
\begin{equation*}
M=\max h(t)\text{, }A_{1}=\min u_{m}(\frac{t}{\alpha }-2m+1)\text{, }%
A_{2}=\max u_{m}(\frac{t}{\alpha }-2m+1),
\end{equation*}%
where in all the above $\max $ and $\min $, the variable $t$ runs over the
closed interval $[(2m-1)\alpha ,2m\alpha ].$ Note that $M<1$. If $%
A_{1}=A_{2} $ (that is, if the function $u_{m}$ is constant on $[0,1]$),
then we can set $\sigma (t)=(1+M)/2$ and easily find a suitable pair of $%
a_{m}$ and $b_{m}$ from (5.18). Let now $A_{1}\neq A_{2}$ and $y=a+bx,$ $%
b\neq 0,$ be a linear function mapping the segment $[A_{1},A_{2}]$ into $%
(M,1).$ Then it is enough to take $a_{m}=a$ and $b_{m}=b.$

At the second stage we define $\sigma $ on the intervals $[2m\alpha
,(2m+1)\alpha ],$ $m=1,2,...,$ so that it is in $C^{\infty }(\mathbb{R})$
and satisfies the inequality (5.20). Finally, in all of $(-\infty ,\alpha )$
we define $\sigma $ while maintaining the $C^{\infty }$ strict monotonicity
property, and also in such a way that $\lim_{t\rightarrow -\infty }\sigma
(t)=0.$ We obtain from the properties of $h$ and the condition (5.20) that $%
\sigma (t)$ is a $\lambda $-strictly increasing function on the interval $%
[\alpha ,+\infty )$ and $\sigma (t)\rightarrow 1$, as $t\rightarrow +\infty
. $

From the above construction of $\sigma $, that is, from (5.19) it follows
that for each $m=1,2,...,$ there exists numbers $A_{m}$,$\ B_{m}$ and $r_{m}$
such that
\begin{equation*}
u_{m}(t)=A_{m}\sigma (\alpha t-r_{m})-B_{m},\eqno(5.21)
\end{equation*}%
where $A_{m}\neq 0.$

Let $f$ be any continuous function on the unit cube $[0,1]^{d}.$ By the
Kolmogorov superposition theorem the expansion (5.16) is valid for $f.$ For
the exterior continuous univariate function $g(t)$ in (5.16) and for any $%
\varepsilon >0$ there exists a polynomial $u_{m}(t)$ of the above form such
that
\begin{equation*}
\left\vert g(t)-u_{m}(t)\right\vert <\frac{\varepsilon }{2(2d+1)},
\end{equation*}%
for all $t\in \lbrack 0,1].$ This together with (5.21) means that
\begin{equation*}
\left\vert g(t)-[a\sigma (\alpha t-r)-b]\right\vert <\frac{\varepsilon }{%
2(2d+1)},\eqno(5.22)
\end{equation*}%
for some $a,b,r\in \mathbb{R}$ and all $t\in \lbrack 0,1].$

Substituting (5.22) in (5.16) we obtain that
\begin{equation*}
\left\vert f(x_{1},...,x_{d})-\sum_{p=1}^{2d+1}\left( a\sigma \left( \alpha
\cdot \sum_{q=1}^{d}\lambda _{q}\phi _{p}(x_{q})-r\right) -b\right)
\right\vert <\frac{\varepsilon }{2}\eqno(5.23)
\end{equation*}%
for all $(x_{1},...,x_{d})\in \lbrack 0,1]^{d}.$

For each $p\in \{1,2,...,2d+1\}$ and $\delta >0$ there exist constants $%
a_{p},b_{p}$ and $r_{p}$ such that
\begin{equation*}
\left\vert \phi _{p}(x_{q})-[a_{p}\sigma (\alpha
x_{q}-r_{p})-b_{p}]\right\vert <\delta ,\eqno(5.24)
\end{equation*}%
for all $x_{q}\in \lbrack 0,1].$ Since $\lambda _{q}>0,$ $q=1,...,d,$ $%
\sum_{q=1}^{d}\lambda _{q}=1,$ it follows from (5.24) that

\begin{equation*}
\left\vert \sum_{q=1}^{d}\lambda _{q}\phi _{p}(x_{q})-\left[
\sum_{q=1}^{d}\lambda _{q}a_{p}\sigma (\alpha x_{q}-r_{p})-b_{p}\right]
\right\vert <\delta ,\eqno(5.25)
\end{equation*}%
for all $(x_{1},...,x_{d})\in \lbrack 0,1]^{d}.$

Now since the function $a\sigma (\alpha t-r)$ is uniformly
continuous on every closed interval, we can choose $\delta $ sufficiently
small and obtain from (5.25) that
\[
\left\vert \sum_{p=1}^{2d+1}a\sigma \left( \alpha
\sum_{q=1}^{d}\lambda _{q}\phi _{p}(x_{q})-r\right) \
-\sum_{p=1}^{2d+1}a\sigma \left( \alpha \left[ \sum_{q=1}^{d}\lambda
_{q}a_{p}\sigma (\alpha x_{q}-r_{p})-b_{p}\right] -r\right) \right\vert
\]
\[
< \frac{\varepsilon }{2}.
\]
This inequality may be rewritten as
\begin{equation*}
\left\vert \sum_{p=1}^{2d+1}a\sigma \left( \alpha
\sum_{q=1}^{d}\lambda _{q}\phi _{p}(x_{q})-r\right)
-\sum_{p=1}^{2d+1}d_{p}\sigma \left( \sum_{q=1}^{d}c_{pq}\sigma (\mathbf{w}%
^{pq}\cdot \mathbf{x}-\theta _{pq})-\gamma _{p}\right) \right\vert <\frac{%
\varepsilon }{2}.\eqno(5.26)
\end{equation*}%
From (5.23) and (5.26) it follows that
\begin{equation*}
\left\vert f(\mathbf{x})-\left[ \sum_{p=1}^{2d+1}d_{p}\sigma \left(
\sum_{q=1}^{d}c_{pq}\sigma (\mathbf{w}^{pq}\cdot \mathbf{x-}\theta
_{pq})-\gamma _{p}\right) -s\right] \right\vert <\varepsilon ,\eqno(5.27)
\end{equation*}%
where $s=(2d+1)b$. Since the constant $s$ can be written in the form
\begin{equation*}
s=d\sigma \left( \sum_{q=1}^{d}c_{q}\sigma (\mathbf{w}^{q}\cdot \mathbf{x-}%
\theta _{q})-\gamma \right) ,
\end{equation*}%
from (5.27) we finally obtain the validity of (5.17).
\end{proof}

\textbf{Remark 5.3.} It is easily seen in the proof of Theorem 5.9 that all
the weights $\mathbf{w}^{ij}$ are fixed (see (5.26)). Namely, $\mathbf{w}%
^{ij}=\alpha \mathbf{e}^{j},$ for all $i=1,...,2d+2,$ $j=1,...,d,$ where $%
\mathbf{e}^{j}$ is the $j$-th coordinate vector of the space $\mathbb{R}^{d}$%
.

\bigskip

The next theorem follows from Theorem 5.9 easily, since the Kolmogorov
superposition theorem is valid for all compact sets of $\mathbb{R}^{d}$.

\bigskip

\textbf{Theorem 5.10.} \textit{Let $Q$ be a compact set in $\mathbb{R}^{d}.$
For any numbers $\alpha \in \mathbb{R}$ and $\lambda >0,$ there exists a
\textit{$C^{\infty }(\mathbb{R}),$ sigmoidal activation function }$\sigma :$
$\mathbb{R\rightarrow R}$ which is strictly increasing on $(-\infty ,\alpha
) $, $\lambda$-strictly increasing on $[\alpha ,+\infty )$, and
satisfies the following property: For any $f\in C(Q)$ and $\varepsilon >0$
there exist real numbers $d_{i},$ $c_{ij},$ $\theta _{ij}$, $\gamma _{i},$
and vectors $\mathbf{w}^{ij}\in \mathbb{R}^{d}$ for which}
\begin{equation*}
\left\vert f(\mathbf{x})-\sum_{i=1}^{2d+2}d_{i}\sigma \left(
\sum_{j=1}^{d}c_{ij}\sigma (\mathbf{w}^{ij}\cdot \mathbf{x-}\theta
_{ij})-\gamma _{i}\right) \right\vert <\varepsilon
\end{equation*}%
\textit{for all $\mathbf{x}=(x_{1},...,x_{d})\in Q.$}

\bigskip

\textbf{Remark 5.4.} In some literature, a single hidden layer perceptron is
defined as the function
\begin{equation*}
\sum_{i=1}^{r}c_{i}\sigma (\mathbf{w}^{i}\mathbf{\cdot x}-\theta _{i})-c_{0}.
\end{equation*}%
A two hidden layer network then takes the form
\begin{equation*}
\sum_{i=1}^{s}d_{i}\sigma \left( \sum_{j=1}^{r}c_{ij}\sigma (\mathbf{w}%
^{ij}\cdot \mathbf{x-}\theta _{ij})-\gamma _{i}\right) -d_{0}.\eqno(5.28)
\end{equation*}%
The proof of Theorem 5.9 shows that for networks of type (5.28) the theorem
is valid if we take $2d+1$ neurons in the second hidden layer (instead of $%
2d+2$ neurons as above). That is, there exist networks of type (5.28) having
the universal approximation property and for which the number of units in
the hidden layers is equal to the number of summands in the Kolmogorov
superposition theorem.

\bigskip

\textbf{Remark 5.5.} It is known that the $2d+1$ in the Kolmogorov
superposition theorem is minimal (see Sternfeld \cite{130}). Thus it is
doubtful if the number of neurons in Theorems 5.9 and 5.10 can be reduced.

\bigskip

\textbf{Remark 5.6.} Inequality (5.22) shows that single hidden layer neural
networks of the form (5.28) with the activation function $\sigma$ and with
only one neuron in the hidden layer can approximate any continuous function
on the interval $[0,1]$ with arbitrary precision. Since the number $b$ in
(5.22) can always be written as $b=a_1\sigma (0 \cdot t-r_1)$ for some $a_1$
and $r_1$, we see that two neurons in the hidden layer are sufficient for
traditional single hidden layer neural networks with the activation function
$\sigma$ to approximate continuous functions on $[0,1]$. Applying the linear
transformation $x=a+(b-a)t$ it can be proven that the same argument holds
for any interval $[a,b]$.

\bigskip

\section{Construction of a universal sigmoidal function}

In the preceding section, we considered two theorems (Theorem 5.7 of Maiorov
and Pinkus, and Theorem 5.9) on the approximation capabilities of the MLP
model of neural networks with a prescribed number of hidden neurons. Note
that both results are more theoretical than practical, as they indicate only
the existence of the corresponding activation functions.

In this section, we construct algorithmically a smooth, sigmoidal, almost
monotone activation function $\sigma $ providing approximation to an
arbitrary continuous function within any degree of accuracy. This algorithm
is implemented in a computer program, which computes the value of $\sigma $
at any reasonable point of the real axis.

\bigskip

\subsection{A construction algorithm}

In this subsection, we construct algorithmically a sigmoidal function $%
\sigma $ which we use in our results in Section 5.3.3.

To start with the construction of $\sigma$, assume that we are given a
closed interval $[a, b]$ and a sufficiently small real number $\lambda$. We
construct $\sigma$ algorithmically, based on two numbers, namely $\lambda$
and $d := b - a$. The following steps describe the algorithm.

\textit{Step 1.} Introduce the function
\begin{equation*}
h(x) := 1 - \frac{\min\{1/2, \lambda\}}{1 + \log(x - d + 1)}.
\end{equation*}
Note that this function is strictly increasing on the real line and
satisfies the following properties:

\begin{enumerate}
\item $0 < h(x) < 1$ for all $x \in [d, +\infty)$;

\item $1 - h(d) \le \lambda$;

\item $h(x) \to 1$, as $x \to +\infty$.
\end{enumerate}

We want to construct $\sigma $ satisfying the inequalities
\begin{equation*}
h(x)<\sigma (x)<1\eqno(5.29)
\end{equation*}%
for $x\in \lbrack d,+\infty )$. Then our $\sigma $ will tend to $1$ as $x$
tends to $+\infty $ and obey the inequality
\begin{equation*}
|\sigma (x)-h(x)|\leq \lambda ,
\end{equation*}%
i.e., it will be a $\lambda $-increasing function.

\textit{Step 2.} Before proceeding to the construction of $\sigma$, we need
to enumerate the monic polynomials with rational coefficients. Let $q_n$ be
the Calkin--Wilf sequence (see~\cite{CW00}). Then we can enumerate all the
rational numbers by setting
\begin{equation*}
r_0 := 0, \quad r_{2n} := q_n, \quad r_{2n-1} := -q_n, \ n = 1, 2, \dots.
\end{equation*}
Note that each monic polynomial with rational coefficients can uniquely be
written as $r_{k_0} + r_{k_1} x + \ldots + r_{k_{l-1}} x^{l-1} + x^l$, and
each positive rational number determines a unique finite continued fraction
\begin{equation*}
[m_0; m_1, \ldots, m_l] := m_0 + \dfrac1{m_1 + \dfrac1{m_2 + \dfrac1{\ddots
+ \dfrac1{m_l}}}}
\end{equation*}
with $m_0 \ge 0$, $m_1, \ldots, m_{l-1} \ge 1$ and $m_l \ge 2$. We now
construct a bijection between the set of all monic polynomials with rational
coefficients and the set of all positive rational numbers as follows. To the
only zeroth-degree monic polynomial 1 we associate the rational number 1, to
each first-degree monic polynomial of the form $r_{k_0} + x$ we associate
the rational number $k_0 + 2$, to each second-degree monic polynomial of the
form $r_{k_0} + r_{k_1} x + x^2$ we associate the rational number $[k_0; k_1
+ 2] = k_0 + 1 / (k_1 + 2)$, and to each monic polynomial
\begin{equation*}
r_{k_0} + r_{k_1} x + \ldots + r_{k_{l-2}} x^{l-2} + r_{k_{l-1}} x^{l-1} +
x^l
\end{equation*}
of degree $l \ge 3$ we associate the rational number $[k_0; k_1 + 1, \ldots,
k_{l-2} + 1, k_{l-1} + 2]$. In other words, we define $u_1(x) := 1$,
\begin{equation*}
u_n(x) := r_{q_n-2} + x
\end{equation*}
if $q_n \in \mathbb{Z}$,
\begin{equation*}
u_n(x) := r_{m_0} + r_{m_1-2} x + x^2
\end{equation*}
if $q_n = [m_0; m_1]$, and
\begin{equation*}
u_n(x) := r_{m_0} + r_{m_1-1} x + \ldots + r_{m_{l-2}-1} x^{l-2} +
r_{m_{l-1}-2} x^{l-1} + x^l
\end{equation*}
if $q_n = [m_0; m_1, \ldots, m_{l-2}, m_{l-1}]$ with $l \ge 3$. For example,
the first few elements of this sequence are
\begin{equation*}
1, \quad x^2, \quad x, \quad x^2 - x, \quad x^2 - 1, \quad x^3, \quad x - 1,
\quad x^2 + x, \quad \ldots.
\end{equation*}

\textit{Step 3.} We start with constructing $\sigma$ on the intervals $%
[(2n-1)d, 2nd]$, $n = 1, 2, \ldots$. For each monic polynomial $u_n(x) =
\alpha_0 + \alpha_1 x + \ldots + \alpha_{l-1} x^{l-1} + x^l$, set
\begin{equation*}
B_1 := \alpha_0 + \frac{\alpha_1-|\alpha_1|}{2} + \ldots + \frac{%
\alpha_{l-1} - |\alpha_{l-1}|}{2}
\end{equation*}
and
\begin{equation*}
B_2 := \alpha_0 + \frac{\alpha_1+|\alpha_1|}{2} + \ldots + \frac{%
\alpha_{l-1} + |\alpha_{l-1}|}{2} + 1.
\end{equation*}
Note that the numbers $B_1$ and $B_2$ depend on $n$. To avoid complication
of symbols, we do not indicate this in the notation.

Introduce the sequence
\begin{equation*}
M_n := h((2n+1)d), \qquad n = 1, 2, \ldots.
\end{equation*}
Clearly, this sequence is strictly increasing and converges to $1$.

Now we define $\sigma $ as the function
\begin{equation*}
\sigma (x):=a_{n}+b_{n}u_{n}\left( \frac{x}{d}-2n+1\right) ,\quad x\in
\lbrack (2n-1)d,2nd],\eqno(5.30)
\end{equation*}%
where
\begin{equation*}
a_{1}:=\frac{1}{2},\qquad b_{1}:=\frac{h(3d)}{2},\eqno(5.31)
\end{equation*}%
and
\begin{equation*}
a_{n}:=\frac{(1+2M_{n})B_{2}-(2+M_{n})B_{1}}{3(B_{2}-B_{1})},\qquad b_{n}:=%
\frac{1-M_{n}}{3(B_{2}-B_{1})},\qquad n=2,3,\ldots .\eqno(5.32)
\end{equation*}

It is not difficult to notice that for $n>2$ the numbers $a_{n}$, $b_{n}$
are the coefficients of the linear function $y=a_{n}+b_{n}x$ mapping the
closed interval $[B_{1},B_{2}]$ onto the closed interval $%
[(1+2M_{n})/3,(2+M_{n})/3]$. Besides, for $n=1$, i.e. on the interval $%
[d,2d] $,
\begin{equation*}
\sigma (x)=\frac{1+M_{1}}{2}.
\end{equation*}%
Therefore, we obtain that
\begin{equation*}
h(x)<M_{n}<\frac{1+2M_{n}}{3}\leq \sigma (x)\leq \frac{2+M_{n}}{3}<1,\eqno%
(5.33)
\end{equation*}%
for all $x\in \lbrack (2n-1)d,2nd]$, $n=1$, $2$, $\ldots $.

\textit{Step 4.} In this step, we construct $\sigma $ on the intervals $%
[2nd,(2n+1)d]$, $n=1,2,\ldots $. For this purpose we use the \textit{smooth
transition function }%
\begin{equation*}
\beta _{a,b}(x):=\frac{\widehat{\beta }(b-x)}{\widehat{\beta }(b-x)+\widehat{%
\beta }(x-a)},
\end{equation*}%
where
\begin{equation*}
\widehat{\beta }(x):=%
\begin{cases}
e^{-1/x}, & x>0, \\
0, & x\leq 0.%
\end{cases}%
\end{equation*}%
Obviously, $\beta _{a,b}(x)=1$ for $x\leq a$, $\beta _{a,b}(x)=0$ for $x\geq
b$, and $0<\beta _{a,b}(x)<1$ for $a<x<b$.

Set
\begin{equation*}
K_n := \frac{\sigma(2nd) + \sigma((2n+1)d)}{2}, \qquad n = 1, 2, \ldots.
\end{equation*}
Note that the numbers $\sigma(2nd)$ and $\sigma((2n+1)d)$ have already been
defined in the previous step. Since both the numbers $\sigma(2nd)$ and $%
\sigma((2n+1)d)$ lie in the interval $(M_n, 1)$, it follows that $K_n \in
(M_n, 1)$.

First we extend $\sigma $ smoothly to the interval $[2nd,2nd+d/2]$. Take $%
\varepsilon :=(1-M_{n})/6$ and choose $\delta \leq d/2$ such that
\begin{equation*}
\left\vert a_{n}+b_{n}u_{n}\left( \frac{x}{d}-2n+1\right) -\left(
a_{n}+b_{n}u_{n}(1)\right) \right\vert \leq \varepsilon ,\quad x\in \lbrack
2nd,2nd+\delta ].\eqno(5.34)
\end{equation*}%
One can choose this $\delta $ as
\begin{equation*}
\delta :=\min \left\{ \frac{\varepsilon d}{b_{n}C},\frac{d}{2}\right\} ,
\end{equation*}%
where $C>0$ is a number satisfying $|u_{n}^{\prime }(x)|\leq C$ for $x\in
(1,1.5)$. For example, for $n=1$, $\delta $ can be chosen as $d/2$. Now
define $\sigma $ on the first half of the interval $[2nd,(2n+1)d]$ as the
function
\begin{equation*}
\sigma(x) :=K_{n}-\beta _{2nd,2nd+\delta }(x)
\end{equation*}
\begin{equation*}
\times \left(
K_{n}-a_{n}-b_{n}u_{n}\left( \frac{x}{d}-2n+1\right) \right) , x\in \left[
2nd,2nd+\frac{d}{2}\right].\eqno(5.35)
\end{equation*}

Let us prove that $\sigma (x)$ satisfies the condition~(5.29). Indeed, if $%
2nd+\delta \leq x\leq 2nd+d/2$, then there is nothing to prove, since $%
\sigma (x)=K_{n}\in (M_{n},1)$. If $2nd\leq x<2nd+\delta $, then $0<\beta
_{2nd,2nd+\delta }(x)\leq 1$ and hence from~(5.35) it follows that for each $%
x\in \lbrack 2nd,2nd+\delta )$, $\sigma (x)$ is between the numbers $K_{n}$
and $A_{n}(x):=a_{n}+b_{n}u_{n}\left( \frac{x}{d}-2n+1\right) $. On the
other hand, from~(5.34) we obtain that
\begin{equation*}
a_{n}+b_{n}u_{n}(1)-\varepsilon \leq A_{n}(x)\leq
a_{n}+b_{n}u_{n}(1)+\varepsilon ,
\end{equation*}%
which together with~(5.30) and~(5.33) yields $A_{n}(x)\in \left[ \frac{%
1+2M_{n}}{3}-\varepsilon ,\frac{2+M_{n}}{3}+\varepsilon \right] $ for $x\in
\lbrack 2nd,2nd+\delta )$. Since $\varepsilon =(1-M_{n})/6$, the inclusion $%
A_{n}(x)\in (M_{n},1)$ is valid. Now since both $K_{n}$ and $A_{n}(x)$
belong to $(M_{n},1)$, we finally conclude that
\begin{equation*}
h(x)<M_{n}<\sigma (x)<1,\quad \text{for }x\in \left[ 2nd,2nd+\frac{d}{2}%
\right] .
\end{equation*}

We define $\sigma $ on the second half of the interval in a similar way:
\begin{equation*}
\begin{split}
\sigma (x)& :=K_{n}-(1-\beta _{(2n+1)d-\overline{\delta },(2n+1)d}(x)) \\
& \times \left( K_{n}-a_{n+1}-b_{n+1}u_{n+1}\left( \frac{x}{d}-2n-1\right)
\right) ,\quad x\in \left[ 2nd+\frac{d}{2},(2n+1)d\right] ,
\end{split}%
\end{equation*}%
where
\begin{equation*}
\overline{\delta }:=\min \left\{ \frac{\overline{\varepsilon }d}{b_{n+1}%
\overline{C}},\frac{d}{2}\right\} ,\qquad \overline{\varepsilon }:=\frac{%
1-M_{n+1}}{6},\qquad \overline{C}\geq \sup_{\lbrack -0.5,0]}|u_{n+1}^{\prime
}(x)|.
\end{equation*}%
One can easily verify, as above, that the constructed $\sigma (x)$ satisfies
the condition~(5.29) on $[2nd+d/2,2nd+d]$ and
\begin{equation*}
\sigma \left( 2nd+\frac{d}{2}\right) =K_{n},\qquad \sigma ^{(i)}\left( 2nd+%
\frac{d}{2}\right) =0,\quad i=1,2,\ldots .
\end{equation*}

Steps 3 and 4 construct $\sigma$ on the interval $[d, +\infty)$.

\textit{Step 5.} On the remaining interval $(-\infty ,d)$, we define $\sigma
$ as
\begin{equation*}
\sigma (x):=\left( 1-\widehat{\beta }(d-x)\right) \frac{1+M_{1}}{2},\quad
x\in (-\infty ,d).
\end{equation*}%
It is not difficult to verify that $\sigma $ is a strictly increasing,
smooth function on $(-\infty ,d)$. Note also that $\sigma (x)\rightarrow
\sigma (d)=(1+M_{1})/2$, as $x$ tends to $d$ from the left and $\sigma
^{(i)}(d)=0$ for $i=1$, $2$, $\ldots $. This final step completes the
construction of $\sigma $ on the whole real line.

\bigskip

\subsection{Properties of the constructed sigmoidal function}

\begin{figure}[tbp]
\includegraphics[width=0.75\textwidth]{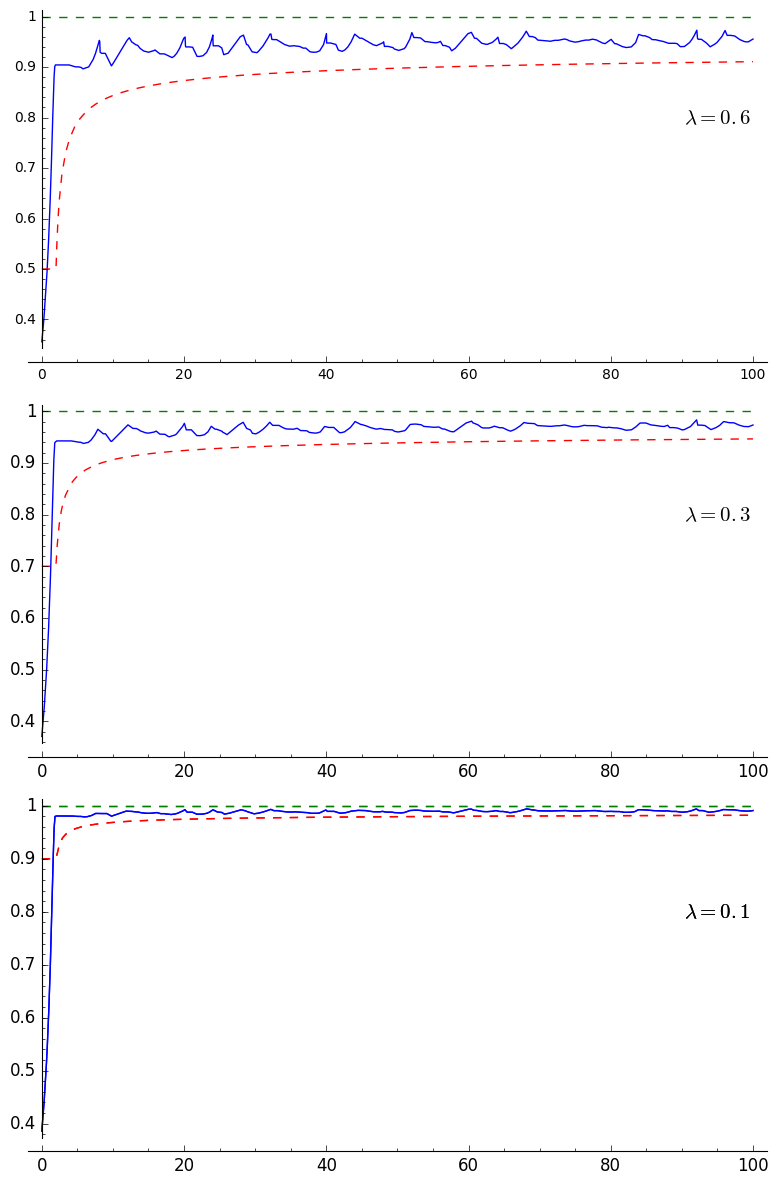}
\caption{The graph of $\protect\sigma$ on $[0, 100]$ ($d = 2$)}
\label{fig:sigma100}
\end{figure}

It should be noted that the above algorithm allows one to compute the
constructed $\sigma$ at any point of the real axis instantly. The code of
this algorithm is available at %
\url{http://sites.google.com/site/njguliyev/papers/monic-sigmoidal}. As a
practical example, we give here the graph of $\sigma$ (see Figure~\ref%
{fig:sigma}) and a numerical table (see Table~\ref{tbl:sigma}) containing
several computed values of this function on the interval $[0, 20]$. Figure~%
\ref{fig:sigma100} shows how the graph of $\lambda$-increasing function $%
\sigma$ changes on the interval $[0,100]$ as the parameter $\lambda$
decreases.

The above $\sigma$ obeys the following properties:

\begin{enumerate}
\item $\sigma$ is sigmoidal;

\item $\sigma \in C^{\infty}(\mathbb{R})$;

\item $\sigma$ is strictly increasing on $(-\infty, d)$ and $\lambda$%
-strictly increasing on $[d, +\infty)$;

\item $\sigma$ is easily computable in practice.
\end{enumerate}

All these properties are easily seen from the above exposition. But the
essential property of our sigmoidal function is its ability to approximate
an arbitrary continuous function using only a fixed number of translations
and scalings of $\sigma$. More precisely, only two translations and scalings
are sufficient. We formulate this important property as a theorem in the
next section.

\begin{figure}[tbp]
\includegraphics[width=1.0\textwidth]{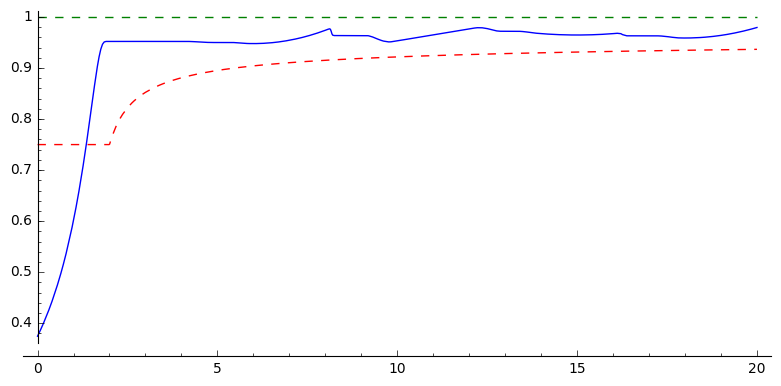}
\caption{The graph of $\protect\sigma$ on $[0, 20]$ ($d = 2$, $\protect%
\lambda = 1/4$)}
\label{fig:sigma}
\end{figure}

\begin{table}[tbp]
\caption{Some computed values of $\protect\sigma$ ($d = 2$, $\protect\lambda %
= 1/4$)}
\label{tbl:sigma}%
\begin{tabular}{|c|c|c|c|c|c|c|c|c|c|}
\hline
$t$ & $\sigma$ & $t$ & $\sigma$ & $t$ & $\sigma$ & $t$ & $\sigma$ & $t$ & $%
\sigma$ \\ \hline
$0.0$ & $0.37462$ & $4.0$ & $0.95210$ & $8.0$ & $0.97394$ & $12.0$ & $%
0.97662 $ & $16.0$ & $0.96739$ \\ \hline
$0.4$ & $0.44248$ & $4.4$ & $0.95146$ & $8.4$ & $0.96359$ & $12.4$ & $%
0.97848 $ & $16.4$ & $0.96309$ \\ \hline
$0.8$ & $0.53832$ & $4.8$ & $0.95003$ & $8.8$ & $0.96359$ & $12.8$ & $%
0.97233 $ & $16.8$ & $0.96309$ \\ \hline
$1.2$ & $0.67932$ & $5.2$ & $0.95003$ & $9.2$ & $0.96314$ & $13.2$ & $%
0.97204 $ & $17.2$ & $0.96307$ \\ \hline
$1.6$ & $0.87394$ & $5.6$ & $0.94924$ & $9.6$ & $0.95312$ & $13.6$ & $%
0.97061 $ & $17.6$ & $0.96067$ \\ \hline
$2.0$ & $0.95210$ & $6.0$ & $0.94787$ & $10.0$ & $0.95325$ & $14.0$ & $%
0.96739$ & $18.0$ & $0.95879$ \\ \hline
$2.4$ & $0.95210$ & $6.4$ & $0.94891$ & $10.4$ & $0.95792$ & $14.4$ & $%
0.96565$ & $18.4$ & $0.95962$ \\ \hline
$2.8$ & $0.95210$ & $6.8$ & $0.95204$ & $10.8$ & $0.96260$ & $14.8$ & $%
0.96478$ & $18.8$ & $0.96209$ \\ \hline
$3.2$ & $0.95210$ & $7.2$ & $0.95725$ & $11.2$ & $0.96727$ & $15.2$ & $%
0.96478$ & $19.2$ & $0.96621$ \\ \hline
$3.6$ & $0.95210$ & $7.6$ & $0.96455$ & $11.6$ & $0.97195$ & $15.6$ & $%
0.96565$ & $19.6$ & $0.97198$ \\ \hline
\end{tabular}%
\end{table}

\bigskip

\subsection{Theoretical results}

The following theorems are valid.

\bigskip

\textbf{Theorem 5.11.} \textit{Assume that $f$ is a continuous function on a
finite segment $[a,b]$ of $\mathbb{R}$ and $\sigma $ is the sigmoidal
function constructed in Section~5.3.1. Then for any sufficiently small $%
\varepsilon >0 $ there exist constants $c_{1}$, $c_{2}$, $\theta _{1}$ and $%
\theta _{2}$ such that}
\begin{equation*}
|f(x)-c_{1}\sigma (x-\theta _{1})-c_{2}\sigma (x-\theta _{2})|<\varepsilon
\end{equation*}%
\textit{for all $x\in \lbrack a,b]$.}

\bigskip

\begin{proof} Set $d:=b-a$ and divide the interval $[d,+\infty )$ into the segments
$[d,2d]$, $[2d,3d]$, $\ldots $. It follows from~(5.30) that
\begin{equation*}
\sigma (dx+(2n-1)d)=a_{n}+b_{n}u_{n}(x),\quad x\in \lbrack 0,1]\eqno(5.36)
\end{equation*}%
for $n=1$, $2$, $\ldots $. Here $a_{n}$ and $b_{n}$ are computed by~(5.31)
and~(5.32) for $n=1$ and $n>1$, respectively.

From~(5.36) it follows that for each $n=1$, $2$, $\ldots $,
\begin{equation*}
u_{n}(x)=\frac{1}{b_{n}}\sigma (dx+(2n-1)d)-\frac{a_{n}}{b_{n}}.\eqno(5.37)
\end{equation*}

Let now $g$ be any continuous function on the unit interval $[0,1]$. By the
density of polynomials with rational coefficients in the space of continuous
functions on any compact subset of $\mathbb{R}$, for any $\varepsilon >0$
there exists a polynomial $p(x)$ of the above form such that
\begin{equation*}
|g(x)-p(x)|<\varepsilon
\end{equation*}%
for all $x\in \lbrack 0,1]$. Denote by $p_{0}$ the leading coefficient of $p$%
. If $p_{0}\neq 0$ (i.e., $p\not\equiv 0$) then we define $u_{n}$ as $%
u_{n}(x):=p(x)/p_{0}$, otherwise we just set $u_{n}(x):=1$. In both cases
\begin{equation*}
|g(x)-p_{0}u_{n}(x)|<\varepsilon ,\qquad x\in \lbrack 0,1].
\end{equation*}%
This together with~(5.37) means that
\begin{equation*}
|g(x)-c_{1}\sigma (dx-s_{1})-c_{0}|<\varepsilon \eqno(5.38)
\end{equation*}%
for some $c_{0}$, $c_{1}$, $s_{1}\in \mathbb{R}$ and all $x\in \lbrack 0,1]$%
. Namely, $c_{1}=p_{0}/b_{n}$, $s_{1}=d-2nd$ and $c_{0}=p_{0}a_{n}/b_{n}$.
On the other hand, we can write $c_{0}=c_{2}\sigma (dx-s_{2})$, where $%
c_{2}:=2c_{0}/(1+h(3d))$ and $s_{2}:=-d$. Hence,
\begin{equation*}
|g(x)-c_{1}\sigma (dx-s_{1})-c_{2}\sigma (dx-s_{2})|<\varepsilon .\eqno(5.39)
\end{equation*}%
Note that~(5.39) is valid for the unit interval $[0,1]$. Using linear
transformation it is not difficult to go from $[0,1]$ to the interval $[a,b]$%
. Indeed, let $f\in C[a,b]$, $\sigma $ be constructed as above, and $%
\varepsilon $ be an arbitrarily small positive number. The transformed
function $g(x)=f(a+(b-a)x)$ is well defined on $[0,1]$ and we can apply the
inequality~(5.39). Now using the inverse transformation $x=(t-a)/(b-a)$, we
can write
\begin{equation*}
|f(t)-c_{1}\sigma (t-\theta _{1})-c_{2}\sigma (t-\theta _{2})|<\varepsilon
\end{equation*}%
for all $t\in \lbrack a,b]$, where $\theta _{1}=a+s_{1}$ and $\theta
_{2}=a+s_{2}$. The last inequality completes the proof.
\end{proof}

Since any compact subset of the real line is contained in a segment $[a,b]$,
the following generalization of Theorem 5.11 holds.

\bigskip

\textbf{Theorem 5.12.} \textit{Let $Q$ be a compact subset of the real line
and $d$ be its diameter. Let $\lambda $ be any positive number. Then one can
algorithmically construct a computable sigmoidal activation function $\sigma
\colon \mathbb{R}\rightarrow \mathbb{R}$, which is infinitely
differentiable, strictly increasing on $(-\infty ,d)$, $\lambda $-strictly
increasing on $[d,+\infty )$, and satisfies the following property: For any $%
f\in C(Q)$ and $\varepsilon >0$ there exist numbers $c_{1}$, $c_{2}$, $%
\theta _{1}$ and $\theta _{2}$ such that}
\begin{equation*}
|f(x)-c_{1}\sigma (x-\theta _{1})-c_{2}\sigma (x-\theta _{2})|<\varepsilon
\end{equation*}%
\textit{for all $x\in Q$.}

\bigskip

\textbf{Remark 5.7.} Theorems 5.11 and 5.12 show that single hidden layer
neural networks with the constructed sigmoidal activation function $\sigma $
and only two neurons in the hidden layer can approximate any continuous
univariate function. Moreover, in this case, one can fix the weights equal
to $1$. For the approximation of continuous multivariate functions two
hidden layer neural networks with $3d+2$ hidden neurons can be taken.
Namely, Theorem 5.9 (and hence Theorem 5.10) is valid with the constructed
in Section 5.3.1 activation function $\sigma$. Indeed, the proof of this
theorem shows that any activation function with the property (5.22)
suffices. But the activation function constructed
in Section 5.3.1 satisfies this property (see (5.38)).

\bigskip

\subsection{Numerical results}

We prove in Theorem 5.11 that any continuous function on $[a,b]$ can be
approximated arbitrarily well by single hidden layer neural networks with
the fixed weight $1$ and with only two neurons in the hidden layer. An
activation function $\sigma $ for such a network is constructed in
Section 5.3.1. We have seen from the proof that our approach is totally
constructive. One can evaluate the value of $\sigma $ at any point of the
real axis and draw its graph instantly, using the programming interface at
the URL shown at the beginning of Section 5.3.2. In the current subsection, we
demonstrate our result in various examples. For different error bounds we
find the parameters $c_{1}$, $c_{2}$, $\theta _{1}$ and $\theta _{2}$ in
Theorem 5.11. All computations were done in SageMath~\cite{Sage}. For
computations, we use the following algorithm, which works well for analytic
functions. Assume $f$ is a function, whose Taylor series around the point $%
(a+b)/2$ converges uniformly to $f$ on $[a,b]$, and $\varepsilon >0$.

\begin{enumerate}
\item Consider the function $g(t) := f(a + (b - a) t)$, which is
well-defined on $[0, 1]$;

\item Find $k$ such that the $k$-th Taylor polynomial
\begin{equation*}
T_k(x) := \sum_{i=0}^k \frac{g^{(i)}(1/2)}{i!} \left( x - \frac{1}{2}
\right)^i
\end{equation*}
satisfies the inequality $|T_k(x) - g(x)| \le \varepsilon / 2$ for all $x
\in [0, 1]$;

\item Find a polynomial $p$ with rational coefficients such that
\begin{equation*}
|p(x) - T_k(x)| \le \frac{\varepsilon}{2}, \qquad x \in [0, 1],
\end{equation*}
and denote by $p_0$ the leading coefficient of this polynomial;

\item If $p_0 \ne 0$, then find $n$ such that $u_n(x) = p(x) / p_0$.
Otherwise, set $n := 1$;

\item For $n=1$ and $n>1$ evaluate $a_{n}$ and $b_{n}$ by~(5.31) and~(5.32),
respectively;

\item Calculate the parameters of the network as
\begin{equation*}
c_1 := \frac{p_0}{b_n}, \qquad c_2 := \frac{2 p_0 a_n}{b_n (1 + h(3d))},
\qquad \theta_1 := b - 2 n (b - a), \qquad \theta_2 := 2 a - b;
\end{equation*}

\item Construct the network $\mathcal{N}=c_{1}\sigma (x-\theta
_{1})+c_{2}\sigma (x-\theta _{2}).$ Then $\mathcal{N}$ gives an $\varepsilon
$-approximation to $f.$
\end{enumerate}

In the sequel, we give four practical examples. To be able to make
comparisons between these examples, all the considered functions are given
on the same interval $[-1,1]$. First we select the polynomial function $%
f(x)=x^{3}+x^{2}-5x+3$ as a target function. We investigate the sigmoidal
neural network approximation to $f(x)$. This function was considered in
\cite{HH04} as well. Note that the authors of \cite{HH04} chose the sigmoidal
function as
\begin{equation*}
\sigma (x)=%
\begin{cases}
1, & \text{if }x\geq 0, \\
0, & \text{if }x<0,%
\end{cases}%
\end{equation*}%
and obtained the numerical results (see Table~\ref{tbl:Hahm}) for single
hidden layer neural networks with $8$, $32$, $128$, $532$ neurons in the
hidden layer (see also~\cite{CX10} for an additional constructive result
concerning the error of approximation in this example).

\begin{table}[tbp]
\caption{The Heaviside function as a sigmoidal function}
\label{tbl:Hahm}%
\begin{tabular}{|c|c|c|}
\hline
$N$ & Number of neurons ($2N^2$) & Maximum error \\ \hline
$2$ & $8$ & $0.666016$ \\ \hline
$4$ & $32$ & $0.165262$ \\ \hline
$8$ & $128$ & $0.041331$ \\ \hline
$16$ & $512$ & $0.010333$ \\ \hline
\end{tabular}%
\end{table}

As it is seen from the table, the number of neurons in the hidden layer
increases as the error bound decreases in value. This phenomenon is no
longer true for our sigmoidal function. Using Theorem 5.11, we can construct
explicitly a single hidden layer neural network model with only two neurons
in the hidden layer, which approximates the above polynomial with
arbitrarily given precision. Here by \textit{explicit construction} we mean
that all the network parameters can be computed directly. Namely, the
calculated values of these parameters are as follows: $c_{1}\approx
2059.373597$, $c_{2}\approx -2120.974727$, $\theta _{1}=-467$, and $\theta
_{2}=-3$. It turns out that for the above polynomial we have an exact
representation. That is, on the interval $[-1,1]$ we have the identity
\begin{equation*}
x^{3}+x^{2}-5x+3\equiv c_{1}\sigma (x-\theta _{1})+c_{2}\sigma (x-\theta
_{2}).
\end{equation*}

Let us now consider the other polynomial function
\begin{equation*}
f(x) = 1 + x + \frac{x^2}{2} + \frac{x^3}{6} + \frac{x^4}{24} + \frac{x^5}{%
120} + \frac{x^6}{720}.
\end{equation*}
For this function we do not have an exact representation as above.
Nevertheless, one can easily construct a $\varepsilon$-approximating network
with two neurons in the hidden layer for any sufficiently small
approximation error $\varepsilon$. Table~\ref{tbl:polynomial} displays
numerical computations of the network parameters for six different
approximation errors.

\begin{table}[tbp]
\caption{Several $\protect\varepsilon$-approximators of the function $1 + x
+ x^2 / 2 + x^3 / 6 + x^4 / 24 + x^5 / 120 + x^6 / 720$}
\label{tbl:polynomial}%
\begin{tabular}{|c|c|c|l|c|c|}
\hline
Number of & \multicolumn{4}{|c|}{Parameters of the network} & Maximum \\
\cline{2-5}
neurons & $c_1$ & $c_2$ & \multicolumn{1}{c|}{$\theta_1$} & $\theta_2$ &
error \\ \hline
$2$ & $2.0619 \times 10^{2}$ & $2.1131 \times 10^{2}$ & $-1979$ & $-3$ & $%
0.95$ \\ \hline
$2$ & $5.9326 \times 10^{2}$ & $6.1734 \times 10^{2}$ & $-1.4260 \times
10^{8}$ & $-3$ & $0.60$ \\ \hline
$2$ & $1.4853 \times 10^{3}$ & $1.5546 \times 10^{3}$ & $-4.0140 \times
10^{22}$ & $-3$ & $0.35$ \\ \hline
$2$ & $5.1231 \times 10^{2}$ & $5.3283 \times 10^{2}$ & $-3.2505 \times
10^{7}$ & $-3$ & $0.10$ \\ \hline
$2$ & $4.2386 \times 10^{3}$ & $4.4466 \times 10^{3}$ & $-2.0403 \times
10^{65}$ & $-3$ & $0.04$ \\ \hline
$2$ & $2.8744 \times 10^{4}$ & $3.0184 \times 10^{4}$ & $-1.7353 \times
10^{442}$ & $-3$ & $0.01$ \\ \hline
\end{tabular}%
\end{table}

\begin{figure}[tbp]
\includegraphics[width=1.0\textwidth]{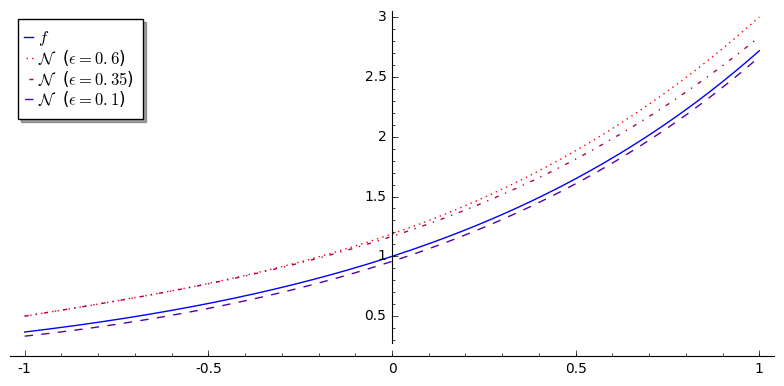}
\caption{The graphs of $f(x) = 1 + x + x^2 / 2 + x^3 / 6 + x^4 / 24 + x^5 /
120 + x^6 / 720$ and some of its approximators ($\protect\lambda = 1/4$)}
\label{fig:polynomial}
\end{figure}

\begin{figure}[tbp]
\includegraphics[width=1.0\textwidth]{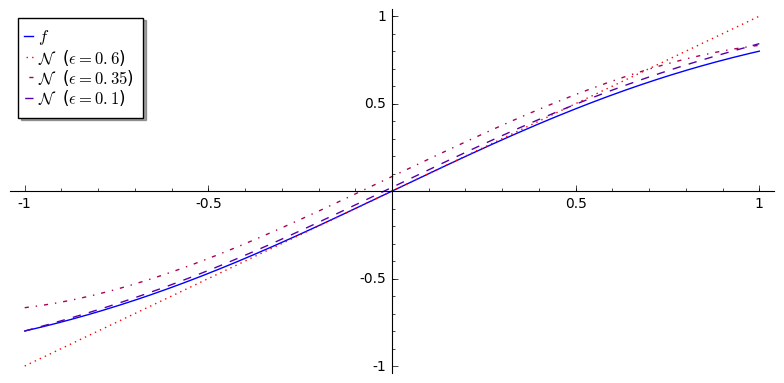}
\caption{The graphs of $f(x) = 4x / (4 + x^2)$ and some of its approximators
($\protect\lambda = 1/4$)}
\label{fig:nonpolynomial}
\end{figure}

\begin{figure}[tbp]
\includegraphics[width=1.0\textwidth]{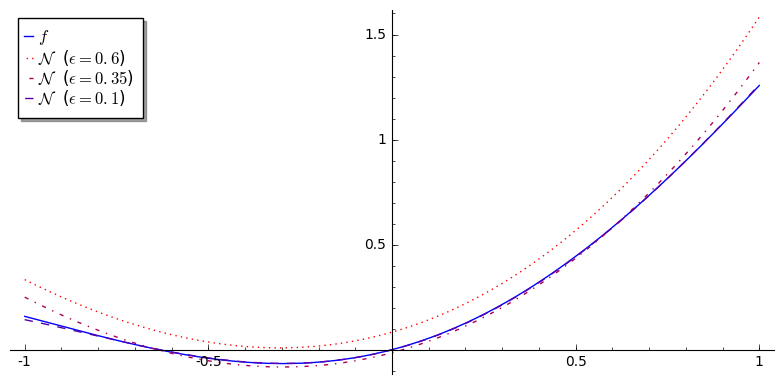}
\caption{The graphs of $f(x)=\sin x-x\cos (x+1)$ and some of its
approximators ($\protect\lambda =1/4$)}
\label{fig:nonpolynomial2}
\end{figure}

At the end we consider the nonpolynomial functions $f(x) = 4x / (4 + x^2)$
and $f(x) = \sin x - x \cos(x + 1)$. Tables~\ref{tbl:nonpolynomial} and \ref%
{tbl:nonpolynomial2} display all the parameters of the $\varepsilon$%
-approximating neural networks for the above six approximation error bounds.
As it is seen from the tables, these bounds do not alter the number of
hidden neurons. Figures~\ref{fig:polynomial}, \ref{fig:nonpolynomial} and %
\ref{fig:nonpolynomial2} show how graphs of some constructed networks $%
\mathcal{N}$ approximate the corresponding target functions $f$.

\begin{table}[tbp]
\caption{Several $\protect\varepsilon$-approximators of the function $4x /
(4 + x^2)$}
\label{tbl:nonpolynomial}%
\begin{tabular}{|@{\hspace{4pt}}c|@{\hspace{4pt}}c|c|l|c|c|}
\hline
Number of & \multicolumn{4}{|c|}{Parameters of the network} & Maximum \\
\cline{2-5}
neurons & $c_1$ & $c_2$ & \multicolumn{1}{c|}{$\theta_1$} & $\theta_2$ &
error \\ \hline
$2$ & $\phantom{-}1.5965 \times 10^{2}$ & $\phantom{-}1.6454 \times 10^{2}$
& $-283$ & $-3$ & $0.95$ \\ \hline
$2$ & $\phantom{-}1.5965 \times 10^{2}$ & $\phantom{-}1.6454 \times 10^{2}$
& $-283$ & $-3$ & $0.60$ \\ \hline
$2$ & $-1.8579 \times 10^{3}$ & $-1.9428 \times 10^{3}$ & $-6.1840 \times
10^{11}$ & $-3$ & $0.35$ \\ \hline
$2$ & $\phantom{-}1.1293 \times 10^{4}$ & $\phantom{-}1.1842 \times 10^{4}$
& $-4.6730 \times 10^{34}$ & $-3$ & $0.10$ \\ \hline
$2$ & $\phantom{-}2.6746 \times 10^{4}$ & $\phantom{-}2.8074 \times 10^{4}$
& $-6.8296 \times 10^{82}$ & $-3$ & $0.04$ \\ \hline
$2$ & $-3.4218 \times 10^{6}$ & $-3.5939 \times 10^{6}$ & $-2.9305 \times
10^{4885}$ & $-3$ & $0.01$ \\ \hline
\end{tabular}%
\end{table}

\begin{table}[tbp]
\caption{Several $\protect\varepsilon$-approximators of the function $\sin x
- x \cos(x + 1)$}
\label{tbl:nonpolynomial2}%
\begin{tabular}{|c|c|c|l|c|c|}
\hline
Number of & \multicolumn{4}{|c|}{Parameters of the network} & Maximum \\
\cline{2-5}
neurons & $c_1$ & $c_2$ & \multicolumn{1}{c|}{$\theta_1$} & $\theta_2$ &
error \\ \hline
$2$ & $\phantom{-}8.950 \times 10^{3}$ & $\phantom{-}9.390 \times 10^{3}$ & $%
-3.591 \times 10^{53}$ & $-3$ & $0.95$ \\ \hline
$2$ & $\phantom{-}3.145 \times 10^{3}$ & $\phantom{-}3.295 \times 10^{3}$ & $%
-3.397 \times 10^{23}$ & $-3$ & $0.60$ \\ \hline
$2$ & $\phantom{-}1.649 \times 10^{5}$ & $\phantom{-}1.732 \times 10^{5}$ & $%
-9.532 \times 10^{1264}$ & $-3$ & $0.35$ \\ \hline
$2$ & $-4.756 \times 10^{7}$ & $-4.995 \times 10^{7}$ & $-1.308 \times
10^{180281}$ & $-3$ & $0.10$ \\ \hline
$2$ & $-1.241 \times 10^{7}$ & $-1.303 \times 10^{7}$ & $-5.813 \times
10^{61963}$ & $-3$ & $0.04$ \\ \hline
$2$ & $\phantom{-}1.083 \times 10^{9}$ & $\phantom{-}1.138 \times 10^{9}$ & $%
-2.620 \times 10^{5556115}$ & $-3$ & $0.01$ \\ \hline
\end{tabular}%
\end{table}

\end{document}